\documentclass[12pt,a4paper]{article}

\usepackage[T2A]{fontenc}
\usepackage[utf8]{inputenc}
\usepackage[russian,english]{babel}
\usepackage{amsmath,amssymb,amsthm}
\usepackage[margin=2.5cm]{geometry}
\usepackage{enumitem}
\usepackage[unicode,colorlinks,linkcolor=blue,citecolor=blue,urlcolor=blue]{hyperref}

\def\phi{\varphi}

\newcommand{\thmTheorem}{Theorem}
\newcommand{\thmProposition}{Proposition}
\newcommand{\thmLemma}{Lemma}
\newcommand{\thmCorollary}{Corollary}
\newcommand{\thmDefinition}{Definition}
\newcommand{\thmRemark}{Remark}
\newcommand{\thmQuestion}{Question}
\newcommand{\thmConjecture}{Conjecture}
\newcommand{\thmExample}{Example}

\theoremstyle{plain}
\newtheorem{theorem}{\thmTheorem}
\newtheorem{proposition}[theorem]{\thmProposition}
\newtheorem{lemma}[theorem]{\thmLemma}
\newtheorem{corollary}[theorem]{\thmCorollary}
\theoremstyle{definition}
\newtheorem{definition}{\thmDefinition}
\newtheorem{remark}{\thmRemark}

\newtheorem{example}{\thmExample}

\newcommand{\PA}{\mathrm{PA}}
\newcommand{\N}{\mathbb{N}}
\newcommand{\gn}[1]{#1}
\newcommand{\Prf}{\mathrm{Prf}}
\newcommand{\Prov}{\mathrm{Pr}}
\newcommand{\Con}{\mathrm{Con}}
\newcommand{\val}[2][*]{[\![#2]\!]_{#1}}
\newcommand{\Tru}{\text{T}}
\newcommand{\Fls}{\text{F}}
\newcommand{\Int}{\mathrm{Int}}
\newcommand{\Lp}{\mathcal L}
\newcommand{\J}{\mathfrak J}
\newcommand{\R}{\mathfrak R}
\newcommand{\Reff}{\R_{\mathrm{eff}}}
\newcommand{\Rcnst}{\R_{\mathrm{const}}}
\newcommand{\GL}{\mathrm{GL}}
\newcommand{\GLS}{\mathrm{GLS}}
\newcommand{\HC}{\mathrm{HC}}
\newcommand{\HCI}{\HC \mathrm I}
\newcommand{\HCB}{\HC  \mathrm B}
\newcommand{\HCE}{\HC \mathrm E}
\newcommand{\HCgl}{\HC \mathrm I}
\newcommand{\HCglm}{\HCI^{-}}
\newcommand{\HCgls}{\HCE}
\newcommand{\HCzero}{\HC^{0}}
\newcommand{\ZF}{\mathrm{ZF}}
\newcommand{\TA}{\mathrm{TA}}

\newcommand{\sI}{\mathsf{I}}
\newcommand{\sE}{\mathsf{E}}
\newcommand{\sB}{\mathsf{B}}
\newcommand{\nf}[1]{{{#1}^{\dag}}}
\newcommand{\sS}{\mathsf{S}}

\newcommand{\ato}{\supset}
\newcommand{\aand}{\wedge}
\newcommand{\aor}{\vee}
\newcommand{\aneg}{\neg}
\newcommand{\aeq}{\equiv}
\newcommand{\arbot}{\bot}

\newcommand{\pto}{\to}
\newcommand{\pand}{\curlywedge}
\newcommand{\por}{\curlyvee}
\newcommand{\pneg}{{\sim}}
\newcommand{\peq}{\leftrightarrow}
\newcommand{\bang}{{\textnormal{\sffamily\bfseries !}}}
\newcommand{\quest}{{\textnormal{\sffamily\bfseries ?}}}

\newcommand{\remove}[1]{}

\makeatletter
\let\saved@maketitle\maketitle
\let\saved@@maketitle\@maketitle
\let\saved@title\title
\let\saved@author\author
\let\saved@date\date
\let\saved@thanks\thanks
\let\saved@and\and
\makeatother

\title{Arithmetical and set-theoretic semantics\\
of the logic of problems and propositions}
\author{Anastasia Onoprienko \qquad Nikolay Vereshchagin\\[6pt]
  \normalsize Lomonosov Moscow State University}
\date{}

\begin{document}
\maketitle

\begin{abstract}
Kolmogorov interpreted intuitionistic logic as a calculus of problems and
later remarked that a complete treatment should deal with problems and
propositions together. Melikhov's logic $\HC$ is such a joint two-sorted
system: it combines classical propositions with intuitionistic problems by
means of two modalities, the problem $\bang A$ (``find a proof of the
proposition $A$'') and the proposition $\quest\varphi$ (``the problem
$\varphi$ has a solution''). We give the language of $\HC$ a formal
semantics in which the modalities are read literally: solutions of
$\bang A$ are proofs of $A$ in a fixed formal theory, $\quest\varphi$
states that $\varphi$ has a solution, and the connectives on problems are
interpreted as in Kleene realizability. Three variants are studied: the
arithmetical semantics $\sI$ (provability in $\PA$), its external version
$\sE$ (truth in the standard model), and the set-theoretic semantics $\sB$
over $\ZF$, where problem variables are ``black boxes'' and a quantifier
$\Pi$ over them is available. Validity is defined as the existence of a
computable uniform solution. We determine which schemes are valid, build
calculi $\HCI$, $\HCE$, $\HCB$ sound with respect to the three semantics,
and derive L\"ob's scheme from the fixed-point lemma inside the common
core of our calculi and $\HC$. The propositional fragments of $\HCI$ and
$\HCE$ are exactly the provability logics $\GL$ and $\GLS$; the positive
modality-free fragment of valid problems is intuitionistic logic, while
the full modality-free fragment lies strictly above it (Rose's formula is
valid, Scott's principle is not). We compare the calculi with $\HC$,
which turns out to be incompatible with fixed points, prove independence
results and list open questions. The Russian original of the paper is
appended.
\end{abstract}

\tableofcontents

\section{Introduction}

In Kolmogorov's paper~\cite{Kolm32} the intuitionistic calculus was
interpreted as a \emph{calculus of problems}. And in his comments on this paper in the
Selected Works~\cite{Kolm85} Kolmogorov wrote that problems and propositions are different things and that a complete exposition
should speak about both sorts at once: intuitionistic logic describes
problems, classical logic---propositions, and one needs a common system containing
both together with the connection between them.

An example of such a system is Melikhov's logic $\HC$~\cite{Melikhov13,Melikhov24}
---a joint logic of problems and propositions. It combines in one two-sorted
language classical propositions and intuitionistic problems, connecting the two sorts
by two modalities: the problem $\bang A$ is read as ``find a proof of
the proposition $A$'', and the proposition $\quest\varphi$---as ``the problem $\varphi$
has a solution'', that is, exactly as Kolmogorov had outlined.
Several kinds of semantics are known for $\HC$: topological ones,
interpretations in plane geometry and in algebraic equations, an informal
semantics going back to Kolmogorov himself (Melikhov,
\cite{Melikhov13,Melikhov15}), and also a Kripke semantics with marked
worlds (Onoprienko, \cite{Onoprienko20}). We, however, are interested in a semantics
in which the modalities are read literally: the solutions of the problem $\bang A$
are the proofs of the proposition $A$ in a fixed formal theory, and the
proposition $\quest\varphi$ asserts that $\varphi$ has a solution.
Such a semantics---\emph{arithmetical}, in the spirit of Kleene realizability---
for the full language of $\HC$ did not exist until now: for the fragment without
the modality $\quest$ it was sketched in~\cite{Melikhov24}, but remained
unfinished.

In the present paper such a semantics is constructed and the notion of validity
it generates is studied; then calculi sound with respect to this semantics
are introduced; the axioms are selected by validity. The order itself,
``from interpretation to calculus'', is not new: the logic
$\HC$ is built in the same way, its axioms being extracted from Kolmogorov's interpretation of problems and
the Orlov--Heyting interpretation of proofs~\cite{Melikhov24}.
The difference is that our interpretation is a formal semantics, and
validity in it is a precise mathematical notion, so that the question
of soundness and completeness of the calculus has a precise meaning. This order---
statement of the problem, the semantics and its study, the calculi and their study
---also determines the plan of the paper.

 More precisely, we define three different semantics. We begin with the first two, the arithmetical ones,
denoted $\sI$ and $\sE$.
In them a proposition $P$ is translated into a closed
formula $P^*$ of Peano arithmetic, and a problem $\varphi$---into a formula $\varphi^*(x)$ with (at most) one free variable;
the connectives on problems are understood in the sense of Kleene~\cite{Kleene45}, the proposition $\quest\varphi$ is understood as
$\exists x\,\varphi^*(x)$, and the solutions of the problem $\bang A$ are the G\"odel
numbers of proofs of the formula $A^*$ in $\PA$ (\S\ref{sec:val}). It turns out that not all
axioms of HC are true in this semantics, that is, HC is not sound with respect to it.
The most important point of our interpretation is the notion of validity of a problem: a solution of the problem
must be exhibited \emph{by an algorithm computing it from the code of the
interpretation} (Definition~\ref{def:valid}). As shown in
Remark~\ref{rem:three}, it is exactly this requirement that decides whether the
fragment of problems remains intuitionistic: under a weaker understanding of validity it
becomes classical.
In the semantics $\sI$ it is required that under every
interpretation the formula be provable in PA, and in the semantics $\sE$---true in the standard interpretation of PA.
In the semantics $\sI$ the set of valid propositions
essentially coincides with the provability logic $\GL$,
while the set of valid problems is a rather complicated
logic whose decidability
is unknown.
In the semantics $\sE$ the set of valid propositions
essentially coincides with the provability logic $\GLS$,
while the set of valid problems is again a rather complicated
logic whose decidability is unknown.

Then we present calculi (logics) sound with respect to these
semantics.
 We did not strive to obtain calculi close to $\HC$: from it we
take only the language and the informal understanding of the modalities, while the axioms
were selected by validity---if a scheme is valid and is not derivable from
those already accepted, it was added. This is how we obtained the calculus
$\HCI$, corresponding to what
$\PA$ proves, and the calculus $\HCE$, corresponding to what is true
(\S\ref{sec:calc}); their soundness is Theorem~\ref{th:sound}.

What are such calculi good for? For instance, using the calculus $\HCI$ one can
obtain a formal derivation of L\"ob's theorem
from the fixed-point theorem for $\PA$.
A peculiarity of these calculi is that almost all their axioms are
problems. There are three exceptions: the axioms of the classical propositional
calculus (A1); the purely structural schemes (A4), asserting
that $\quest$ translates connectives of problems into connectives of propositions (a Boolean
homomorphism); and the reflection (A8) of the calculus $\HCE$---it is not
structural but substantive, yet it is a proposition, not a problem.
All the other axioms are problems. They either speak of the possibility of an algorithmic transformation
of some arithmetical proofs into others: for example, (A5)---``from a solution of the problem $\phi$ (for problems $\phi$ of a special form) one can construct
an arithmetical proof that $\phi$ has a solution'', (A6)---
``an arithmetical proof of the formula $\square A \ato A$
can be transformed into an arithmetical proof of the formula $A$'' (L\"ob's theorem), or
of the realizability of certain problems without the modality $\quest$---the axiom (A7).
On the level of propositions the calculi $\HCI$ and $\HCE$ are complete and coincide with GL and GLS, while on the level of problems
one again obtains
certain  
logics about which we do not know whether they are decidable.

In the third semantics $\sB$ we slightly extend the language by introducing the universal quantifier $\Pi$
over problem variables. The problem $\Pi\alpha\,\varphi$ is understood as ``exhibit a program that solves
the problem $\varphi$ under every interpretation of the problem $\alpha$''.
Recall that without the quantifier over $\alpha$ this program may depend computably
on the interpretation $\alpha^*$ of the variable $\alpha$.
Why do we need such a quantifier? The reason is a strong asymmetry between the two fragments, the problem one and the propositional one.
It might seem that we have translations in both directions: the modality $\bang$ turns a problem
$A$ into the proposition $\bang A$, and the modality $\quest$ acts in the opposite direction.
But it turns out that the first translation is faithful: validity of the problem $\bang A$
is equivalent to validity of the proposition $A$.
The translation $\phi\mapsto \quest\varphi$, however,
does not have this property. Often, proving in PA the proposition
$(\quest\varphi)^*$ under every fixed interpretation $*$
is much easier than proving the solvability of the problem $\varphi^*$ uniformly in $*$.
This is the source of the
asymmetry: the set of valid propositions is simple, while the set of
valid problems is not.

For example,
the law of excluded middle $\alpha\por\pneg\alpha$ is not valid (there is no algorithm that,
given an arithmetical formula $\alpha^*(x)$, would find out whether the proposition
$\exists x \alpha^*(x)$ is true). At the same time the proposition $?(\alpha\por\pneg\alpha)$ is
valid, because PA is a classical theory and can
use the law of excluded middle to prove
the formula $\exists x (\alpha\por\pneg\alpha)^*(x)$ for every fixed
formula $\alpha^*(x)$.  Another example: Rose's formula $\rho$ is
a valid formula, and the algorithm realizing it is not trivial.
But after prefixing the modality $\quest$ one obtains the trivially
valid formula $? \rho$: under every interpretation $*$ the formula  $\exists x \rho^*(x)$
is easily provable in $\PA$ by means of the law of excluded middle.

Unfortunately, to write down the assertion that the problem $\quest\Pi\alpha\,\varphi$ is solvable,
one needs a universal quantifier over arithmetical interpretations. In Peano arithmetic this is apparently inexpressible, since by Tarski's theorem it is impossible to express the truth of the formula with a given G\"odel number.
Therefore instead of PA we pass to the set theory ZF. In the semantics $\sB$ (from the words ``black box'')
the interpretations of problem variables are \emph{arbitrary} subsets of $\N$---``black boxes'',
which the solving program does not see; the universal quantifier over them is a first-order formula of the language of ZF,
and the proposition ``the program $p$ solves the problem $\varphi$ under every interpretation of the problem $\alpha$''
is written as a set-theoretic formula. Validity in $\sB$
is defined via provability in ZF.
A sound calculus for it is arranged as $\HCI$ plus
a quantifier layer: the quantifier-free part of its axiomatics is exactly
$\HCI$. The purpose of the extension was
to enlarge the stock
of valid propositions, and in the extended language this certainly
succeeded: for example, the proposition
$\quest\Pi\beta\Pi\gamma\rho$---``Rose's formula is uniformly
solvable''---is valid, and it is inexpressible without the quantifier. The extension also yields new \emph{quantifier-free}
theorems---for example, $\square\quest\alpha\ato\square\bot$---but
only because in the semantics $\sB$ the modality $\bang$ binds
problem variables and is read with closure. For \emph{clean} formulas,
in which $\bang$ is applied only to closed propositions,
the picture is different: on the propositional side the extension is conservative over
$\HCI$ (assuming $\Sigma_1$-soundness of $\ZF$), while for clean problems the question of
conservativity is open.
We also do not know whether in the new semantics the fragment of problems has the same complexity as
the fragment of propositions.\footnote{Unfortunately, the translation $\phi \mapsto ?(\text{closure of } \phi)$, which could prove that the complexity of the problem fragment is at most that of the propositional one,
may translate valid problems into non-valid propositions.}

Validity of all axioms of our calculi is summarized in
Table~\ref{tab:axioms} (\S\ref{sec:tables}), from which the structure of the
calculi can be seen directly.


\paragraph{The place of our calculi among related logics.} Among the systems with which the result
is worth comparing, $\HC$ stands out in one respect: it has the same syntax. Everything else
diverges. From $\HC$ we take the language and the reading of the modalities, but not the list of theorems:
our calculi are in general position with respect to HC.

On the propositional side we land in classical provability logic:
the $\square$-fragments of the calculi $\HCI$ and $\HCE$ are exactly Solovay's $\GL$ and $\GLS$~\cite{Solovay76} (Corollary~\ref{cor:solovay}), and two-sortedness gives
nothing new there. Intuitionistic provability logics describe something else:
$\mathrm{iGL}$~\cite{VZ19,DvdGM24} is strictly contained in our
$\triangledown$-fragment, while $\mathrm{iGLC}$ is not a fragment of ours
---its completeness principle $A\ato\square A$ is co-reflection. Still,
there is a similarity: there completeness is required only on $\Sigma_1$, and for us
co-reflection holds exactly on $\Sigma$-problems
(axiom (A5)).

In Artemov's school~\cite{Artemov01,Sidon97,Nogina14,Protopopescu16} proofs
are also present explicitly, but differently: there the language is enriched with terms $t{:}F$,
``$t$ is a proof of $F$'', whereas for us proofs are the
solutions of the problem $\bang A$, and the language remains modal. The two-sortedness ``a problem is
a set, a solution is a number'' is absent there.

Finally, problems as sets are a long-standing tradition: the mass problems of Medvedev and
Muchnik and the logics of realizability~\cite{Rose53,Plisko09}. The difference is that there
is one sort and no modalities: there is no place there for proofs and for the modality $\bang$. On the other hand, our notion of validity coincides exactly with one of the
known ones: it is $\Reff$---``a realization can be found effectively from the substituted
formulas''~\cite{VSSCh}. Hence the results on the fragment of problems
(\S\ref{sec:modalfreeval}): its positive part is intuitionistic, and the fragment itself lies
strictly between $\Int$ and the logic of weak excluded middle.

The combination itself---two sorts, problems as arbitrary arithmetically
definable sets with Kleene's operations, $\bang A$ as the set of proofs,
$\quest$ as non-emptiness, and uniformity in the code of the interpretation---we have
not found in the literature.

The plan of the paper follows the order ``problem---semantics---calculi''.
The semantics and its study occupy \S\ref{sec:val}--\S\ref{sec:validity}:
the interpretations of all three semantics, the notion of validity and everything
we know about valid formulas: validity of the future axioms,
a description of the modality-free fragment of problems, examples of non-valid formulas,
the laws of the operator $\Pi$ (\S\ref{sec:piatoms}) and the summary table
of validity (\S\ref{sec:tables}). The calculi and their study
---\S\ref{sec:calc}: three calculi, one per semantics,
derived rules and examples of derivations, soundness, what the quantifier
$\Pi$ gives, normal form and
completeness on the propositional side together with the characterization of valid
propositions by the logics $\GL$ and $\GLS$ (\S\ref{sec:asscompl}), normal
form with $\Pi$-atoms (\S\ref{sec:pinf}),
independence of the axioms and the little that is known about derivable problems. Comparison with $\HC$ and other related
logics---\S\ref{sec:hc}, open questions---\S\ref{sec:todo}. The
appendices contain the deferred proofs, an analysis of the classes of problems
in the schemes (A5) and (A7), the independence of the scheme (A7) and an analysis of what is new
in our derivation of L\"ob's scheme.

\section{Language and notation}\label{sec:lang}
We first define the syntax for the arithmetical semantics.
We use the language of $\HC$, but change the letters so that the sort
of a formula is visible at once:
\begin{itemize}[leftmargin=1.6em,itemsep=1pt]
\item \textbf{capital Latin} letters $A,B,C,\dots$ denote formulas
  of \emph{sort proposition};
\item \textbf{Greek} letters denote formulas of \emph{sort problem}.
\end{itemize}
In doing so we adhere to the following convention. The letters $\alpha$, $\beta$,
$\gamma$ (with indices if necessary) denote only \emph{variables}
of sort problem, whereas the letters $\varphi$, $\psi$, $\chi$, $\theta$ denote arbitrary
problems. For the sort of propositions we do not reserve separate letters for variables:
capital Latin letters denote, generally speaking, arbitrary
propositions, and where it matters that a letter denotes a variable,
this is stated explicitly.
Formulas are built as follows:
\[
  A ::= A_i \mid \bot \mid (A\aand A) \mid (A\aor A)
             \mid (A\ato A)\mid \quest\varphi ,
\]
\[
  \varphi ::= \alpha_i \mid \varnothing \mid (\varphi\pand\varphi) \mid
              (\varphi\por\varphi) \mid (\varphi\pto\varphi) \mid \bang A,
\]
where $A_i$ are variables of sort proposition and $\alpha_i$ are variables of sort
problem.

The constant $\bot$ is the false proposition, the constant $\varnothing$ is
the unsolvable problem (the problem with the empty set of solutions). Negation and equivalence are not part of
the language but are defined:
\[
  \aneg A:=A\ato\bot ,\ \pneg\varphi:=\varphi\pto\varnothing, \
 A\aeq  B :=(A\ato B)\ \aand\ (B\ato A) ,\ \phi\peq  \psi :=(\phi\pto \psi)\ \pand\  (\psi\pto \phi).
\]

\medskip
\noindent\textbf{The symbols of the connectives differ between the sorts.} The language is two-sorted, and one
and the same connective means different things on the two sorts: the implication of propositions is
the usual classical one, the implication of problems is Kleene's ``from a solution of the premise
compute a solution of the conclusion''. Therefore each sort gets its own symbols:
\[
\begin{array}{l|c|c}
   & \text{propositions and arithmetic} & \text{problems}\\[2pt]\hline
  \text{implication}      & \ato   & \pto\\
  \text{conjunction}      & \aand  & \pand\\
  \text{disjunction}      & \aor   & \por\\
  \text{negation}       & \aneg  & \pneg\\
  \text{equivalence} & \aeq   & \peq\\
  \text{falsity constant}   & \bot   & \varnothing
\end{array}
\]
Propositions received the classical pair $\ato$, $\aeq$, problems
the arrows $\pto$, $\peq$; conjunction, disjunction and negation on problems
are marked by a curl and a tilde.

The point of the distinction is already visible in the scheme (A3), which is written as
$\bang(A\ato B)\pto(\bang A\pto\bang B)$: inside the parentheses stands the implication of propositions,
outside stand two implications of problems, and these are different things. The symbols $\ato$, $\aand$, $\aor$, $\bot$,
$\aneg$, $\aeq$ are also used inside \emph{arithmetical} formulas (in $A^*$,
$\varphi^*$ and the like): the language of arithmetic is classical, and there is no reason to introduce
a third set of symbols for it. Connectives of sort problem never occur in arithmetical
formulas, and the asterisk helps to distinguish an object formula from an arithmetical one.
Derived modalities: $\square A=\quest\bang A$ (a proposition) and
$\triangledown\varphi=\bang\quest\varphi$ (a problem).

\section{Arithmetical preliminaries}

We fix the language of arithmetic and the standard model $\N$; ``true'' everywhere means
``true in $\N$''. Truth values are denoted by $\Tru$ and $\Fls$.
We fix:
\begin{itemize}[leftmargin=1.6em,itemsep=2pt]
\item a G\"odel numbering of unary computable (in general, partial)
  functions
  \[
    U_0,\ U_1,\ U_2,\ \dots
  \]
\item a computable bijection $c\colon\N^2\to\N$ (pairing);
\item an arithmetical formula $\Prf(x,y)$ expressing ``$x$ is (the number of) a proof
  in $\PA$ of the formula with number $y$'', and $\Prov(y)=\exists x\,\Prf(x,y)$. As
  usual, $\Con=\aneg\Prov\arbot $.
\end{itemize}

\begin{remark}[lightening the notation]\label{rem:notation}
We \emph{identify a formula with its G\"odel number}: where formally
one should write $\Prf(x,\ulcorner\sigma\urcorner)$, we write
$\Prf(x,\sigma)$. Likewise, where the value $n$ of a computable function
is substituted into an arithmetical formula, formally one should
substitute the numeral $\bar n$, the standard term denoting $n$; we write
simply $n$. Finally, the expression $U_i(u)$ inside arithmetical formulas is read
via the graph: $\psi(U_i(u))$ is an abbreviation for $\exists w\,\bigl(U_i(u)
\simeq w\aand\psi(w)\bigr)$, where $U_i(u)\simeq w$ is a $\Sigma_1$-formula
expressing ``$U_i$ is defined on $u$ and equals $w$''; in particular, such an expression
implies definedness. All these provisos are made only for readability and are nowhere
used in an essential way.
\end{remark}

\section{Interpretations}\label{sec:val}

\subsection{Arithmetical interpretations}

\begin{definition}[arithmetical interpretation]\label{def:star}
An \emph{arithmetical interpretation} 
is a map $*$ assigning
\begin{itemize}[leftmargin=1.6em,itemsep=1pt]
\item to each variable $A$ of sort proposition a closed arithmetical
  formula $A^*$;
\item to each variable $\alpha$ of sort problem an arithmetical formula
  $\alpha^*(x)$ with (at most) one free variable $x$.
\end{itemize}
\end{definition}

The interpretation is extended to all formulas of the language of $\HC$ by \emph{simultaneous
recursion}: to a proposition $A$ a closed formula $A^*$ is assigned, to a problem
$\varphi$ a formula $\varphi^*(x)$ with (at most) one free variable.

\begin{definition}[extension of $*$]\label{def:tr}
For propositions:
\[
\begin{aligned}
  \bot^*&:=\arbot ,\\
  (A\aand B)^*&:=A^*\aand B^*, &
  (A\aor B)^*&:=A^*\aor B^*,\\
  (A\ato B)^*&:=A^*\ato B^*, &
  (\quest\varphi)^*&:=\exists x\,\varphi^*(x). &&
\end{aligned}
\]
For problems:
\[
\begin{aligned}
  \varnothing^*(x)&:=\arbot ,\\
  (\varphi\pand\psi)^*(x)&:=\exists u\,\exists v\,
     \bigl(x=c(u,v)\aand\varphi^*(u)\aand\psi^*(v)\bigr), &&\\
  (\varphi\por\psi)^*(x)&:=\exists u\,
     \bigl((x=c(0,u)\aand\varphi^*(u))\aor(x=c(1,u)\aand\psi^*(u))\bigr), &&\\
  (\varphi\pto\psi)^*(x)&:=\forall v\,\bigl(\varphi^*(v)\ato
     \psi^*(U_x(v))\bigr), &&\\
  (\bang A)^*(x)&:=\Prf\bigl(x,\gn{A^*}\bigr). &&
\end{aligned}
\]
\end{definition}

We stress that $A^*$ is always a closed formula: by induction on the construction it is easy
to prove that the translation of any problem is a formula with (at most) one free variable,
and the translation of any proposition is a closed formula.

\begin{definition}[value]\label{def:val}
For a proposition $A$ we put $\val{A}:=\Tru$ if the formula $A^*$ is true in
$\N$, and $\val{A}:=\Fls$ otherwise. For a problem $\varphi$ we put
\[
  \val{\varphi}:=\{n\in\N\ :\ \N\models\varphi^*(n)\}.
\]
\end{definition}

Thus the value of a proposition is a truth value, and the value of a problem is
a set of natural numbers (the ``set of solutions'', which is also the set of
realizations). The index $*$ in the notation $\val{\cdot}$ reminds us that the value
depends on the interpretation.

\begin{lemma}[the connectives behave as promised]\label{lem:clauses}
For every interpretation $*$:
\begin{enumerate}[leftmargin=2.2em,itemsep=1pt,label=\textup{(\arabic*)}]
\item $\val{\bot}=\Fls$, and the connectives $\aand,\aor,\ato$ on propositions
  act classically;
\item $\val{\quest\varphi}=\Tru\iff\val{\varphi}\neq\varnothing$;
\item $\val{\varnothing}=\varnothing$;
\item $\val{\varphi\pand\psi}=\{c(u,v):u\in\val{\varphi},\,v\in\val{\psi}\}$;
\item $\val{\varphi\por\psi}=\{c(0,u):u\in\val{\varphi}\}\cup
      \{c(1,v):v\in\val{\psi}\}$;
\item $\val{\varphi\pto\psi}=\{i:\ \text{for all}\ u\in\val{\varphi}\
      \text{the value}\ U_i(u)\ \text{is defined and lies in}\ \val{\psi}\}$;
      in particular $\val{\pneg\varphi}=\N$ if $\val{\varphi}=\varnothing$, and
      $\val{\pneg\varphi}=\varnothing$ otherwise;
\item $\val{\bang A}=\{n:\ n\ \text{is the number of a proof of the formula}\ A^*\
      \text{in}\ \PA\}$.
\end{enumerate}
\end{lemma}

\begin{lemma}[double negation is $\quest$ on the problem side]\label{lem:dn}
For every problem $\varphi$ and every interpretation $*$
\[
  \val{\pneg\pneg\varphi}=
  \begin{cases}
    \N, & \text{if } \val{\varphi}\neq\varnothing,\\
    \varnothing, & \text{otherwise,}
  \end{cases}
  \qquad\text{that is,}\qquad
  \val{\pneg\pneg\varphi}\neq\varnothing\iff\val{\quest\varphi}=\Tru .
\]
Moreover, $\PA$ proves $\forall x((\pneg\pneg\varphi)^*(x)\aeq
(\quest\varphi)^*)$,
and a proof of this equivalence is constructed from
$\gn{*}$.
\end{lemma}

The proof is in the appendix.

\remove{
\begin{remark}\label{rem:dn}
Итак, $\pneg\pneg$ --- это тот же оператор <<решение существует>>, только
оставшийся на стороне задач. Отсюда видно, почему $\val{\pneg\pneg\varphi}$
двузначно ($\N$ или $\varnothing$) и почему никакая схема не может извлечь из
$\pneg\pneg\varphi$ решение самой $\varphi$: в $\N$ выбирать нечего.
\end{remark}

\begin{remark}[почему рекурсия должна тащить формулу, а не только значение]
Значение $\val{\bang A}$ зависит не от значения $\val{A}$, а от самой формулы
$A^*$: два высказывания с одним и тем же истинностным значением (скажем, две
разные истинные формулы, из которых одна доказуема, а другая нет) дают разные
$\val{\bang A}$. Поэтому семантика не композициональна на уровне значений, и
определение~\ref{def:tr} обязано определять перевод $*$ параллельно со
значением.
\end{remark}
}

\subsection{Set-theoretic semantics: black boxes}\label{sec:valBS}

Alongside the arithmetical interpretation we consider one more,
a set-theoretic one. Together with the two levels of validity
(\S\ref{sec:validity}) this gives three semantics in total: $\sI$ and $\sE$
(arithmetical) and $\sB$ (black boxes).

\paragraph{Language.} For the new semantics we add to the grammar generating
formulas two more rules---one for each sort:
\[
\varphi \leftarrow  \Pi\alpha_i\,\varphi, \qquad
A \leftarrow  \Pi\alpha_i\,A .
\]
The operator $\Pi\alpha$ (``uniformly in $\alpha$'') applies to
formulas of both sorts and \emph{binds}
the occurrences of the variable $\alpha$ within its scope. Moreover---
and this is essential---in the semantics $\sB$ the modality $\bang$ is
also regarded as binding: the problem $\bang A$ binds all problem
variables free in $A$. The reason will become clear from the definitions:
the solutions of the problem $\bang A$ are proofs of a single formula, and a proof cannot
depend on the contents of a box.

Free and bound occurrences, renaming of bound variables and
the proviso ``$\theta$ is free for $\alpha$ in $\varphi$'' are understood as
usual for quantifiers. By $FV(\varphi)$, $FV(A)$ we denote the set of free
problem variables of a formula (of either sort).

\begin{definition}[the semantics $\sB$]
\label{def:starB}
In the semantics $\sB$ an \emph{interpretation} $*$ assigns to each
variable of sort proposition a closed formula of the language of $\ZF$ (for instance,
an arithmetical one); problem variables are not interpreted---they are
\emph{black boxes}. The interpretation is the ``transparent'' part of the semantics
(it has a code), the boxes are the ``black'' part (a box has no code). The translation
is defined by the same clauses as in Definition~\ref{def:tr}, with
four changes:
\[
  \alpha^*(x):=(x\in X_\alpha),\qquad
  (\bang A)^*(x):=\Prf_{\ZF}\bigl(x,\gn{\forall\vec X\,A^*}\bigr),
\]
\[
  (\Pi\alpha\,\varphi)^*(x):=\forall X_\alpha{\subseteq}\,\omega\
  \varphi^*(x),\qquad
  (\Pi\alpha\,A)^*:=\forall X_\alpha{\subseteq}\,\omega\ A^* .
\]
Here $X_\alpha$ is a set variable of the language of $\ZF$, a separate one for
each problem variable $\alpha$: the boxes are represented in the translation
by the free variables $X_\alpha$. In the clause for $\bang$, $\vec X$
denotes the variables $X_\beta$ for all problem variables $\beta$
free in $A$, that is, those bound by this occurrence of $\bang$; the G\"odel number
of the closure is a numerical term, and the variables $X_\beta$ do not occur in the formula
$(\bang A)^*$. The closure $\forall\vec X\,A^*$ is a
block of quantifiers $\forall X_\beta{\subseteq}\,\omega$ over all such
$\beta$ in a fixed order (say, in increasing order of the indices of the
variables). In other words, $\forall\vec X\,A^*$ is the translation of the
proposition $\Pi\vec\beta\,A$, where $\vec\beta$ are the free problem
variables of $A$ listed in the same order, and the translations of the problems
$\bang A$ and $\bang\Pi\vec\beta\,A$ coincide. We shall use this when
verifying the scheme (A6). The \emph{code} of an interpretation is the code
of the transparent part alone (Definition~\ref{def:code}): a box has
no code.

\emph{Values.} We speak of values of formulas in the semantics $\sB$ only
inside $\ZF$. For given sets $X_\alpha\subseteq\omega$
the value of a problem $\varphi$ is the set
$\val{\varphi}:=\{n\in\omega:\varphi^*(n)\}$, which exists by the separation
axiom; the words ``for every contents of the boxes'' mean a universal
quantifier over $X_\alpha$ inside $\ZF$. One may imagine that the
reasoning takes place in some model of $\ZF$. But $\ZF$ has no canonical
model, and therefore values do not enter the definition of validity
(Definition~\ref{def:valid}).
\end{definition}

\begin{remark}[why $\bang$ binds]\label{rem:bangbinds}
By induction on construction it is easy to prove that the translation $A^*$ of any
proposition $A$ is a formula whose free variables can only
be $X_\alpha$ for free problem variables $\alpha$
of the proposition $A$, and the translation $\varphi^*$ of any problem may have in addition
at most one free numerical variable. A solution of the problem
$\Pi\alpha\,\varphi$ is a \emph{single} solution of $\varphi$ that works for
any contents of the box $\alpha$; the proposition $\Pi\alpha\,A$ reads
``$A$ for any contents of the box $\alpha$''. 

The clause for $\bang$ could also be written more briefly,
$\Prf_{\ZF}(x,\gn{A^*})$, with the open formula $A^*$ under the
number sign: the proofs of an open formula are the proofs of its closure, so that
both forms define the same set of solutions up to a
computable transformation. But the second form hides the quantifier, and this is
dangerous.

Here is an example illustrating this phenomenon. Consider the problem
$\phi:=\bang\quest\alpha$. Its interpretation is the formula
$\phi^*=\Prf(y,\forall X \exists x\in X)$. The formula $\forall X\,\exists x\,(x\in X)$ is refutable in $\ZF$ (the box
$\varnothing$), and from any proof of it, together with a fixed
refutation, one effectively constructs a proof of the formula $\arbot $,
and the correctness of this transformation is provable in $\ZF$; conversely,
from a proof of $\arbot $ one constructs a proof of anything. Thus our problem
$\phi$ is equivalent to the problem $\bang\bot$: the formula
$\phi\peq\bang\bot$ is valid. 
\remove{
Если бы $\bang$ не считалась связывающей, то задача
$\Pi\alpha(\bang\quest\alpha\pto\varnothing)$ была бы общезначима:
посылка $\bang\quest\alpha$ --- выводы опровержимого замыкания
$\forall X\,\exists y\,(y\in X)$, то есть при $\Con(\ZF)$ пустое
множество. А её экземпляр по закону $(\Pi\mathrm{E})$ с
$\theta:=\pneg\varnothing$ --- задача
$\bang\quest\pneg\varnothing\pto\varnothing$ --- решений не имеет:
формула $\exists x\,(\pneg\varnothing)^*(x)$ доказуема, и выводы у неё
есть. Подстановка внутрь гёделева номера не достаёт, и закон
$(\Pi\mathrm{E})$ ломался бы. Считая $\bang$ связывающей, мы получаем,
что в записи $\Pi\alpha(\bang\quest\alpha\pto\varnothing)$ переменная
$\alpha$ связана знаком $\bang$, квантор $\Pi\alpha$ пуст, и
экземпляров с подстановкой у этой задачи нет.

По той же причине высказывание $\alpha\pto\bang\quest\beta$ не имеет
ничего общего с корефлексией: $\bang\quest\beta$ --- выводы замыкания
$\forall X_\beta\,\exists y\,(y\in X_\beta)$, то есть, по существу,
задача $\bang\bot$, и вся формула равносильна $\pneg\alpha$.}
(End of remark.)
\end{remark}

\begin{remark}[on the rejected semantics $\sS$]\label{rem:semS}
Alongside $\sB$ we also considered one more set-theoretic
semantics---the \emph{substitutional} one, $\sS$ (from the word ``substitute''),
in two versions. In the first, problem variables were interpreted, as
in $\sI$, by arithmetical formulas $\alpha^*(x)$, and the solver received
their code, while the quantifier $\Pi$ was read substitutionally: ``$\varphi^*$
holds when any arithmetical formula is substituted for $\alpha$''. This can be written as a formula of $\ZF$ only via
relativization of the quantifier over boxes: with nested quantifiers
the substitution instances of $\varphi^*$ themselves contain the truth
predicate, and by Tarski's theorem a direct formulation does not go through; therefore
the box $X_\alpha$ ranged over arithmetically definable sets, which
is expressed by a formula $\mathrm{Def}(X)$ of the language of $\ZF$ (via the truth
predicate for arithmetic definable in $\ZF$). We did not manage to prove anything essentially
new about this semantics compared with $\sI$ and $\sB$, while the combination of the code of the substituted formula with the binding
modality $\bang$ required provisos in the scheme $(\Pi\mathrm{E})$ and in
the generalization rule. In the second version the code was removed altogether: $\sS$ is
$\sB$ in which the boxes range only over arithmetically definable
sets (the quantifiers over $X_\alpha$ in the clauses for $\Pi$ and for $\bang$
are relativized by the formula $\mathrm{Def}$). Such a semantics hardly
differs from $\sB$: all positive results carry over by
relativizing the quantifiers, and all our refutations use
arithmetically definable boxes. The only point of possible divergence known to us
is predicativity: in $\sS$ the values of
$\Pi$-containing problems are merely definable in $\ZF$ and do not fall into the range of the quantifier,
so the scheme $(\Pi\mathrm{E})$ could be justified
only for $\Pi$-free $\theta$, and closure under
substitution only for $\Pi$-free substituted problems; in $\sB$
the quantifier is not predicative, and both provisos are superfluous. We chose not to pay for a
semantics that, as far as what has been proved goes, does not differ from $\sB$ with provisos in
all statements, and excluded it from the paper. Whether
restricting the boxes to definable sets changes anything remains
an open question (\S\ref{sec:todo}); an incidental argument in favour of
definable boxes is that the translations of the propositions $\quest\Pi\vec\alpha\,\varphi$
become arithmetical relative to $0^{(\omega)}$ rather than
$\Pi^1_1$-formulas.
(End of remark.)
\end{remark}

\begin{remark}\label{rem:clean}
Whenever possible, we shall try to satisfy
one syntactic restriction on formulas for this semantics:
\emph{no variable has both free and
bound occurrences in a formula, different binding operators bind different
variables, and the modality $\bang$ is applied only to propositions
without free problem variables}.
Such formulas are called \emph{clean}. Every formula is equivalent to a
clean one: the translations of the problems $\bang A$ and $\bang\Pi\vec\alpha\,A$ coincide
(Definition~\ref{def:starB}), and it suffices to make the quantifier explicit and
then rename the bound variables.
We do this to avoid misunderstandings. For example,
in the semantics $\sB$
the co-reflection $\alpha\pto\bang\quest\alpha$, which contains both a free and a bound
occurrence of the variable $\alpha$, has an entirely different meaning than in the arithmetical
semantics. In the semantics $\sB$ co-reflection is equivalent to the clean formula
 $\alpha\pto\bang\Pi\beta\,\quest\beta$, which in turn is equivalent to the
formula $\alpha\pto\varnothing$. The latter is by no means equivalent to co-reflection in the arithmetical semantics. The meaning of clean formulas, on the other hand, is roughly the same in all semantics.
\end{remark}

Lemma~\ref{lem:clauses} carries over to the semantics $\sB$ (its statements
are provable in $\ZF$); in addition, $\ZF$ proves
\[
  \val{\Pi\alpha\,\varphi}
   =\bigcap_{X\subseteq\omega}\val{\varphi}[X_\alpha:=X]
  .
\]

The parameters of the three semantics
are collected in the table:
\begin{center}\footnotesize
\begin{tabular}{l|l|l|l}
 & $\sI$ & $\sB$ & $\sE$\\
\hline
problem variables & arith.\ formulas & arbitrary boxes & arith.\ formulas\\
propositional variables & arith.\ sentences & $\ZF$-sentences & arith.\ sentences\\
operator $\Pi$ & no & $\forall X_\alpha$ & no\\
$\bang$ binds problem variables & no & yes & no\\
internal theory & $\PA$ & $\ZF$ & $\PA$\\
external theory & $\PA$ & $\ZF$ & $\TA$\\
\end{tabular}
\end{center}
The last two rows introduce two terms that we shall
use constantly. The \emph{internal theory} of a semantics is the one
whose proofs serve as solutions of the problems $\bang A$; provability in it
is what the modality $\square=\quest\bang$ expresses, and $\Prov$, $\Prf$ and $\Con$ without
further specification refer to it. The \emph{external theory} of a semantics is the one in
which validity is verified (\S\ref{sec:validity}): a formula
is considered established if it is provable in the external theory. For
$\sE$ the external theory is $\TA$---its axioms are all arithmetical sentences true in
$\N$, so that
``$\TA\vdash\sigma$'' simply means the truth of $\sigma$ in $\N$. In the
internal semantics $\sI$, $\sB$ the internal and the external theory
coincide (whence the name of these semantics); in $\sE$ the external theory
is strictly stronger than the internal one. In all three cases the internal theory
is contained in the external one. 

\section{Validity}\label{sec:validity}

\subsection{Definition}

\begin{definition}[code of an interpretation]\label{def:code}
Let a given formula $\varphi$ 
contain the variables $A_1,\dots,A_k$ of sort
proposition and $\alpha_1,\dots,\alpha_m$ of sort problem. The \emph{code}
of an arithmetical interpretation $*$ (relative to the formula $\varphi$) is the tuple 
\[
  \gn{*}_\varphi\ :=\ \bigl\langle\, A_1^*,\dots,A_k^*,\
        \alpha_1^*,\dots,\alpha_m^*\,\bigr\rangle .
\]
This is a finite object: the values of $*$ on variables not occurring in $\varphi$
are not included in it. For the semantics $\sB$ we omit the second part of the
tuple: there the problem variables are boxes, and they have no code.
\end{definition}

\begin{definition}[validity]\label{def:valid}
A proposition $A$ is \emph{valid} in a given semantics if for
every interpretation~$*$ the formula $A^*$ is provable in the external theory of this
semantics. 
A problem $\varphi$ is \emph{valid} in a given semantics if there exists
a total computable function $S$ such that for every
interpretation $*$ the value $S(\gn{*}_\varphi)$ is defined and the formula
$\varphi^*\bigl(S(\gn{*})\bigr)$ is provable in the external theory.\footnote{More precisely, instead of the value of the function one should
substitute into the formula the standard term denoting this value in the standard interpretation, but
we do not want to overload the notation.}
In the semantics $\sB$ the translations $A^*$ and
$\varphi^*\bigl(S(\gn{*})\bigr)$ may contain free
set variables $\vec X$---the boxes of the free problem
variables of the formula. By provability of such a formula in $\sB$ we mean
provability of its closure $\forall\vec X\,\varphi^*\bigl(S(\gn{*})\bigr)$.
\end{definition}

Validity in a semantics $Z$ is also called \emph{$Z$-validity};
validity in several semantics at once---for instance,
\emph{$\sI,\sB$-validity}.
For $\sE$ provability in $\TA$ is truth in $\N$, that is,
membership $S(\gn{*})\in\val{\varphi}$; for $\sI$ this
membership follows from provability, since $\PA$ is contained in
$\TA$. For $\sB$ values do not enter the definition: $\ZF$
has no canonical model. If one imagines some model
$M$ of $\ZF$ and the values in it (Definition~\ref{def:starB}), then
a valid problem is solved by the function $S$ in $M$ under every
interpretation and every contents of the boxes; but this is not part of the definition, only its informal reading.

\begin{remark}[metatheoretic assumptions]\label{rem:assumptions}
All proofs of \emph{non}-validity in the semantics $\sB$ use
the assumption that $\ZF$ is consistent, which is necessary,
since otherwise all formulas are valid.

More precisely, to prove that a problem $\phi$ is not valid in these
semantics, we must, for every total computable function
$S$, exhibit an interpretation $*$ for which
the formula $\phi^*(S(*))$ is not provable in $\ZF$.
Instead, we shall prove that the negation of this formula
is provable in $\ZF$. Under the assumption of consistency
this guarantees that the formula itself is not provable.
Similarly, instead of proving that
the formula $A^*$ (for a proposition $A$) is not provable, we shall prove that
$\aneg A^*$ is provable. We proceed in exactly the same way for the arithmetical semantics,
except that in this case no assumptions need to be made.
In a few places consistency is not enough, and we use the stronger
assumption of \emph{$\Sigma_1$-soundness} of the external theory:
every arithmetical $\Sigma_1$-sentence provable in it is
true; this is indicated explicitly each time. There are no other metatheoretic
assumptions in the paper: everywhere we speak either of the consistency
of the external theory or of its $\Sigma_1$-soundness, and for the external
theories $\PA$ and $\TA$ both properties hold.

Both assumptions are a consequence of the fact that our metatheory is not
fixed and is kept as weak as possible. In a sufficiently strong
metatheory, for instance when formalizing in the system Lean, they become
theorems: the standard model $\N$ and the model of $\mathrm{ZFC}$ from
the Mathlib library (the cumulative hierarchy up to an inaccessible cardinal
$\kappa$) are constructed there, so that the consistency and
$\Sigma_1$-soundness of $\PA$ and $\ZF$ are provable. At the same time an
analogue of the semantics $\sE$ for $\ZF$ appears: the external theory is truth in
this model $V_\kappa$, while the internal theory is still $\ZF$.
(End of remark.)
\end{remark}

\begin{remark}[another internal theory]\label{rem:otherth}
In the semantics $\sI$ Peano arithmetic plays both roles---that of the internal and
that of the external theory. Instead of it one can take any computably enumerable
$\Sigma_1$-sound theory $T$ containing $\PA$ (or
interpreting arithmetic), for instance $\ZF$: propositional variables
are interpreted by closed formulas of $T$, problem variables
by formulas of $T$ with one free variable ranging over the natural
numbers, solutions of the problem $\bang A$ are proofs of the formula $A^*$ in $T$, and
validity is provability in $T$. Denote this semantics by $\sI_T$,
so that $\sI=\sI_{\PA}$. Everything proved below about the semantics $\sI$
carries over to $\sI_T$ together with the proofs: they use
only $\Sigma_1$-completeness, the fixed-point lemma, the formalizability
of realizability and the Rosser sentence---properties common to all
such theories; for $T=\ZF$ the corresponding computations are in fact carried out below
in the treatment of the semantics $\sB$. Note that $\sI_{\ZF}$ is not the same
as $\sB$ without the quantifier: in $\sB$ the modality $\bang$ binds
problem variables, and, for instance, the proposition
$\square\quest\alpha\ato\square\bot$ is $\sB$-valid
(Remark~\ref{rem:boxqatom}), but not valid in $\sI_{\ZF}$ (take
$\alpha^*(x):=(x=x)$). Whether the sets of valid problems in
$\sI_{\PA}$ and $\sI_{\ZF}$ coincide, we do not know: transferring the proofs is
not enough for this, and there is no completeness on the problem side. (End of remark.)
\end{remark}

\begin{remark}[a proof of the formula is found by exhaustive search]\label{rem:whycomp}
Note that in the definition of $\sI,\sB$-validity
the proofs of the formulas $A^*$ and $\varphi^*\bigl(S(\gn{*})\bigr)$ can be found by exhaustive search.
For example, for the formula $A^*$: from
the code $\gn{*}$ one computes the formula $A^*$ itself (the translation $\cdot\mapsto\cdot^*$ of
Definition~\ref{def:tr} is a computable operation on codes),
the relation $\Prf$ is decidable, and by assumption the search terminates for
every $*$; hence the function $P\colon\gn{*}\mapsto{}$``the first proof found
of the formula $A^*$'' is total and computable.
\end{remark}


\begin{remark}[three possible notions and why the middle one is chosen]\label{rem:three}
In the definition of validity of a problem $\varphi$ one could require
\[
  \text{(a)}\ \ \forall *\ \exists e\ \ e\in\val{\varphi};
  \qquad
  \text{(b)}\ \ \exists\, S\ \forall *\ \ S(\gn{*})\in\val{\varphi};
  \qquad
  \text{(c)}\ \ \exists e\ \forall *\ \ e\in\val{\varphi},
\]
where in~(b) the function $S$ is computable. Obviously, (c) implies (b), and (b) implies~(a).
Neither of the two extreme notions is suitable.

\emph{Under definition \textup{(a)}}
the choice of a solution is not constrained in any way,
and the problem fragment becomes classical:
$\alpha\por\pneg\alpha$, Peirce's law and in general all classical tautologies are valid
(Proposition~\ref{prop:noem} below shows that under~(b) this
is no longer the case).

\emph{The notion \textup{(c)}} is too strong. Already the problem $\bang(B\ato B)$ will not
be valid: its solutions under an interpretation $*$ are the proofs of the formula
$B^*\ato B^*$, and there is no single number suitable for all $*$, since different
formulas have different proofs. Meanwhile, we would like all intuitionistically derivable formulas
to be valid.
Under~(b) everything is fine: a proof of the formula
$B^*\ato B^*$ is effectively constructed from $B^*$.
(End of remark.)
\end{remark}

\remove{
\begin{remark}[семантика интенсиональна]\label{rem:intens}
Функция $S$ получает \emph{синтаксис} интерпретации, а не сами множества:
множество решений --- бесконечный объект, подать его на вход нельзя. Поэтому
на двух кодах логически эквивалентных формул $S$ вправе вести себя
по-разному. Для модальности $\bang$, которая и так смотрит на формулу $A^*$, а
не на её истинностное значение, это естественно.
\end{remark}
}

\begin{remark}[the choice of numbering, pairing and coding does not matter]\label{rem:indep}
In Definitions~\ref{def:tr}--\ref{def:code} three arbitrary
choices were made: the G\"odel numbering $U_0,U_1,\dots$, the pairing function $c$ and the way
of coding an interpretation by $\gn{*}$. None of them affects which
formulas are valid---in any of the semantics. The justification is given in
the appendix.
\end{remark}

We conclude this section by proving that validity is preserved
under substitution of formulas for variables, as well as under application of the rule
Modus Ponens (MP).

\begin{lemma}[closure]\label{lem:closure}
In the arithmetical semantics the set of valid formulas
is closed under the rule \textup{(MP)} and under
substitution: of problems for free problem variables, of propositions
for propositional variables.
\end{lemma}

For the semantics $\sB$ the analogue of this lemma holds not for all substitutions. Here is a counterexample:
the proposition $\quest\Pi\alpha(\beta\por\pneg\beta)$
is valid in $\sB$: 
the theory $\ZF$ proves the formula $\exists x\,(\beta\por\pneg\beta)^*(x)$
by case analysis: if the value of $\beta$ is nonempty, say $\beta^*(u)$ holds,
then
$c(0,u)$ is a solution,
and otherwise $c(1,0)$ is. If
we substitute $\alpha$ for $\beta$
, we obtain the
proposition $\quest\Pi\alpha(\alpha\por\pneg\alpha)$, which $\ZF$
\emph{refutes}: a solution of the form $c(0,u)$ is rejected by the box
$\varnothing$, and a solution of the form $c(1,e)$ by the box $\{0\}$, since
the label $1$ promises that the value is empty.

For the semantics $\sB$ it will suffice for us to have this lemma for
the so-called \emph{safe substitutions}. By definition, these are substitutions
with the following two properties:
\begin{enumerate}[leftmargin=2.2em,itemsep=1pt,label=\textup{(\arabic*)}]
\item The free problem variables of the substituted formula
  must not fall under the binding quantifier $\Pi$. This phenomenon is called
\emph{variable capture};
the necessity of forbidding capture we have just demonstrated by the counterexample.
\item Similarly, the free problem variables of the substituted formula
  must not fall under the binding modality $\bang$. This
  restriction concerns only substitutions of propositions for
  propositional variables: when a problem is substituted for a problem
  variable, it is satisfied automatically, since free
  occurrences of a problem variable do not lie under the symbol $\bang$. An example
of a violation: $\bang A\mapsto \bang\quest \alpha$. We cannot demonstrate the necessity of this restriction by an example:
we do not know whether it is really necessary.
\end{enumerate}



\begin{lemma}[closure]\label{lem:closurebs}
In the semantics $\sB$ the set of valid formulas
is closed under the rule \textup{(MP)} and under
safe substitutions.
\end{lemma}

The proofs of the lemmas are given in the appendix.

\subsection{
The scheme (A4)}\label{sec:assval}

The only scheme about propositions that will become an axiom is the
group (A4), asserting that the modality $\quest$ is a Boolean
homomorphism; here we verify its validity.
\remove{Полное описание
общезначимых $\Pi$-свободных высказываний --- с точностью до
нормальной формы это логика доказуемости $\GL$ Соловея для внутренних
семантик $\sI$, $\sB$ и её внешний вариант $\GLS$ для $\sE$ --- мы
получим в \S\ref{sec:asscompl}, после теоремы о корректности: там оно
достаётся почти даром из того же описания \emph{выводимых}
высказываний. Для арифметических семантик $\sI$ и $\sE$, в языке
которых квантора $\Pi$ нет вовсе, это описание исчерпывает \emph{все}
высказывания; в полном же языке семантики $\sB$ за его пределами
остаются высказывания с квантором $\Pi$, и логика этого
фрагмента --- открытый вопрос (\S\ref{sec:piatoms}, \S\ref{sec:todo}).
}
In this and the two following subsections all statements are formulated and
proved for all three semantics at once.
\remove{
; специфика семантики $\sB$
--- оператор $\Pi$ --- разобрана
отдельно в \S\ref{sec:piatoms}, поэтому здесь рассматриваются
$\Pi$-свободные формулы. Соглашения: внутренняя и внешняя теории каждой
семантики --- как в таблице \S\ref{sec:valBS} (внутренняя: $\PA$ для
$\sI$, $\sE$ и $\ZF$ для $\sB$; внешняя отличается от неё только
в $\sE$, где это $\TA$); $\Con$ обозначает непротиворечивость внутренней
теории. Кодом интерпретации в $\sB$ служит код прозрачной
части (определение~\ref{def:code}); переводы в этой семантике могут
содержать свободные множественные переменные $\vec X$, и выводимость
таких формул понимается как выводимость замыкания
(определение~\ref{def:valid}).
Леммы~\ref{lem:clauses} и~\ref{lem:dn}
используются во всех семантиках: их доказательства схемны и одинаково
проходят в $\PA$ и в $\ZF$, в том числе для переводов со свободными
множественными переменными. Отрицательные утверждения для $\sB$,
как всегда, доказываются в допущениях
замечания~\ref{rem:assumptions}.
}

\begin{lemma}[$\quest$ is a Boolean homomorphism at the level of translations]\label{lem:qsem}
For any problems $\varphi,\psi$ the following equivalences are provable in
the internal theory (and hence in the external one),
and the proofs are computably constructed from the code of the interpretation:
\[
\begin{aligned}
  (\quest(\varphi\pto\psi))^* &\aeq\bigl((\quest\varphi)^*\ato(\quest\psi)^*\bigr),&
  (\quest(\varphi\pand\psi))^* &\aeq\bigl((\quest\varphi)^*\aand(\quest\psi)^*\bigr),\\
  (\quest(\varphi\por\psi))^* &\aeq\bigl((\quest\varphi)^*\aor(\quest\psi)^*\bigr),&
  (\quest\varnothing)^* &\aeq\arbot .
\end{aligned}
\]
\end{lemma}

\begin{proof}
For conjunction: the formula $(\quest(\varphi\pand\psi))^*$ is $\exists x\,
\exists u\exists v\,(x=c(u,v)\aand\varphi^*(u)\aand\psi^*(v))$, and the internal
theory proves its equivalence to the conjunction $\exists u\,\varphi^*(u)\aand
\exists v\,\psi^*(v)$: from left to right by the projections of the pair, from right to left
by forming the pair. The same for disjunction, with a case distinction on which of the labels
$c(0,u)$, $c(1,v)$ is present. For $\varnothing$: $(\quest\varnothing)^*$ is
$\exists x\,\arbot $.

For implication, from left to right the internal theory proves $\exists
x\,(\varphi\pto\psi)^*(x)\ato(\exists y\,\varphi^*(y)\ato\exists z\,\psi^*(z))$---this
is pure logic. From right to left some work is needed: an index must be exhibited.
The internal theory proves $(\exists u\,\varphi^*(u)\ato\exists v\,\psi^*(v))\ato\exists
x\,(\varphi\pto\psi)^*(x)$ by distinguishing two cases. If $\aneg\exists
u\,\varphi^*(u)$, then the condition in $(\varphi\pto\psi)^*(x)$ is vacuous and any
$x$ will do. If, on the other hand, $\exists v\,\psi^*(v)$, then, taking the least such $v$,
$x:=k(v)$ will do, where $k$ is a total computable function that yields, given $v$, an index
of the constant function $u\mapsto v$ (the parametrization theorem). The case distinction
here is classical---this is the \emph{residual classicality} of our semantics,
to which we shall return (Remark~\ref{rem:qhom}). All computations
are the same in $\PA$ and $\ZF$; in $\sB$ they are carried out on top of the free
set variables $X_\alpha$, which take part neither in the case distinction nor in
the formation of pairs.
\end{proof}

\subsection{Modality-free problems: the schemes of $\Int$}\label{sec:intval}

We turn to problems and begin with the simplest class of them---problems without
modalities, that is, built from variables and the constant $\varnothing$
by the connectives $\pto$, $\pand$, $\por$. The future axioms (A2)---the schemes
of the intuitionistic calculus $\Int$---are valid.

\begin{proposition}[the schemes of $\Int$ are valid]\label{prop:intval}
Every problem derivable in the intuitionistic calculus is valid in
all three semantics, moreover with a constant solving function.
\end{proposition}

\begin{proof}
The Kleene realizers of the axiom schemes of $\Int$ are built from indices supplied
by the parametrization theorem and do not depend on the interpretation; the correctness
of each realizer is provable in the internal theory, and the proof is constructed from the code:
the computations are the same in $\PA$ and $\ZF$, and in $\sB$ they are carried out on top of
free set variables. The rule (MP) preserves
validity by Lemmas~\ref{lem:closure} and~\ref{lem:closurebs}.
\end{proof}

Everything else that we know about valid modality-free problems---in
particular, that there are no classical schemes among them, and that they are not
exhausted by the intuitionistic ones---is more convenient to prove after the soundness
theorem; this is done in~\S\ref{sec:modalfreeval}.

\subsection{Valid problems with modalities}\label{sec:taskval}

About the set of all valid problems---with modalities---we know much
less. It is (in each of the three semantics) closed under (MP) and
substitution---safe substitution in $\sB$ (Lemmas~\ref{lem:closure}
and~\ref{lem:closurebs}). Here we verify the validity of
those schemes with modalities that will become axioms of the calculi of
\S\ref{sec:calc} for the semantics $\sI$ and $\sB$: the scheme (A3) and the schemes (A5) and
(A7) about problems of a special form. The axioms of the calculus for $\sE$---reflection
(A8) and its analogue for problems (A9)---will be verified when it is
defined. 

\begin{proposition}[the scheme (A3)]\label{prop:a3}
For any propositions $A$, $B$ the problem
$\bang(A\ato B)\pto(\bang A\pto\bang B)$ is valid in all three
semantics.
\end{proposition}

\begin{proof}
In $\sI$ and $\sE$: the function that, given proofs of the formulas $A^*\ato B^*$ and
$A^*$ in the internal theory, constructs a proof of the formula $B^*$ by the rule modus ponens
is computable; the parametrization theorem yields a constant
solving function $S$, and its correctness is provable in the external theory.
In $\sB$ the same, with a correction for binding: solutions of the problems
$\bang(A\ato B)$, $\bang A$, $\bang B$ are proofs of the closures of the formulas
$A^*\ato B^*$, $A^*$, $B^*$ with respect to the free variables of the corresponding
propositions. From a proof of
the formula $\forall\vec X\,(A^*\ato B^*)$ and a proof of the formula
$\forall\vec X_A\,A^*$ one constructs a proof of the formula $\forall\vec X_B\,B^*$:
the variables free in $A$ but not in $B$ are instantiated by the set
$\varnothing$, after which only modus ponens and
generalization remain. The transformation is computable, the solver is constant,
and its correctness is provable in $\ZF$.
\end{proof}

\subsubsection{$\Sigma$-problems}

We single out a class of problems defined grammatically; it will appear in the axioms
(A5) and (A7) of the calculi (\S\ref{sec:calc}).

\begin{definition}[$\Sigma$-problems and $\Sigma$-propositions]\label{def:sigma}
\emph{$\Sigma$-problems} are the problems built by the grammar
\[
  \varphi ::= \varnothing \mid \bang A \mid (\varphi\pand\varphi) \mid
             (\varphi\por\varphi) ,
\]
where $A$ is an arbitrary proposition. In other words, a $\Sigma$-problem
is assembled from problems of the form $\bang A$ and the constant $\varnothing$ by means of
$\pand$ and $\por$ alone: variables of sort problem, implication and, hence, negation
do not enter here. For example, $\bang A$, $\bang A\por\bang B$ and
$(\bang A\pand\bang B)\por\bang C$ are $\Sigma$-problems, whereas $\alpha$, $\pneg\bang A$ and
$\bang A\pto\bang B$ are not. In the semantics $\sB$ a $\Sigma$-problem
has no free problem variables: all its problem variables occur
under $\bang$ and are bound by this symbol.

\emph{$\Sigma$-propositions} are the propositions built by the
grammar
\[
  A ::= \bot \mid \quest\phi \mid (A\aand A) \mid (A\aor A) ,
\]
where $\phi$ ranges over $\Sigma$-problems. Again there are no propositional variables,
no implication and no negation here. In our calculus
every such
proposition will be equivalent
to a combination of $\bot$ and formulas $\square
B$ by means of $\aand$ and $\aor$: $\quest\varnothing\aeq\bot$,
$\quest\bang B\aeq\square B$, and $\quest$ passes through $\pand$ and $\por$.
In particular, $\square B$ is a $\Sigma$-proposition for any $B$.
\end{definition}

We verify two schemes about $\Sigma$-problems, which will become the axioms (A5) and (A7).

\begin{lemma}\label{lem:sigmaenum}
There is a total computable function which, given a $\Sigma$-problem
$\phi$ and a code $\gn{*}$, outputs an index of a program enumerating the set
$\val{\phi}$. Moreover, the external theory proves that this program
enumerates exactly the set defined by the formula $\phi^*$.
\end{lemma}

\begin{proof}
Induction on the construction of $\phi$, relying on Lemma~\ref{lem:clauses}.
For $\phi=\varnothing$---the empty enumeration. For $\phi=\bang A$ the set
$\val{\bang A}$ consists of the proofs of the formula $A^*$ (in $\sB$, of its
closure), and ``$n$ is a proof of the given formula'' is a decidable relation,
so an exhaustive search over all $n$ with
checking will do. For $\phi=\psi\pand\chi$ we enumerate the pairs $c(u,v)$ from
enumerations of $\val{\psi}$ and $\val{\chi}$; for $\phi=\psi\por\chi$,
the elements $c(0,u)$ and $c(1,v)$, alternating the two enumerations. In $\sB$
the translation $\phi^*$ contains no free set variables, and the argument
is the same.
\end{proof}

\begin{proposition}[the schemes (A5) and (A7)]\label{prop:a5a7}
For every $\Sigma$-problem $\phi$ the problems
\[
  \text{(A5)}\ \ \phi\pto\bang\quest\phi
  \quad\text{(co-reflection)}\qquad\text{and}\qquad
  \text{(A7)}\ \ \pneg\pneg\phi\pto\phi
  \quad\text{(double negation elimination)}.
\]
are valid in all three semantics.
\end{proposition}

\begin{proof}
The translation $\phi^*$ is a $\Sigma_1$-formula: the clause for $\bang A$ yields
the decidable relation $\Prf$, and the clauses for $\pand$ and $\por$ add
only $\aand$, $\aor$ and $\exists$.

(A5). A solution of the problem $(\bang\quest\phi)^*$ is a proof of the formula
$\exists x\,\phi^*(x)$ in the internal theory. Let $e$ be an index of the
program which, given anything as input, runs through all numbers and
outputs the first number found that is a proof of the formula
$\exists x\,\phi^*(x)$ in the internal theory; put $S(\gn{*}):=e$.
If $\val{\phi}=\varnothing$, there is nothing to check. If, on the other hand,
$u\in\val{\phi}$, then $\exists x\,\phi^*(x)$ is a true
$\Sigma_1$-sentence; by the $\Sigma_1$-completeness of the internal theory it is
provable, the search halts, and $U_e(u)$ is a proof of it.
Provability in the external theory: it proves
$\phi^*(u)\ato\Prov(\exists x\,\phi^*(x))$ (formalized
$\Sigma_1$-completeness), and given $\Prov(\exists x\,\phi^*(x))$ the search
halts at the first proof found (formalization
of unbounded search, as below in (A7)). Thus the external theory
proves
$\forall u\,\bigl(\phi^*(u)\ato\Prf(U_e(u),\exists x\,\phi^*(x))\bigr)$.

(A7). By Lemma~\ref{lem:sigmaenum}, from $\gn{*}$ one constructs an algorithm
enumerating the set $\val{\phi}$. Let $e$ be an index of the program
which, given anything as input, launches this enumeration and outputs
the first element found; put $S(\gn{*}):=e$. If
$\val{\pneg\pneg\phi}=\varnothing$, there is nothing to check. If, on the other hand,
$\val{\pneg\pneg\phi}=\N$, then by Lemma~\ref{lem:dn} the set
$\val{\phi}$ is nonempty, the search halts, and
$U_e(u)\in\val{\phi}$ for all $u$. Provability in the external theory: by
Lemma~\ref{lem:dn} the formula $(\pneg\pneg\phi)^*(u)$ is provably
equivalent to $\exists v\,\phi^*(v)$, and by Lemma~\ref{lem:sigmaenum}
the external theory proves that the enumeration corresponds to the formula $\phi^*$;
hence it also proves that, given $\exists v\,\phi^*(v)$, the search
halts at the first witness found (the usual formalization
of bounded search). Thus the external theory proves
$\forall u\,\bigl((\pneg\pneg\phi)^*(u)\ato\phi^*(U_e(u))\bigr)$, and
the proof is constructed from $\gn{*}$.

In all cases the computations are the same in $\PA$ and $\ZF$; in $\sB$
the translations of $\Sigma$-problems contain no free set variables.
\end{proof}


\subsection{The semantics $\sB$: the operator $\Pi$}
\label{sec:piatoms}

All facts of Subsections \ref{sec:assval}--\ref{sec:taskval} have already
been proved for all three semantics at once; what is new in this
subsection is the operator $\Pi$ itself. Let us recall the motivation from the introduction. Under the symbol $\quest$
the constructive content of a problem is lost: the proposition
$\quest(\alpha\por\pneg\alpha)$ is valid trivially ($\PA$ is a
classical theory and proves $\exists
x\,(\alpha\por\pneg\alpha)^*(x)$ by a case analysis for every
fixed interpretation), and just as trivially valid is
$\quest\rho$ for Rose's formula (\S\ref{sec:modalfreeval}), although the
problem $\rho$ itself is solved by a nontrivial algorithm. The quantifier $\Pi$ restores what has been lost: later
we shall prove (\S\ref{sec:pigain}) that the proposition
$\quest\Pi\alpha(\alpha\por\pneg\alpha)$ is refutable---its negation is
valid---while $\quest\Pi\beta\Pi\gamma\rho$ is a nontrivial
valid proposition. Here we verify the validity of the laws of the operator itself;
they will become axioms of the calculus for $\sB$.

\begin{proposition}[laws of $\Pi$]\label{prop:pilaws}
The following formulas are valid in $\sB$ (in parentheses, the realizer of the problems;
all realizers are constant).
\begin{enumerate}[leftmargin=1.8em,itemsep=1pt]
\item $(\Pi\mathrm{E})$:
  $\Pi\alpha\,\varphi\pto\varphi[\theta/\alpha]$ and
  $\Pi\alpha\,A\ato A[\theta/\alpha]$, if $\theta$
  is free for $\alpha$ \ (for the problem the realizer is the identity: the translation of
  $\varphi[\theta/\alpha]$ is obtained from $\varphi^*$ by substituting
  the set $\{x:\theta^*(x)\}$ for $X_\alpha$---the symbols $\bang$
  are not affected by the substitution, since no free occurrences of $\alpha$ lie under
  them---and the intersection over all boxes is contained in the value
  at this box; for the proposition, $\ZF$ proves the translation by the same
  instantiation of the quantifier $\forall X_\alpha$);
\item $(\Pi\por)$:
  $\Pi\alpha(\varphi\por\psi)\pto\Pi\alpha\,\varphi\por\Pi\alpha\,\psi$
  \ (the identity; it is essential that the tag of the solution is fixed: a solution
  of the premise is a pair consisting of a tag and a solution $e$, and the same tag with the same $e$
  works for every box);
\item $(\Pi\bang)$: $\bang A\peq\bang\Pi\alpha\,A$ \ (the realizer is
  a transformation of proofs: solutions of the premise are proofs of the closure
  $\forall\vec X\,A^*$, solutions of the conclusion are proofs of the closure
  $\forall\vec X'\,\forall X_\alpha\,A^*$, where $\vec X'$ are the boxes
  of the remaining free problem variables of $A$; these formulas differ
  only in the order of the quantifiers (and, when $\alpha\notin FV(A)$, by a superfluous
  quantifier), and if $\alpha$ is the last in $\vec X$, they coincide;
  permutation of quantifiers in a proof is a computable
  transformation whose correctness is provable in $\ZF$).
\end{enumerate}
All arguments are carried out in $\ZF$ over free set
variables.
\end{proposition}

Besides these axioms, the calculus for $\sB$ will need one rule---%
\emph{Bernays' rule}: from $\varphi\pto\psi$ infer
$\varphi\pto\Pi\alpha\,\psi$ provided $\alpha\notin FV(\varphi)$ (and likewise
for propositions). It preserves validity---this is verified
together with the other rules in Theorem~\ref{th:sound}. The remaining
customary laws of the quantifier---generalization, $\Pi\alpha\,\varphi\peq\varphi$
for $\alpha\notin FV(\varphi)$, distributivity over
implication and conjunction, the converse of $(\Pi\por)$---will be
derived from $(\Pi\mathrm{E})$ and Bernays' rule
(Lemma~\ref{lem:piderived}).

\remove{
Лемма~\ref{lem:qsem} верна и для задач с квантором $\Pi$: её
доказательство не заглядывает внутрь задач $\varphi$ и $\psi$, а
классический разбор случаев проводится в $\ZF$ поверх свободных
множественных переменных. Протолкнуть $\quest$ сквозь сам квантор
$\Pi$ нельзя (замечание после леммы~\ref{lem:semnf}). Исчислениям для
$\sB$ понадобится ещё одно правило --- обобщение по задачной
переменной; его корректность, как и корректность остальных правил, ---
в \S\ref{sec:calc}.
}

\subsection{Summary: the table of validity}\label{sec:tables}

Table~\ref{tab:axioms} summarizes the validity of all schemes which in
\S\ref{sec:calc} will become axioms of the calculi (the schemes (A8) and (A9)
of the calculus for $\sE$ are treated in~\S\ref{sec:calcE}). Notation: $+$ means
valid; $-$ means not valid (under the assumptions of
Remark~\ref{rem:assumptions}); $\cdot$ means not expressible in the given
semantics (no $\Pi$).
The last
column indicates where this is proved; the validity of the scheme (A6) will be
obtained in \S\ref{sec:calc} by a derivation in the calculus (Theorem~\ref{th:sound}).

\begin{table}[ht]\centering\small
\caption{Axioms of the calculi and their validity in the three semantics}
\label{tab:axioms}
\medskip
\begin{tabular}{l l c c c l}
 & scheme & $\sI$ & $\sB$ & $\sE$ & where proved\\
\hline
(A1) & classical tautologies & $+$ & $+$ & $+$ & by definition\\
(A2) & schemes of $\Int$ & $+$ & $+$ & $+$ & Prop.~\ref{prop:intval}\\
(A3) & $\bang(A\ato B)\pto(\bang A\pto\bang B)$ & $+$ & $+$ & $+$ & Prop.~\ref{prop:a3}\\
(A4) & $\quest$ is a Boolean homomorphism & $+$ & $+$ & $+$ & Lemma~\ref{lem:qsem}\\
(A5) & $\phi\pto\bang\quest\phi$, $\phi\in\Sigma$ & $+$ & $+$ & $+$ & Prop.~\ref{prop:a5a7}\\
(A6) & $\bang(\square A\ato A)\pto\bang A$ & $+$ & $+$ & $+$ & Theorem~\ref{th:sound}\\
(A7) & $\pneg\pneg\phi\pto\phi$, $\phi\in\Sigma$
  & $+$ & $+$ & $+$ & Prop.~\ref{prop:a5a7}\\
(A8) & reflection $\square A\ato A$ & $-$ & $-$ & $+$ & Prop.~\ref{prop:hce}(1),(2)\\
(A9) & $\bang\quest\varphi\pto\pneg\pneg\varphi$
  & $-$ & $-$ & $+$ & Prop.~\ref{prop:hce}(1),(2)\\
($\Pi$E) & $\Pi\alpha\,\varphi\pto\varphi[\theta/\alpha]$, both sorts
  & $\cdot$ & $+$ & $\cdot$ & Prop.~\ref{prop:pilaws}\\
($\Pi\por$) & $\Pi\alpha(\varphi\por\psi)\pto\Pi\alpha\,\varphi\por\Pi\alpha\,\psi$
  & $\cdot$ & $+$ & $\cdot$ & Prop.~\ref{prop:pilaws}\\
($\Pi\bang$) & $\bang A\peq\bang\Pi\alpha\,A$
  & $\cdot$ & $+$ & $\cdot$ & Prop.~\ref{prop:pilaws}\\
\hline
\end{tabular}

\end{table}

\remove{
\begin{remark}[что видно из таблицы]\label{rem:read}
Семантики $\sI$, $\sB$ задают \emph{одно и то же} исчисление ---
с двумя оговорками: в $\sI$ $\Pi$-строки отсутствуют в языке, а в
$\sS$ схема ($\Pi$E) ограничена $\Pi$-свободными $\theta$. Семантика
$\sE$ отличается ровно тем, что теряет правило (AtoP) --- оно не
сохраняет $\sE$-общезначимость (теорема~\ref{prop:not}(2)) --- и
приобретает две аксиомы: рефлексию (A8) и рефлексию для задач (A9).
В таблице нет отдельных строк для $\square$-рефлексии
$\square A\ato A$ (это та же (A8), поскольку $\square=\square$),
для формулы Роуза (частный случай (A8)) и для задачи $\pneg\bang\bot$
(частный случай (A9), теорема~\ref{prop:not}(7)).
\end{remark}
}

\section{Calculi}\label{sec:calc}

For each of the three
semantics we define a calculus sound with respect to it:
$\HCI$, $\HCB$ and $\HCE$. The axioms are exactly the schemes of
Table~\ref{tab:axioms} marked with a plus in the corresponding column.
The validity of all axioms of the calculi for $\sI$ and $\sB$, except (A6), has been proved in
\S\ref{sec:validity}, that of the axioms (A8) and (A9) of the calculus for $\sE$ is proved in~\S\ref{sec:calcE};
the scheme (A6) we shall not verify directly---it will be derived from
the other axioms and the fixed-point lemma (Theorem~\ref{th:sound}).

\subsection{The calculus for the semantics $\sI$}\label{sec:calcI}

\begin{definition}[the system $\HCI$]\label{def:hck}
Axiom schemes:
\[
\begin{array}{ll}
  \text{(A1)} & \text{classical tautologies (sort proposition);}\\[2pt]
  \text{(A2)} & \text{axiom schemes of intuitionistic logic (sort problem);}\\[2pt]
  \text{(A3)} & \bang(A\ato B)\pto(\bang A\pto\bang B);\\[2pt]
  \text{(A4)} & \text{$\quest$ is a Boolean homomorphism:}\\[2pt]
   & \qquad\text{(A4a)}\ \ \quest(\varphi\pto\psi)\aeq(\quest\varphi\ato\quest\psi),\\[2pt]
   & \qquad\text{(A4b)}\ \ \quest(\varphi\pand\psi)\aeq(\quest\varphi\aand\quest\psi),\\[2pt]
   & \qquad\text{(A4c)}\ \ \quest(\varphi\por\psi)\aeq(\quest\varphi\aor\quest\psi),\\[2pt]
   & \qquad\text{(A4d)}\ \ \quest\varnothing\aeq\bot;\\[2pt]
  \text{(A5)} & \phi\pto\bang\quest\phi\quad\text{for $\Sigma$-problems $\phi$ (co-reflection on the class $\Sigma$);}\\[2pt]
  \text{(A6)}& \bang(\square A\ato A)\pto\bang A\ \ \text{(a constructive version of L\"ob's theorem);}\\[2pt]
  \text{(A7)} & \pneg\pneg\varphi\pto\varphi\quad\text{for $\Sigma$-problems $\varphi$.}
\end{array}
\]
In all schemes $A$, $B$ are arbitrary
propositions, $\varphi$, $\psi$ are arbitrary problems (and in (A7),
arbitrary problems of the indicated class); a scheme denotes the set of all its
instances. All formulas here are $\Pi$-free. Inference rules:
\[
\begin{array}{lll}
  \text{(MP)} & \dfrac{A,\ A\ato B}{B}\ \ \text{and}\ \
                \dfrac{\varphi,\ \varphi\pto\psi}{\psi}
              & \text{(for both sorts);}\\[12pt]
  \text{(AtoP)} & \dfrac{A}{\bang A}
              & \text{from a proposition, a problem;}\\[12pt]
  \text{(PtoA)} & \dfrac{\varphi}{\quest\varphi}
              & \text{from a problem, a proposition.}
\end{array}
\]
The names of the rules read as ``assertion to problem'' and ``problem to
assertion'': (AtoP) produces from a proposition $A$ the problem $\bang A$ ``find
a proof of $A$'', and (PtoA) produces from a problem $\varphi$ the proposition
$\quest\varphi$ ``$\varphi$ has a solution''.

The special case $\sigma:=\bang A$ of the scheme (A5) is the problem
$\bang A\pto\bang\square A$: from a proof of the formula $A$ one can construct
a proof of the formula $\square A$.

By $\HCglm$ we denote the same calculus \emph{without} the schemes (A6) and (A7). It
will be needed because we shall not verify L\"ob's scheme directly: it will be
derived in $\HCglm$ from the fixed-point lemma. By $\HCzero$
(the \emph{core}) we denote a still smaller system: the schemes (A1), (A2), (A3),
the \emph{left} half of (A4a), (A4b)--(A4d), (A5) and all three inference rules. The core is
exactly that part of our axiomatics which is contained in Melikhov's logic
$\HC$. 
\end{definition}

\subsection{The calculus for the semantics $\sB$}\label{sec:picalc}

\begin{definition}[the system $\HCB$]\label{def:hczf}
Formulas are the formulas of the full language with the quantifier $\Pi$
(\S\ref{sec:valBS}).
The axioms of the system $\HCB$ are the schemes (A1)--(A7) and
the axioms of the operator $\Pi$:
\[
  (\Pi\mathrm{E})\quad\text{(for both sorts)},\qquad
  (\Pi\por),\qquad (\Pi\bang)
\]
(for the formulations see Proposition~\ref{prop:pilaws}; in
$(\Pi\mathrm{E})$ it is understood that $\theta$ is free for
$\alpha$). Inference rules: (MP), (AtoP), (PtoA) and Bernays' rule
for both sorts
\[
\begin{array}{lll}
  \text{(Bern)} & \dfrac{\varphi\pto\psi}{\varphi\pto\Pi\alpha\,\psi}
    \ \ \text{and}\ \ \dfrac{A\ato B}{A\ato\Pi\alpha\,B}
  & \text{(for both sorts),}
\end{array}
\]
where $\alpha$ does not occur free in the premise $\varphi$ (respectively $A$).
\end{definition}

Thus the quantifier-free part of the axiomatics of $\HCB$ coincides with
$\HCI$; what the quantifier layer gives is discussed in
\S\ref{sec:pigain}.

\subsection{The calculus for the semantics $\sE$}\label{sec:calcE}

\begin{definition}[the system $\HCE$]\label{def:hcks}
The axioms of the system $\HCgls$ are all theorems of $\HCgl$, together with the schemes
\[
  \text{(A8)}\ \ \square A\ato A\ \ \text{(reflection)},
  \qquad
  \text{(A9)}\ \ \bang\quest\varphi\pto\pneg\pneg\varphi\ \
  \text{(reflection for problems)}
\]
($A$ is an arbitrary proposition, $\varphi$ is an arbitrary problem).
The inference rules are (MP) and (PtoA); there is \textbf{no} rule (AtoP).
\end{definition}

\begin{proposition}[the schemes (A8) and (A9) and the semantics]\label{prop:hce}
Let $A$ be a variable of sort proposition and $\alpha$ a variable of sort
problem.
\begin{enumerate}[leftmargin=2.2em,itemsep=2pt,label=\textup{(\arabic*)}]
\item The schemes (A8) and (A9) are $\sE$-valid.
\item The schemes (A8) and (A9) are neither $\sI$- nor $\sB$-valid: already their
  instances $\square\bot\ato\bot$ and
  $\bang\quest\varnothing\pto\pneg\pneg\varnothing$ are not valid.
\item The problem $\bang(\square\bot\ato\bot)$ is not valid in any of the three
  semantics.
\item Neither the problem $\bang\quest\alpha\pto\alpha$ nor the converse
  co-reflection $\alpha\pto\bang\quest\alpha$ is valid in any of the
  three semantics.
\end{enumerate}
\end{proposition}

\begin{proof}
(1) The scheme (A8). In $\sE$: if $\val{\quest\bang A}=\Tru$, then $\PA\vdash A^*$, and the formula
$A^*$ is true, since the internal theory $\PA$ is contained in the external
theory $\TA$.

The scheme (A9). In $\sE$ the function $S(\gn{*}):=e$ works, where $e$ is the number of the constant function
$u\mapsto0$; it is the same for all $*$, so the algorithm makes no case
distinction at all. Let us check that this number is suitable. If
$\val{\bang\quest\varphi}=\varnothing$, then
$\val{\bang\quest\varphi\pto\pneg\pneg\varphi}=\N$ and any number is suitable.
If, on the other hand, $\val{\bang\quest\varphi}\neq\varnothing$, then
$\PA\vdash\exists x\,\varphi^*(x)$, hence $\val{\varphi}\neq\varnothing$, and by Lemma~\ref{lem:dn}
$\val{\pneg\pneg\varphi}=\N$; therefore $U_e$ lands in
$\val{\pneg\pneg\varphi}$ on any input.

(2) The scheme (A8). 
The translation of the reflection formula $\square\bot\ato\bot$ is $\Con$, and
the consistency of the internal theory is not provable in that theory itself (which is
also the external one)---G\"odel's second theorem.

The scheme (A9). In the internal semantics already the formula
$\bang\quest\varnothing\pto\pneg\pneg\varnothing$ is not valid. 
Indeed,
the problem $\pneg\pneg\varnothing$ is provably
equivalent to $\varnothing$, and $(\quest\varnothing)^*$ is provably equivalent to
the formula $\arbot $ (Lemma~\ref{lem:qsem}), so the translation of the problem
$\bang\quest\varnothing\pto\pneg\pneg\varnothing$ is provably equivalent to
the translation of the problem $\pneg\bang\bot$, that is, to the formula
$\forall v\,(\Prf(v,\arbot )\ato \arbot )$, which is provably equivalent to $\Con$;
and the consistency of the internal theory is not provable in that theory itself (which is
also the external one)---G\"odel's second theorem.

(3) Suppose that this problem is valid, and denote
its solver by $S$.
Applying it to the formula $\arbot $, we obtain a natural number $n$
for which the external theory proves
the formula $\Prf(n, \square\bot\ato\bot)$. This is a $\Sigma_0$-sentence;
were it false, its negation would be provable ($\Sigma_1$-completeness),
contradicting the consistency of the external theory. Hence it is true: $n$ is a proof
of the formula $\square\bot\ato\bot$ in the internal theory, that is, the internal
theory proves its own consistency---contradicting G\"odel's second
theorem.

(4) For the arithmetical semantics:
fix a computably enumerable undecidable set
$X\subset\N$ and express the property $n\in X$ by an arithmetical
sigma formula $\exists x A(x,n)$.
Put
$\alpha^{*_n}(x):=(x=0\aand\aneg A(x,n))\aor(x=1\aand A(x,n)))$.
Then $\PA$ proves
$\exists x\,\alpha^{*_n}(x)$ for every natural number $n$, so that
$\val[*_n]{\bang\quest\alpha}\neq\varnothing$, while $\val[*_n]{\alpha}$ is $\{0\}$ or
$\{1\}$ according to whether $n$ belongs to $X$; a solution would yield
a decision algorithm for $X$.

This formula is not
clean, so in the set-theoretic semantics the proof is quite different.
The clean formula $\bang\quest\beta\pto\alpha$ equivalent to it is equivalent to $\bang\bot\pto\alpha$ and is not valid: its instance at
$\alpha:=\varnothing$ is the problem $\pneg\bang\bot$, whose translation is
provably equivalent to $\Con$, so that, assuming the consistency of $\ZF$,
it has no solutions.

Co-reflection. For the arithmetical semantics,
suppose that there is an algorithm $S$ which provably produces,
from $*$, a solution of the problem $(\alpha\pto\bang\quest\alpha)^*$.
Let $\alpha^*$ be $\Con$ (with no free variables), that is, the arithmetical formula
expressing the consistency of the internal theory.
Put $n=S(*)$.
Then the external theory proves the formula
$$
(\alpha\pto\bang\quest\alpha)^*(n)=
\forall x(\Con\ato \Prf(U_n(x),\exists y\,\Con)),
$$
which provably implies
$$
\Con\ato  \exists z \Prf(z,\Con).
$$
Meanwhile,
by G\"odel's second theorem
the internal (and hence also the external) theory
proves
 $$
\Con\ato \aneg \exists z \Prf(z,\Con).
$$
Consequently, the external theory proves the false formula
$\aneg\Con$.

For the semantics $\sB$, first note that the variable $\alpha$ in the conclusion of the formula $\alpha\pto\bang\quest\alpha$
is bound, so it means something quite different from what the variable
$\alpha$ in the premise means.
%
The clean formula $\alpha\pto\bang\quest\beta$ equivalent to it is equivalent to
the formula $\alpha\pto\bang\bot$ (Remark~\ref{rem:bangbinds}), and the latter is not
valid: its instance at $\alpha:=\pneg\varnothing$ is equivalent to
the problem $\bang\bot$, which, assuming the consistency of $\ZF$, has no solutions.
\end{proof}

Now it is clear why $\HCE$ has no rule (AtoP): the proposition
$\square\bot\ato\bot$ is $\sE$-valid (a special case of (A8), item~(1)), while
applying the rule (AtoP) to it yields the problem $\bang(\square\bot\ato\bot)$,
which is not valid in any semantics (item~(3)).

\begin{remark}[why (A9), and why it cannot be strengthened]\label{rem:a10}
The scheme (A9) is exactly reflection transferred to the sort of problems:
the proposition $\quest\varphi$ and the problem $\pneg\pneg\varphi$ have, as
Lemma~\ref{lem:dn} shows, one and the same meaning ``$\varphi$ has a solution'', and (A9)
says: if it is proved that a solution exists, then it exists. It cannot be strengthened to
$\bang\quest\varphi\pto\varphi$: already for a variable this problem is not
valid (Proposition~\ref{prop:hce}(4)). The double negation on the right-hand side is
exactly the weakening that erases this information: by
Lemma~\ref{lem:clauses}(6) the set $\val{\pneg\pneg\varphi}$ is either $\N$ or
$\varnothing$, and there is nothing to choose in it.
\end{remark}

The soundness of all three calculi
is Theorem~\ref{th:sound} (\S\ref{sec:sound}); first come the derived
rules and examples of derivations.

One more natural question: why exactly these classes of
problems are taken in (A5) and (A7). The short answer: we took the largest natural syntactically definable
classes of problems for which these schemes are valid. A more detailed answer is in Section~\ref{sec:whyclasses} in the Appendix.

\subsection{Derived inference rules}\label{sec:derived}

\begin{proposition}\label{prop:derived}
In the calculi $\HCglm$, $\HCgl$, $\HCB$ the inference
rules~(1), (4), (5) are admissible, and the deduction lemma~(2) and the
lemma on transferring derivations~(3) hold. In $\HCE$, where there is no rule (AtoP), the
deduction lemma~(2), the second half of lemma~(3) and rule~(4) hold; rules (1),
(5) and the first half of lemma~(3) rely on (AtoP).
\begin{enumerate}[leftmargin=2.2em,itemsep=2pt,label=\textup{(\arabic*)}]
\item the necessitation rule for $\square$: from $A$ one can derive $\square A$.
\item Deduction lemma. By a \emph{derivation from the hypotheses} $\varphi_1,\dots,\varphi_n$
we mean a derivation that uses only the rule (MP) and in which the formulas
$\varphi_1,\dots,\varphi_n$, together with any formulas derivable in
$\HCglm$, are used as axioms. Then, if there is a derivation of the problem $\chi$ from the hypotheses $\varphi_1,\dots,\varphi_n$, the problem
\[
\varphi_1\pto(\varphi_2\pto(\dots(\varphi_n\pto\chi)\dots))
\]
is derivable in $\HCglm$.
The letters here are Greek: the hypotheses and the conclusion are problems. The analogous statement holds
for propositions and for $\HCgl$.
\item Lemma on transferring derivations.
If a proposition $A$ is derivable from the hypotheses $B_1,\dots,B_n$ in $\HCglm$,
then the problem $\bang A$ is derivable from the hypotheses $\bang B_1,\dots,\bang B_n$ in $\HCglm$.
Similarly, if a problem $\chi$ is derivable from the hypotheses
$\varphi_1,\dots,\varphi_n$ in $\HCglm$, then the proposition $\quest\chi$ is derivable from the
hypotheses $\quest\varphi_1,\dots,\quest\varphi_n$ in $\HCglm$. The analogous statement holds
for $\HCgl$ as well.
\item the rule ($\quest\to$):  
From $\varphi\pto\psi$ one can derive
$\quest\varphi\ato\quest\psi$.

\item the rule ($\bang\to$):  
From $A\ato B$ one can derive
$\bang A\pto\bang B$.

\end{enumerate}
\end{proposition}
\begin{proof}

(1) (AtoP), then (PtoA).

(2) The usual proof of the deduction lemma works.

(3) Suppose we are given a derivation of a proposition $A$ from the hypotheses $B_1,\dots,B_n$ in $\HCglm$.
Replace in this derivation every proposition $C$ by the problem $\bang C$.
Under this replacement all formulas derivable in $\HCglm$
(taken as axioms) turn, by
the rule (AtoP), into derivable formulas again. The hypotheses $B_1,\dots,B_n$ turn
into the given hypotheses $\bang B_1,\dots,\bang B_n$, and the last formula $A$---into
the desired formula $\bang A$.
However, the resulting sequence of formulas ceases to be a derivation
by the rule MP. But this is easy to repair, since if a formula $U$ is obtained from formulas $V$ and
$V\ato U$ by the rule MP, then from the formulas $\bang V$ and
$\bang(V\ato U)$ and the axiom (A3) one can derive $\bang U$ by applying MP twice.
The second transfer is justified in the same way, except that instead of (AtoP) one applies
(PtoA), and instead of (A3)---the left half of (A4a): from $\quest V$ and
$\quest(V\pto U)$ it yields $\quest U$.

(4) By (PtoA) and the left half of (A4a).

(5) By (AtoP) and (A3).
\end{proof}


Next we prove the lemma on replacement by an equivalent formula. It requires some care: our
formulas are two-sorted, and the subformula being replaced may occur under modalities.
We stress in advance that no axioms are thereby carried ``inside a modality'':
replacement under the signs $\bang$ and $\quest$ is performed by lifting a
\emph{derivable} equivalence through the modality by the rules (AtoP) and
(PtoA), and these are applicable precisely because the equivalence is a theorem.

\begin{lemma}[on replacement by an equivalent formula]\label{lem:repl}
Suppose that an equivalence $C\aeq D$ (propositions) or
$\chi\peq\chi'$ (problems) is derivable in $\HCglm$, and a formula $E'$ (of either sort) is obtained from a formula
$E$ by replacing some occurrences of the subformula $C$ by $D$ (respectively, $\chi$
by $\chi'$). Then the equivalence $E\aeq E'$ or
$E\peq E'$, according to the sort of $E$, is derivable in $\HCglm$. The same holds in $\HCgl$, in $\HCgls$ and
---for formulas with the quantifier $\Pi$---in $\HCB$.
\end{lemma}

\begin{proof}
Induction on the construction of the formula $E$. If $E$ is a replaced occurrence, the equivalence is given by
hypothesis; if there are no such occurrences inside $E$, then $E=E'$.

The inductive step is treated according to the main connective of $E$.
Propositional connectives: if $\vdash E_1\aeq E_1'$ and $\vdash E_2\aeq E_2'$, then
$\vdash(E_1\star E_2)\aeq(E_1'\star E_2')$ for
$\star\in\{\aand,\aor,\ato\}$, since
$(X\aeq X')\ato\bigl((Y\aeq Y')\ato\bigl((X\star Y)\aeq(X'\star
Y')\bigr)\bigr)$ is a classical tautology, scheme (A1); it remains to apply (MP).

Problem connectives: the analogous congruence schemes
$(\varphi\peq\varphi')\pto\bigl((\psi\peq\psi')\pto\bigl((\varphi\star\psi)
\peq(\varphi'\star\psi')\bigr)\bigr)$ for
$\star\in\{\pand,\por,\pto\}$ are derivable in $\Int$, that is, they belong to (A2).

The modality $\bang$: let $\vdash A\aeq A'$. From both implications, by (A1),
we obtain $\vdash A\ato A'$ and $\vdash A'\ato A$; the rule ($\bang\to$) yields
$\bang A\pto\bang A'$ and $\bang A'\pto\bang A$. 

The modality $\quest$: similarly, except that instead of ($\bang\to$)
we use ($\quest\to$)


The quantifier $\Pi$ (in $\HCB$): let $\vdash\varphi\peq\varphi'$. By
$(\Pi\mathrm{E})$ the formula $\Pi\alpha\,\varphi\pto\varphi$ is derivable, by (A2)---so is
$\varphi\pto\varphi'$, and their composition is $\Pi\alpha\,\varphi\pto\varphi'$;
the variable $\alpha$ does not occur free in $\Pi\alpha\,\varphi$, and Bernays'
rule yields $\Pi\alpha\,\varphi\pto\Pi\alpha\,\varphi'$; the converse
implication is symmetric. 
For propositions
$\Pi\alpha\,A$ the argument is the same, with the propositional $(\Pi\mathrm{E})$, Bernays'
rule and the schemes (A1).
%
\end{proof}

For the calculi with the quantifier we also need derived laws of $\Pi$ itself.

\begin{lemma}[derived laws of $\Pi$]\label{lem:piderived}
The following problems are derivable in $\HCB$.
\begin{enumerate}[leftmargin=2.2em,itemsep=1pt,label=\textup{(\arabic*)}]
\item $\Pi\alpha\,\varphi\peq\varphi$, if $\alpha\notin FV(\varphi)$.
\item $(\Pi\pto)$: $\Pi\alpha(\varphi\pto\psi)\peq(\varphi\pto\Pi\alpha\,\psi)$,
  if $\alpha\notin FV(\varphi)$.
\item $\Pi\alpha\Pi\beta\,\varphi\peq\Pi\beta\Pi\alpha\,\varphi$.
\item $(\Pi\pand)$:
  $\Pi\alpha(\varphi\pand\psi)\peq\Pi\alpha\,\varphi\pand\Pi\alpha\,\psi$,
  as well as the implication converse to the axiom $(\Pi\por)$,
  $\Pi\alpha\,\varphi\por\Pi\alpha\,\psi\pto\Pi\alpha(\varphi\por\psi)$.
\end{enumerate}
The same statements hold for propositions:
$\Pi\alpha\,A\aeq A$ for $\alpha\notin FV(A)$;
$\Pi\alpha(A\ato B)\aeq(A\ato\Pi\alpha\,B)$ for $\alpha\notin FV(A)$;
$\Pi\alpha\Pi\beta\,A\aeq\Pi\beta\Pi\alpha\,A$;
$\Pi\alpha(A\aand B)\aeq\Pi\alpha\,A\aand\Pi\alpha\,B$ and
$\Pi\alpha\,A\aor\Pi\alpha\,B\ato\Pi\alpha(A\aor B)$.
\end{lemma}

\begin{proof}
Everywhere below ``Bernays' rule'' is applied to an implication in whose premise
$\alpha$ does not occur free; this is each time evident from the form of the
premise.

(1) The implication $\Pi\alpha\,\varphi\pto\varphi$ is $(\Pi\mathrm{E})$
with $\theta:=\alpha$; the converse is obtained by Bernays' rule from the
theorem $\varphi\pto\varphi$.

(2) Left to right: from $(\Pi\mathrm{E})$
$\Pi\alpha(\varphi\pto\psi)\pto(\varphi\pto\psi)$ the formula
$\varphi\pand\Pi\alpha(\varphi\pto\psi)\pto\psi$ is derivable in $\Int$; Bernays' rule yields
$\varphi\pand\Pi\alpha(\varphi\pto\psi)\pto\Pi\alpha\,\psi$, and it remains
to pass back to an implication in $\Int$. Right to left: from
$(\Pi\mathrm{E})$ $\Pi\alpha\,\psi\pto\psi$ the formula
$(\varphi\pto\Pi\alpha\,\psi)\pto(\varphi\pto\psi)$ is derivable in $\Int$; Bernays' rule
yields $(\varphi\pto\Pi\alpha\,\psi)\pto\Pi\alpha(\varphi\pto\psi)$.

(3) From $(\Pi\mathrm{E})$ applied twice,
$\Pi\alpha\Pi\beta\,\varphi\pto\varphi$ is derivable; Bernays' rule with respect to $\alpha$,
then with respect to $\beta$, yields
$\Pi\alpha\Pi\beta\,\varphi\pto\Pi\beta\Pi\alpha\,\varphi$; the converse
implication is symmetric.

(4) $(\Pi\pand)$ left to right: from $(\Pi\mathrm{E})$ and the projections
of $\Int$ the formulas $\Pi\alpha(\varphi\pand\psi)\pto\varphi$ and
$\Pi\alpha(\varphi\pand\psi)\pto\psi$ are derivable; Bernays' rule to each of them and
$\pand$-introduction. Right to left: from $(\Pi\mathrm{E})$ applied twice,
$\Pi\alpha\,\varphi\pand\Pi\alpha\,\psi\pto\varphi\pand\psi$, then
Bernays' rule.

The implication converse to $(\Pi\por)$: from
$(\Pi\mathrm{E})$ and $\Int$ the formula $\Pi\alpha\,\varphi\pto\varphi\por\psi$ is derivable,
Bernays' rule yields $\Pi\alpha\,\varphi\pto\Pi\alpha(\varphi\por\psi)$;
likewise for $\Pi\alpha\,\psi$; it remains to apply $\por$-elimination in $\Int$.

For propositions the derivations are the same, with the propositional axiom
$(\Pi\mathrm{E})$, Bernays' rule for propositions and classical
tautologies in place of theorems of $\Int$.
\end{proof}

\remove{
Выводимы также правило обобщения и $(\Pi\mathrm{K})$: из $\varphi$
выводима $\Pi\alpha\,\varphi$, и выводима
$\Pi\alpha(\varphi\pto\psi)\pto(\Pi\alpha\,\varphi\pto\Pi\alpha\,\psi)$
(так же для высказываний). Обобщение: из $\varphi$ в $\Int$ выводима
$(\varnothing\pto\varnothing)\pto\varphi$; правило Бернайса даёт
$(\varnothing\pto\varnothing)\pto\Pi\alpha\,\varphi$, и (MP) с теоремой
$\varnothing\pto\varnothing$ --- $\Pi\alpha\,\varphi$ (для высказываний
вместо $\varnothing\pto\varnothing$ берётся $\aneg\bot$).
$(\Pi\mathrm{K})$: из $(\Pi\mathrm{E})$ дважды в $\Int$ выводима
$\Pi\alpha(\varphi\pto\psi)\pand\Pi\alpha\,\varphi\pto\psi$; правило
Бернайса даёт
$\Pi\alpha(\varphi\pto\psi)\pand\Pi\alpha\,\varphi\pto\Pi\alpha\,\psi$,
и остаётся перейти к импликации.
}

\subsection{Examples of derivations}\label{sec:examples}

\newcommand{\provgl}{\vdash_{\HCgl}}
\newcommand{\provglm}{\vdash_{\HCglm}}
\newcommand{\provgls}{\vdash_{\HCgls}}
\newcommand{\provglb}{\vdash_{\HCB}}

\begin{example}\label{ex:proofs}
We give examples of derivations in the calculi $\HCglm$, $\HCI$, $\HCE$.
Some of them will be needed later; in each item we indicate in which
calculus exactly the derivation is carried out.
\begin{enumerate}[leftmargin=2.2em,itemsep=2pt,label=\textup{(\arabic*)}]
\item In $\HCglm$ the $\square$-version of (A3) is derivable:
  $\square(A\ato B)\ato(\square A\ato\square B)$. To derive it, apply ($\quest\to$)
to (A3), then (A4a) and the Lemma on replacement by an equivalent formula.

\item In $\HCglm$ the formula $A\ato\square A$ is derivable for every
$\Sigma$-proposition $A$ ($\Sigma_1$-completeness on the propositional side);
in particular, the $\square$-version of (A5), the scheme
$\square A\ato\square\square A$, is derivable, since $\square A$ is a
$\Sigma$-proposition. First we note that every $\Sigma$-proposition
$A$ is derivably equivalent to the proposition $\quest\sigma_A$ for some
$\Sigma$-problem $\sigma_A$: by (A4d) $\bot\aeq\quest\varnothing$, by (A4b) and
(A4c) $\quest\sigma\aand\quest\tau\aeq\quest(\sigma\pand\tau)$ and
$\quest\sigma\aor\quest\tau\aeq\quest(\sigma\por\tau)$, and the problems
$\varnothing$, $\sigma\pand\tau$, $\sigma\por\tau$ again lie in the class
$\Sigma$; it remains to do induction on the construction of $A$ with Lemma~\ref{lem:repl}.
Now axiom (A5) gives $\sigma_A\pto\bang\quest\sigma_A$, the rule
$(\quest\to)$ gives $\quest\sigma_A\ato\square\quest\sigma_A$, and
Lemma~\ref{lem:repl} replaces $\quest\sigma_A$ by $A$ on both sides.


\item In $\HCgl$ the following is derivable:
\[
 \text{(D4)}\  \square (\square A\ato A)\ato\square A\ \text{(L\"ob's theorem)}
  \]
Apply ($\quest\to$) to (A6).

\item  In $\HCE$ the formulas $\pneg\bang\bot$ and $\aneg\square\bot$ are derivable.
First we apply the rule $(\bang\to)$ to the formula $\bot\ato\quest\varnothing$, a classical tautology (A1).
We obtain the formula $\bang\bot\pto\bang\quest\varnothing$.

 Next, the scheme (A9) with
$\varphi=\varnothing$ gives $\bang\quest\varnothing\pto\pneg\pneg\varnothing$, and
$\pneg\pneg\varnothing\pto\varnothing$ is an intuitionistic theorem
. Altogether $\bang\bot\pto\varnothing$, that
is, $\pneg\bang\bot$. Now we apply ($\quest\to$) to this formula,
obtaining $\square\bot\pto\quest\varnothing$. Finally, using (A4d),
we obtain $\square\bot\pto\bot$.
(We note in passing that the last formula is a special case of axiom A8,
so it could have been derived more briefly.)

\item This example of derivability in $\HCzero$ will be needed in the proof of validity of axiom (A6).
 Let $P$ and $S$ be arbitrary propositions. Then in $\HCglm$ the following problem is derivable:
\begin{equation}\tag{L}\label{eq:lobimpl}
  \bang\bigl(P\aeq(\square P\ato S)\bigr)\ \pto\
  \bigl(\bang(\square S\ato S)\pto\bang S\bigr).
\end{equation}
First we use the lemma on transferring derivations, and then the deduction lemma. To this end,
we take as hypotheses the propositions
$$
\square S\ato S, \qquad
\square (P\ato (\square P\ato S)), \qquad (\square P\ato S)\ato P
$$
and derive the formulas $\square P\ato S$ and $P$.

\begin{enumerate}[leftmargin=2.6em,itemsep=5pt,label=\textup{(\arabic*)}]

\item $\square P\ato\square(\square P\ato S)$\\
  \emph{this formula is derived by MP from the second hypothesis and the first example of derivable formulas
}
\item $\square(\square P\ato S)\ato(\square\square P\ato\square S)$\\
  \emph{this is the formula from our first example of derivable formulas with $A:=\square P$, $B:=S$.}
\item $\square P\ato\square\square P$\\
  \emph{this is the formula from the second example of derivable formulas.}
\item $\square P\ato\square S$\\
  \emph{Derived from (1), (2), (3) in the classical propositional calculus.}
\item $\square P\ato S$\\
  \emph{Derived from (4) and the first hypothesis in the classical propositional calculus.}
\item $P$\\
\emph{from the third hypothesis and (5) by MP.}

\item  Note that the whole derivation (1)--(6) is a derivation from the three adopted hypotheses.
Therefore we may apply the lemma on transferring derivations. By this lemma
there is a derivation of the formulas $\bang(\square P\ato S)$ and $\bang P$ from the hypotheses
$$
\bang(\square S\ato S),
\quad \bang\square (P\ato(\square P\ato S)),
\quad \bang((\square P\ato S)\ato P).
$$
We show that this derivation can be continued to a derivation of the problem $\bang S$.
\item  From the same hypotheses one can also derive the formula $\bang\square P$
(by appending to the derivation axiom (A5) with $\sigma:=\bang P$ and applying MP).
\item By appending to the derivation axiom (A3)
and applying MP, one can further derive the problem $\bang\square P\pto\bang S$.
\item  Finally, one can derive the formula $\bang S$ by applying MP
once more.

Thus the formula $\bang S$ is derivable from these three hypotheses.

\item Now we show that the second and the third hypotheses can be replaced by a single hypothesis
$$
F:=\bang\bigl(P\aeq(\square P\ato S)\bigr),
$$
which is the first premise of the original formula.
First, the third hypothesis is derivable from $F$: the proposition
$(\square P\ato S)\ato P$ is derivable from $P\aeq(\square P\ato S)$ in
the classical calculus, and the lemma on transferring derivations carries this derivation
under the sign~$\bang$.
In the same way the formula
$\bang(P\ato(\square P\ato S))$ is derivable from $F$. It differs from the second hypothesis only
by the presence of a box after the exclamation mark, and axiom (A5) with
$\sigma:=\bang(P\ato(\square P\ato S))$ is exactly what allows us to derive the second
hypothesis from it by MP.

Thus, from the hypotheses $F$ and $\bang(\square S\ato S)$
one can derive the formula $\bang S$, and moreover by applications of MP alone. It remains
to apply the deduction lemma and obtain formula~\eqref{eq:lobimpl}.

We note that the derivation of this item has been formalized and checked in the
Metamath system, in fully expanded form: the deduction lemma and the lemma on
transferring derivations are not used in the formalization but are unfolded
into an explicit derivation from the axioms (about seven hundred elementary steps). The
derivations of items (1)--(4) have been checked as well.
\end{enumerate}
\end{enumerate}
\end{example}

\begin{remark}[on the class in (A5)]
Axiom (A5) is stated at once for the whole class $\Sigma$, although its
special case $\sigma:=\bang A$ would suffice. Indeed,
by induction on the construction of a $\Sigma$-problem $\sigma$ one derives
$\sigma\pto\bang\quest\sigma$: for $\varnothing$, by $\Int$; for
$\bang A$, this is the special case; if $\sigma=\sigma_1\pand\sigma_2$, then
from $\sigma_i\pto\bang\quest\sigma_i$ we obtain
$\sigma\pto\bang\quest\sigma_1\pand\bang\quest\sigma_2$, then
$\sigma\pto\bang(\quest\sigma_1\aand\quest\sigma_2)$ (from the tautology
$A\ato(B\ato A\aand B)$ by the rule (AtoP) and the scheme (A3)), and
$\quest\sigma_1\aand\quest\sigma_2\aeq\quest\sigma$ by (A4b), and
Lemma~\ref{lem:repl} finishes the argument; for $\por$ likewise, with the rule
$(\bang\to)$ from the tautologies $A\ato A\aor B$, $B\ato A\aor B$ and with (A4c).
This is the same observation as in item~(2) of Example~\ref{ex:proofs}: the schemes
(A4) glue the $\quest$-images of $\Sigma$-problems into the $\quest$-image of a single
$\Sigma$-problem. Stating (A5) at once for the class $\Sigma$ is
more natural: semantically it is effective $\Sigma_1$-completeness
\emph{with a witness on the input} (Proposition~\ref{prop:a5a7}), and
the restriction to $\bang A$ in it is accidental.

Axiom (A5) and item~(2) of Example~\ref{ex:proofs} are one and the same
statement seen from two sides: effective $\Sigma_1$-completeness in
$\bang$-form generates both co-reflection on $\Sigma$-problems and
$\Sigma_1$-completeness on $\Sigma$-propositions, while $\Sigma$-propositions
are, up to derivable equivalence, the propositions
$\quest\sigma$ for $\Sigma$-problems $\sigma$. Therefore the boundaries of both classes
are the same. A variable cannot be admitted into $\Sigma$-propositions: with
$A^*:=\Con$ the scheme $A\ato\square A$ gives $\Con\ato\Prov(\Con)$, which is false by
G\"odel's second theorem. Neither can implication: $\square\bot\ato\bot$ is $\Con$, and
the same thing results.
(End of remark.)
\end{remark}

\remove{
\begin{remark}[почему в $\HCgls$ нет (AtoP)]\label{rem:lobfork}
Рефлексия и необходимость несовместимы. Из (A8) применением (AtoP) и
(PtoA) получается формула $\square(\square A\ato A)$, схема Лёба (D4) даёт $\square A$, а рефлексия --- формулу $A$. При этом
$A$ здесь может быть любым высказыванием. Именно
поэтому в $\HCE$ остаются только (MP) и (PtoA); непротиворечивость
$\HCE$ следует из теоремы~\ref{th:sound}.
\end{remark}

\begin{remark}[о $\square$-фрагментах]\label{rem:names}
В $\HCI$ модальность $\square$ удовлетворяет схеме K, схеме
$\square A\ato\square\square A$ (пример~\ref{ex:proofs}(1),(2)) и
правилу необходимости (предложение~\ref{prop:derived}(1)), то есть
аксиоматике $\GL$ без схемы Лёба, --- это $\mathrm{K4}$; добавление
(A6) даёт $\GL$. Система $\HCE$ добавляет рефлексию при единственных
правилах (MP) и (PtoA) --- так же, как $\GLS$ получается из $\GL$.
$\square$-фрагменты и совпадают с этими логиками в точности ---
теорема~\ref{th:assval} ниже.
\end{remark}
}

\subsection{Soundness of the calculi}\label{sec:sound}

\begin{theorem}[soundness]\label{th:sound}
For each of the three semantics $Z\in\{\sI,\sB,\sE\}$, every theorem
of the calculus $\HC_Z$ is $Z$-valid.
\end{theorem}

\begin{proof}
Induction on the derivation: we have to check that every axiom is valid in
the corresponding semantics and that every rule preserves
validity.

(1) \emph{Axioms other than (A6).} The validity in all three
semantics of the schemes (A1) is seen directly: the translation is a substitution instance
of a classical tautology, and its proof is computably constructed from the code
of the interpretation. The validity of the schemes (A2) is
Proposition~\ref{prop:intval}, of the scheme (A3)
Proposition~\ref{prop:a3}, of the group (A4) Lemma~\ref{lem:qsem}
(the translations are provable in the internal theory, the proofs are constructed from the code),
of the schemes (A5) and (A7) Proposition~\ref{prop:a5a7}. The $\sE$-validity of the schemes (A8) and (A9) is
Proposition~\ref{prop:hce}(1). The $\Pi$-axioms in $\sB$ are Proposition~\ref{prop:pilaws}. All this is also collected in
Table~\ref{tab:axioms}.

(2) \emph{Rules.} Everywhere below, when speaking of provability in the external
theory, in $\sB$ we mean provability of the closure from
Definition~\ref{def:valid}.

(MP) preserves validity by Lemmas~\ref{lem:closure}
and~\ref{lem:closurebs}.

(PtoA): let $S$ solve the problem $\varphi$, that is, the external theory
proves $\varphi^*(S(\gn{*}))$ for every interpretation $*$
(its closure $\overline{\varphi^*}(S(\gn{*}))$ for interpretations in $\sB$).
Append to its proof one axiom $\varphi^*(S(\gn{*}))\ato \exists x\,\varphi^*(x)$
of the predicate calculus and the formula $\exists x\,\varphi^*(x)$---we obtain a proof
of the formula $(\quest\varphi)^*$.  In $\sB$
we derive the formula $\exists x\,\overline{\varphi^*}(x)$, from which the
required formula $\overline{\exists x\,\varphi^*(x)}$ is derived by standard means of the predicate
calculus. Thus the proposition $\quest\varphi$
is valid.

(AtoP) in $\sI$, $\sB$: let the proposition $A$ be valid, that
is, the external theory, which here coincides with the internal one, proves $A^*$
(in $\sB$, the closure $\overline{A^*}$). This is exactly the formula
whose proofs are the solutions of the problem $\bang A$ (Definitions~\ref{def:tr},
\ref{def:starB}). From the code of the interpretation
some proof $P(*)$ of this formula is found by exhaustive search (Remark~\ref{rem:whycomp}).
The statement $\Prf(P(\gn{*}),A^*)$  (or $\Prf(P(\gn{*}),\overline{A^*})$) is a true
$\Sigma_0$-sentence, and the external theory proves it. Hence $P$
solves the problem $\bang A$.

Bernays' rule in $\sB$: let $S$ solve $\varphi\pto\psi$, where
$\alpha\notin FV(\varphi)$, that is, $\ZF$ proves the closure
$\forall\vec X\,\forall X_\alpha\,\forall u\,\bigl(\varphi^*(u)\ato
\psi^*(U_{S(\gn{*})}(u))\bigr)$ over all boxes, including
$X_\alpha$. The translation $(\varphi\pto\Pi\alpha\,\psi)^*(S(\gn{*}))$ is
$\forall u\,\bigl(\varphi^*(u)\ato\forall X_\alpha\,
\psi^*(U_{S(\gn{*})}(u))\bigr)$; the variable $X_\alpha$ does not
occur in $\varphi^*$, and the closure of this formula is obtained from the previous one by pure
logic---by moving the quantifier $\forall X_\alpha$ across the premise. The same
$S$ solves $\varphi\pto\Pi\alpha\,\psi$. For propositions likewise: from
$\ZF\vdash\forall\vec X\,\forall X_\alpha\,(A^*\ato B^*)$ with
$X_\alpha$ not occurring in $A^*$ it follows that
$\ZF\vdash\forall\vec X\,(A^*\ato\forall X_\alpha\,B^*)$.

(3) \emph{The scheme (A6).}
First we note that in the semantics $\sB$ it suffices to treat
propositions $A$ without free problem variables: the instance of the scheme
for an arbitrary $A$ is \emph{derived} in $\HCB$ from the instance for the
proposition $\bar A:=\Pi\vec\alpha\,A$, where $\vec\alpha$ are the free
problem variables of $A$, using the axioms and rules already checked in
items~(1) and~(2); therefore the validity of the instance for $\bar A$
implies the validity of the instance for $A$.
For this it suffices to prove that the problems
$$
\bang(\square A\ato A)\pto\bang A,  \qquad
\bang(\square \bar A\ato \bar A)\pto\bang \bar A
$$
are derivably equivalent in $\HCB$ (without (A6)).
By Lemma~\ref{lem:repl} it suffices to establish the equivalence of their premises
and of their conclusions.

By the axiom $(\Pi\bang)$
the equivalence of the conclusions
$\bang A\peq\bang\bar A$ is derivable, and from it by the rule $(\quest\to)$
$\square A\aeq\square\bar A$.
For the premises the same axiom $(\Pi\bang)$ gives
$$
\bang(\square A\ato A)\peq\bang\Pi\vec\alpha(\square A\ato A)
$$
(the free problem variables of the proposition $\square A\ato A$ are the same
$\vec\alpha$, since $\square A=\quest\bang A$ has none).
By Lemma~\ref{lem:piderived}(2) for propositions
($\vec\alpha$ do not occur in $\square A$) the equivalence
$$
\Pi\vec\alpha(\square A\ato A)\aeq(\square A\ato\Pi\vec\alpha\,A)
$$
is derivable, and by Lemma~\ref{lem:repl} in its right-hand side $\square A$ can be replaced
by $\square\bar A$; we obtain
$\Pi\vec\alpha(\square A\ato A)\aeq(\square\bar A\ato\bar A)$.
Finally, we apply Lemma~\ref{lem:repl} once more, now under the sign
$\bang$, and obtain $\bang\Pi\vec\alpha(\square A\ato A)\peq
\bang(\square\bar A\ato\bar A)$, and with it the required equivalence
of the premises $\bang(\square A\ato A)\peq\bang(\square\bar A\ato\bar A)$.

Now that the validity of the remaining axioms and
the soundness of the rules have been established, Löb's scheme need not be checked
directly; instead it can be \emph{derived} from them and the fixed-point lemma for
the internal theory. We need the fixed-point lemma in an effective
form:
\begin{quote}
from a closed formula $\sigma$ (arithmetical for $\PA$; of the language of
$\ZF$ for $\ZF$) one can compute a closed formula $\pi$ and a proof in the
internal theory of the equivalence $\pi\aeq(\Prov(\pi)\ato\sigma)$.
\end{quote}

The derivation
splits into two parts. First, a derivation in $\HCglm$ of the formula
\begin{equation}\tag{L$'$}
  \bang\bigl(P_0\aeq(\square P_0\ato A)\bigr)\ \pto\
  \bigl(\bang(\square A\ato A)\pto\bang A\bigr),
\end{equation}
where $P_0$ is any propositional variable not occurring in $A$:
this is Example~\ref{ex:proofs}(5) with $P:=P_0$, $S:=A$. Since
(L$'$) is derived without (A6), by what has already been proved it is valid---in
$\sI$, and also in $\sB$, where $\HCglm$ is contained in $\HCB$. Let $A$ be a proposition (in $\sB$, by what has been proved, without
free problem variables) for which we need to prove the
validity of the problem $\bang(\square A\ato A)\pto\bang A$. We have to
construct an algorithm which, given an interpretation $*$, finds an $x$
for which the external theory proves
$\bigl(\bang(\square A\ato A)\pto\bang A\bigr)^*(x)$. Here we may
assume that $*$ assigns values only to the variables occurring in
our formula; in particular, it assigns no value to the variable $P_0$.
Since the formula (L$'$) is valid, we have an algorithm that
delivers the required $x$ if it is given as input any $y,*'$ for
which the external theory proves
$\bigl(\bang(P_0\aeq(\square P_0\ato A))\bigr)^{*'}(y)$ and $*'$
extends $*$ to the variable $P_0$. Such $y,*'$ are delivered by the
fixed-point lemma: put $*'(P_0)=\pi$, for which the
internal theory proves the equivalence
$\pi\aeq(\Prov(\pi)\ato A^*)$; $y$ is its proof, found by exhaustive search
through all proofs (the formula $\Prf(y,\dots)$ is a true
$\Sigma_0$-sentence, and the external theory proves it). It is exactly here
that the closedness of the formula $A^*$ is needed: in $\sB$ a solution of the problem
$\bang(P_0\aeq(\square P_0\ato A))$ is a proof of the closure of the equivalence
over the free problem variables of $A$, whereas the fixed-point lemma gives
the equivalence only for closed $A^*$; the formula (L$'$) itself is
valid for any $A$.

(4) \emph{The semantics $\sE$.} The axioms of $\HCE$ are the theorems of $\HCI$ and the schemes
(A8), (A9). Every theorem of $\HCI$ is $\sI$-valid by what has been proved, and
$\sI$-validity implies $\sE$-validity, since the internal
theory $\PA$ is contained in the external theory $\TA$. The schemes (A8) and (A9)
are $\sE$-valid by Proposition~\ref{prop:hce}(1), and the rules (MP)
and (PtoA) preserve $\sE$-validity by item~(2).
\end{proof}

\subsection{What the quantifier $\Pi$ gives}\label{sec:pigain}

Recall the motivation for the quantifier (\S\ref{sec:piatoms}): under the sign
$\quest$ the constructive content of a problem is lost, and the propositions
$\quest(\alpha\por\pneg\alpha)$, $\quest\rho$ are valid trivially.
Now that we have the calculi and the soundness theorem, what was promised
can be proved.

\begin{theorem}[three problems with the quantifier $\Pi$ equivalent to $\varnothing$]\label{prop:piatoms}
In $\sB$ the equivalences
\[
  \Pi\alpha\,\alpha\peq\varnothing,\qquad
  \Pi\alpha(\alpha\por\pneg\alpha)\peq\varnothing,\qquad
  \Pi\alpha(\pneg\pneg\alpha\pto\alpha)\peq\varnothing
\]
are valid, and so are the negations of the propositions
\[
  \quest\Pi\alpha\,\alpha,\qquad
  \quest\Pi\alpha(\alpha\por\pneg\alpha),\qquad
  \quest\Pi\alpha(\pneg\pneg\alpha\pto\alpha).
\]
These three propositions themselves are therefore not valid (under the
assumption that $\ZF$ is consistent).
\end{theorem}

\begin{proof}We derive all six formulas in the calculus $\HCB$.
The implications from right to left for the first three formulas
are axioms of $\Int$. Now we derive the implications in the opposite direction.

%
(1) 
The formula $\Pi\alpha\,\alpha\pto\varnothing$ is an axiom
$(\Pi\mathrm{E})$ with $\theta:=\varnothing$.

(2) 
Indeed, $(\Pi\mathrm{E})$ with
$\theta:=\varnothing$ gives $\Pi\alpha\,\alpha\pto\varnothing$.
Similarly, $(\Pi\mathrm{E})$ with
$\theta:=\pneg\varnothing$ gives $\Pi\alpha\,\pneg\alpha\pto\pneg\pneg\varnothing$,
and $\pneg\pneg\varnothing\pto\varnothing$ is derivable in $\Int$.
Hence $\Pi\alpha\,\pneg\alpha\pto\varnothing$ is derivable.
Together with
$(\Pi\por)$ this gives $\Pi\alpha(\alpha\por\pneg\alpha)\pto\varnothing$.

(3) 
Indeed,
for $H:=\Pi\alpha(\pneg\pneg\alpha\pto\alpha)$: $(\Pi\mathrm{E})$ with
$\theta:=\alpha\por\pneg\alpha$
gives
$$
H\pto  (\pneg\pneg  (\alpha\por\pneg\alpha) \pto (\alpha\por\pneg\alpha)).
$$
The theorem $\pneg\pneg(\alpha\por\pneg\alpha)$ of $\Int$
gives $H\pto(\alpha\por\pneg\alpha)$. Then
Bernays' rule ($\alpha$ does not occur free in $H$) gives
$H\pto\Pi\alpha(\alpha\por\pneg\alpha)$,
and from this $H\pto\varnothing$ by the previous item.

To prove the second assertion (on the derivability of the negations)
 we apply to each of the derived implications ($\quest\to$) and then (A4d).
 %
%
%
\end{proof}

Thus $\quest\Pi\alpha(\alpha\por\pneg\alpha)$ is refutable, whereas
without the quantifier $\quest(\alpha\por\pneg\alpha)$ is a theorem (by (A4) it
is equivalent to the tautology $\quest\alpha\aor\aneg\quest\alpha$)
: the classical case analysis does not pass through the quantifier, and now the
difference between the law of excluded middle and its constructive version
becomes visible at the level of propositions as well.

\begin{example}\label{ex:pider}
Let $\rho$ be Rose's formula, that is, the problem \eqref{rose-formula} for
$\varphi=\pneg\beta\por\pneg\gamma$, where $\beta,\gamma$ are problem
variables (\S\ref{sec:modalfreeval} below). In $\sB$ the proposition
$\quest\Pi\beta\Pi\gamma\rho$ is valid: the solver of Rose's formula is constant and
does not look at the contents of the boxes (Proposition~\ref{cor:rose} below), so that
$\ZF$ proves $\forall X_\beta\forall X_\gamma\,\rho^*(r)$ for its
code $r$, that is, $(\Pi\beta\Pi\gamma\rho)^*(r)$, and with it
$(\quest\Pi\beta\Pi\gamma\rho)^*$.
Without the quantifier the proposition $\quest\rho$ is by (A4) equivalent to a tautology and
says nothing; with the quantifier it asserts that Rose's formula is
solved by a single algorithm for any contents of the boxes. This is the first
case where a discovery on the problem side produces a new
valid proposition. Whether it is derivable in $\HCB$ we do not know:
$\rho$ itself is not derivable in $\HCI$ (\S\ref{sec:indep}), and
the conservativity of $\HCB$ over $\HCI$ on clean problems is open
(Remark~\ref{rem:qfree}).
\end{example}

The second question is how new the problems and propositions with the quantifier are.
In all nontrivial cases that we managed to analyze, the quantifier
is eliminated: a problem with $\Pi$ turns out to be equivalent to a problem without
$\Pi$. Let us collect these cases together.

\begin{proposition}[quantifier elimination]\label{prop:pienc}
Let $\varphi$ be a problem in which the variable $\beta$ does not occur
free, and let $\alpha,\beta,\gamma$ be variables. The following
equivalences are $\sB$-valid:
\begin{enumerate}[leftmargin=2.2em,itemsep=2pt,label=\textup{(\arabic*)}]
\item $\Pi\beta\,(\varphi\pto\beta)\ \peq\ \pneg\varphi$;
\item $\Pi\beta\,\bigl((\varphi\pto\beta)\pto\beta\bigr)\ \peq\ \varphi$;
\item $\Pi\beta\,\bigl(((\alpha\pto\beta)\pto\alpha)\pto\alpha\bigr)
  \ \peq\ (\pneg\pneg\alpha\pto\alpha)$;
\item $\Pi\alpha\,\bigl(((\alpha\pto\beta)\pto\alpha)\pto\alpha\bigr)
  \ \peq\ \beta$;
\item $\Pi\beta\,\bigl(((\alpha\pto\beta)\pto\beta)\pto(\alpha\por\beta)\bigr)
  \ \peq\ \alpha\por\pneg\alpha$;
\item $\Pi\alpha\,\bigl(((\alpha\pto\beta)\pto\beta)\pto(\alpha\por\beta)\bigr)
  \ \peq\ \beta$;
\item $\Pi\alpha\,\bigl((\alpha\pto\beta\por\gamma)\pto
  ((\alpha\pto\beta)\por(\alpha\pto\gamma))\bigr)
  \ \peq\ (\gamma\pto\beta)\por(\beta\pto\gamma)$;
\item $\Pi\alpha\,\bigl((\beta\pand\gamma\pto\alpha)\pto
  ((\beta\pto\alpha)\por(\gamma\pto\alpha))\bigr)
  \ \peq\ (\beta\pto\gamma)\por(\gamma\pto\beta)$.
\end{enumerate}
All of them, except items~(4) and~(6), are derivable in $\HCB$; in items~(4)
and~(6) the implication from right to left is derivable. Together with each
equivalence of problems $\Pi\vec\alpha\,\chi\peq\psi$, the equivalence
of propositions $\quest\Pi\vec\alpha\,\chi\aeq\quest\psi$ is derivable
(respectively, valid) as well.
\end{proposition}

\begin{proof}
The last assertion: from the problem $\Pi\vec\alpha\,\chi\peq\psi$ the rule
$(\quest\to)$ gives both implications between $\quest\Pi\vec\alpha\,\chi$ and
$\quest\psi$; the semantic variant is Lemma~\ref{lem:qsem} and
Theorem~\ref{th:sound}.

In all the problems the implication from right to left is proved in the same way.
Let $\Pi\delta\,\lambda$ denote the left-hand side and $\mu$ the right-hand side of the formula.
Then the implication $\mu\pto\lambda$ is derivable in Int, which is checked directly
in each of the items. After that we apply Bernays' rule.

The implication from left to right is also proved in the same way in all items, except
items (4) and (6). Namely, we apply the axiom $(\Pi\mathrm{E})$, substituting into the formula $\lambda$
some formula
$\theta$ for the bound variable in such a way that for the resulting formula $\lambda(\theta)$ the formula
$\lambda(\theta)\pto\mu$ is derivable in Int. We indicate how to choose the formula $\theta$ in each item:
(1) $\theta:=\varnothing$,
(2) $\theta:=\varphi$, (3) $\theta:=\varnothing$,
(5) $\theta:=\pneg\alpha$, (7)  $\theta:=\beta\por\gamma$;
(8) $\theta:=\beta\pand\gamma$.

Items (4) and (6) are treated separately.
\remove{
Направление справа налево --- правило Бернайса из
$\Int$-теорем $\beta\pto(((\alpha\pto\beta)\pto\alpha)\pto\alpha)$ и
$\beta\pto(((\alpha\pto\beta)\pto\beta)\pto(\alpha\por\beta))$.
Направление слева направо докажем семантически,}
We begin with (4). We write
$Y$ for $X_\beta$; we argue in $\ZF$ with $Y$ free. For
sets $X,Y$ put $X{\to}Y:=\{q:\forall x\in X\ U_q(x){\downarrow}\in Y\}$
and
$$
\Phi(X):
=\{p:\forall q\in X{\to}Y\ U_p(q){\downarrow}\in X\}.
$$
 Then
$e\in\val{\Pi\alpha(((\alpha\pto\beta)\pto\alpha)\pto\alpha)}$ means:
for all $X$ and all $p\in\Phi(X)$ we have $U_e(p){\downarrow}\in X$.
Fix a code $I$ of the identity program ($U_I(q)=q$) and a code $a$
of the program identically equal to zero. The solver is as follows:
\[
  f(e):=U_n(a),\qquad\text{where } n:=U_e(I).
\]
We show that for every $Y$ and every such $e$ the value $f(e)$
is defined and lies in $Y$.

\emph{Step 1: $n$ is defined.} For $X:=\omega$ we have $I\in\Phi(\omega)$,
since $U_I(q)=q\in\omega$ for every $q$; hence $n=U_e(I)$
is defined.

\emph{Step 2: $U_n$ maps $\omega\setminus\{n\}$ into $Y$.} Suppose
the contrary: there is $x\neq n$ for which $U_n(x)$ is undefined or does not
lie in $Y$. Put $X:=\omega\setminus\{n\}$. Then $I\in\Phi(X)$:
if $q\in X{\to}Y$, then $q\neq n$ (because of the point $x\in X$ the number $n$ does not
belong to $X{\to}Y$), that is, $U_I(q)=q\in X$. By the definition of $e$
this gives $n=U_e(I)\in X=\omega\setminus\{n\}$---a contradiction.

\emph{Step 3: $f(e)\in Y$.} If $a\neq n$, then by Step~2 the value
$U_n(a)$ is defined and lies in $Y$. If $a=n$, then $U_n$ is identically
zero, by Step~2 the number $0=U_n(x)$ for any $x\neq n$ lies in $Y$,
and $f(e)=0\in Y$.

All of this is an argument in $\ZF$ with $Y$ free. 

For (6) the argument is similar. We use the notation from the argument
for case (4). Put $k:=\lambda q.\,U_q(q)$, that is, $U_k(q)=U_q(q)$. The number $k$ belongs to
$\val{(\omega\pto Y)\pto Y}$,
so $U_e(k)$ is defined and lies in
$\omega\sqcup Y$. If it is $\langle1,y\rangle$, then $y\in Y$. If it is
$\langle0,n\rangle$, then, as in
Step~2, $U_n$
maps $\omega\setminus\{n\}$ into $Y$. Indeed, otherwise
$k$ is a solution of the problem $(\omega\setminus\{n\} \pto Y)\pto Y$
 (if $q$ is a solution of the problem $\omega\setminus\{n\}\pto Y$, then $q\ne n$ by assumption, hence
 $U_q(q)\in Y$, since $U_q(x)\in Y$ for all $x\ne n$). Consequently $U_n(a)\in Y$, which is proved in the same way as in Step~3 above.

 In the first case the solver outputs the $y$ for which $U_e(k)=c(1,y)$, and in the second case $U_n(a)$.
\end{proof}

Equivalence~(2) is the well-known impredicative representation of a
problem by a second-order quantifier: ``to solve $\varphi$'' is equivalent to
``to be able to solve every problem $\beta$ reducible to $\varphi$''.
\remove{Заметим, что
равносильности $\Pi\beta\bigl((\varphi\pto\beta)\pto\beta\bigr)\peq
\pneg\pneg\varphi$, которую подсказывает подстановка $\theta:=\varnothing$,
в общем случае нет: она дала бы $\pneg\pneg\varphi\pto\varphi$, а при
переменной $\varphi:=\alpha$ эта задача не общезначима
(предложение~\ref{prop:noem}). Для $\Sigma$-задач $\varphi$ она есть ---
по схеме (A7).}

\remove{
Первое тождество при $A:=\bot$ даёт $\quest\Pi\beta(\bang\bot\pto\beta)
\aeq\aneg\square\bot$: высказывание $C:=\quest\Pi\beta(\bang\bot\pto\beta)$
равновыводимо с $\aneg\square\bot$, а потому $\ZF$ доказывает
$C^*\aeq\Con(\ZF)$ (теорема~\ref{th:sound}).

\begin{remark}\label{rem:nosigma}
Схема $\quest\Pi\vec\alpha\,\varphi\ato\square\,\quest\Pi\vec\alpha\,\varphi$
не общезначима. Например, при $\varphi=(\bang\bot\pto\beta)$ перевод
высказывания $\quest\Pi\beta\,\varphi\ato\square\,\quest\Pi\beta\,\varphi$
$\ZF$-эквивалентен формуле
$\Con(\ZF)\ato\Prov(\Con(\ZF))$; если бы $\ZF$ её доказывала, то вместе
с формализованной второй теоремой Гёделя это дало бы
$\ZF\vdash\aneg\Con(\ZF)$, что невозможно при $\Sigma_1$-корректности
$\ZF$. Этим высказывания с квантором $\Pi$ отличаются от высказываний
$\quest\sigma$ с $\Sigma$-задачей $\sigma$, для которых схема
$\quest\sigma\ato\square\quest\sigma$ выводима из аксиомы (A5) правилом
(PtoA) и схемой (A4a): переводы $\Sigma$-задач --- $\Sigma_1$-формулы,
и работает
формализованная $\Sigma_1$-полнота, а перевод высказывания
$\quest\Pi\vec\alpha\,\varphi$ --- $\Pi^1_1$-формула.
\end{remark}
}

\begin{remark}[on items (4) and (6)]\label{rem:peirce}
Whether the left-to-right implications in items~(4) and~(6) of
Proposition~\ref{prop:pienc} are derivable in $\HCB$, we do not know. Their weakening with
$\pneg\pneg\beta$ in the conclusion is derivable: under the hypothesis $\pneg\beta$ the implication
$\alpha\pto\beta$ is equivalent to $\pneg\alpha$, Peirce's law turns into
$\pneg\pneg\alpha\pto\alpha$, Bernays' rule gives
$\Pi\alpha(\pneg\pneg\alpha\pto\alpha)$, and this problem is equivalent to
$\varnothing$ (Theorem~\ref{prop:piatoms}).

As for the non-derivability of the right-to-left implications as they stand (with the formula $\beta$), we have managed to prove it only
in the calculus HCI minus $(\Pi\por)$.
Indeed, take the two-point Kripke frame
$w_0<w_1$ in which $\beta$ is true only in $w_1$, problems are cones,
$\Pi$ is the intersection over all cones, $\bang A$ is the top, and
$\quest\varphi$ is true when $\varphi$ is true in $w_1$. All axioms and
rules of $\HCB$ except $(\Pi\por)$ hold in it (cf.\
Theorem~\ref{prop:modalfree}), and in $w_0$ the problem
$\Pi\alpha(((\alpha\pto\beta)\pto\alpha)\pto\alpha)$ holds: for any cone
$X$ the premise $X\pto\beta$ is true in $w_0$ exactly when $X$ does not contain
$w_0$, and then $(X\pto\beta)\pto X$ is false in $w_0$. And $\beta$ is false in $w_0$.
The axiom $(\Pi\por)$ fails in this frame: in $w_1$ the premise
$\Pi\alpha(\alpha\por\pneg\alpha)$ is true, while the conclusion
$\Pi\alpha\,\alpha\por\Pi\alpha\,\pneg\alpha$ is false (both disjuncts are empty:
among the cones there are the empty one and the full one), so this instance of
$(\Pi\por)$ is false both in $w_1$ and in $w_0$. Hence a possible derivation must
use $(\Pi\por)$ essentially. In the semantic
proof the decisive step was instantiating the quantifier by the box
$\omega\setminus\{n\}$, which depends on the answer of the solver itself; such a set is never the value
of a formula. Items~(4) and~(6) are candidates for
examples of incompleteness of $\HCB$ (Corollary~\ref{cor:incomplete}).
(End of remark.)
\end{remark}

An example of \emph{ineliminability} of the quantifier---a problem with $\Pi$ not
equivalent to any problem without $\Pi$---we have not managed to find. This,
of course, does not mean that the quantifier is unnecessary: the elimination in
Proposition~\ref{prop:pienc} itself requires the quantifier, and in items~(4)
and~(6) it moreover goes beyond the calculus.

\begin{remark}[the quantifier-free part]\label{rem:qfree}
The axiomatics of $\HCB$ is arranged as $\HCI$ plus a quantifier layer:
discarding the axioms for $\Pi$ and Bernays' rule, we obtain exactly $\HCI$.
The aim
of the extension is \emph{to increase the stock of provable formulas}. In the extended
language the growth is certain: new substantive propositions are derivable, such as
$\quest\Pi\beta\Pi\gamma\rho$ (Example~\ref{ex:pider}), not expressible without
the quantifier. There is growth on \emph{quantifier-free} formulas as well:
$\square\quest\alpha\ato\square\bot$ is derivable in $\HCB$ and not derivable in
$\HCI$ (Remark~\ref{rem:boxqatom}). But this growth is not new
knowledge about problems; it is a consequence of the fact that in $\sB$ the symbol $\bang$ binds
problem variables: in $\HCB$ the formula means
$\square\Pi\alpha\,\quest\alpha\ato\square\bot$ and is not clean. Therefore
the question of conservativity should properly be posed for \emph{clean}
quantifier-free formulas (Remark~\ref{rem:clean}): does the
set of clean quantifier-free theorems of $\HCB$ coincide with the set of theorems of
$\HCI$? For propositions the answer is known: in the normal form of a clean
$\Pi$-free proposition the atoms $\quest\alpha$ stand outside the symbols
$\square$, and under $\Sigma_1$-soundness of $\ZF$ its derivability in $\HCB$
is equivalent to its derivability in $\HCI$ (Proposition~\ref{prop:nf},
Theorem~\ref{th:assval}(1),(3)). For clean problems the question is open
(\S\ref{sec:todo}).
\end{remark}

\remove{
\begin{remark}[лёбовская механика не задета]\label{rem:lobpi}
Вывод задачи (L), пример~\ref{ex:proofs} о $\square$-версиях, вывод
схемы (A6) из эффективной леммы о неподвижной точке и
теорема~\ref{th:lob} переносятся в $\Pi$-расширение дословно: в них не
участвуют ни ящики, ни $\Pi$, а условия выводимости и диагональная
лемма верны для обеих внутренних теорий, $\PA$ и $\ZF$.
\end{remark}
}

\subsection{Normal form and completeness on the propositional side}
\label{sec:asscompl}

Here we describe two families of $\Pi$-free propositions: those
\emph{derivable} in each of the calculi and those \emph{valid} in
each of the semantics. Both descriptions are obtained via one and the same
\emph{normal form}, in which the modality $\quest$ is eliminated down to
the combinations $\quest\alpha$ and $\square B$, and both refer to well-known
objects: up to normal form, this is Solovay's provability logic
$\GL$ (for $\HCI$ and for the internal semantics $\sI$, $\sB$) and its
external variant $\GLS$ (for $\HCE$ and for $\sE$). The proof goes around
a circle: from derivability in our calculus to validity (the soundness
theorem), from validity to derivability in $\GL$ or $\GLS$
(Solovay's theorem), and from derivability in $\GL$ or $\GLS$ to derivability
in our calculus (direct verification). As a result, on the propositional side
the calculi $\HCI$ and $\HCE$ are complete.

The normal form is built by pushing $\quest$ inside according to the schemes (A4):
down to variables and to problems of the form $\bang B$.

\begin{definition}[translation into normal form]\label{def:nf}
To every $\Pi$-free formula $\chi$ we assign its normal form
$\nf{\chi}$. For propositions: variables and the constant $\bot$ are mapped
to themselves, the classical connectives are preserved
($\nf{(A\aand B)}=\nf{A}\aand\nf{B}$ and so on for $\aor$, $\ato$),
and $\nf{(\quest\varphi)}:=\varphi^\quest$, where the auxiliary translation
$\varphi^\quest$ of problems into propositions is given by the rules
\[
\begin{aligned}
  \alpha^\quest&:=\quest\alpha, &
  \varnothing^\quest&:=\bot, &
  (\bang B)^\quest&:=\square\,\nf{B},\\
  (\varphi\pto\psi)^\quest&:=\varphi^\quest\ato\psi^\quest, &
  (\varphi\pand\psi)^\quest&:=\varphi^\quest\aand\psi^\quest, &
  (\varphi\por\psi)^\quest&:=\varphi^\quest\aor\psi^\quest.
\end{aligned}
\]
In what follows we call $\varphi^\quest$ the \emph{Boolean collapse} of the problem $\varphi$.
For problems: $\nf{\alpha}=\alpha$, $\nf{\varnothing}=\varnothing$,
the problem connectives are preserved, and $\nf{(\bang B)}:=\bang\,\nf{B}$ (where
the proposition $B$ is normalized according to the rules above).
\end{definition}

It is easy to see that in the normal form of every $\Pi$-free
proposition $A$ the modality occurs only in the combinations
$\quest\alpha$ ($\alpha$ a variable of sort problem) and
$\quest\bang B=\square B$.

\begin{proposition}[normal form]\label{prop:nf}
For every $\Pi$-free formula, the equivalence
\[
  A\aeq\nf{A}\quad\text{for propositions},\qquad
  \varphi\peq\nf{\varphi}\quad\text{for problems}
\]
is derivable in $\HCglm$ (and hence also in $\HCgl$,
$\HCB$, $\HCgls$).
Consequently, under every interpretation the external theory of each of the three
semantics proves $A^*\aeq\nf{A}^*$ (in $\sB$, the closure of this
equivalence), and a proposition $A$ is valid in a given semantics
if and only if $\nf{A}$ is valid in it.
\end{proposition}

\begin{proof}
Simultaneous induction on the construction of the formula. For a variable and for $\bot$
this is obvious; the connectives go through by Lemma~\ref{lem:repl}. For $A=\quest\varphi$
we have to check $\provglm\quest\varphi\aeq\varphi^\quest$ by induction on
$\varphi$: for $\varphi=\alpha$ this is an identity
($\alpha^\quest=\quest\alpha$); for $\varphi=\varnothing$ it is the scheme
(A4d); for $\varphi=\bang B$ we need
$\provglm\square B\aeq\square\,\nf{B}$, which is obtained from the outer
induction hypothesis $\provglm B\aeq\nf{B}$ by Lemma~\ref{lem:repl};
the cases $\pto$, $\pand$, $\por$ are exactly the three equivalences
(A4a), (A4b), (A4c) together with the induction hypothesis and
Lemma~\ref{lem:repl}. This is where the Boolean nature of $\quest$ does its work. For problems:
a variable and $\varnothing$ are obvious, the connectives go through by
Lemma~\ref{lem:repl}, and $\nf{(\bang B)}=\bang\,\nf{B}$ is obtained from
$\provglm B\aeq\nf{B}$ by the same lemma. Only (A1), (A2),
(A3), the group (A4) and the three rules have been used.

The second claim: by Theorem~\ref{th:sound} the derived
equivalence is $\sI$-valid, hence also $\sE$-valid
($\PA\subseteq\TA$), and $\sB$-valid, since
$\HCglm\subseteq\HCB$. The set of theorems of the external theory is closed under
(MP), so $A^*$ is provable in the external theory if and only if
$\nf{A}^*$ is provable.
\end{proof}

The proviso of $\Pi$-freeness is essential: the schemes (A4) push
$\quest$ through the connectives $\pto,\pand,\por$ and the constant $\varnothing$, but not
through the quantifier $\Pi$, and in the full language of the semantics $\sB$ normalization
stops at propositions of the form $\quest\Pi\alpha\,\varphi$. There is no analogue
of the schemes (A4) for $\Pi$, and there cannot be one: the proposition
$\quest\Pi\alpha(\alpha\por\pneg\alpha)$ is \emph{refutable}---its
negation is valid (Theorem~\ref{prop:piatoms}), whereas without the
quantifier $\quest(\alpha\por\pneg\alpha)\aeq\quest\alpha\aor\aneg\quest\alpha$
is a valid classical tautology: the classical case distinction does not pass
through $\Pi$. The laws of the quantifier $\Pi$ allow one to
simplify such propositions as well; what comes out of this is described in
\S\ref{sec:pinf}.

\begin{corollary}[propositional fragment]\label{cor:assfrag}
In normal form, propositions are built by the classical connectives and
the modality $\square$ from the variables $A_i$ and from the formulas $\quest\alpha_j$. Moreover, the
formulas $\quest\alpha_j$ behave like \emph{additional independent variables
of sort proposition}: taking $\alpha^*:=\sigma$ for an arbitrary closed
$\sigma$, we get $(\quest\alpha)^*=\exists x\,\sigma$, which is logically equivalent
to $\sigma$. Therefore the question about the propositional fragment reduces entirely to the question about the
$\square$-fragment, which we take up next.
\end{corollary}

A \emph{$\square$-formula} is a proposition built by the
classical connectives $\bot,\aor,\aand,\ato$ from \emph{atoms} of two
kinds---propositional variables and propositions of the form
$\quest\beta$, where $\beta$ is a variable of sort problem---and the modality
$\square=\quest\bang$, applied to $\square$-formulas; there are no other occurrences
of $\bang$ and $\quest$. For derivability in $\GL$ and $\GLS$, atoms of both kinds
play the role of propositional variables of the modal language.
In normal form, a $\Pi$-free proposition is built by the classical connectives from
the variables $A_i$, the propositions $\quest\alpha_j$ and the propositions $\square B$, where
$B$ are again propositions in normal form. In other words,
\emph{the normal form of a $\Pi$-free proposition is a $\square$-formula}, and no separate translation into the modal language is needed. The $\quest$-collapse $\varphi^\quest$ of every $\Pi$-free
problem $\varphi$ is also a $\square$-formula (Definition~\ref{def:nf}).

Let us recall the definitions of the logics $\GL$ and $\GLS$
(see~\cite[Ch.~1]{Boolos93}); both are formulated in the language of $\square$-formulas. 
\begin{itemize}[leftmargin=1.6em,itemsep=1pt]
\item $\GL$: the axioms are all classical tautologies (in
  the modal language), the scheme K $\square(M\ato N)\ato(\square M\ato\square N)$,
  the scheme $\square M\ato\square\square M$ (traditionally called
  scheme~4) and L\"ob's scheme $\square(\square M\ato M)\ato\square M$;
  the rules are (MP) and the \emph{necessitation rule} $M/\square M$.
  (Scheme~4 is derivable from the remaining axioms, but we keep it in the list.)
\item $\GLS$: the set of formulas derivable by the rule (MP) from the theorems
  of $\GL$ and all instances of the reflection scheme $\square M\ato M$.
  There is no necessitation rule in $\GLS$, and the set of its theorems
  is not closed under it: $\GLS\vdash\square\bot\ato\bot$, but
  $\GLS\nvdash\square(\square\bot\ato\bot)$, since this formula
  is $\GL$-equivalent to $\square\bot$.
\end{itemize}
Solovay's theorems say that $\GL$ is exactly the set of $\square$-formulas all of whose arithmetical realizations are provable in $\PA$, and
$\GLS$ is exactly the set of formulas all of whose realizations are true; in our
terms this is validity in $\sI$ and in $\sE$.

For the semantics $\sI$ and $\sE$ we shall use Solovay's
theorems~\cite{Solovay76} for $\PA$; for $\sB$ their
variant for $\ZF$ is needed, and here the second assumption of
Remark~\ref{rem:assumptions} is required---the $\Sigma_1$-soundness of $\ZF$.
One cannot get by with
consistency alone already for closed formulas:
the non-validity of the formula $\square \bot$ means that
$\ZF$ does not prove $\Prov_{\ZF}(\bot)$.
Consistency does not guarantee this, while $\Sigma_1$-soundness does.
To simplify the formulation we use the facts that
(1) the language of $\GL$ is a part of our language and (2) for
the semantics $\sI$ and $\sE$ our definition of validity of
$\square$-formulas is the definition of validity in the sense of Solovay: the interpretations of
atoms of both kinds are arbitrary closed arithmetical
formulas. In the semantics $\sB$ only atoms of the first
kind are read this way (they correspond to arbitrary closed formulas of the language of $\ZF$); an atom
$\quest\beta$ there is the sentence $\exists y\,(y\in X_\beta)$ about a box, and
under the sign $\square$ it is moreover read as a closure
(Remark~\ref{rem:bangbinds}). Therefore for $\sB$ the lemma
is formulated for $\square$-formulas without atoms of the second kind.

\begin{lemma}[Solovay's theorem]\label{lem:solzf}
\begin{enumerate}[leftmargin=2.2em,itemsep=1pt,label=\textup{(\arabic*)}]
\item For every $\square$-formula $M$: $\GL\vdash M$ if and only
  if $M$ is valid in $\sI$; $\GLS\vdash M$ if and only
  if $M$ is valid in $\sE$.
\item Let $\ZF$ be $\Sigma_1$-sound. For every $\square$-formula
  $M$ containing no atoms of the form $\quest\beta$, the following statements
  are equivalent: $\GL\vdash M$ and $M$ is valid in $\sB$.
\end{enumerate}
\end{lemma}

\begin{remark}
Item (1) is Solovay's first and second theorems for $\PA$; item
(2) is Solovay's theorem for $\ZF$. The metatheoretic assumption is needed
for the direction ``valid $\Rightarrow$ derivable'': in Solovay's
construction, from a formula not derivable in $\GL$ one builds a realization whose
unprovability in $\ZF$ is established under the assumption of
$1$-consistency, that is, $\Sigma_1$-soundness. In item~(2)
the interpretations of propositional variables are arbitrary closed
formulas of the language of $\ZF$, while Solovay's theorem is proved for arithmetical
realizations; this is no obstacle: the direction ``$\GL\vdash M\Rightarrow
\ZF\vdash M^*$'' holds for any substituted sentences (the derivability
conditions and the formalized L\"ob theorem do not depend on them), and
for the converse it suffices to have it for arithmetical realizations, which are
a special case. In $\sB$ the translation of a formula without atoms
$\quest\beta$ contains no set variables, and both closures in
Definition~\ref{def:valid} are empty. For details of Solovay's construction
see~\cite{Solovay76} and~\cite[Ch.~9]{Boolos93}.
\end{remark}

\begin{theorem}[derivable and valid $\square$-formulas are $\GL$ and $\GLS$]
\label{th:assval}
Let $M$ be a $\square$-formula.
\begin{enumerate}[leftmargin=2.2em,itemsep=1pt,label=\textup{(\arabic*)}]
\item The following three statements are equivalent: $\provgl M$;\
  $M$ is valid in $\sI$;\ $\GL\vdash M$.
\item The following three statements are equivalent: $\provgls M$;\
  $M$ is valid in $\sE$;\ $\GLS\vdash M$.
\item Suppose that atoms of the form $\quest\alpha$ occur in $M$ only outside the signs
  $\square$ (under the sign $\square$ in the semantics $\sB$ a problem variable
  would be bound and read as a closure---
  Remark~\ref{rem:bangbinds}). Then $\provglb M$ implies
  the validity of $M$ in $\sB$, $\GL\vdash M$ implies $\provglb M$, and
  under the assumption of $\Sigma_1$-soundness of $\ZF$ validity in
  $\sB$ implies $\GL\vdash M$, so that all three statements are equivalent.
\end{enumerate}
\end{theorem}

\begin{proof}
Each item is proved by going around a circle: from derivability in our calculus
to validity (the soundness theorem), from validity to
derivability in $\GL$ or $\GLS$ (Solovay's theorem), from derivability in
$\GL$ or $\GLS$ to derivability in our calculus (direct verification).
The last step is common to all items, so we begin with it.

\emph{$\GL\vdash M$ implies $\provgl M$.} Every $\GL$-derivation,
read on our atoms (propositional variables and formulas
$\quest\alpha$), turns into a derivation in $\HCgl$: the scheme K is
item~(1) of Example~\ref{ex:proofs}, scheme~4 is item~(2), L\"ob's scheme
is item~(3); classical tautologies are given by the schemes (A1); the necessitation
rule is admissible, since from $\provgl A$ by (AtoP) we get
$\provgl\bang A$, and by (PtoA) $\provgl\quest\bang A$, that is,
$\provgl\square A$. \emph{$\GLS\vdash M$ implies $\provgls M$:}
$\GLS$ consists of the theorems of $\GL$ together with all special cases of
reflection $\square N\ato N$, closed under (MP); the theorems of $\GL$
are derivable in $\HCgl\subseteq\HCgls$, reflection is the scheme (A8), and the necessitation
rule is not used here. \emph{$\GL\vdash M$ implies
$\provglb M$}, since $\HCgl\subseteq\HCB$.

\emph{Items (1) and (2).} $\provgl M$ implies the $\sI$-validity of
$M$, and $\provgls M$ implies $\sE$-validity (Theorem~\ref{th:sound}).
$\sI$-validity implies $\GL\vdash M$, and
$\sE$-validity implies $\GLS\vdash M$ (Lemma~\ref{lem:solzf}(1)).
The circle is closed.

\emph{Item (3).} $\provglb M$ implies $\sB$-validity by
Theorem~\ref{th:sound}. It remains to derive $\GL\vdash M$ from the
$\sB$-validity of $M$ under the assumption of $\Sigma_1$-soundness.
The formula $M$ is a Boolean combination of the variables $A_i$, the atoms
$\quest\alpha_1,\dots,\quest\alpha_m$ and propositions $\square B$, where
by hypothesis the atoms $\quest\alpha_j$ stand \emph{outside} the signs $\square$, and
the arguments $B$ contain no problem variables. Fix an
interpretation of the propositional variables. The translation $M^*$ is a Boolean
combination of the sentences $A_i^*$, the atoms $\exists y\,(y\in X_j)$ and
the sentences $\Prov(B^*)$ (the formulas $B^*$ are closed). The validity of $M$
means $\ZF\vdash\forall\vec X\,M^*$.

We show that this condition is equivalent to the following one:
\[
  \text{($\dagger$)}\qquad \ZF\vdash M^*[\vec\sigma]\ \ \text{for every
  tuple}\ \vec\sigma\in\{0{=}0,\ \arbot \}^m,
\]
where $M^*[\vec\sigma]$ is obtained from $M^*$ by replacing each atom
$\exists y\,(y\in X_j)$ with the sentence $\sigma_j$. Validity
implies ($\dagger$): substitute for $X_j$ the set $\varnothing$
when $\sigma_j=(\arbot )$ and $\omega$ when $\sigma_j=(0{=}0)$; after
the substitution the atom
$\exists y\,(y\in X_j)$ is $\ZF$-equivalent to $\sigma_j$, and replacement of
equivalents in the Boolean skeleton yields $M^*[\vec\sigma]$. Conversely, suppose
($\dagger$) holds. We argue in $\ZF$: by classical logic
$\bigwedge_j\bigl(\exists y\,(y\in X_j)\aor\aneg\exists y\,(y\in X_j)\bigr)$,
which gives a case distinction with $2^m$ cases; in each case every atom
is $\ZF$-equivalent to its $\sigma_j$, replacement of equivalents in the Boolean
skeleton yields $M^*[\vec\sigma]$, and the latter is provable by ($\dagger$). Thus
$\ZF\vdash M^*$, and with it the closure.

It remains to relate ($\dagger$) to $\GL$. Let $M[\vec\sigma]$ be
the result of replacing in $M$ each atom $\quest\alpha_j$ with $\aneg\bot$ when
$\sigma_j=(0{=}0)$ and with $\bot$ when $\sigma_j=(\arbot )$. This is
a $\square$-formula without atoms $\quest\beta$, and $M^*[\vec\sigma]$ is $\ZF$-equivalent
to its translation. By Lemma~\ref{lem:solzf}(2) the condition ($\dagger$)
is equivalent to $\GL\vdash M[\vec\sigma]$ for all $\vec\sigma$.
And this is equivalent to $\GL\vdash M$: the atoms $\quest\alpha_j$ stand outside
$\square$, so the formula
$M\aeq\bigvee_{\vec\sigma}\bigl(\bigwedge_j(\quest\alpha_j\aeq\sigma_j)
\aand M[\vec\sigma]\bigr)$ is a classical tautology, 
and from all $\GL\vdash M[\vec\sigma]$ we get $\GL\vdash M$; while
conversely the $M[\vec\sigma]$ are substitution instances of $M$.
\end{proof}

\begin{remark}[atoms $\quest\alpha$ under the sign $\square$]\label{rem:boxqatom}
The proviso of item~(3) is essential: if atoms $\quest\alpha$ stand under
the sign $\square$, validity in $\sB$ does not imply derivability in $\GL$.
Example: $\square\quest\alpha\ato\square\bot$
is valid in $\sB$, since $\ZF$ proves the implication
$(\square\quest\alpha)^*\ato(\square\bot)^*$: the argument of the inner
$\square$ is read as the closure $\forall X\,\exists y\,(y\in X)$,
which is refutable, and from a
proof of it together with a fixed refutation one obtains a proof of $\bot$. In
$\GL$, however, the formula $\square p\ato\square\bot$ is not derivable, and in $\sI$, $\sE$
the proposition $\square\quest\alpha\ato\square\bot$ is not valid: for
$\alpha^*:=(0{=}0)$ it turns into a formula provably equivalent to
the unprovable formula $\square\bot$. In $\HCB$, on the other hand, this proposition is
derivable, and the derivation shows how the calculus takes the hidden
quantifier into account: the propositional scheme $(\Pi\mathrm{E})$ with
$\theta:=\varnothing$ gives $\Pi\alpha\,\quest\alpha\ato\quest\varnothing$, that
is, $\Pi\alpha\,\quest\alpha\ato\bot$ by (A4d); the rule $(\bang\to)$
gives the problem $\bang\Pi\alpha\,\quest\alpha\pto\bang\bot$; the axiom
$(\Pi\bang)$ gives $\bang\quest\alpha\pto\bang\Pi\alpha\,\quest\alpha$; altogether
the problem $\bang\quest\alpha\pto\bang\bot$ is derivable, and the rule $(\quest\to)$
yields $\square\quest\alpha\ato\square\bot$. In $\HCI$ this $\Pi$-free
formula is not derivable: by item~(1) of the theorem it would be a theorem
of $\GL$. The formula is, of course, not clean (Remark~\ref{rem:clean}):
the variable $\alpha$ stands under the sign $\bang$; we return to the question of
conservativity of the $\Pi$-extension in
Remark~\ref{rem:qfree}.
\end{remark}

\begin{corollary}[completeness on the propositional side]\label{cor:solovay}
For every $\Pi$-free proposition $A$ the following statements are equivalent:
\begin{itemize}
\item $\provgl A$,
\item $A$ is $\sI$-valid,
\item $\GL\vdash\nf{A}$.
\end{itemize}
The following statements are also equivalent:
\begin{itemize}
\item
$\provgls A$,
\item  $A$
is $\sE$-valid,
\item $\GLS\vdash\nf{A}$.
\end{itemize}
In particular, the $\square$-fragment of $\HCgl$ is exactly $\GL$, and
the $\square$-fragment of $\HCgls$ is exactly $\GLS$.
\end{corollary}

\begin{proof}
By Proposition~\ref{prop:nf} the derivability and the validity of $A$
are equivalent to the derivability and the validity of the $\square$-formula $\nf{A}$;
Theorem~\ref{th:assval} does the rest.
\end{proof}

Thus, on the side of $\Pi$-free propositions our semantics gives
nothing new compared with provability logic: two-sortedness and
the modalities add not a single new scheme to $\GL$ and $\GLS$. All the
novelty is concentrated on the problem side, and for $\sB$ also in
propositions with the quantifier $\Pi$ (\S\ref{sec:pigain}, \S\ref{sec:pinf}).

\subsection{Normal form in the $\Pi$-extension}\label{sec:pinf}

In the $\Pi$-free language the normal form
(Proposition~\ref{prop:nf}) is built from the variables $A_i$ and
the propositions $\quest\alpha$ and $\square B$. In the language with the quantifier $\Pi$
the reduction by (A4) stops at propositions
$\quest\Pi\vec\alpha\,\chi$ (the remark after Proposition~\ref{prop:nf}), but
the laws of $\Pi$ allow one to simplify these as well. By an \emph{atom} we mean a
formula at which such simplifications come to an end. This notion is
syntactic: it is defined by the set of equivalences listed below,
and an atom may turn out to be equivalent to a simpler
formula. In this subsection we consider formulas in which the quantifier
$\Pi$ is applied only to problems; propositions of the form $\Pi\alpha\,A$ are
not simplified here.

\begin{definition}[$\Pi$-atom]\label{def:piatom}
A \emph{$\Pi$-atom} is a proposition
$\quest\Pi\alpha_1\ldots\Pi\alpha_k\,\chi$, $k\geqslant1$, in which every
variable $\alpha_i$ occurs free in the problem $\chi$, and the problem $\chi$ is either
a variable, or a problem of the form $\bang A$, or an implication
$\varphi\pto\psi$ in whose premise $\varphi$ every
$\alpha_i$ occurs free.
\end{definition}

\begin{proposition}[normal form in the $\Pi$-extension]\label{prop:pinf}
Every formula in which the quantifier $\Pi$ is applied only to
problems is equivalent in $\HCB$ to a formula in
which $\quest$ occurs only in the combinations $\quest\alpha$ and
$\square B$ and in $\Pi$-atoms.
\end{proposition}

\begin{proof}
As in Proposition~\ref{prop:nf}, but to the equivalences (A4a)--(A4d)
we add the equivalences $(\Pi\pand)$, $(\Pi\por)$ (an axiom and
Lemma~\ref{lem:piderived}(4)) and the equivalences of
Lemma~\ref{lem:piderived}; all replacements are legitimate by Lemma~\ref{lem:repl}. Let
$\quest\Pi\vec\alpha\,\chi$ be an occurrence that is not a $\Pi$-atom.
A quantifier over a variable that does not occur free in $\chi$ is dropped
(Lemma~\ref{lem:piderived}(1)). If $\chi$ is a conjunction or a disjunction,
we move the quantifiers into the members by $(\Pi\pand)$, $(\Pi\por)$, and $\quest$
passes through the connective by (A4b) or (A4c). If $\chi=\varphi\pto\psi$ and
some $\alpha_i$ does not occur in the premise, we permute its quantifier
to the inside of the block (Lemma~\ref{lem:piderived}(3)) and move it into the conclusion
(Lemma~\ref{lem:piderived}(2)); if no quantifiers over the implication remain,
$\quest$ passes through it by (A4a). The process terminates: every step, except
the permutation of quantifiers, decreases in the lexicographic order the triple
consisting of the number of problem connectives standing under the symbols $\quest$, the total
number of problem connectives in the scopes of the quantifiers $\Pi$ standing under
the symbols $\quest$, and the number of such quantifiers; and a permutation is performed only
immediately before a step that decreases the triple.
\end{proof}

The atoms of the $\Pi$-free normal form, the variables $A_i$ and
the propositions $\quest\alpha$, are unsimplifiable also in the semantic sense: they
can independently be assigned arbitrary truth values
(Corollary~\ref{cor:assfrag}). For $\Pi$-atoms this is not so.

\begin{definition}[eliminable $\Pi$-atom]\label{def:elim}
A $\Pi$-atom $\theta$ is \emph{eliminable} if for some
$\Pi$-free proposition $A$ the equivalence $\theta\aeq A$
is $\sB$-valid.
\end{definition}

Let us note at once that quantifier elimination in problems
(Proposition~\ref{prop:pienc}) also yields eliminability of atoms: if
the problem $\Pi\vec\alpha\,\chi$ is $\sB$-equivalent to a $\Pi$-free problem
$\psi$, then the proposition $\quest\Pi\vec\alpha\,\chi$ is equivalent to
$\quest\psi$, and the latter is reduced by (A4) to a $\Pi$-free normal
form; in this way many atoms that look irreducible at first glance are eliminable,
for example $\quest\Pi\alpha((\alpha\pto\beta\por\gamma)\pto
((\alpha\pto\beta)\por(\alpha\pto\gamma)))$, which is equivalent to
$(\quest\gamma\ato\quest\beta)\aor(\quest\beta\ato\quest\gamma)$, that is,
to a classical tautology.

All $\Pi$-atoms encountered so far are eliminable:
$\quest\Pi\alpha\,\alpha$ and $\quest\Pi\alpha(\pneg\pneg\alpha\pto\alpha)$
are equivalent to $\bot$ (Theorem~\ref{prop:piatoms}),
$\quest\Pi\beta\bigl((\bang A\pto\beta)\pto\beta\bigr)$ is equivalent to
$\square A$ (Proposition~\ref{prop:pienc}(2) with $\varphi:=\bang A$), and $\quest\Pi\beta\Pi\gamma\rho$
is valid (Example~\ref{ex:pider}) and hence equivalent to $\aneg\bot$.
The atom $\quest\Pi\alpha(((\alpha\pto\beta)\pto\alpha)\pto\alpha)$ is eliminable too
---it is $\sB$-equivalent to $\quest\beta$---but this equivalence has been proved
only semantically, and whether it is derivable in $\HCB$ is unknown
(Proposition~\ref{prop:pienc}(4), Remark~\ref{rem:peirce}). Whether
ineliminable $\Pi$-atoms exist, we do not know (\S\ref{sec:todo}).

\subsection{Modality-free $\Pi$-free valid problems}\label{sec:modalfreeval}

Let us return to the simplest class of problems---problems without modalities,
that is, those built from variables and the constant $\varnothing$
by the connectives $\pto$, $\pand$, $\por$. For each semantics
$Z\in\{\sI,\sB,\sE\}$ we denote by $\Lp_Z$ the set of modality-free
problems valid in $Z$ (note that by their grammar
modality-free problems do not contain the quantifier $\Pi$ either).
By
Lemmas~\ref{lem:closure} and~\ref{lem:closurebs} all three sets are
closed under (MP) and under substitution of modality-free problems (such
substitutions affect neither $\bang$ nor $\Pi$ and are therefore safe), that
is, they are
\emph{superintuitionistic logics} in the language
$\pto,\pand,\por,\varnothing$ (the constant $\varnothing$ plays the role of falsity):
they contain $\Int$ by Proposition~\ref{prop:intval}. Now that the
soundness theorem has been proved, every fact of non-validity also gives a fact of
non-derivability; we shall refer to this in~\S\ref{sec:indep}.

At the same time, none of the three logics contains classical schemes.

\begin{proposition}[the problems have remained intuitionistic]\label{prop:noem}
Let $\alpha$, $\beta$ be variables of sort problem. The problems
$\alpha\por\pneg\alpha$, $\pneg\alpha\por\pneg\pneg\alpha$,
$\pneg\pneg\alpha\pto\alpha$ and $((\alpha\pto\beta)\pto\alpha)\pto\alpha$ are not
valid in any of the three semantics.
\end{proposition}

\begin{proof}
It suffices to prove the non-validity of the second problem---the weak law of
excluded middle $\pneg\alpha\por\pneg\pneg\alpha$. Indeed,
the first problem implies the second in $\Int$, and the third and the fourth are schemes each
of which together with $\Int$ gives the whole of classical logic; since
validity is preserved under substitution of problems for problem variables and
under the rule (MP) (Lemmas~\ref{lem:closure} and~\ref{lem:closurebs}), the
validity of any of them would imply the validity of all classical
tautologies, including the second problem.

In the semantics $\sI,\sE$: we fix a computably enumerable undecidable set
$X\subset\N$ and express the property $n\in X$ by an arithmetical
$\Sigma$-formula $\exists x A(x,n)$.
Arguing by contradiction,
suppose there exists a solution
$S$ of the problem $\pneg\alpha\por\pneg\pneg\alpha$.
For an arbitrary $n\in\N$  consider the interpretation $*_n$ that interprets the variable $\alpha$ by the formula $A(n,x)$.
Then the solutions of the problem $(\pneg\alpha\por\pneg\pneg\alpha)^{*_n}$ are all numbers of the form
$c(0,m)$ (and only they)
if $n\notin X$, and all numbers $c(1,m)$ (and only they)
if $n\in X$. And this fact is provable in the external theory.
Thus the external theory proves
$n\in X\aeq \exists m ( S(*_n)=c(1,m))$, and hence this statement is true for all $n$.
Therefore, by means of the algorithm $S$ applied to the code of the interpretation
$*_n$, one can decide whether $n$ belongs to $X$. A contradiction.

\remove{
--- её решение. По
лемме~\ref{lem:clauses}(5),(6) элемент множества
$\val{\alpha\por\pneg\alpha}$ имеет вид $c(0,u)$ с $u\in\val{\alpha}$, если
$\val{\alpha}\neq\varnothing$, и вид $c(1,v)$, если $\val{\alpha}=\varnothing$
(при непустом $\val{\alpha}$ множество $\val{\pneg\alpha}$ пусто). В
семантиках $\sI$, $\sE$ по
метке числа $S(\gn{\alpha^*})$ вычисляется истинностное значение предложения
$\exists x\,\alpha^*(x)$. Но уже для $\Sigma_1$-формул $\alpha^*$ это
невозможно: получился бы разрешающий алгоритм для проблемы остановки.}
In
the semantics $\sB$ it is even simpler: there is no code, the solution is a single fixed
number $c(i,u)$, and for $i=0$ it is rejected by the box $\varnothing$, while for
$i=1$ by the box $\{0\}$, and $\ZF$ proves the negation
of the translation for such a box.
\end{proof}

Thus all three logics $\Lp_Z$ lie between $\Int$ and
classical logic, without coinciding with the latter. Neither do they coincide with $\Int$,
and a separating example was found by Rose~\cite{Rose53}---as
a propositional formula realizable under every substitution but not
intuitionistically derivable; it was the first to show that the propositional
logic of realizability is strictly wider than $\Int$. In our semantics it is obtained as
a special case of a general statement about \emph{surveyable} problems (Theorem~\ref{prop:obz}
below).

\subsubsection{Surveyable problems}
Let us fix an interpretation $*$ (in the semantics $\sB$ we argue
inside $\ZF$ for given contents of the boxes). By \emph{choice programs}
for a problem $\varphi$ under a given interpretation we mean the elements of the set
$\val{\pneg\pneg\varphi\pto\varphi}$, that is, the solutions of the problem
$\pneg\pneg\varphi\pto\varphi$ under this interpretation. By
Lemma~\ref{lem:clauses}(6) they are structured as follows:
if $\val{\varphi}\neq\varnothing$, then $\val{\pneg\pneg\varphi}=\N$, the input
conveys nothing, and a choice program is a total program with
values in $\val{\varphi}$; if, however, $\val{\varphi}=\varnothing$, then
$\val{\pneg\pneg\varphi}=\varnothing$ and a choice program is
\emph{any} program, including the nowhere defined one.

\begin{definition}[surveyable problem]\label{def:obz}
A problem $\varphi$ is called \emph{surveyable} if from the code $\gn{*}$ one can compute
a finite list of programs $q_1,\dots,q_k$ (both their number and the programs themselves may
depend on the code) with the property: if $\val{\varphi}\neq\varnothing$, then at
least one of $q_1,\dots,q_k$ is total and all its values lie in
$\val{\varphi}$. (Recall that in $\sB$ the code is the code of the transparent
part, and the property of the list is required for arbitrary contents of the boxes.)
\end{definition}

In other words, surveyability means that there are finitely many candidates for choice programs
and they are written out in advance from the code of the interpretation. For example,
the problem $\pneg\beta\por\pneg\gamma$
is surveyable with
a constant list of two programs---``output $c(0,0)$'' and ``output
$c(1,0)$''.

\begin{theorem}[Rose's scheme for surveyable problems]\label{prop:obz}
Let a problem $\varphi$ be surveyable and $\psi:=\pneg\pneg\varphi\por\pneg\varphi$.
Then the problem
\begin{align}\tag{R}\label{rose-formula}
  \bigl((\pneg\pneg\varphi\pto\varphi)\pto \psi\bigr)\pto \psi
\end{align}
is valid in the given semantics as soon as the defining property
of the list from Definition~\ref{def:obz} is, for every $*$, provable in its
external theory. In $\sE$ this condition is always fulfilled (the external theory
is $\TA$, and the property of the list is true by definition), so that there
the problem~\eqref{rose-formula} is valid unconditionally.
\end{theorem}

Before the proof we recall that 
Rose's formula
$\rho$ is the problem \eqref{rose-formula} for the (surveyable) formula $\phi=\pneg\beta\por\pneg\gamma$.

\begin{proof}
Let us first work out what a solution of the problem $\psi$ is. By Lemmas~\ref{lem:clauses}(5)
and~\ref{lem:dn} the set $\val{\pneg\pneg\varphi}$ equals $\N$ when
$\val{\varphi}\neq\varnothing$ and is empty otherwise, while $\val{\pneg\varphi}$ is
the other way round. Therefore a solution of the problem $\psi$ is \emph{a label alone}, and moreover
a correct one: when $\val{\varphi}\neq\varnothing$ the solutions are exactly the numbers with
label~$0$, and when $\val{\varphi}=\varnothing$, exactly the numbers with
label~$1$. In particular, $c(0,0)$ will do in the first case and $c(1,0)$ in the second.

A solution of the premise is a program $p$ which, on every choice program for
$\varphi$, halts and outputs a solution of the problem~$\psi$. From it we construct a solution
of the conclusion.

From the code $\gn{*}$ we compute the list $q_1,\dots,q_k$ of
Definition~\ref{def:obz} and run $p$ in parallel on all $q_i$. Then:
\begin{itemize}[leftmargin=1.6em,itemsep=1pt]
\item as soon as some run has halted with an answer having label~$0$,
  we halt and output $c(0,0)$;
\item if all $k$ runs have halted with answers having label~$1$,
  we output $c(1,0)$.
\end{itemize}

\emph{Verification.} Let $\val{\varphi}=\varnothing$. Then every program
is a choice program, in particular all the $q_i$; hence $p$ halts on each of them,
and all answers have label~$1$---there are no other solutions of the problem $\psi$
when $\val{\varphi}$ is empty. Therefore the first clause will not fire, all $k$
runs will terminate, and we shall output $c(1,0)$---a correct solution.

Now let $\val{\varphi}\neq\varnothing$. By the definition of surveyability
some $q_i$ is a choice program, hence $p$ halts on it and
outputs a solution of the problem~$\psi$; when $\val{\varphi}$ is nonempty this is a number with label~$0$.
Consequently the first clause will fire sooner or later, and we shall output $c(0,0)$---
again a correct solution.

The list $q_1,\dots,q_k$ is computed from the code, and the rest of the construction does not
depend on the code, so the solving function is computable. The provability of the formula
$\varphi^*(S(\gn{*}))$ in the external theory
is obtained by the same argument carried out inside it: the case
distinction on
$\exists x\,\varphi^*(x)$ is classical, and everything else is
a $\Sigma_1$-argument about the work of a machine; in $\sB$ the argument is carried out
over free set variables, and the list is computed from the code of the
transparent part. The only place where an
additional condition is required is the step ``if $\exists x\,\varphi^*(x)$, then
some $q_i$ is a choice program'': inside the external theory it can be
used only if it has been proved there (for $\TA$ this is just
truth; for $\PA$ and $\ZF$ truth is not enough).
\end{proof}

The class of surveyable problems is fairly wide; to describe it explicitly, we single out
one more class, given grammatically.

\begin{definition}[$\mathcal O$-problems]\label{def:obzgram}
\emph{$\mathcal O$-problems} are the problems built by the grammar
\[
  \varphi ::= \sigma \mid \pneg\psi \mid (\varphi\pand\varphi) \mid
              (\varphi\por\varphi) ,
\]
where $\sigma$ ranges over $\Sigma$-problems and $\psi$ over arbitrary problems.
In other words, all possible negations are added to the $\Sigma$-problems, and all this
is closed under $\pand$ and $\por$. Every $\Sigma$-problem is an
$\mathcal O$-problem; the converse fails already for $\pneg\bang A$.
\end{definition}

\begin{lemma}[the stock of surveyable problems]\label{lem:obz}
\begin{enumerate}[leftmargin=2.2em,itemsep=1pt,label=\textup{(\arabic*)}]
\item Every problem $\phi$ for which the problem $\pneg\pneg\varphi\pto\varphi$ is valid
  is surveyable. In particular, all $\Sigma$-problems and all problems of the form
  $\pneg\psi$ are surveyable.
\item If $\varphi$ and $\psi$ are surveyable, then so are $\varphi\pand\psi$ and
  $\varphi\por\psi$.
\end{enumerate}
All of this applies to each of the three semantics, and in all cases
the defining property of the list is provable in the external theory.
\end{lemma}

\begin{proof}
(1) Let $S$ be a computable function solving
$\pneg\pneg\varphi\pto\varphi$. Take the list consisting of a single program:
$q:=$ ``on any input output $U_{S(\gn{*})}(0)$''. If
$\val{\varphi}\neq\varnothing$, then $U_{S(\gn{*})}$ is total with
values in $\val{\varphi}$, hence $q$ is total and its values lie
in $\val{\varphi}$. For $\Sigma$-problems the scheme
is valid in all three semantics---this is the scheme (A7),
Proposition~\ref{prop:a5a7};
for $\pneg\psi$ it is derivable already in $\Int$ and valid by
Proposition~\ref{prop:intval}.

Provability in the external theory: for $\Sigma$-problems it is part of
the validity of (A7) by Definition~\ref{def:valid}---the external theory
proves
$\forall u\,((\pneg\pneg\varphi)^*(u)\ato\varphi^*(U_{S(\gn{*})}(u)))$, and
the premise is equivalent to $\exists x\,\varphi^*(x)$ by Lemma~\ref{lem:dn};
for $\pneg\psi$ it is part of the validity of a scheme of $\Int$.

(2) Let $q_1,\dots,q_k$ be a list for $\varphi$, and $q'_1,\dots,q'_m$
a list for $\psi$. For $\varphi\por\psi$ take the list of $k+m$ programs: ``run
$q_i$ and output $c(0,\text{result})$'' and ``run $q'_j$ and output
$c(1,\text{result})$''. If $\val{\varphi\por\psi}\neq\varnothing$, then
$\val{\varphi}$ or $\val{\psi}$ is nonempty, and the corresponding wrapped
program is total with values in $\val{\varphi\por\psi}$. For
$\varphi\pand\psi$ take the $km$ programs ``run $q_i$ and $q'_j$ and output
$c(\text{results})$''; if $\val{\varphi\pand\psi}\neq\varnothing$, then
both factors are nonempty, and a suitable pair exists. Both arguments are formalized
in the external theory word for word, so the provability of the property of the list
is preserved.
\end{proof}

\begin{proposition}[Rose's scheme is valid on the class $\mathcal O$]\label{cor:rose}
For every $\mathcal O$-problem $\varphi$ the problem~\eqref{rose-formula}
is valid  in all three semantics.
\end{proposition}

\begin{proof}
Induction on the grammar of Definition~\ref{def:obzgram} shows that
every $\mathcal O$-problem is surveyable, and moreover the defining property of the list
is provable in the external theory: for $\Sigma$-problems and problems of the form $\pneg\psi$ this is
clause~(1) of Lemma~\ref{lem:obz}, and the connectives $\pand$ and $\por$ preserve
surveyability by its clause~(2). It remains to apply
Theorem~\ref{prop:obz}.
\end{proof}

We did \emph{not} include Rose's scheme among the axioms of the calculi,
although it is
valid in all three semantics. The calculi are composed so as
to facilitate proofs of validity, and for this the axioms
must be natural principles with a simple formulation---such are the axioms
(A1)--(A7).
Rose's scheme is not a
principle but a ready-made result: it is long, the class of $\mathcal O$-problems
is not simple to define, it apparently does not give completeness on the problem side,
and it is hard to imagine anyone using it for
proving the validity of anything else. Therefore we
record it among the examples of \emph{incompleteness}: Rose's formula is valid
but not derivable in $\Int$, and hence not in $\HCI$ either  (\S\ref{sec:indep}).

Surveyability in Theorem~\ref{prop:obz} is essential: if $\varphi$ is
a variable, the scheme~\eqref{rose-formula} is not valid. In $\Int$ it is
equivalent to the scheme
\[
  \bigl((\pneg\pneg\varphi\pto\varphi)\pto(\varphi\por\pneg\varphi)\bigr)\pto
  (\pneg\pneg\varphi\por\pneg\varphi),
\]
which differs in that the conclusion of the premise is
$\varphi\por\pneg\varphi$ rather than $\pneg\pneg\varphi\por\pneg\varphi$:
the premises are equivalent. Indeed,
$(\varphi\por\pneg\varphi)\pto(\pneg\pneg\varphi\por\pneg\varphi)$, since
$\varphi\pto\pneg\pneg\varphi$; conversely, under the assumption
$\pneg\pneg\varphi\pto\varphi$, from $\pneg\pneg\varphi\por\pneg\varphi$
we get $\varphi\por\pneg\varphi$. This scheme is called \emph{Scott's
axiom}: Kreisel and Putnam~\cite{KP57} mention that D.~Scott proved
the disjunction property of the logic obtained by adding it to $\Int$.
Rose's formula is its special case for $\varphi:=\pneg\beta\por\pneg\gamma$.
We shall call this problem \emph{Scott's principle}. Yankov~\cite{Yankov63}
proved that Scott's principle with a variable is not Kleene realizable;
the following theorem is an analogue of this result for our semantics (the
proof, as in Yankov's, uses diagonalization over the number of
computation steps).

\begin{theorem}[Scott's principle is not valid]\label{prop:scott}
The problem
\[\phi:=
  \bigl((\pneg\pneg\alpha\pto\alpha)\pto(\pneg\pneg\alpha\por\pneg\alpha)\bigr)\pto
  (\pneg\pneg\alpha\por\pneg\alpha),
\]
where $\alpha$ is a variable of sort problem (that is, the scheme~\eqref{rose-formula}
for $\varphi:=\alpha$), is not valid in any of the three semantics.
\end{theorem}


\begin{proof}
Let $S$ be an arbitrary algorithm. We need to define
an interpretation $*$ for which the formula $\phi^*(S(*))$ is unprovable or false.
As was already said in Remark~\ref{rem:assumptions},
we shall prove the provability of the negation of this formula in the external theory.
By the letter $X$ we denote the value $\val{\alpha}$: in the semantics $\sB$ this is
the contents of the box, in the others it is the set defined by the formula
$\alpha^*$.

\emph{The semantics $\sB$.}
Here the proof is simpler: there are no propositional variables in the problem $\phi$,
the code of the interpretation is empty, and the solver is a single number. So we
begin with it.

Let us describe what the solutions of the subproblems occurring in the problem are.

By Lemma~\ref{lem:clauses}(6), a solution of the problem $\pneg\pneg\alpha\pto\alpha$ is,
for $X\neq\varnothing$, the number of a total function with values
in~$X$ (the set $\val{\pneg\pneg\alpha}$ then equals $\N$, so the input
conveys nothing), and for $X=\varnothing$ it is an arbitrary number. We call such a
number $i$ \emph{$X$-correct}.

A solution of the problem $\pneg\pneg\alpha\por\pneg\alpha$ is a number $c(0,v)$
(if $X\neq\varnothing$) or a number $c(1,v)$ (if $X=\varnothing$), where $v$
is arbitrary: the value of each of the disjuncts, when it is nonempty, is
$\N$. Thus a solution is essentially just a tag, that is, the answer
to the question whether $X$ is empty; we shall say that the program answers ``nonempty''
or ``empty''.

So, a solution of the premise is a program $p$ which, on every
$X$-correct number $i$, halts and correctly answers whether
$X$ is empty.

The code 
of the interpretation $*$ is empty, hence
$S(\gn{*})$ is some
fixed program not depending
on the box. Put $r:=S(\gn{*})$

\medskip
\emph{The diagonal construction.}
By
Kleene's fixed-point theorem there is a program $p$ which on input
$i$ does the following:
\begin{enumerate}[leftmargin=2.2em,itemsep=1pt,label=\textup{(\arabic*)}]
\item it simulates the run of the program $r$ on input $p$, that is, on its
  own number, for exactly $i$ steps;
\item if during this time $r$ has halted with the answer ``empty'', it answers
  ``nonempty'';
\item otherwise it answers ``empty''.
\end{enumerate}
The bound on the number of steps is essential in two respects. First, the simulation is finite,
and no looping arises, although $r$ may along the way simulate
$p$ itself. Second---and this is the main point---if $r$ halts with the
answer ``empty'', then the inputs on which $p$ answers ``empty'' form
a \emph{finite} initial segment. (The program $p$ uses the input $i$ itself
only as the number of steps.)

With $r$ we associate a set $X_r$: if $r$ on input $p$ halts at step
$s$ with the answer ``empty'', we put $X_r:=\{\nu(s)\}$, where $\nu(s)$ is
the least number not occurring among $U_i(0)$ for $i<s$ (there are
finitely many such numbers, so $\nu(s)$ is defined; computability of $\nu$ will not be
needed, arithmetical definability suffices); in all other
cases we put $X_r:=\varnothing$.

\medskip
\emph{Case analysis.} Let us check: if $X=X_r$, then $p$ is a solution of the premise,
while $r$ on input $p$ diverges or outputs a wrong tag.

\emph{Case 1: $r$ on input $p$ diverges.} Then for every $i$ step~(2) does not
fire, and $p(i)$ answers ``empty''. Hence $p$ is total and always
answers ``empty''. Moreover $X=\varnothing$, every number is $X$-correct, and
the answer ``empty'' is a correct solution of the problem $\pneg\pneg\alpha\por\pneg\alpha$. Therefore
$p$ is a solution of the premise, and $r$ diverges on it.

\emph{Case 2: $r$ on input $p$ halts with the answer ``nonempty''.}
Step~(2) again does not fire for any $i$, so $p$ still everywhere
answers ``empty'' and is a solution of the premise for $X=\varnothing$. But then
the correct answer is ``empty'', whereas $r$ answered otherwise.

\emph{Case 3: $r$ on input $p$ halts at step $s$ with the answer
``empty''.} Then $X=\{\nu(s)\}$ is a nonempty set. Let us check that $p$ is
a solution of the premise. Let $i$ be an $X$-correct number; then $U_i$ is
total and $U_i(0)=\nu(s)$. If we had $i<s$, then $\nu(s)$ would occur
among $U_i(0)$ for $i<s$, which is impossible by the choice of $\nu(s)$. Hence
$i\geq s$; when simulating $i\geq s$ steps the program $p$ sees the halting of $r$
with the answer ``empty'', step~(2) fires, and $p$ answers ``nonempty''
---correctly, since $X$ is nonempty. So,
$p$ is a solution of the premise, the set $X$ is nonempty, and $r$ answered ``empty''.

In all three cases $p$ is a solution of the premise for $X=X_r$, and $r$ on it
errs or diverges.


Take the set $X_r$ as the contents of the box: by
what has been proved, $p$ solves the premise, and $r$ errs or diverges on it, that
is, $r$ does not belong to the value of the problem for these contents of the box.
This completes the proof of the negation of the formula
$\phi^*(S(*))$ in set theory.

\emph{The semantics $\sE$ and $\sI$.}
This time we need to start by constructing a formula $\alpha^*(x)$,
which we define by means of the fixed-point theorem for PA.
Again
put $r:=S(\gn{\alpha^*})$. Now the program $r$ is built by the algorithm $S$ from the 
formula $\alpha^*$ which we must exhibit, and we must exhibit it so
that the value $\val{\alpha}$ coincides with $X_r$. We get a vicious circle,
and the existence of $\alpha^*$ is ensured by the fixed-point lemma for
$\PA$: there is a formula $\alpha^*(x)$ for which $\PA$ proves
\[
  \alpha^*(x)\ \aeq\ \exists s\,\bigl(\text{``$
  S(\gn{\alpha^*})
  $ on input
  $p$ halts at step $s$ with the answer `empty'\,''}\ \aand\ x=\nu(s)\bigr),
\]
where $p$ is built from $r=S(\gn{\alpha^*})$ by the Diagonal construction given above. The lemma is applicable: the condition in
parentheses is arithmetical, and $\nu$ is arithmetically definable.
By construction the formula $\alpha^*(x)$ is equivalent in PA
to the formula $x\in X_r$, so the case
analysis refutes the formula $\phi^*(S(*))$ in the external theory.
\end{proof}

Let us note in passing that the class $\mathcal O$ in Rose's scheme can neither be extended to
all problems nor substantially narrowed: with a variable in place of
$\varphi$ the scheme implies Scott's principle, which is not valid
(Theorem~\ref{prop:scott}); on $\Sigma$-problems it is derivable from (A7)
by a single (MP), so that all its content rests on the negations.



\subsubsection{The positive modality-free fragment}\label{sec:posfrag}

A modality-free problem is called \emph{positive} if it does not contain
$\varnothing$ (and hence does not contain $\pneg$ either). The \emph{positive fragment} of a logic is the set of its positive
formulas; a logic \emph{has an intuitionistic positive fragment} if its
positive fragment coincides with the positive fragment of $\Int$. We shall prove
that the positive fragments of all three logics $\Lp_Z$ are
intuitionistic: the whole difference between these logics and $\Int$ lives on formulas
with~$\varnothing$.

\begin{definition}[Medvedev]\label{def:crit}
A \emph{critical implication} is a positive formula of the form
$\bigwedge_i((P_i\pto Q_i)\pto Q_i)\pto R$, where the $P_i$ are conjunctions
of variables, the $Q_i$ and $R$ are disjunctions of variables, and moreover for each $i$
the formulas $P_i$ and $Q_i$ have no variables in common, while each of $P_i$, $Q_i$, $R$
contains at least one variable. For each $n>0$ put
\[
  J_n:=\bigwedge_{\varnothing\neq E\subsetneq\{1,\dots,n\}}
    \Bigl(\bigl(\bigwedge_{j\notin E}\alpha_j\pto\bigvee_{i\in E}\alpha_i\bigr)
      \pto\bigvee_{i\in E}\alpha_i\Bigr)\pto\bigvee_{i=1}^n \alpha_i
\]
---the weakest critical implication in the problem variables
$\alpha_1,\dots,\alpha_n$; in particular $J_1=\alpha_1$.
\end{definition}

\begin{theorem}[Medvedev, 1962;~\cite{Medvedev62}]\label{th:medv}
A superintuitionistic logic has an intuitionistic positive fragment if and
only if no critical implication is derivable in it,
and this is equivalent to the condition that none of the formulas $J_n$ is derivable in it.
\end{theorem}

\begin{theorem}[Yankov, 1968;~\cite{Yankov68}]\label{th:jan}
A superintuitionistic logic has an intuitionistic positive fragment if and
only if it is contained in $\J:=\Int+(\pneg\alpha\por
\pneg\pneg\alpha)$---the logic of weak excluded middle.
\end{theorem}

Medvedev's criterion is convenient for logics given semantically (and ours is given
exactly in this way), Yankov's criterion---for logics given axiomatically. We use
them in tandem: the first yields an intuitionistic positive fragment, the second
turns this into the inclusion $\Lp_Z\subseteq\J$.

The logic of realizable propositional formulas is arranged as follows: variables
are replaced by arithmetical formulas, connectives---by Kleene's operations, and
a formula is considered to belong to the logic if under every substitution it has
a realization. Three variants are distinguished~\cite{Plisko09}: $\R$---a realization merely
exists; $\Reff$---a realization can be \emph{found effectively} from the
substituted formulas; $\Rcnst$---one and the same realization works under all
substitutions. Obviously $\Rcnst\subseteq\Reff\subseteq\R$.

The three levels of Remark~\ref{rem:three} are exactly these. Our
notion of validity is the middle one: a solution is produced by a computable function from
the code of the interpretation. The only difference from $\Reff$ is that for us the value of a
variable is an arbitrary arithmetically definable set, whereas there it is
the set of realizations of an arithmetical formula.

For the realizability logics $\R$, $\Reff$, $\Rcnst$ the intuitionisticity of the
positive fragment is known (Medvedev, Plisko;~\cite{Plisko09}), but
we cannot refer to it here. It might seem that the
inclusion $\Lp_{\sE}\subseteq\R$ suffices---the condition $\exists S\,\forall{*}$ is stronger than
$\forall{*}\,\exists$---and then the positive fragment of $\Lp_{\sE}$ would be squeezed
between the coinciding positive fragments of $\Int$ and $\R$. But the semantics do not
coincide: for us the value of a variable is an arbitrary arithmetically
definable set, whereas in the definition of $\R$ it is the set of realizations of an
arithmetical formula, and \emph{not every arithmetical set is a
set of realizations}. Therefore we take the direct route: by
Theorem~\ref{th:medv} it suffices to check that no formula $J_n$ is
valid. The required construction is in~\cite{VSSCh}; the operations on sets
there are exactly the same as in our Definition~\ref{def:val} (labels $0$ and~$1$ for
disjunction, pairs for conjunction, $e(x)\in Y$ for implication), so that
the transfer reduces to closing the circular dependence of the solver on the code
of the interpretation.

We say that a program \emph{co-enumerates} a set if it enumerates its
complement.

\begin{lemma}[\cite{VSSCh}]\label{lem:vssch}
Let there be given a natural number $m\ge 1$, a program $q$ enumerating a set $M$ of
at most $2^m$ numbers, and a critical implication $J_n$ in $n$ variables.
Then one can effectively (that is, by some total computable
function of $m$, $n$ and $q$) find programs $a_1,\dots,a_n$ co-enumerating
nonempty sets $A_1,\dots,A_n$ for which
\[
  M\cap J_n(A_1,\dots,A_n)=\varnothing .
\]
Moreover, all elements of the sets $A_j$ do not exceed $2^{Cm}$, where the constant $C$
depends only on $n$.
\end{lemma}

The proof of this lemma does not go beyond PA.
The construction in a nutshell: the sets $A_j$ are initially taken to be
finite segments of suitable sizes, and then elements are removed from them step by step---every
time the next program $e\in M$ reveals itself
on suitable arguments, its answer is removed from the corresponding $A_j$.
Simultaneously, by the fixed-point theorem, auxiliary
programs are constructed witnessing membership in the premises
$((P_i\pto Q_i)\pto Q_i)$, and the main work is the estimate of the number of removals, which
shows that the process terminates and the sets do not become empty. It is exactly
here that the condition from the definition of a critical implication, that $P_i$ and
$Q_i$ have no variables in common, is used.

\begin{theorem}[the positive fragment is intuitionistic]\label{th:posfrag}
A positive problem is valid (in any of the three semantics) if and
only if it is derivable in~$\Int$.
\end{theorem}

\begin{proof}
One direction is Proposition~\ref{prop:intval}. Conversely, by
Theorem~\ref{th:medv}
it suffices to check that no formula $J_n$ is valid---in
none of the three semantics.

\emph{The semantics $\sB$---the simplest case.} There are no propositional
variables in $J_n$, so the code of the transparent part is empty and the solver
reduces to a single number $e_0$. Apply Lemma~\ref{lem:vssch} with
$m:=1$ to the program ``print $e_0$ and halt'': we obtain nonempty
co-enumerable sets $A_1,\dots,A_n$ with
$e_0\notin J_n(A_1,\dots,A_n)$. Let the contents of the boxes be
$X_{\alpha_j}:=A_j$; by Lemma~\ref{lem:clauses} the values of the connectives are
exactly the operations from Lemma~\ref{lem:vssch}, so that
$\val{J_n}=J_n(A_1,\dots,A_n)$ and the statement $e_0\notin\val{J_n}$ is provable in PA, and all the more so in ZF.

\emph{The semantics $\sE$ and $\sI$.} Here the solver receives a code, and
self-reference is needed.
Let $S$ be a total computable function, presumably solving $J_n$.
This
means: the formula $J_n^*(S(\gn{*}))$ is provable or true under every interpretation~$*$.

\emph{Construction from the code.} Take $m:=1$. For an arbitrary tuple $c$ of $n$ formulas
denote by $q_c$ the program ``print $S(c)$ and halt''; its number
is computable from $c$. Let $a_1(c),\dots,a_n(c)$ be the programs that the algorithm of Lemma~\ref{lem:vssch} outputs on
$q_c$. This algorithm is a total
computable function of the number of the input program, so $a_j(c)$ are defined
for \emph{every} $c$; put
\[
  A_j(c):=\{x:\ \text{the program } a_j(c)\ \text{does not print } x\}.
\]
The condition on the right is $\forall i\ \aneg\bigl(U(a_j(c),i)\!\downarrow=x\bigr)$,
that is, a $\Pi_1$-condition, and uniformly in~$c$. Hence there are arithmetical
formulas $\theta_j(c,x)$ defining $A_j(c)$ for each~$c$.

The statement of the lemma guarantees that for all $c$
the sets $A_j(c)$ are nonempty and that
\[
  S(c)\notin J_n\bigl(A_1(c),\dots,A_n(c)\bigr) \tag{$\ast$},
\]
and moreover this statement is provable in PA. Note that
this statement is the negation of the formula $(J_n)^{*}(S(c))$ for the interpretation
$*$ (depending on $c$) defined by the equalities $\alpha_j^*:=\theta_j(c,x)$.
If we managed to choose $c$ so that $c$ coincides with $\langle \theta_1(c,x),\dots,\theta_n(c,x)\rangle$,
we would have reached our goal: the formula $\aneg(J_n)^{*}(S(c))$ would be
provable in PA for the corresponding interpretation $*$. But clearly this is hopeless.


\emph{
The vicious circle is broken}
by the joint fixed-point lemma (see, e.g.,
\cite[Ch.~1]{Boolos93}): for any formulas $\theta_1(c,x),\dots,\theta_n(c,x)$
there exist formulas $\chi_1(x),\dots,\chi_n(x)$ for which
\[
  \PA\vdash \chi_j(x)\aeq\theta_j\bigl(\langle
   \gn{\chi_1},\dots,\gn{\chi_n}\rangle,\,x\bigr),\qquad j=1,\dots,n .
\]
Now we can define the interpretation $*$ by putting $\alpha_j^*:=\chi_j$ and $c:=\langle
   \gn{\chi_1},\dots,\gn{\chi_n}\rangle$.
By construction we can prove in PA the equivalence of $\theta_j(c,x)$ and $\chi_j(x)$,
and hence also the equivalence of $J_n\bigl(A_1(c),\dots,A_n(c)\bigr)$
and $J_n\bigl(\chi_1,\dots,\chi_n\bigr)$. Hence the formula $S(c)\notin J_n\bigl(\chi_1,\dots,\chi_n\bigr)$
is likewise provable in PA. And this formula is exactly $\aneg(J_n)^{*}(S(c))$.

As usual, for the semantics $\sE$ we can simply say that provability implies truth,
so $J_n$ is not valid. And for the semantics $\sI$, as above for
$\sB$, we use the consistency of the external theory.
\end{proof}

Let us put both bounds together.

\begin{corollary}\label{cor:sandwich}
For each of the three semantics $Z$
\[
  \Int\ \subsetneq\ \Lp_Z\ \subsetneq\ \J .
\]
Moreover, $\Lp_{\sI}\subseteq\Lp_{\sE}$.
\end{corollary}

\begin{proof}
The first inclusion and its strictness are Propositions~\ref{prop:intval}
and~\ref{cor:rose} (Rose's formula is valid in all three semantics,
and with a constant solver at that). Further: by Theorem~\ref{th:posfrag}
each logic $\Lp_Z$ has an intuitionistic positive fragment, and
then by Yankov's criterion (Theorem~\ref{th:jan}) $\Lp_Z\subseteq\J$.
Strictness: weak excluded middle belongs to $\J$ by
definition, but is not valid in any of the semantics
(Proposition~\ref{prop:noem}). The inclusion
$\Lp_{\sI}\subseteq\Lp_{\sE}$: the internal theory $\PA$ is contained in
the external theory $\TA$.
\end{proof}

\subsection{Independence of the axioms}\label{sec:indep}

Nontrivial facts about derivability are, first of all, facts about
\emph{non}-derivability. We collect them here.

\paragraph{Incompleteness on the problem side.} The families of \emph{valid} modality-free problems---the logics $\Lp_Z$---are
described in~\S\ref{sec:modalfreeval} up to a sandwich. Here we deal with
\emph{derivable} modality-free problems, namely with derivability in $\HCI$.
About $\HCE$ we say right away the
little we know: its modality-free fragment contains the fragment of $\HCI$
(every theorem of $\HCI$ is an axiom of $\HCE$) and is contained in $\Lp_{\sE}$ (soundness,
Theorem~\ref{th:sound}), but whether it coincides with the fragment of
$\HCI$ is unknown: from (A8) and (A9) by (MP) new modality-free problems could
in principle be derived, and the model of
Theorem~\ref{prop:modalfree} below does not rule this out---the schemes (A8) and
(A9) are false in it (see also \S\ref{sec:todo}).

About $\HCI$ the following is known.

\begin{theorem}[the modality-free fragment is exactly $\Int$]
\label{prop:modalfree}
Let $\varphi$ be a problem containing neither $\bang$ nor $\quest$. Then
$\varphi$ is derivable in $\HCI$ if and only if it is
derivable in $\Int$.
\end{theorem}

The proof is in the appendix.

\begin{corollary}[incompleteness on the problem side]\label{cor:incomplete}
The calculus $\HCI$ is incomplete with respect to each of the semantics $\sI$, $\sE$:
Rose's formula is valid in them (Proposition~\ref{cor:rose}), but it is
modality-free and not derivable in $\Int$ by Rose's theorem~\cite{Rose53}, and
hence not derivable in $\HCI$ either.
\end{corollary}
Whether Rose's formula is derivable in $\HCB$, we do not know
(Remark~\ref{rem:qfree}); two more candidates for an example of incompleteness of
$\HCB$ are the problems
\[
  \Pi\alpha\bigl(((\alpha\pto\beta)\pto\alpha)\pto\alpha\bigr)\pto\beta
  \qquad\text{and}\qquad
  \Pi\alpha\bigl(((\alpha\pto\beta)\pto\beta)\pto(\alpha\por\beta)\bigr)\pto\beta,
\]
which are $\sB$-valid by Proposition~\ref{prop:pienc}(4),(6); their
derivability in $\HCB$ is unknown (Remark~\ref{rem:peirce}).

\paragraph{The scheme (A7).}

\begin{theorem}\label{prop:a7main}
The scheme (A7) is not derivable in $\HCI$ from the remaining axioms.
\end{theorem}

The proof is a two-point Kripke model in which all the schemes
(A1)--(A6) and all three inference rules hold, while the problem
$\pneg\pneg\bang\bot\pto\bang\bot$ does not; it is constructed in
Appendix~\ref{sec:a7indep} (Proposition~\ref{prop:a7indep}). For $\HCE$
the question is open, and it cannot be settled by models in this way (see the remark below).

\paragraph{The schemes (A8) and (A9).} They are not merely non-derivable in $\HCI$, but
incompatible with it.

\begin{remark}
\label{rem:a10fork}
Neither (A8) nor (A9) can be added to $\HCgl$: \emph{both calculi
$\HCgl+$(A8) and $\HCgl+$(A9) are inconsistent}. A single instance of
(A8) suffices---the proposition $\square\bot\ato\bot$: the rule (AtoP) gives
$\bang(\square\bot\ato\bot)$, L\"ob's scheme (A6) with $A:=\bot$ gives the problem
$\bang\bot$, the rule (PtoA) gives the proposition $\quest\bang\bot$, that is,
$\square\bot$, and $\square\bot\ato\bot$ yields $\bot$.

The scheme (A9) supplies this instance: in Example~\ref{ex:proofs}(4) the problem
$\pneg\bang\bot$ is derived from (A9), and from it, by the rule (PtoA) and
the schemes (A4a), (A4d), the proposition $\quest\bang\bot\ato\bot$, that is,
$\square\bot\ato\bot$.

In other words, (A9) is related to the sort of problems exactly as reflection (A8)
is related to the sort of propositions: both are
incompatible with the pair ``necessitation $+$ L\"ob''. This is why both stand in
$\HCE$, where the rule (AtoP) is absent, and this is also why $\HCE$, despite its appearance,
is \emph{not} closed under (AtoP): its consistency is given only by
Theorem~\ref{th:sound}.

A practical consequence: no single algebraic model can
satisfy all the axioms and rules of $\HCgl$ and in addition all the
instances of (A9). Independence of anything from $\HCgls$ cannot be
proved by models in this way---one needs a technique like the one used for $\GLS$.
(End of remark.)
\end{remark}

The converse question---whether (A8) is derivable from (A9)---is settled
negatively, although special cases of (A8) are obtained from (A9): thus in
Example~\ref{ex:proofs}(4) the proposition $\square\bot\ato\bot$ is derived.

\begin{proposition}[(A8) is not derivable from (A9)]\label{prop:a8indep}
Let $A$ be a variable of sort proposition. The proposition
$\square A\ato A$ is not derivable from the theorems of $\HCI$ and the scheme (A9) by the rules
(MP) and (PtoA). On the other hand, all instances of reflection
$\square\quest\varphi\ato\quest\varphi$ are derivable from them.
\end{proposition}

\begin{proof}
Derivability: apply to an instance of (A9) the rule $(\quest\to)$ and the schemes
(A4a), (A4d): $\quest\bang\quest\varphi\ato\quest\pneg\pneg\varphi$, and
$\quest\pneg\pneg\varphi\aeq\aneg\aneg\quest\varphi\aeq\quest\varphi$.

Non-derivability. Let $C_1,\dots,C_m$ be a derivation of the proposition $C_m$ from
theorems of $\HCI$ and instances of (A9) by the rules (MP) and (PtoA). To
each formula of the derivation we assign a proposition $C_i'$: if $C_i$ is a proposition,
then $C_i':=C_i$, and if it is a problem, then $C_i':=\quest C_i$
(we do not look inside the formulas).

We show that
$\nf{C_m}$ is derivable from theorems of $\GL$ and finitely many propositions
of the form $\square N\ato N$ by the rule (MP) alone, where the $N$ are
$\square$-formulas in which all propositional variables occur under
$\square$ (only the atoms $\quest\alpha$ occur freely); in particular,
$N=A$ is not among them. If $C_i$ is a theorem of
$\HCI$, then $C_i'$ is a theorem of $\HCI$ as well (applying
(PtoA) if necessary), and $\GL\vdash\nf{C_i'}$ by Corollary~\ref{cor:solovay}. If
$C_i$ is an instance of (A9) $\bang\quest\varphi\pto\pneg\pneg\varphi$, then
$\nf{C_i'}=\square N\ato\aneg\aneg N$, where $N=\varphi^\quest$; it is
obtained by (MP) from $\square N\ato N$ and a theorem of classical logic.
The formula $N$ is of the indicated form, since $\varphi$ is a problem:
a proposition enters a problem only in the form $\bang B$, and
$(\bang B)^\quest=\square\,\nf{B}$. If $C_i$ is obtained by (MP) from the problems $C_j$ and
$C_l=C_j\pto C_i$, then $\nf{C_i'}$ is obtained by (MP) from $\nf{C_j'}$ and
$\nf{C_l'}=\nf{C_j'}\ato\nf{C_i'}$; the rule (MP) on propositions
is carried to itself. Finally, if $C_i=\quest C_j$ is obtained by (PtoA),
then $C_i'=C_j'$.

Thus $\nf{C_m}$ is obtained from theorems of $\GL$ and finitely
many propositions $\square N_i\ato N_i$ of the indicated form
by applications of (MP) alone, and therefore
\[
  \GL\vdash\textstyle\bigwedge_i(\square N_i\ato N_i)\ato\nf{C_m}.
\]
For $C_m=\square A\ato A$ (here $\nf{C_m}=C_m$) this is impossible. Let $k$
be greater than the modal depth of all the $N_i$, and consider the $\GL$-model on the chain
$w_0Rw_1R\dots Rw_k$ (transitive and irreflexive) in which $A$
is false at $w_0$ and true at the remaining points, and all other atoms
are constant. Then $w_0\Vdash\square A$ and $w_0\nVdash A$. On the other
hand, at the points $w_1,\dots,w_k$ all atoms are constant, and by induction on
the modal depth $d$ of an arbitrary $\square$-formula $X$ (in which $A$
may also occur outside modalities) one checks that at all points
$w_1,\dots,w_{k-d}$ the truth value of $X$ is the same. Base ($d=0$): $X$
is a Boolean combination of atoms, which are constant at these points. Induction step: let
$X=\square Y$ have depth $d+1$ and $1\le i\le k-d-1$; the successors of the point
$w_i$ are $w_{i+1},\dots,w_{k-d}$, at which $Y$ has the same value by the induction
hypothesis, and $w_{k-d+1},\dots,w_k$, which are common
to all such $w_i$. Similarly, for a formula $N_i$ in which $A$
occurs only under $\square$, the truth value of $N_i$ at $w_0$ coincides with
its truth value at $w_1$. Indeed, $N_i$ is a Boolean combination of atoms
other than $A$ and of formulas of the form $\square X$; the atoms other than $A$
have the same values at $w_0$ and $w_1$, while $\square X$ at $w_0$
differs from $\square X$ at $w_1$ only by the value of $X$ at $w_1$,
which coincides with the value of $X$ at $w_2$ by what has been proved (the depth of
$X$ does not exceed $k-2$). Hence $w_0\Vdash\square N_i$ implies $w_1\Vdash N_i$, and with
it $w_0\Vdash N_i$; all the premises $\square N_i\ato N_i$ are true at $w_0$,
while $\square A\ato A$ is not.
\end{proof}

Thus the schemes (A8) and (A9) in $\HCE$ are independent at least in one
direction; the converse independence---of (A9) from $\HCI+$(A8)---remains
open (\S\ref{sec:todo}).

\paragraph{Schemes that did not become axioms.} The axioms of the calculi were selected by
validity, and some natural candidates were rejected along the way.
Two of them are worth naming---as traces of the search for axioms:
the distribution of $\bang$ over disjunction (and its $\square$-version) and the attempt
to invert the scheme (A3), obtaining from a transformation of proofs a proof of an implication
between provabilities.

\begin{theorem}\label{prop:notaxioms}
Let $A$, $B$ be variables of sort proposition. The following formulas are not
valid in any of the three semantics, and hence (Theorem~\ref{th:sound})
not derivable in any of the three calculi:
\begin{enumerate}[leftmargin=2.2em,itemsep=2pt,label=\textup{(\arabic*)}]
\item 
$\square(A\aor B)\ato(\square A\aor  \square B)$;
\item $(\bang A\pto\bang B)\pto\bang(\square A\ato\square B)$.
\end{enumerate}
\end{theorem}

\begin{proof}
(1) 
Let $\sigma$ be a Rosser sentence for the internal theory:
the consistency of the theory implies that neither $\sigma$ nor $\aneg\sigma$ is
provable in it, and this argument can be formalized in the theory itself
(see, e.g.,~\cite[Ch.~2]{Boolos93}). Take $A^*:=\sigma$,
$B^*:=\aneg\sigma$. The disjunction $A^*\aor B^*$ is a tautology and hence
provable, while neither $A^*$ nor $B^*$ is provable.
In the semantics $\sE$ everything is already clear:
the premise is true and the conclusion is false. In the
internal semantics we must refute provability in the external theory
$T$ (which is also the internal one): the premise is a true $\Sigma_1$-sentence and
hence provable, so provability of the proposition would imply
$T\vdash\Prov(\sigma)\aor\Prov(\aneg\sigma)$, and together with
the formalized Rosser theorem
\[
  T\vdash\Con\ato\aneg\Prov(\sigma)\aand\aneg\Prov(\aneg\sigma)
\]
this would give $T\vdash\aneg\Con$. The latter is impossible under
$\Sigma_1$-soundness of $T$: $\aneg\Con$ is a false $\Sigma_1$-sentence.

(2) Take $A^*:=\aneg\Con$ (that is, $\square\bot$) and $B^*:=\arbot $;
denote by $L$ the translation of the proposition $\square A\ato\square B$, that
is, the formula $\square\square\bot\ato\square\bot$, and by $T$ the external
theory. Here we need $\Sigma_1$-soundness of $T$
(Remark~\ref{rem:assumptions});
since the internal theory is contained in the external one, it is then
$\Sigma_1$-sound as well. We must show that $T\nvdash\varphi^*(n)$ for
every $n$, where
$$
\varphi=(\bang\square\bot\pto\bang\bot)\ato \bang(\square\square\bot\pto\square\bot)
$$
 (for $\sE$ this means
that the value of the problem is empty). The formula $\varphi^*(n)$ has
the form $$\forall v\,\bigl(P(v)\ato\Prf(U_n(v),L)\bigr),$$ where $P(v)$ is
the translation of the premise $\bang A\pto\bang B$, that is,
$$\forall w\,\bigl(\Prf(w,\aneg\Con)\ato\Prf(U_v(w),\arbot )\bigr).
$$
Suppose
$T\vdash\varphi^*(n)$. We argue in $T$: if $\aneg\square\square\bot$,
that is, there are no proofs of the formula $\aneg\Con$, then $P(v)$ holds for every
$v$, since its premise is false, hence $\Prf(U_n(0),L)$ and all the more $\square L$. Thus
$T\vdash\aneg\square\square\bot\ato\square L$. On the other hand,
the formalized L\"ob theorem $\square(\square C\ato C)\ato\square C$ is a
theorem of $\PA$, and hence of $T$; with
$C:=\square\bot$ it gives $T\vdash\square L\ato\square\square\bot$. Together
we obtain $T\vdash\square\square\bot$. This is a $\Sigma_1$-sentence, and by
$\Sigma_1$-soundness of $T$ it is true: the internal theory proves
$\square\bot$, that is, $\aneg\Con$---contrary to its
$\Sigma_1$-soundness.
\end{proof}


\subsection{Derivable problems with modalities: what is known}\label{sec:taskder}

We have no complete description of the family of derivable problems---just as we have no description
of the family of valid ones (\S\ref{sec:taskval}); it is not even known whether the
system $\HCgl$ is complete on the problem side (\S\ref{sec:todo}). Here we collect what is known about
problems of a special form---about the $\triangledown$-fragment.

For $\triangledown=\bang\quest$ the schemes K, 4 and L\"ob are derivable. The first scheme
$\triangledown(\varphi\pto\psi)\pto(\triangledown\varphi\pto\triangledown\psi)$
is obtained from (A3) with $A:=\quest\varphi$, $B:=\quest\psi$, since
$\bang\quest(\varphi\pto\psi)\peq\bang(\quest\varphi\ato\quest\psi)$ by (A4a) and
Lemma~\ref{lem:repl}. The second scheme $\triangledown\varphi\pto\triangledown\triangledown\varphi$
is (A5) with $\sigma:=\bang\quest\varphi$. Finally,
$\triangledown(\triangledown\varphi\pto\varphi)\pto\triangledown\varphi$
is obtained from (A6) with $A:=\quest\varphi$, since
$\bang\quest(\triangledown\varphi\pto\varphi)\peq\bang(\square\quest\varphi\ato\quest\varphi)$
by (A4a) and Lemma~\ref{lem:repl}. The necessitation rule for
$\triangledown$ is (PtoA) followed by (AtoP). Co-reflection
$\varphi\pto\triangledown\varphi$ is not derivable: it is not valid
(Proposition~\ref{prop:hce}(4)). The converse scheme
$\triangledown\triangledown\varphi\pto\triangledown\varphi$
is $\sE$-valid, but not $\sI$-valid and, consequently, not derivable in $\HCgl$:
it would require $\PA\vdash\square\square A\ato\square A$.

Separately we note what has already been proved: the scheme $\pneg\pneg\triangledown\varphi\pto\triangledown
\varphi$ is derivable as a special case of (A7), and co-reflection
$\varphi\pto\triangledown\varphi$ is derivable on the whole class of $\Sigma$-problems and
only on it ($\triangledown\varphi=\bang\quest\varphi$ is a
$\Sigma$-problem by Definition~\ref{def:sigma}; the axiom (A5) and
Proposition~\ref{prop:hce}(4)). The main point is the following.

\begin{theorem}\label{prop:nabla}
If the Boolean collapse of a problem $\chi$ by the schemes (A4) is a classical
tautology, then the problem $\triangledown\chi$ is derivable in $\HCgl$. In particular,
the following are derivable:
\[
  \triangledown(\varphi\por\pneg\varphi),\qquad
  \triangledown(\pneg\pneg\varphi\pto\varphi),\qquad
  \triangledown\bigl(((\varphi\pto\psi)\pto\varphi)\pto\varphi\bigr).
\]
\end{theorem}

\begin{proof}
By the schemes (A4) the proposition $\quest\chi$ is derivably equivalent to
its Boolean collapse (Proposition~\ref{prop:nf}), and the latter is a
classical tautology, the scheme (A1); hence $\quest\chi$ is a theorem of
$\HCgl$, and, applying (AtoP), we obtain the problem
$\bang\quest\chi=\triangledown\chi$.
\end{proof}

\remove{
\begin{remark}\label{rem:nabla}
Причина расхождения структурная, а не случайная. Модальность
$\triangledown=\bang\quest$ пропускает свой аргумент через $\quest$, а $\quest$ по (A4)
--- булев гомоморфизм; поэтому под знаком $\triangledown$ интуиционистская
природа задач стирается, и туда попадает вся классическая логика.
\end{remark}
}

\section{Comparison with related logics}\label{sec:hc}

Of the systems listed in the introduction, only $\HC$ has the same syntax, and
only the comparison with it requires separate work; most of this section is
devoted to it. The remaining comparisons rest on what has already been proved
and are only named here---at the end of the section.

\subsection{Melikhov's logic $\HC$}

Recall the axiomatics of $\HC$: the classical schemes (propositions),
the intuitionistic schemes (problems), the schemes
\[
  \bang(A\ato B)\pto(\bang A\pto\bang B),\quad
  \quest(\varphi\pto\psi)\ato(\quest\varphi\ato\quest\psi),\quad
  \quest\bang A\ato A,\quad
  \varphi\pto\bang\quest\varphi,\quad
  \pneg\bang\bot
\]
and the rules (MP), (AtoP), (PtoA).

The comparison is two-sided: which theorems of $\HC$ are valid in our semantics and
which of our axioms are derivable in $\HC$. The first half is given by
Proposition~\ref{prop:hce} and the following theorem.

\begin{proposition}\label{prop:hcval}
The problem $\pneg\bang\bot$ is $\sE$-valid, but neither $\sI$- nor
$\sB$-valid.
\end{proposition}

\begin{proof}
In $\sE$: this problem is derivable in $\HCE$ (Example 1(4)), hence valid.
In the internal
semantics: the formula $(\pneg\bang\bot)^*(x)$ is provably equivalent to $\Con$
---again G\"odel's second theorem.
\end{proof}

The first two schemes displayed above are our (A3) and
the left half of (A4a); they are $\sI$-valid (Theorem~\ref{th:sound}).
Reflection $\quest\bang A\ato A$ and the scheme $\pneg\bang\bot$ are only $\sE$-valid
(Proposition~\ref{prop:hce}(1),(2) and Theorem~\ref{prop:hcval}), while co-reflection
$\varphi\pto\bang\quest\varphi$ is not valid at all
(Proposition~\ref{prop:hce}(4)). The rules: (MP) and (PtoA) preserve both $\sI$- and
$\sE$-validity (Lemma~\ref{lem:closure} and Theorem~\ref{th:sound}), whereas (AtoP)
preserves only $\sI$-validity: by Proposition~\ref{prop:hce}(3), applying (AtoP) to
the $\sE$-valid proposition $\square\bot\ato\bot$ yields a problem that is not
valid. The rule $\quest\varphi/\varphi$, the converse of (PtoA), is not admissible
in our semantics either: the proposition $\quest(\alpha\por\pneg\alpha)$ is valid,
since $\val{\alpha\por\pneg\alpha}\neq\varnothing$ under every
interpretation, whereas the problem $\alpha\por\pneg\alpha$ is not valid
(Proposition~\ref{prop:noem}).

It is convenient to name the outcome in advance: the common part of the two calculi is the \emph{core}
$\HCzero$, the axioms (A1), (A2), (A3), the left half of (A4a), (A4b)--(A4d) and (A5) with all
three rules; everything beyond the core diverges, and in both directions.

The second half remains---what of our axiomatics is derivable in $\HC$.
We begin with the group (A4): it is examined item by item.

\begin{remark}[the group (A4) in $\HC$; residual classicality]\label{rem:qhom}
The scheme (A4d) is derivable in $\HC$ in full. The direction $\bot\ato\quest\varnothing$
is classical logic; for the converse we take the intuitionistic theorem
$\varnothing\pto\bang\bot$, apply (PtoA) and the left half of (A4a), which gives
$\quest\varnothing\ato\square\bot$, and reflection $\square\bot\ato\bot$ remains.

The scheme (A4c) is also derivable in full. The easy direction
$(\quest\varphi\aor\quest\psi)\ato\quest(\varphi\por\psi)$ is obtained from the embeddings
$\varphi\pto\varphi\por\psi$ by the same device: (PtoA), then the left half of
(A4a). The substantive direction goes through co-reflection: from
$\varphi\pto\bang\quest\varphi$ and the $\bang$-version of the tautology
$\quest\varphi\ato(\quest\varphi\aor\quest\psi)$ (the rule (AtoP)), by (A3) we obtain
$\varphi\pto\bang(\quest\varphi\aor\quest\psi)$, likewise for $\psi$, and hence also
$\varphi\por\psi\pto\bang(\quest\varphi\aor\quest\psi)$; applying (PtoA) and the left
half of (A4a), we obtain
$\quest(\varphi\por\psi)\ato\square(\quest\varphi\aor\quest\psi)$, and reflection
remains. Conceptually this is to be expected: in $\HC$ the operator $\quest$ is the left
adjoint of $\bang$, and left adjoints preserve disjunctions.

Of the scheme (A4b), the half
$\quest(\varphi\pand\psi)\ato(\quest\varphi\aand\quest\psi)$ is derivable in $\HC$---from the projections, again
(PtoA) plus the left half of (A4a). In the converse direction: to the intuitionistic axiom
$\phi\ato(\psi\ato(\phi\aand\psi))$ we apply the left half of (A4a), obtaining
$\quest\phi\ato\quest(\psi\pto(\phi\pand\psi))$. Using the axiom
$\quest(\psi\ato(\phi\aand\psi))\ato(\quest\psi\ato\quest(\phi\aand\psi))$
and MP we derive $\quest\phi\ato(\quest\psi\ato\quest(\phi\aand\psi))$, which classically is the same
as the desired formula.


Thus what remains essentially new is the right half of (A4a),
\[
  (\quest\varphi\ato\quest\psi)\ato\quest(\varphi\pto\psi),
\]
and its derivative $\aneg\quest\varphi\ato\quest\pneg\varphi$ (obtained from
(A4a) and (A4d)). $\HC$ does not have them: in Kripke frames, a quantifier
over future worlds arises in the premise. [TODO Add a reference to the work of Melikhov where this is proved.]
(End of remark.)
\end{remark}

\begin{proposition}\label{prop:incl}\phantom{l}
\begin{enumerate}[leftmargin=2.2em,itemsep=2pt,label=\textup{(\arabic*)}]
\item $\HCzero\subseteq\HC$.
\item $\HC\nsubseteq\HCgls$; a fortiori $\HC\nsubseteq\HCgl$.
\item $\HCglm\nsubseteq\HC$; a fortiori $\HCgl\nsubseteq\HC$ and
  $\HCgls\nsubseteq\HC$.
\item The calculus $\HC+$(A6) is inconsistent.
\end{enumerate}
\end{proposition}

\begin{proof}
(1) The schemes (A1), (A2), (A3), the left half of (A4a) and the rules are exactly
a part of the axiomatics of $\HC$, the schemes (A4b)--(A4d) are derivable in $\HC$
(Remark~\ref{rem:qhom}), and (A5) is a special case of co-reflection
$\varphi\pto\bang\quest\varphi$ (for $\Sigma$-problems $\varphi$).

(2) Co-reflection is derivable in $\HC$, and by Proposition~\ref{prop:hce}(4) it
is not $\sE$-valid, while $\HCE$ is sound with respect to $\sE$ (Theorem~\ref{th:sound});
the theorems of $\HCgl$ are among those of $\HCgls$.

(3) Already the right half of (A4a) is not derivable in $\HC$ (Remark~\ref{rem:qhom});
nor is (A7): its special case $\pneg\pneg\bang A\pto\bang A$ is Markov's
principle for $\bang$, while the problems in $\HC$ are intuitionistic, without Markov-style additions.

(4) $\HC$ has reflection $\square A\ato A$ and the necessitation rule for
$\square$ (apply (AtoP), then (PtoA)); from (A6) we would obtain
$\square(\square A\ato A)$, then $\square A$ for every $A$, and for $A=\bot$
a contradiction with reflection. Meanwhile $\HC$ is consistent: it is sound
with respect to Kripke frames with marked worlds~\cite{Onoprienko20}.
\end{proof}

Thus, of our calculus exactly the core lives inside $\HC$, and the divergence runs
along several independent lines: the classicality of the modality $\quest$ (the right
half of (A4a) with its derivative $\aneg\quest\varphi\ato\quest\pneg\varphi$),
Markov's principle (A7) and the L\"ob line---the scheme (A6), which is not merely
absent from $\HC$ but incompatible with it.
In the opposite direction there is a single line---
co-reflection.
It is essential that the derivation of the formula (L) goes through already in the core
(Theorem~\ref{th:lob}): therefore Corollary~\ref{cor:hcfp} on the
incompatibility of $\HC$ with fixed points relies only on the common part of
the two calculi.

The incompatibility of $\HC$ with the reading of $\bang$ as formal provability
is visible even without fixed points. In $\HC$ the problem
$\bang(\square\bot\ato\bot)$ is derivable: we apply the rule
(AtoP) to the special case $\square\bot\ato\bot$ of reflection.
If the solutions of the problem $\bang A$ are proofs of $A$ in a consistent theory
$T$ proving the Hilbert--Bernays derivability conditions,
then a solution of this problem is a proof of
$\mathrm{Con}(T)$ in $T$, which does not exist by G\"odel's second theorem. This is
a two-sorted version of G\"odel's observation~\cite{Goedel33}: the $\square$-fragment
of $\HC$ is $\mathrm{S4}$~\cite{Melikhov13}, and in $\mathrm{S4}$
$\square\aneg\square\bot$ is derivable, so $\square$ in $\mathrm{S4}$ cannot be
read as formal provability. What is excluded is precisely provability in
a single fixed formal theory; a semantics with explicit
proofs in the spirit of Artemov's logic of
proofs~\cite{Artemov01} is not affected by this argument: there
$\aneg(t{:}\bot)$ is derivable for every term $t$, but the statement
``there is no proof of $\bot$ at all'' is not in the language.
\remove{
Наконец, всё это сравнение держится на выбранном понятии общезначимости. При
слабом варианте~(а) замечания~\ref{rem:three} фрагмент задач классичен, схему
(A2) пришлось бы усиливать до классической, и связь с интуиционистской частью
$\HC$ пропала бы вовсе; содержательным сопоставление двух исчислений делает
именно требование предъявлять решение вычислимо.
}

\subsubsection*{Fixed points are incompatible with $\HC$}

From the derivation of the formula (L) in Example~\ref{ex:proofs}(5) one can extract more than
the soundness of (A6). A careful inspection shows that the derivation goes through already in the
common part of the two calculi---in the core $\HCzero$.

\begin{theorem}\label{th:lob}
Already in the core $\HCzero$ (and hence in $\HCglm$ and in $\HC$) the problem
\[
  \bang\bigl(P\aeq(\square P\ato S)\bigr)\ \pto\
  \bigl(\bang(\square S\ato S)\pto\bang S\bigr)
\]
is derivable for all propositions $P$ and $S$.
\end{theorem}

\begin{proof}
This is item~(5) of Example~\ref{ex:proofs}. It suffices to note that it
uses only (A1), (A2), (A3), (A5), the rules (MP), (AtoP), (PtoA),
the deduction lemma and the lemma on transferring derivations
(Proposition~\ref{prop:derived}(2),(3)) and the left half of (A4a)---the latter
only in items~(1) and~(2) of the same example, in deriving the $\square$-versions of (A3)
and (A5). Neither the right half of (A4a), nor (A6), nor (A7) is needed.
\end{proof}
\remove{
\begin{remark}[как устроен вывод и почему посылку удаётся снять]
Посылкой служит не сама эквивалентность $P\aeq(\square P\ato S)$, а
задача $\bang(\dots)$ <<найти её вывод>>: (A3) умеет проносить $\bang$ через
импликацию, и потому работать можно только с выводами.

Сам вывод устроен в два этажа. Внизу, на сорте высказывание, всё делается
классически: шаги~(1)--(6) пункта~(5) примера~\ref{ex:proofs} --- это обычный
модальный вывод, использующий лишь $\square$-версии (A3) и (A5). Наверх, под
знак~$\bang$, он переносится целиком леммой о переводе выводов
(предложение~\ref{prop:derived}(3)); правило (AtoP) при этом применяется
только к уже доказанным формулам, а к гипотезам --- никогда, и потому лемму о
дедукции в конце применить можно.

Схема (A5) работает дважды, и оба раза для того, чтобы приписать бокс под
знаком~$\bang$: она даёт $\bang\square P$ из $\bang P$ и вторую гипотезу из
$\bang(P\ato(\square P\ato S))$. Она занимает то самое место, которое в
стандартном доказательстве теоремы Лёба занимает схема
$\square A\ato\square\square A$.
(Конец замечания.)
\end{remark}
}
\begin{corollary}[fixed points are incompatible with $\HC$]\label{cor:hcfp}
Suppose that a fixed point for falsity is added to $\HC$: that is, the axiom
$\bang\bigl(P\aeq(\square P\ato\bot)\bigr)$,
where $P$ is a variable of sort proposition.
Then the resulting
system is inconsistent---the problem $\varnothing$ is derivable in it, and hence
every problem.
\end{corollary}

\begin{proof}
By Proposition~\ref{prop:incl}(1) the core $\HCzero$ is contained in
$\HC$, so Theorem~\ref{th:lob} gives in $\HC$ the implication
\[
  \bang\bigl(P\aeq(\square P\ato\bot)\bigr)\pto
  \bigl(\bang(\square\bot\ato\bot)\pto\bang\bot\bigr) .
\]
The premise is derivable by assumption. Further, $\square\bot\ato\bot$ is a special case
of reflection, that is, a theorem of $\HC$; by (AtoP),
$\bang(\square\bot\ato\bot)$ is derivable, and by two applications of (MP) we obtain $\bang\bot$.
Together with the axiom $\pneg\bang\bot$ this gives $\varnothing$.
\end{proof}
\remove{
Это синтаксическая форма развилки Лёба: $\HC$ несовместима не с какой-то одной
арифметической конструкцией, а с самим наличием неподвижных точек --- а в
арифметической интерпретации они есть всегда. Система $\HCgl$, наоборот,
неподвижные точки спокойно выдерживает и даже извлекает из них схему Лёба.
}

\subsection{Other systems}

Here only a list, with an indication of the places where the corresponding statements
are proved.

\emph{Provability logic.} The $\square$-fragment of $\HCgl$ is exactly $\GL$,
the $\square$-fragment of $\HCgls$ is exactly $\GLS$; together with Solovay's
theorems~\cite{Solovay76} this yields completeness on the propositional side:
Corollary~\ref{cor:solovay}.

\emph{Intuitionistic provability logics.} Recall the
definitions~\cite{VZ19,DvdGM24}: $\mathrm{iK}$ is the intuitionistic
tautologies plus the scheme $\square(A\ato B)\ato(\square A\ato\square B)$ with
the rules (MP) and necessitation; $\mathrm{iGL}=\mathrm{iK}+{}$L\"ob's scheme
$\square(\square A\ato A)\ato\square A$; finally
$\mathrm{iGLC}=\mathrm{iGL}+{}$the \emph{completeness principle} $A\ato\square A$.
The Kripke semantics for $\mathrm{iGL}$ consists of birelational frames
$\langle W,\preccurlyeq,\sqsubset\rangle$, where $\preccurlyeq$ is a
partial order (the intuitionistic relation) and $\sqsubset$ is a
transitive conversely well-founded modal relation, compatible
with the order in both directions: $x\preccurlyeq y\sqsubset z$ implies
$x\sqsubset z$, and $x\sqsubset y\preccurlyeq z$ implies $x\sqsubset z$
(Iemhoff calls such frames \emph{brilliant}). Valuations of variables are
monotone with respect to $\preccurlyeq$, and $\square A$ is forced at a point $x$ if
$A$ is forced at all its $\sqsubset$-successors. The logic $\mathrm{iGL}$
is complete with respect to such frames~\cite[Proposition~9]{Iemhoff99}; see\
also~\cite{DvdGM24}. We compare them with
the $\triangledown$-fragment of $\HCgl$ (\S\ref{sec:taskder}), reading $\square$ as
$\triangledown$.

\begin{proposition}\label{prop:igl}
The logic $\mathrm{iGL}$ is contained in the $\triangledown$-fragment of $\HCgl$,
and strictly so: every theorem of $\mathrm{iGL}$, read with
$\square:=\triangledown$, is derivable in $\HCgl$, and the problem
$\triangledown(\varphi\por\pneg\varphi)$ is derivable in $\HCgl$, whereas
the formula $\square(p\aor\aneg p)$ is not derivable in $\mathrm{iGL}$.
\end{proposition}

\begin{proof}
Inclusion: the scheme K, L\"ob's scheme and the necessitation rule for
$\triangledown$ are derived in \S\ref{sec:taskder}. Strictness is given by the problem
$\triangledown(\varphi\por\pneg\varphi)$ (Theorem~\ref{prop:nabla}): we show
that $\square(p\aor\aneg p)$ is not derivable in $\mathrm{iGL}$.

Take the frame $\langle W,\preccurlyeq,\sqsubset\rangle$ with $W=\{w,v,u\}$, where
$\preccurlyeq$ is the reflexive relation together with $v\preccurlyeq u$, and
$\sqsubset$ consists of the pairs $(w,v)$ and $(w,u)$. This is a frame for
$\mathrm{iGL}$: the relation $\sqsubset$ is transitive and conversely well-founded
(there are no chains of length two), and both compatibility conditions are checked
by exhaustive inspection---the only nontrivial case $w\sqsubset v\preccurlyeq u$
gives $w\sqsubset u$. Let the variable $p$ be true only at $u$;
the valuation is monotone, since $u$ is maximal. Then $v\nVdash p$ and
$v\nVdash\aneg p$ (since $v\preccurlyeq u$ and $u\Vdash p$), hence
$v\nVdash p\aor\aneg p$, and from $w\sqsubset v$ it follows that
$w\nVdash\square(p\aor\aneg p)$. By the completeness of $\mathrm{iGL}$ with respect to
such frames, the formula $\square(p\aor\aneg p)$ is not derivable in $\mathrm{iGL}$.
\end{proof}

\begin{remark}\label{rem:igl}
The reason for the divergence is that 
under the sign
$\triangledown$ all of classical logic appears, which no
intuitionistic modal logic can provide, and $\mathrm{iGL}$ is no
exception. The same shows where to look for an exact characterization:
the $\triangledown$-fragment is, apparently, not some
intuitionistic modal logic, but $\mathrm{iGL}$ together with all
formulas $\triangledown\chi$ whose Boolean collapse $\chi$ is a
tautology. Whether the fragment coincides exactly with this completion is
an open question (\S\ref{sec:todo}).

The completeness principle of the logic $\mathrm{iGLC}$ is co-reflection,
refuted by Proposition~\ref{prop:hce}(4); but the restriction of completeness to
$\Sigma_1$, adopted in~\cite{VZ19}, corresponds to our co-reflection on
$\Sigma$-problems (the axiom (A5)). Finally, in Melikhov's logic $\HC$, for
$\triangledown$ the scheme K, the scheme~4 and co-reflection hold, but not L\"ob's scheme
---this is the intuitionistic modal logic $\mathrm{H4}$, also known as
$\mathrm{IEL}^+$ of Artemov and Protopopescu~\cite{Protopopescu16}; in our setting
co-reflection is displaced by L\"ob's scheme, and the $\triangledown$-fragment of $\HCgl$ is
a ``L\"obian'' analogue of $\mathrm{H4}$ over the intuitionistic logic of problems.
(End of remark.)
\end{remark}

\emph{Realizability logics.} Our notion of validity is $\Reff$ from
the classification $\Rcnst\subseteq\Reff\subseteq\R$~\cite{VSSCh}; the three levels of
Remark~\ref{rem:three} are exactly these three logics. Hence the results of
\S\ref{sec:modalfreeval} on the fragment of problems: the positive part is intuitionistic
(Theorem~\ref{th:posfrag}), and the fragment itself lies strictly between $\Int$ and the logic
of weak excluded middle (Corollary~\ref{cor:sandwich}).

\section{Open questions}\label{sec:todo}

\begin{enumerate}[leftmargin=1.8em,itemsep=3pt]
\item \textbf{What logic has arisen on the problem side?} It is known that
  $\Int\subsetneq\Lp_Z\subsetneq\J$ for all three semantics
  (Corollary~\ref{cor:sandwich}) and that
  the positive fragment is intuitionistic (Theorem~\ref{th:posfrag}), so that the whole
  difference from $\Int$ lives on formulas containing $\varnothing$. On the left we are separated by
  Rose's formula (Proposition~\ref{cor:rose}), on the right by Scott's principle
  (Theorem~\ref{prop:scott}). Hence
  the scheme (A2) on the problem side is certainly insufficient, and $\HCI$ is
  incomplete on the problem side (Corollary~\ref{cor:incomplete}); we did not close the gap
  with the ready-made Rose scheme (\S\ref{sec:modalfreeval}), and what a
  complete calculus looks like is unknown. The setting
  is close to the propositional logic of realizability of Rose and
  Plisko~\cite{Rose53,Plisko09}---there, too, variables are replaced by
  arithmetical formulas---and it is from there that the next trial
  schemes should be taken.

  One should be aware of how different the two sides are. On
  the propositional side we deal with a known object: the fragments are
  exactly $\GL$ and $\GLS$, completeness follows
  from Solovay's theorems (Corollary~\ref{cor:solovay}), and this side (in
  the $\Pi$-free language) will no longer change, whatever valid
  schemes about problems we may add: new axioms yield only
  valid theorems, and all valid propositions are already derivable. On the problem side there is no known object: the logics
  $\Lp_Z$ lie strictly between $\Int$ and $\J$ and coincide with neither of them.
  The modality-free fragment of the calculus $\HCI$ is exactly
  $\Int$, and the calculus itself reads as a gluing of two ready-made systems---
  $\Int$ on problems and $\GL$ on propositions, joined by the schemes
  (A3)--(A5) and (A7). The semantics, however, is not a gluing: on one side
  stands the classical provability logic, on the other an intermediate
  logic that is yet to be identified.

  It may well be that it cannot be identified as a calculus at all. The closest
  relatives of the logics $\Lp_Z$ are the logics of realizability, and there the situation is as follows: for
  the logic of realizable propositional formulas neither finite nor even
  recursive axiomatizability is known, and Medvedev's logic of finite problems,
  as was shown recently, is $\Pi^0_1$-complete and hence not computably enumerable:
  there is no computably enumerable complete calculus for it
  at all~\cite{Pawlowski26,AlmeidaKnudstorp26}. Our notion of
  validity is $\Reff$, that is, it belongs to exactly this family. Therefore
  the fact that the valid Rose scheme (Proposition~\ref{cor:rose}) is a series
  indexed by the number of disjuncts is more likely a property of the subject than an accident. And
  for the same reason, selecting axioms by validity cannot
  serve as a sign of completeness: new schemes may stop being found not
  because the list is complete.

\item \textbf{How strong are the $\bang$-forms.} From (A5) and (A6) their
  $\square$-versions are derivable (Example~\ref{ex:proofs}(2),(3)); is the
  converse true? The $\square$-versions speak only of the existence of proofs, the $\bang$-versions
  of a computable transformation of proofs, so that the converse passage
  is not obvious. The answer determines how substantive
  Theorem~\ref{th:lob} is in comparison with its modal analogue
  (Appendix~\ref{sec:novelty}).

\item \textbf{Is (A5) needed for the derivation of L\"ob's theorem}, that is, is
  (A6) derivable from $\HCglm$ without (A5) (more precisely, without its special case
  $\bang A\pto\bang\square A$) together with the premise
  $\bang(P\aeq(\square P\ato S))$? In the proof of
  Theorem~\ref{th:sound} the scheme (A5) takes the place of the scheme
  $\square A\ato\square\square A$, without which the fixed-point argument
  does not go through in the modal case either; but there is no proof of non-derivability.
  (That (A5) is not derivable from the remaining axioms of $\HCgls$ together with the rule
  ``from $\bang A\pto\bang B$ infer $\bang(\square A\ato\square B)$'' and
  all problems of the form $\bang T$ for classical tautologies $T$
  is seen from a reflexive non-transitive three-world Kripke frame in which
  $\bang$ is read as $\square$ and $\quest$ as the identity.)

\item \textbf{Completeness on the problem side.} On the propositional side the question is closed
  (Corollary~\ref{cor:solovay}), and Proposition~\ref{prop:nf} shows that
  as far as schemes about $\quest$ are concerned the system (A4) is exhaustive: $\quest$ is eliminated down to $\quest\alpha$ and
  $\square B$. Hence everything unresolved has gathered on the problem side. The known
  gaps are Rose's formula (Corollary~\ref{cor:incomplete}) and, possibly,
  Peirce's law with a quantifier (Proposition~\ref{prop:pienc}(4),
  Remark~\ref{rem:peirce}). The nearest check:
  whether there remain valid schemes
  of the same kind, where the solver has to guess the label of a disjunction while the list
  of candidates is not finite; see Section~\ref{sec:whyclasses}.

\item \textbf{Is the rule ``from $\bang A$ infer
  $A$'' admissible in $\HCB$?} This rule preserves $\sB$-validity without any
  assumptions: if $\ZF\vdash\mathrm{Prf}(n,\gn{\forall\vec X A^*})$, then,
  provided $\ZF$ is consistent, this $\Delta_0$-formula is true and
  $\ZF\vdash\forall\vec X A^*$, while if $\ZF$ is inconsistent everything is valid.
  In $\HCE$ the rule is derived ((PtoA), then (A8)); in $\HCI$ it is admissible:
  from $\provgl\bang A$, by Theorem~\ref{th:sound}, the
  $\sI$-validity of $A$ follows, and by Corollary~\ref{cor:solovay} its
  derivability. Is it true that $\HCB\vdash\bang A$ implies
  $\HCB\vdash A$? Adding the rule to $\HCB$ yields no new problems
  (by induction on the derivation: for every proposition $A$ derivable with this
  rule, the problem $\bang A$ is derivable in $\HCB$), so the question concerns
  only propositions with the quantifier $\Pi$. A negative answer would mean
  that $\HCB$ lacks some law: under the sign $\bang$
  the calculus derives more than outside it.

\remove{
\item \textbf{Внутренний уровень:} совпадают ли внутренняя и внешняя
  общезначимость на каких-нибудь естественных фрагментах (для задач вида
  $\bang A$ совпадают, замечание~\ref{rem:whycomp}).
}

\item \textbf{The modality-free fragment of $\HCgls$.} Does the set of
  modality-free theorem-problems of $\HCE$ coincide with that of $\HCI$? All that is known is
  $\HCI\subseteq\HCE$ and that both fragments lie in $\Lp_{\sE}$
  (\S\ref{sec:indep}); the frame with a marked maximal point does not settle the
  question---the schemes (A8) and (A9) are false in it.

\item \textbf{Independence of (A7) in $\HCE$.} For $\HCI$
  the question is settled: Proposition~\ref{prop:a7indep} (Appendix~\ref{sec:a7indep})
  gives a two-point Kripke model in
  which (A1)--(A6) and all the rules hold, while
  $\pneg\pneg\bang\bot\pto\bang\bot$ does not. For $\HCgls$ one cannot proceed by models in this way
  (Remark~\ref{rem:a10fork}). Note that syntactically the problem reduces to
  the following: the theorem-problems of $\HCgls$ are the closure under (MP) of the set of
  theorem-problems of $\HCgl$ together with all instances of (A9), since without (AtoP) there is
  no passage from a proposition to a problem. Hence one needs a structure in
  which all theorem-problems of $\HCgl$ receive $1$ (but the rule (AtoP) as a whole
  is not preserved), and all instances of (A9) do as well.

\item \textbf{Independence of (A9).} In the opposite direction the question is settled:
  (A8) is not derivable from (A9) (Proposition~\ref{prop:a8indep}), although
  (A9) yields reflection for all $\square$-formulas with modalized
  propositional variables. The difficulty with (A9) itself is that its
  $\quest$-image is a special case of (A8) (Proposition~\ref{prop:nf}):
  only a construction that distinguishes a problem from
  its $\quest$-image can tell (A9) apart from (A8).

\item \textbf{How far does (A7) reach.} It is adopted for $\Sigma$-problems;
  Appendix~\ref{sec:whyclasses} shows that neither variables nor
  disjunctions of arbitrary problems can be added. The exact class of problems $\varphi$
  for which $\pneg\pneg\varphi\pto\varphi$ is valid has not been described. Note that
  this class depends on the syntax: from $\gn{*}$ one must be able to construct
  an enumeration of $\val{\varphi}$, and for an implication there is nowhere to take it from.

\item \textbf{Does the $\triangledown$-fragment coincide with $\mathrm{iGL}$?}
  (The definitions of the logics $\mathrm{iK}$, $\mathrm{iGL}$, $\mathrm{iGLC}$ are in
  \S\ref{sec:hc}.) The question of coincidence is settled negatively: by
  Proposition~\ref{prop:igl}, $\mathrm{iGL}$ is contained in the
  $\triangledown$-fragment of $\HCgl$ \emph{strictly}---already
  $\triangledown(\varphi\por\pneg\varphi)$ is derivable in our system and not derivable in
  $\mathrm{iGL}$. It remains to find an exact characterization; a plausible conjecture
  is written out in Remark~\ref{rem:igl}.

  On the other hand, $\mathrm{iGLC}$ is \emph{not} our fragment: the completeness principle
  $\varphi\pto\triangledown\varphi$ is co-reflection, and it has been refuted
  (Proposition~\ref{prop:hce}(4)). Nevertheless there is a similarity with the arithmetical
  part of~\cite{VZ19}, and a substantive one: there completeness is required not of all
  formulas but of $\Sigma_1$ ones (the condition $\square\triangle S\ato\square S$ for
  $S\in\Sigma_1$)---in our setting, in exactly the same way, co-reflection holds precisely on
  $\Sigma$-problems (axiom (A5)). It is worth understanding
  whether this coincidence is accidental.

  \emph{Proviso.} We are speaking of \emph{derivable} formulas. The set of
  \emph{valid} $\triangledown$-formulas is strictly wider than $\mathrm{iGL}$ already
  because Rose's formula is valid (Proposition~\ref{cor:rose}), and it is
  not derivable even in $\Int$.

\remove{
\item \textbf{Не следует ли (A6) из теоремы реализации.} Схема (A6) говорит о
  вычислимом преобразовании выводов; в логике доказательств такому отвечает
  явный аналог $\GL$ с теоремой реализации~\cite{Artemov01,Sidon97,Nogina14}.
  Стоит проверить, нельзя ли получить (A6) прямо оттуда, а не выводом из
  неподвижной точки, как при выводе (A6) в теореме~\ref{th:sound}.
}

\item \textbf{An intuitionistic version of the semantics.} Since the divergence from
  $\mathrm{iGL}$ is created by the Boolean nature of $\quest$,  
  it is tempting to replace $\PA$ by $\mathrm{HA}$. It is useful to see in advance that
  this is not enough: the Boolean nature has two independent sources. The \emph{internal} one is
  the classicality of $\PA$: the right half of (A4a) was derived by a case analysis
  ``either $\aneg\exists u\,\varphi^*(u)$ or $\exists v\,\psi^*(v)$'', and
  $\mathrm{HA}$ provides no such case analysis, so here the replacement will work.
  The \emph{external} one is the two-valuedness of the standard model itself: external
  validity is defined via truth in~$\N$, and
  $\val{\quest\varphi}=\Tru$ is equivalent to $\val{\varphi}\neq\varnothing$---a
  two-valued fact about~$\N$ under any theory under the sign~$\bang$. To remove
  this source as well, one would have to compute values not in the standard model but in Kripke frames
  or in Heyting-valued models; but then the sort of propositions would cease to be
  classical, that is, we would be dealing no longer with $\HC$ but with its intuitionistic
  relatives. One should also keep the price in mind: with $\mathrm{HA}$ under the sign~$\bang$
  the modality $\square$ ceases to correspond to $\GL$ and becomes the provability logic
  of $\mathrm{HA}$~\cite{AM19}, and from the disjunction property of
  $\mathrm{HA}$ the scheme $\square(A\aor B)\ato\square A\aor\square B$ is added
  for closed $A$, $B$.

\item \textbf{Definable boxes.} Does anything change if the boxes
  range only over arithmetically definable sets
  (Remark~\ref{rem:semS})? In all the facts examined, no;
  the only candidate for a divergence is the full scheme
  $(\Pi\mathrm{E})$, which with definable boxes we are able to
  justify only for $\Pi$-free $\theta$, and the related
  question of substituting $\Pi$-containing problems
  (Lemma~\ref{lem:closurebs}). Here belongs also
  the conservativity of the $\Pi$-extension over $\HCI$ on clean
  quantifier-free problems: for clean propositions it is proved assuming
  $\Sigma_1$-soundness of $\ZF$, while on non-clean formulas it fails
  (Remark~\ref{rem:qfree}). Do
  \emph{ineliminable} $\Pi$-atoms exist (Definition~\ref{def:elim})? All
  those we have encountered reduce to $\Pi$-free propositions
  (\S\ref{sec:pinf}), but for one of them---Peirce's law with the
  quantifier over $\alpha$---the reduction is proved only semantically, and
  the derivability of the problem
  $\Pi\alpha(((\alpha\pto\beta)\pto\alpha)\pto\alpha)\pto\beta$ in
  $\HCB$ is a separate open question (Proposition~\ref{prop:pienc}(4),
  Remark~\ref{rem:peirce}). For ineliminable atoms to exist
  it suffices that the modality-free
  black-box logic be not computably enumerable; its complexity
  is a separate question. Finally, \emph{returning to $\PA$}: for which
  $\varphi$ is the translation $(\quest\Pi\vec\alpha\,\varphi)^*$ equivalent to
  an arithmetical formula? This is proved for the class of problems in which all
  occurrences of variables lie under $\pneg$ (the dependence on the boxes
  factors through a finite vector of emptiness bits); the general case is
  open; with definable boxes the translation is arithmetical relative to
  $0^{(\omega)}$.

\item \textbf{The two-theory level.} Carry out a full audit of the pair
  $(\PA,\ZF)$ (modality---provability in $\PA$, validity---
  derivability in $\ZF$)---it returns to
  the calculus $\HCE$ a model-independent syntactic semantics---
  compare it with Solovay's set-theoretic interpretation
  (the logic $SOL$,~\cite[Ch.~13]{Boolos93},~\cite{EG08}) and choose a
  verification theory for $\sB$ (for example,
  $\ZF$ plus an inaccessible cardinal). Adjacent to this is the question of the
  dual operator $\Sigma\alpha$ (value
  $\bigcup_{A}$): is it needed for the prenex normalization of
  $\Pi$-atoms, and what is its logic.

\item \textbf{The logic of $\Pi$-problems without modalities.} The most
  unexpected of what we have found happens on the problem side without any
  participation of propositions: the fragment of the language of $\sB$ consisting of problems built from
  variables and $\varnothing$ by the connectives $\pto,\pand,\por$ and the quantifier
  $\Pi$, together with its validity in $\sB$, form a logic that we have not
  identified. Syntactically it is second-order intuitionistic
  propositional logic with a single universal quantifier;
  semantically it is Kleene realizability with uniform quantification
  over sets, the same as the one used to interpret second-order arithmetic
  $\mathrm{HAS}$; the axiom $(\Pi\por)$ is the propositional form of
  Troelstra's uniformity principle. On $\Pi$-free problems it is
  $\Lp_{\sB}$. In all the cases we have examined (all problems
  of depth at most three in two variables) the quantifier is eliminable: for example,
  $\Pi\alpha\bigl((\alpha\pto\beta\por\gamma)\pto((\alpha\pto\beta)\por
  (\alpha\pto\gamma))\bigr)$ is equivalent to the linearity axiom
  $(\gamma\pto\beta)\por(\beta\pto\gamma)$, and
  $\Pi\beta\bigl(((\alpha\pto\beta)\pto\beta)\pto(\alpha\por\beta)\bigr)$
  to the law $\alpha\por\pneg\alpha$; this is akin to Pitts' uniform
  interpolation for $\Int$, but the elimination sometimes relies on
  boxes like $\omega\setminus\{n\}$, which never occur as the value of a formula
  (Proposition~\ref{prop:pienc}(4), Remark~\ref{rem:peirce}). Questions: is the quantifier
  always eliminable; is this logic computably enumerable (for the predicate
  logic of realizability Plisko proved non-arithmeticity); which
  principles beyond $(\Pi\mathrm{E})$, Bernays' rule and $(\Pi\por)$ hold
  in it---the first candidate is
  $\Pi\alpha(((\alpha\pto\beta)\pto\alpha)\pto\alpha)\peq\beta$.
\end{enumerate}

\appendix

\section{Deferred proofs}\label{sec:app}

\begin{proof}[Proof of Lemma~\ref{lem:dn}]
By Lemma~\ref{lem:clauses}(6) the set $\val{\pneg\varphi}$ is $\N$ when
$\val{\varphi}=\varnothing$ and $\varnothing$ otherwise; applying the same once more,
we obtain what is required. For the second claim we unfold the definitions:
$(\pneg\varphi)^*(u)$ is $\forall v\,(\varphi^*(v)\ato \arbot )$---a
formula not depending on $u$ and equivalent to $\aneg\exists v\,\varphi^*(v)$; hence
$(\pneg\pneg\varphi)^*(x)$ is equivalent to $\aneg\aneg\exists v\,\varphi^*(v)$, and
$\PA$ is classical.
\end{proof}

\begin{proof}[Justification of Remark~\ref{rem:indep}]
For $\gn{*}$ this is obvious: any two effective codings of one and the same
finite collection of formulas are translated into each other computably, and the translator can
be prepended to $S$.

For $U$ and $c$ the argument is slightly longer. Let $\val{\cdot}$ and
$\val[*]{\cdot}'$ be the values computed with respect to two sets of choices, and let
$\varphi^*$ and $\varphi^{*\prime}$ be the corresponding translations. By simultaneous
induction on the construction of the formula we build a \emph{pair} of computable
(uniformly in $\varphi$ and in $\gn{*}$) functions $h_\varphi$ and $g_\varphi$
translating solutions in both directions:
\[
  u\in\val{\varphi}\ \Longrightarrow\ h_\varphi(u)\in\val[*]{\varphi}' ,
  \qquad
  v\in\val[*]{\varphi}'\ \Longrightarrow\ g_\varphi(v)\in\val{\varphi} ,
\]
and moreover this very claim is provable in $\PA$, and its proof is constructed from
$\gn{*}$. For conjunction and disjunction both functions permute the labels, replacing
$c$ by~$c'$ and back; for implication $h_{\varphi\pto\psi}$, given an index $i$,
outputs the index of the program $v\mapsto h_\psi\bigl(U_{i}(g_\varphi(v))\bigr)$ in
the new numbering (the parametrization theorem), and symmetrically for
$g_{\varphi\pto\psi}$. Here one sees why both functions are needed at once: for
implication the direction of translation is reversed on the premise. The case $\bang A$ is handled as follows: by
the induction hypothesis $\PA$ proves
$\exists x\,\varphi^*(x)\aeq\exists x\,\varphi^{*\prime}(x)$ for
all problems $\varphi$ occurring in $A$ under the sign~$\quest$, and hence also
$A^*\aeq A^{*\prime}$; a proof of this equivalence is computed from
$\gn{*}$, and a proof of the formula $A^*$ is transformed into a proof of the formula
$A^{*\prime}$.

Hence if $S$ is a computable section for the first set of choices, then
$\gn{*}\mapsto h_\varphi(S(\gn{*}))$ is a computable section for the second. The same
argument, carried out inside $\PA$, gives independence for
$\sI$-validity as well.

The same argument goes through when the formula $\Prf$ is replaced by any other standard
provability formula for $\PA$: proofs of one system are transformed into proofs of the
other computably, and the correctness of the translation is provable in $\PA$. What matters
is only that $\Prf$ defines a \emph{decidable} relation and that $\PA$ proves
the formalized derivability conditions.
\end{proof}

\begin{proof}[Proof of Lemmas~\ref{lem:closure} and~\ref{lem:closurebs}]
\emph{(MP)}. For propositions the rule is obvious: provability in $\PA$ or $\ZF$ and
truth in $\N$ are closed under modus ponens; in $\sB$---together with
the closure (from $\ZF\vdash\forall\vec X\,A^*$ and
$\ZF\vdash\forall\vec X\,(A^*\ato B^*)$ it follows that
$\ZF\vdash\forall\vec X\,B^*$). For problems: if $S_1$ solves $\varphi$ and
$S_2$ solves $\varphi\pto\psi$, then $S(c):=U_{S_2(c)}(S_1(c))$ is computable and
by Lemma~\ref{lem:clauses}(6) yields an element of $\val{\psi}$. In $\sB$ the solutions
are assembled in the same way (they do not depend on the boxes), and a proof of the closure for
$\psi$ is obtained from the proofs of the closures for the premises by the same scheme.

\emph{Substitution.} The semantics $\sI$ and $\sE$:
 let $\varphi(\alpha)$ be valid with the function $S$, and let us
substitute a formula $\chi$ for $\alpha$. The translation $(\varphi(\chi))^*$
coincides with the translation $\varphi(\alpha)^{*'}$, where the interpretation $*'$ coincides with $*$ on
the other variables and assigns to the variable $\alpha$ the formula $\chi^*$. The code
$\gn{*'}$ is computed from $\gn{*}$, since the translation $\chi\mapsto\chi^*$ is
a computable operation on codes (Definition~\ref{def:tr}). Hence
$S'(\gn{*}):=S(\gn{*'})$ is the required function.

Substitution of a proposition $C$ for a propositional variable $A$ in
the semantics $\sI$ and $\sE$ is handled by the same composition of interpretations
($A^{*'}:=C^*$).

In the semantics $\sB$, substitution of a problem $\chi$ for a free
variable $\alpha$ is an instantiation of the quantifier. Free
occurrences of $\alpha$ do not lie under the sign $\bang$ (they would be bound),
so the translation refers to $\alpha$ only through atoms $x\in X_\alpha$.
From $\ZF\vdash\forall\vec X\,\forall X_\alpha\,
\varphi^*\bigl(\overline{S(\gn{*})}\bigr)$, substituting the set
$\{x:\chi^*(x)\}$ for $X_\alpha$, we obtain the closure of the translation of
$\varphi[\chi/\alpha]$: the free set variables of the formula
$\chi^*$ go under the common closure $\forall\vec X$, and there is no capture by
property (1) of safe substitution. The solver need not be changed, since the boxes do not
enter the code. In $\sB$ the set $\{x:\chi^*(x)\}$ is a legitimate content of a
box for any $\chi$, including one with the quantifier $\Pi$: the boxes
range over \emph{all} sets (impredicativity).

Substitution of a proposition $C$ for a propositional variable $A$
in a formula
is handled in two ways. If $C$ has no free
problem variables, then $C^*$ is a closed formula of the language of $\ZF$, that
is, an admissible interpretation, and the composition of interpretations
$A^{*\prime}:=C^*$ works, including for occurrences of $A$ under the sign $\bang$.

If, however, $C$ has free problem variables, then by property (2)
all occurrences of $A$ lie outside $\bang$. This is possible only if the formula
is a proposition $B(A)$ (propositions enter problems only through
$\bang$), and for the same reason the occurrences of $A$ do not lie under the sign
$\quest$ either. Therefore $(B(C))^*$ is obtained from $B^*$ by replacing the occurrences of $A^*$
by $C^*$, and none of them stands inside a G\"odel number. Since
$B(A)$ is valid, and $(0=0)$ and $\arbot $ are admissible values of the
variable $A$, $\ZF$ proves the closures of $B^*(0=0)$ and $B^*\arbot $.
Reasoning in $\ZF$ with arbitrary contents of the boxes, we distinguish cases:
if $C^*$ is true, the formula $C^*$ is equivalent to $(0=0)$; if $C^*$ is false, to the formula
$\arbot $. By replacement of equivalents, $B^*(C^*)$ is equivalent to $B^*(0=0)$ or to
$B^*\arbot $. The replacement is legitimate also under the quantifiers $\forall X_\beta$ corresponding to
$\Pi\beta$: by property (1), $X_\beta$ does not occur in $C^*$. Consequently,
$\ZF$ proves the closure of $B^*(C^*)=(B(C))^*$.
\end{proof}

\begin{proof}[Proof of Theorem~\ref{prop:modalfree}]
One direction is the scheme (A2). The converse is obtained by a model.

Take a finite Kripke frame, that is, a finite partially ordered
set $\langle P,\le\rangle$, and fix in it a \emph{maximal}
point~$p$. We build a model $K$. Problem variables are evaluated by
arbitrary upward closed subsets of $P$ (as in Kripke models
for $\Int$), propositional variables by arbitrary truth
values. The forcing $x\Vdash\varphi$ of problems at points $x\in P$
is defined by the usual Kripke clauses for $\varnothing$, $\pand$, $\por$,
$\pto$; propositions receive truth values, their connectives
are understood classically. The modalities:
\[
  x\Vdash\bang A\ \ \text{for all $x$ and all $A$},
  \qquad
  \quest\varphi\ \text{is true}\iff p\Vdash\varphi .
\]
The definition is correct: $\bang A$ is forced everywhere independently of $A$,
so there is no circular dependence between the sorts; forcing is
monotone with respect to $\le$ (induction on the construction, trivial for $\bang A$).
A problem is \emph{true in $K$} if it is forced at all points under every
valuation of the variables; a proposition, if it is true under every valuation.

We check the axioms. (A1): propositions are two-valued and the connectives are classical;
(A2): the schemes of $\Int$ are true in every Kripke model. The schemes (A3), (A5) and (A6)
are trivial: the conclusion of each of them is a problem of the form $\bang X$ (for (A3),
after discharging both premises), forced everywhere, and an implication with a
conclusion forced everywhere is itself forced everywhere.

The scheme (A4) is the only place where the choice of the point is at work. Since $p$
is maximal, the only point above it is itself, so forcing
at~$p$ is classical:
\[
  p\Vdash\varphi\pto\psi\iff
  \bigl(p\Vdash\varphi\Rightarrow p\Vdash\psi\bigr),
\]
and likewise $p\Vdash\varphi\pand\psi$ is equivalent to the conjunction, and
$p\Vdash\varphi\por\psi$ to the disjunction; finally $p\nVdash\varnothing$.
This is exactly (A4a)--(A4d): the map $\varphi\mapsto\quest\varphi$ is a Boolean
homomorphism.

The scheme (A7): a $\Sigma$-problem is assembled from $\varnothing$ and problems of the form
$\bang X$ by the connectives $\pand$ and $\por$, hence (induction on the construction)
it is either forced everywhere or forced nowhere. For such problems
$\pneg\pneg\sigma\pto\sigma$ is forced everywhere: if $\sigma$ is forced
everywhere, this is clear; if nowhere, then $\pneg\sigma$ is forced everywhere,
$\pneg\pneg\sigma$ nowhere, and the implication is forced vacuously.

The rules. (MP): if $\varphi$ and $\varphi\pto\psi$ are forced at a point, then
so is $\psi$; for propositions this is classical logic. (AtoP): $\bang A$
is forced everywhere for any $A$, all the more for a derivable one. (PtoA): if
$\varphi$ is forced everywhere, then in particular at $p$, hence $\quest\varphi$ is true.

Thus all axioms and rules of $\HCI$ are true in this model. (Rose's formula
is, of course, not satisfied by the model---it is not derivable in $\Int$; this is precisely
what makes the model useful.) Now let a problem
$\varphi$ contain no modalities. Then its forcing is the usual
forcing in the Kripke frame~$P$.
If $\varphi$ is derivable in $\HCI$, then it is true in every such model, that is,
in every finite Kripke frame, and hence derivable in $\Int$ by the completeness theorem
for intuitionistic logic with respect to finite frames.
\end{proof}
\section
{Why exactly these classes of problems in (A5) and (A7)}
\label{sec:whyclasses}

\label{rem:whysigma}
For the scheme (A5) the answer is short: beyond the class $\Sigma$, co-reflection
$\varphi\pto\bang\quest\varphi$ fails already for a variable
(Proposition~\ref{prop:hce}(4), where $\alpha^*:=\Con$), and semantically
(A5) is effective $\Sigma_1$-completeness with the witness on the input
(Proposition~\ref{prop:a5a7}), which rests on the fact that the translations
of $\Sigma$-problems are $\Sigma_1$-formulas. With the scheme (A7) matters are more complicated.
The class of $\Sigma$-problems cannot be extended in a natural way. For example, it
cannot be replaced by the class $\mathcal O$ from Rose's scheme: there the scheme (A7) fails
already on $\bang A\por\pneg\bang A$, as we shall see below in this section.

First of all, the scheme (A7) ceases to be valid on variables: the problem
$\pneg\pneg\alpha\pto\alpha$ with a variable $\alpha$ is not valid
(Proposition~\ref{prop:noem}).

For disjunctions, even the scheme
$\pneg\pneg(\varphi\por\psi)\pto(\pneg\pneg\varphi\por\pneg\pneg\psi)$ is not valid.
The word ``even'' is justified as follows. A direct extension of (A7) to disjunctions
of arbitrary problems would be the scheme $\pneg\pneg(\varphi\por\psi)\pto
(\varphi\por\psi)$, whereas the conclusion of the displayed one is weaker: in $\Int$ we have
$\varphi\pto\pneg\pneg\varphi$, hence $(\varphi\por\psi)\pto
(\pneg\pneg\varphi\por\pneg\pneg\psi)$. A solution of the displayed scheme is required
to provide not a solution of one of the disjuncts, but merely a \emph{label}---an indication of which
of them has a nonempty value (Lemma~\ref{lem:dn}). Nevertheless, this scheme is not valid
either, and then by Lemmas~\ref{lem:closure} and~\ref{lem:closurebs} the direct extension
is not valid as well.
Take variables $\alpha$, $\beta$ and put
$\alpha^*(x):=\Prf(x,\sigma)$, $\beta^*:=\aneg\Prov(\sigma)$. Then
$\val{\alpha\por\beta}\neq\varnothing$ for every $\sigma$, hence
$\val{\pneg\pneg(\alpha\por\beta)}=\N$ by Lemma~\ref{lem:dn}; while by
Lemma~\ref{lem:clauses}(5),(6) an element of the set
$\val{\pneg\pneg\alpha\por\pneg\pneg\beta}$ carries a label distinguishing the cases
$\val{\alpha}\neq\varnothing$ and $\val{\beta}\neq\varnothing$, that is,
deciding the provability of~$\sigma$.

One can also see exactly where the proof breaks down: the enumeration on which
Lemma~\ref{lem:sigmaenum} rests is nowhere to be found for a variable and for an implication.
For a variable
$\alpha$ the set $\val{\alpha}=\{n:\N\models\alpha^*(n)\}$ is arbitrary: the
formula $\alpha^*$ can be anything whatsoever, and there is nothing to enumerate. For
an implication the condition ``$x\in\val{\psi\pto\chi}$'' is, by
Lemma~\ref{lem:clauses}(6), the condition ``$U_x$ is defined everywhere on $\val{\psi}$ and
lands in $\val{\chi}$''---a condition with a universal quantifier, again not
enumerable.

Now about the promised class $\mathcal O$ from Rose's scheme
(Proposition~\ref{cor:rose}): on it the scheme (A7)
also ceases to be valid. It suffices to take
$\varphi:=\bang A\por\pneg\bang A$---the first disjunct is a $\Sigma$-problem, the second
has the form $\pneg\psi$, so that $\varphi$ is an $\mathcal O$-problem. Under every
interpretation $\val{\varphi}\neq\varnothing$: if $A^*$ is provable, then
$\val{\bang A}$ is nonempty, and if not, then $\val{\pneg\bang A}=\N$. Hence
$\val{\pneg\pneg\varphi}=\N$ always, and the scheme (A7) for this $\varphi$
would require computing from the code a solution of $\varphi$---that is, a label answering
the question whether $A^*$ is provable. This is exactly the counterexample analyzed above with
a disjunction whose label would decide provability---written inside the
class $\mathcal O$. So it is no accident that the scheme (A7) and Rose's scheme are stated on different
classes: Rose's scheme requires only a \emph{finite list of candidates}---they
are written out in advance, and which of them is suitable is learned from the interlocutor who solves
the premise---whereas (A7) requires a \emph{search} that finds a witness, that
is, enumerability of the value. The class $\mathcal O$ provides the former and not
the latter: as soon as a negation enters a disjunction, enumerability is lost.
Inside $\Sigma$, on the other hand, disjunction preserves enumerability---this is why (A7)
works there, while Rose's scheme, conversely, degenerates: on $\Sigma$-problems
it is derivable from (A7) by a single (MP).

Let us now ask the converse question: can the class be narrowed without losing anything? In other
words, which of the two admitted connectives really strengthen the scheme (A7)?
It turns out that only disjunction does.

A problem $\varphi$ is called \emph{stable} if the double negation elimination scheme
$\pneg\pneg\varphi\pto\varphi$ holds for it; in these words the scheme (A7)
says that all $\Sigma$-problems are stable. A connective is said to
\emph{inherit} stability if stability of the arguments implies
stability of the result; in $\Int$ stability is inherited by $\pand$ and $\pto$,
but not by $\por$.

Conjunction therefore does not strengthen the scheme (A7): the case
$\varphi=\psi\pand\chi$ is derived from the cases $\psi$ and
$\chi$, since already in $\Int$ from $\pneg\pneg(\psi\pand\chi)$ follows
$\pneg\pneg\psi\pand\pneg\pneg\chi$. Disjunction does strengthen it: the scheme
$\pneg\pneg(\bang A\por\bang B)\pto(\bang A\por\bang B)$ is no longer derivable in $\Int$ from $\pneg\pneg\bang A\pto\bang A$ and
$\pneg\pneg\bang B\pto\bang B$: even
$\pneg\pneg\alpha\pto\alpha$ and $\pneg\pneg\beta\pto\beta$ do not yield
$\pneg\pneg(\alpha\por\beta)\pto(\alpha\por\beta)$ (the Kripke frame with a root and
two incomparable successors, in which $\alpha$
and~$\beta$ respectively are true). Meanwhile
Lemma~\ref{lem:sigmaenum} does give this scheme: a parallel exhaustive search for proofs of both
formulas. Thus the whole content of (A7), beyond Markov's principle
$\pneg\pneg\bang A\pto\bang A$, consists in that \emph{inside the class $\Sigma$
stability is declared to be inherited by disjunction}.

Conjunction, however, still cannot be removed from the grammar: without it one
cannot build disjunctions like $(\bang A\pand\bang B)\por\bang C$, for which the scheme (A7)
is, by what has just been said, no longer derivable.

%

\section{Independence of the scheme (A7)}\label{sec:a7indep}

The scheme (A7) is not derivable from the remaining axioms of $\HCgl$. The proof is
a Kripke model with two points.

\begin{definition}[two-point model $M$]\label{def:model3}
The frame is the two-element chain $a<b$. Problem variables are evaluated by
arbitrary upward closed sets ($\varnothing$, $\{b\}$,
$\{a,b\}$), propositional variables by arbitrary truth
values; the connectives of problems are understood according to the Kripke clauses, the connectives
of propositions classically, $\varnothing$ is forced nowhere, $\bot$
is false. The modalities:
\[
  x\Vdash\bang A\iff\bigl(x=b\ \text{ or }\ A\ \text{is true}\bigr),
  \qquad
  \quest\varphi\ \text{is true}\iff b\Vdash\varphi .
\]
In other words, $\bang A$ is always forced at the upper point, and at the lower one
if and only if $A$ is true (forcing is monotone); $\quest\varphi$
means that $\varphi$ is forced at least at one point, since every
nonempty upward closed set contains $b$. A formula is \emph{true in
$M$} if under every valuation of the variables it is forced at both points (for
problems) or is true (for propositions).
\end{definition}

\begin{proposition}\label{prop:a7indep}
All the schemes (A1)--(A6) are true in $M$, and all three rules (MP), (AtoP), (PtoA)
preserve truth in $M$. Meanwhile the problem $\pneg\pneg\bang\bot\pto\bang\bot$---a
special case of (A7)---is not true in $M$: it is not forced at the point $a$.
Consequently, the scheme \textup{(A7)} is not derivable in $\HCgl$ from the remaining axioms.
\end{proposition}

\begin{proof}
Note at once: $\bang A$ is forced at $b$ for every $A$, hence
\[
  \square A=\quest\bang A\quad\text{is true for all }A .
\]
In this model ``everything is provable''; this is precisely why Löb's scheme holds in it for free.

Two general observations. First, the point $b$ is maximal, so
forcing at it is classical: $b\Vdash\varphi\pto\psi$ is equivalent to
$(b\Vdash\varphi\Rightarrow b\Vdash\psi)$, and likewise for $\pand$, $\por$.
Second, for problems of the form $\bang X$ forcing at $a$ is the truth
of $X$, while at $b$ it is automatic; therefore for problems built from
problems $\bang X$ by implications, forcing at $a$ reduces to an implication
between truth values: $a\Vdash\bang X\pto\bang Y$ is equivalent to
$(X\ \text{is true}\Rightarrow Y\ \text{is true})$, and likewise for
nested implications.

(A1) Propositions are two-valued, the connectives are classical.

(A2) The schemes of $\Int$ are true in every Kripke model.

(A3) At $b$ everything is forced. At $a$, by the second observation, we need: if
$A\ato B$ and $A$ are true, then $B$ is true.

(A4) By the first observation $\quest\varphi$ is true if and only if
$b\Vdash\varphi$, and forcing at $b$ is classical; this gives (A4a)--(A4c),
while (A4d) is $b\nVdash\varnothing$.

(A5) If $\quest\sigma$ is true, then $\bang\quest\sigma$ is forced everywhere
and the implication $\sigma\pto\bang\quest\sigma$ is forced everywhere. If
$\quest\sigma$ is false, then $b\nVdash\sigma$, and then by monotonicity
$\sigma$ is forced nowhere, and the implication is forced vacuously.

(A6) $\square A$ is true, so $\square A\ato A$ is true if and
only if $A$ is true; by the second observation
$a\Vdash\bang(\square A\ato A)\pto\bang A$, while at $b$ everything is forced.

The rules. (MP) is the usual property of forcing and of classical implication.
(AtoP): if $A$ is true, then $\bang A$ is forced at both points. (PtoA):
if $\varphi$ is forced at both points, then in particular at $b$, and $\quest\varphi$
is true.

Finally, $\bot$ is false, so $\bang\bot$ is forced only at $b$. Then
$\pneg\bang\bot$ is forced nowhere: at $b$ because $b\Vdash\bang\bot$,
at $a$ because $b$ lies above $a$. Hence $\pneg\pneg\bang\bot$
is forced everywhere, and at the point $a$ the premise of the problem
$\pneg\pneg\bang\bot\pto\bang\bot$ is forced while the conclusion is not.
\end{proof}

\section{What is new in the derivation of L\"ob's scheme}\label{sec:novelty}

Here we examine whether there is anything new
in our proof of L\"ob's scheme (item~(3) of the proof of
Theorem~\ref{th:sound}).
In the modal language such a derivation is known and standard: in
$\GL$ without L\"ob's scheme the fixed point $\psi\aeq(\square\psi\ato
\sigma)$ yields L\"ob's formula. Indeed, from $\psi\ato(\square\psi\ato\sigma)$
by the necessitation rule and the scheme K we obtain
$\square\psi\ato\square(\square\psi\ato\sigma)$, then by K and the scheme
$\square\psi\ato\square\square\psi$ the implication $\square\psi\ato\square\sigma$;
hence $(\square\sigma\ato\sigma)\ato(\square\psi\ato\sigma)$, and therefore also
$\square(\square\sigma\ato\sigma)\ato\square(\square\psi\ato\sigma)$, which together with
$\square(\square\psi\ato\sigma)\ato\square\psi$ (the second half of the fixed
point) gives $\square(\square\sigma\ato\sigma)\ato\square\sigma$.

Therefore what is new in Theorem~\ref{th:lob} is not that L\"ob's scheme is derived from
a fixed point, but \emph{in what form} it is derived. The scheme (A6)
is written via $\bang$: it asserts the existence of a \emph{computable
transformation of proofs}, not the existence of proofs, and it is
inexpressible in the modal language. From (A6) the $\square$-version is obtained by applying (PtoA) and (A4a);
the converse passage is not obvious and remains open (see~\S\ref{sec:todo}).

In substance, what we get is this: \emph{the constructive L\"ob theorem is derived
from effective $\Sigma_1$-completeness}. Usually it is proved by going through the
proof of L\"ob's theorem step by step and keeping track that all transformations of proofs
are computable; here it suffices to have effective $\Sigma_1$-completeness---the scheme (A5)---and
the effectiveness of the fixed-point lemma, and the rest is done by the formal derivation.
The function converting a proof of the formula $\Prov(\sigma)\ato\sigma$ into a proof of
the formula $\sigma$ need not be constructed by hand: it is read off from this derivation.

\clearpage
\selectlanguage{russian}
\makeatletter
\let\maketitle\saved@maketitle
\let\@maketitle\saved@@maketitle
\let\title\saved@title
\let\author\saved@author
\let\date\saved@date
\let\thanks\saved@thanks
\let\and\saved@and
\renewcommand\tableofcontents{\section*{\contentsname}\@starttoc{rtoc}}
\let\orig@addcontentsline\addcontentsline
\renewcommand\addcontentsline[3]{\orig@addcontentsline{r#1}{#2}{#3}}
\@ifpackageloaded{hyperref}{%
  \def\Hy@bookmarkstype{rtoc}%
  \def\Hy@chapapp{section}%
  \renewcommand{\theHsection}{ru.\arabic{section}}%
  \renewcommand{\theHsubsection}{ru.\arabic{section}.\arabic{subsection}}%
  \renewcommand{\theHsubsubsection}{ru.\arabic{section}.\arabic{subsection}.\arabic{subsubsection}}%
  \renewcommand{\theHtheorem}{ru.\arabic{theorem}}%
  \renewcommand{\theHdefinition}{ru.\arabic{definition}}%
  \renewcommand{\theHremark}{ru.\arabic{remark}}%
  \renewcommand{\theHquestion}{ru.\arabic{question}}%
  \renewcommand{\theHexample}{ru.\arabic{example}}%
  \renewcommand{\theHequation}{ru.\arabic{equation}}%
  \renewcommand{\theHtable}{ru.\arabic{table}}%
  \renewcommand{\theHfootnote}{ru.\arabic{footnote}}%
  \let\orig@appendix\appendix
  \renewcommand\appendix{\orig@appendix\renewcommand{\theHsection}{ru.\Alph{section}}}%
}{}
\makeatother
\renewcommand{\thmTheorem}{Теорема}
\renewcommand{\thmProposition}{Предложение}
\renewcommand{\thmLemma}{Лемма}
\renewcommand{\thmCorollary}{Следствие}
\renewcommand{\thmDefinition}{Определение}
\renewcommand{\thmRemark}{Замечание}
\renewcommand{\thmQuestion}{Вопрос}
\renewcommand{\thmConjecture}{Гипотеза}
\renewcommand{\thmExample}{Пример}
\renewcommand{\Tru}{\text{И}}
\renewcommand{\Fls}{\text{Л}}
\setcounter{section}{0}
\renewcommand{\thesection}{\arabic{section}}
\setcounter{theorem}{0}
\setcounter{definition}{0}
\setcounter{remark}{0}
\setcounter{question}{0}
\setcounter{example}{0}
\setcounter{equation}{0}
\setcounter{table}{0}
\setcounter{footnote}{0}
\title{Арифметическая и теоретико-множественная семантики\\
логики задач и высказываний}
\author{Н.\,К.~Верещагин \qquad А.\,А.~Оноприенко\\[6pt]
  \normalsize Московский государственный университет
  имени М.\,В.~Ломоносова}
\date{}

\maketitle

\tableofcontents

\section{Введение}

В работе Колмогорова~\cite{ru:Kolm32} интуиционистское исчисление было 
истолковано как \emph{исчисление задач}. А в комментарии к этой работе в составе 
Избранных трудов~\cite{ru:Kolm85} Колмогоров написал, что задачи и высказывания --- вещи разные и что полное изложение
должно говорить об обоих сортах сразу: интуиционистская логика описывает
задачи, классическая --- высказывания, и нужна общая система, в которой есть
и то и другое вместе со связью между ними.

Примером такой системы является логика $\HC$ Мелихова~\cite{ru:Melikhov13,ru:Melikhov24}
--- совместная логика задач и высказываний. Она соединяет в одном двусортном
языке классические высказывания и интуиционистские задачи, связывая два сорта
двумя модальностями: задача $\bang A$ читается как <<найти доказательство
высказывания $A$>>, а высказывание $\quest\varphi$ --- как <<у задачи $\varphi$
есть решение>>, то есть, именно так, как наметил Колмогоров. 
Для $\HC$ известны семантики нескольких видов: топологические,
интерпретации в планиметрии и в алгебраических уравнениях, неформальная
семантика, восходящая к самому Колмогорову (Мелихов,
\cite{ru:Melikhov13,ru:Melikhov15}), а также семантика Крипке с отмеченными
мирами (Оноприенко, \cite{ru:Onoprienko20}). Нас же интересует семантика,
в которой модальности читаются буквально: решения задачи $\bang A$ ---
это выводы высказывания $A$ в фиксированной формальной теории, а
высказывание $\quest\varphi$ утверждает, что у $\varphi$ есть решение.
Такой --- \emph{арифметической}, в духе реализуемости Клини ---
семантики для полного языка $\HC$ до сих пор не было: для фрагмента без
модальности $\quest$ она намечена в~\cite{ru:Melikhov24}, но осталась
незавершённой.

В настоящей работе такая семантика строится и изучается порождаемое ею
понятие общезначимости, затем вводятся исчисления, корректные
относительно этой семантики; аксиомы отбираются по общезначимости. Сам
порядок <<от толкования к исчислению>> не нов: так построена и логика
$\HC$, аксиомы которой извлечены из колмогоровского толкования задач и
орловско-гейтинговского толкования доказательств~\cite{ru:Melikhov24}.
Разница в том, что у нас толкование --- формальная семантика, и
общезначимость в ней --- точное математическое понятие, а потому вопрос
о корректности и полноте исчисления имеет точный смысл. Этот порядок ---
постановка задачи, семантика и её изучение, исчисления и их изучение
--- определяет и план статьи.

 Точнее, мы определяем три разных семантики. Начнем с первых двух, арифметических,
обозначаемых $\sI$ и $\sE$.
В них высказывание $P$ переводится в замкнутую
формулу $P^*$ арифметики Пеано, задача $\varphi$ --- в формулу $\varphi^*(x)$ с (не более чем) одной свободной переменной;
связки на задачах понимаются по Клини~\cite{ru:Kleene45}, высказывание $\quest\varphi$ понимается как
$\exists x\,\varphi^*(x)$, а решения задачи $\bang A$ --- гёделевы
номера выводов формулы $A^*$ в $\PA$ (\S\ref{ru:sec:val}). При этом оказывается, что не все 
аксиомы HC истинны в этой семантике, то есть, HC некорректна относительно нее.
Самое главное в нашей интерпретации --- понятие общезначимости задачи: решение задачи
требуется предъявлять \emph{алгоритмом, вычисляющим его по коду
интерпретации} (определение~\ref{ru:def:valid}). Как показано в
замечании~\ref{ru:rem:three}, именно от этого требования зависит, останется ли
фрагмент задач интуиционистским: при более слабом понимании общезначимости он
становится классическим. 
В семантике $\sI$ требуется, чтобы при любой
интерпретации формула была выводима в PA, а в семантике $\sE$ --- истинна в стандартной интерпретации PA.
В семантике $\sI$ множество общезначимых высказываний  
по существу совпадает с логикой доказуемости $\GL$,
а множество общезначимых задач является довольно сложной 
логикой, разрешимость которой
неизвестна. 
В семантике $\sE$ множество общезначимых высказываний
по существу совпадает с логикой доказуемости $\GLS$,
а множество общезначимых задач опять является довольно сложной 
логикой, разрешимость которой неизвестна.

Затем мы приводим исчисления (логики), корректные относительно этих 
семантик. 
 Мы не стремились получить исчисления, близкие к $\HC$: от неё
берутся только язык и неформальное понимание модальностей, а аксиомы
отбирались по общезначимости --- если схема общезначима и не выводится из
уже принятых, она добавлялась. Так получились исчисление
$\HCI$, отвечающее тому, что
$\PA$ доказывает, и исчисление $\HCE$, отвечающее тому, что истинно
(\S\ref{ru:sec:calc}); их корректность --- теорема~\ref{ru:th:sound}. 

Зачем нужны такие исчисления? Например, используя исчисление $\HCI$, можно
получить формальный вывод теоремы Лёба 
из теоремы о неподвижной точке для $\PA$.
Особенность этих исчислений в том, что почти все их аксиомы являются
задачами. Исключений три: аксиомы классического исчисления
высказываний (A1); чисто структурные схемы (A4), утверждающие,
что $\quest$ переводит связки задач в связки высказываний (булев
гомоморфизм); и рефлексия (A8) исчисления $\HCE$ --- она не
структурная, а содержательная, но это высказывание, а не задача.
Все остальные аксиомы --- задачи. Они либо говорят о возможности алгоритмического преобразования
одних арифметических выводов в другие: например, (A5) --- ``по решению задачи $\phi$ (для задач $\phi$ специального вида) можно построить
арифметическое доказательство того, что у $\phi$ есть решение'', (A6) ---
``арифметическое доказательство формулы $\square A \ato A$ 
можно преобразовать в арифметическое доказательство формулы $A$'' (теорема Лёба), либо
о реализуемости некоторых задач без модальности $\quest$ --- аксиома (A7).
На уровне высказываний исчисления $\HCI$ и $\HCE$ полны и совпадают с GL и GLS, а на уровне задач 
опять получаются 
некоторые  
логики, про которые мы не знаем, разрешимы ли они. 

В третьей семантике $\sB$ мы немного расширяем язык, вводя в него квантор всеобщности $\Pi$
по задачным переменным. Задача $\Pi\alpha\,\varphi$ понимается как <<укажи программу, которая решает
задачу $\varphi$ при любой интерпретации задачи $\alpha$>>.
Напомним, что без квантора по $\alpha$ эта программа может зависеть вычислимым образом
от  интерпретации $\alpha^*$ переменной $\alpha$.  
Зачем нам нужен такой квантор? Дело в сильной асимметрии двух фрагментов, задачного и высказывательного.
Казалось бы, у нас есть переводы в обе стороны: модальность $\bang$ из задачи 
$A$ дает высказывание $\bang A$, а модальность $\quest$ действует в обратную сторону.
Но оказывается, что первый перевод является честным: общезначимость задачи $\bang A$
равносильна общезначимости высказывания $A$. 
А вот перевод $\phi\mapsto \quest\varphi$
этим свойством не обладает. Часто доказать в PA высказывание   
$(\quest\varphi)^*$ при любой фиксированной интерпретации $*$
значительно проще, чем доказать разрешимость задачи $\varphi^*$ равномерно по $*$. 
Из-за этого и возникает
асимметрия: множество общезначимых высказываний устроено просто, а множество
общезначимых задач --- нет. 

Например, 
закон исключённого третьего $\alpha\por\pneg\alpha$ не является общезначимым (нет алгоритма, который
по арифметической формуле $\alpha^*(x)$ выяснял бы, истинно ли высказывание
$\exists x \alpha^*(x)$). Но при этом  высказывание $?(\alpha\por\pneg\alpha)$ является
общезначимым, потому что PA является классической теорией и может
использовать закон исключённого третьего для доказательства
формулы $\exists x (\alpha\por\pneg\alpha)^*(x)$ для любой фиксированной
формулы $\alpha^*(x)$.  Другой пример: формула Роуза $\rho$ задает
общезначимую формулу, при этом реализующий ее алгоритм не тривиален.
Но при навешивании модальности $\quest$ получается тривиально
общезначимая формула $? \rho$: при любой интерпретации $*$ формула  $\exists x \rho^*(x)$
легко доказывается в $\PA$ с помощью закона исключённого третьего. 

К сожалению, чтобы записать утверждение о разрешимости задачи $\quest\Pi\alpha\,\varphi$,
нужен квантор всеобщности по арифметическим интерпретациям. В арифметике Пеано такое, видимо, невыразимо, поскольку по теореме Тарского невозможно выразить истинность формулы с данным гёделевым номером.
Поэтому вместо PA мы переходим к  теории множеств ZF. В семантике $\sB$ (от слов <<black box>>)
интерпретациями задачных переменных служат \emph{любые} подмножества $\N$ --- <<чёрные ящики>>,
которые решающая программа не видит; квантор всеобщности по ним --- формула первого порядка языка ZF,
и высказывание <<программа $p$ решает задачу $\varphi$ при любой интерпретации задачи $\alpha$>>
записывается теоретико-множественной формулой. Общезначимость в $\sB$
определяется через выводимость в ZF.
Корректное исчисление для неё устроено как $\HCI$ плюс
кванторный слой: бескванторная часть его аксиоматики --- это в точности
$\HCI$. Целью расширения было
увеличить запас 
общезначимых высказываний, и в расширенном языке это заведомо
удалось: например, общезначимым является
высказывание
$\quest\Pi\beta\Pi\gamma\rho$ --- <<формула Роуза равномерно
разрешима>>, --- невыразимое без квантора. Новые \emph{бескванторные} теоремы
расширение тоже даёт --- например, $\square\quest\alpha\ato\square\bot$, --- но
лишь за счёт того, что в семантике $\sB$ модальность $\bang$ связывает
задачные переменные и читается с замыканием. Для \emph{чистых} формул,
в которых $\bang$ применяется только к замкнутым высказываниям,
картина другая: на стороне высказываний расширение консервативно над
$\HCI$ (при $\Sigma_1$-корректности $\ZF$), а для чистых задач вопрос о
консервативности открыт. 
Мы также не знаем, имеет ли в новой семантике фрагмент задач ту же сложность, что 
и фрагмент высказываний.\footnote{К сожалению, перевод $\phi \mapsto ?(\text{замыкание } \phi)$, который мог бы доказать, что сложность задачного фрагмента не больше сложности высказывательного,  
может переводить общезначимые задачи в необщезначимые высказывания.}

Общезначимость всех аксиом наших исчислений сведена в
таблицу~\ref{ru:tab:axioms} (\S\ref{ru:sec:tables}), из которой устройство
исчислений видно непосредственно.


\paragraph{Место наших исчислений среди близких логик.} Из систем, с которыми стоит сравнивать
получившееся, $\HC$ выделяется одним: у неё тот же синтаксис. Всё остальное
расходится. От $\HC$ мы берём язык и чтение модальностей, но не список теорем:
наши исчисления находятся с HC в общем положении.

На стороне высказываний мы попадаем в классическую логику доказуемости:
$\square$-фрагменты исчислений $\HCI$ и $\HCE$ суть в точности $\GL$ и $\GLS$
Соловея~\cite{ru:Solovay76} (следствие~\ref{ru:cor:solovay}), и двусортность там
ничего нового не даёт. Интуиционистские логики доказуемости описывают другое:
$\mathrm{iGL}$~\cite{ru:VZ19,ru:DvdGM24} содержится в нашем
$\triangledown$-фрагменте строго, а $\mathrm{iGLC}$ нашим фрагментом не
является --- её принцип полноты $A\ato\square A$ есть ко-рефлексия. Впрочем,
сходство есть: там полнота требуется только на $\Sigma_1$, а у нас
ко-рефлексия верна в точности на $\Sigma$-задачах
(аксиома (A5)).

В школе Артёмова~\cite{ru:Artemov01,ru:Sidon97,ru:Nogina14,ru:Protopopescu16} доказательства
тоже присутствуют явно, но иначе: там в язык вводятся термы $t{:}F$
<<$t$ есть доказательство $F$>>, тогда как у нас доказательства --- это
решения задачи $\bang A$, а язык остаётся модальным. Двусортности <<задача есть
множество, решение есть число>> там нет.

Наконец, задачи как множества --- давняя традиция: массовые задачи Медведева и
Мучника и логики реализуемости~\cite{ru:Rose53,ru:Plisko09}. Отличие в том, что там
один сорт и нет модальностей: доказательствам и модальности $\bang$ там места
нет. Зато понятие общезначимости у нас в точности совпадает с одним из
известных: это $\Reff$ --- <<реализацию можно эффективно найти по подставляемым
формулам>>~\cite{ru:VSSCh}. Отсюда и результаты о фрагменте задач
(\S\ref{ru:sec:modalfreeval}): его позитивная часть интуиционистская, а сам он лежит
строго между $\Int$ и логикой слабого закона исключённого третьего.

Самого сочетания --- два сорта, задачи как произвольные арифметически
определимые множества с клиниевскими операциями, $\bang A$ как множество выводов,
$\quest$ как непустота и равномерность по коду интерпретации --- мы в литературе
не нашли.

План статьи повторяет порядок <<задача --- семантика --- исчисления>>.
Семантике и её изучению посвящены \S\ref{ru:sec:val}--\S\ref{ru:sec:validity}:
интерпретации всех трёх семантик, понятие общезначимости и всё, что
мы знаем об общезначимых формулах: общезначимость будущих аксиом,
описание безмодального фрагмента задач, примеры необщезначимых формул,
законы оператора $\Pi$ (\S\ref{ru:sec:piatoms}) и сводная таблица
общезначимости (\S\ref{ru:sec:tables}). Исчисления и их изучение
--- \S\ref{ru:sec:calc}: три исчисления, по одному на семантику,
производные правила и примеры выводов, корректность, что даёт квантор
$\Pi$, нормальная форма и
полнота на стороне высказываний вместе с характеризацией общезначимых
высказываний логиками $\GL$ и $\GLS$ (\S\ref{ru:sec:asscompl}), нормальная
форма с $\Pi$-атомами (\S\ref{ru:sec:pinf}),
независимость аксиом и то немногое, что известно о выводимых задачах. Сравнение с $\HC$ и другими близкими
логиками --- \S\ref{ru:sec:hc}, открытые вопросы --- \S\ref{ru:sec:todo}. В
приложения вынесены отложенные доказательства, разбор классов задач
в схемах (A5) и (A7), независимость схемы (A7) и разбор того, что ново
в нашем выводе схемы Лёба.

\section{Язык и обозначения}\label{ru:sec:lang}
Сначала определим синтаксис для арифметических семантик.
Мы пользуемся языком $\HC$, но меняем буквы, чтобы сорт
формулы был виден сразу:
\begin{itemize}[leftmargin=1.6em,itemsep=1pt]
\item \textbf{большие латинские} буквы $A,B,C,\dots$ --- формулы
  \emph{сорта высказывание};
\item \textbf{греческие} буквы --- формулы \emph{сорта задача}.
\end{itemize}
При этом мы придерживаемся следующего соглашения. Буквы $\alpha$, $\beta$,
$\gamma$ (при необходимости с индексами) обозначают только \emph{переменные}
сорта задача, а буквы $\varphi$, $\psi$, $\chi$, $\theta$ --- произвольные
задачи. Для сорта высказывание отдельных букв под переменные мы не отводим:
большие латинские буквы обозначают, вообще говоря, произвольные
высказывания, а там, где существенно, что буква обозначает именно переменную,
это оговаривается явно.
Формулы строятся так:
\[
  A ::= A_i \mid \bot \mid (A\aand A) \mid (A\aor A)
             \mid (A\ato A)\mid \quest\varphi ,
\]
\[
  \varphi ::= \alpha_i \mid \varnothing \mid (\varphi\pand\varphi) \mid
              (\varphi\por\varphi) \mid (\varphi\pto\varphi) \mid \bang A,
\]
где $A_i$ --- переменные сорта высказывание, $\alpha_i$ --- переменные сорта
задача. 

Константа $\bot$ --- ложное высказывание, константа $\varnothing$ ---
неразрешимая задача (задача с пустым множеством решений). Отрицание и эквиваленция не входят
в язык, а определяются:
\[
  \aneg A:=A\ato\bot ,\ \pneg\varphi:=\varphi\pto\varnothing, \  
 A\aeq  B :=(A\ato B)\ \aand\ (B\ato A) ,\ \phi\peq  \psi :=(\phi\pto \psi)\ \pand\  (\psi\pto \phi).
\]

\medskip
\noindent\textbf{Знаки связок различаются по сортам.} Язык двусортный, и одна
и та же связка на двух сортах значит разное: импликация высказываний ---
обычная классическая, импликация задач --- клиниевская <<по решению посылки
вычислить решение заключения>>. Поэтому каждый сорт получает свои знаки:
\[
\begin{array}{l|c|c}
   & \text{высказывания и арифметика} & \text{задачи}\\[2pt]\hline
  \text{импликация}      & \ato   & \pto\\
  \text{конъюнкция}      & \aand  & \pand\\
  \text{дизъюнкция}      & \aor   & \por\\
  \text{отрицание}       & \aneg  & \pneg\\
  \text{эквивалентность} & \aeq   & \peq\\
  \text{константа лжи}   & \bot   & \varnothing
\end{array}
\]
Высказывания получили классическую пару $\ato$, $\aeq$, задачи ---
стрелки $\pto$, $\peq$; конъюнкция, дизъюнкция и отрицание на задачах
помечены кочергой и тильдой.

Смысл различения виден уже на схеме (A3), которая записывается как
$\bang(A\ato B)\pto(\bang A\pto\bang B)$: внутри скобок стоит импликация высказываний,
снаружи --- две импликации задач, и это разные вещи. Знаки $\ato$, $\aand$, $\aor$, $\bot$,
$\aneg$, $\aeq$ употребляются и внутри \emph{арифметических} формул (в $A^*$,
$\varphi^*$ и им подобных): язык арифметики классический, и заводить для него
третий набор знаков незачем. Связки сорта задача в арифметических формулах не
встречаются никогда, а отличить объектную формулу от арифметической помогает
звёздочка.
Производные модальности: $\square A=\quest\bang A$ (высказывание) и
$\triangledown\varphi=\bang\quest\varphi$ (задача).

\section{Арифметический аппарат}

Фиксируем язык арифметики и стандартную модель $\N$; <<истинно>> всюду значит
<<истинно в $\N$>>. Истинностные значения обозначаем $\Tru$ и $\Fls$.
Фиксируем:
\begin{itemize}[leftmargin=1.6em,itemsep=2pt]
\item гёделеву нумерацию одноместных вычислимых (вообще говоря, частичных)
  функций
  \[
    U_0,\ U_1,\ U_2,\ \dots
  \]
\item вычислимую биекцию $c\colon\N^2\to\N$ (спаривание);
\item арифметическую формулу $\Prf(x,y)$, выражающую <<$x$ есть (номер) вывода
  в $\PA$ формулы с номером $y$>>, и $\Prov(y)=\exists x\,\Prf(x,y)$. Как
  обычно, $\Con=\aneg\Prov\arbot $.
\end{itemize}

\begin{remark}[облегчение обозначений]\label{ru:rem:notation}
Мы \emph{отождествляем формулу с её гёделевым номером}: там, где формально
следовало бы писать $\Prf(x,\ulcorner\sigma\urcorner)$, мы пишем
$\Prf(x,\sigma)$. Точно так же там, где в арифметическую формулу
подставляется значение $n$ вычислимой функции, формально следовало бы
подставлять нумерал $\bar n$ --- стандартный терм, задающий $n$; мы пишем
просто $n$. Наконец, запись $U_i(u)$ внутри арифметических формул читается
через график: $\psi(U_i(u))$ --- это сокращение для $\exists w\,\bigl(U_i(u)
\simeq w\aand\psi(w)\bigr)$, где $U_i(u)\simeq w$ --- $\Sigma_1$-формула,
выражающая <<$U_i$ определена на $u$ и равна $w$>>; в частности, такая запись
влечёт определённость. Все оговорки сделаны только ради читаемости и нигде не
используются по существу.
\end{remark}

\section{Интерпретации}\label{ru:sec:val}

\subsection{Арифметические интерпретации}

\begin{definition}[арифметической интерпретации]\label{ru:def:star}
\emph{Арифметической интерпретацией} 
называется отображение $*$, сопоставляющее
\begin{itemize}[leftmargin=1.6em,itemsep=1pt]
\item каждой переменной $A$ сорта высказывание --- замкнутую арифметическую
  формулу $A^*$;
\item каждой переменной $\alpha$ сорта задача --- арифметическую формулу
  $\alpha^*(x)$ с (не более чем) одной свободной переменной $x$.
\end{itemize}
\end{definition}

Интерпретация продолжается на все формулы языка $\HC$ \emph{одновременной
рекурсией}: высказыванию $A$ сопоставляется замкнутая формула $A^*$, задаче
$\varphi$ --- формула $\varphi^*(x)$ с (не более чем) одной свободной переменной.

\begin{definition}[продолжение $*$]\label{ru:def:tr}
Для высказываний:
\[
\begin{aligned}
  \bot^*&:=\arbot ,\\
  (A\aand B)^*&:=A^*\aand B^*, &
  (A\aor B)^*&:=A^*\aor B^*,\\
  (A\ato B)^*&:=A^*\ato B^*, &
  (\quest\varphi)^*&:=\exists x\,\varphi^*(x). &&
\end{aligned}
\]
Для задач:
\[
\begin{aligned}
  \varnothing^*(x)&:=\arbot ,\\
  (\varphi\pand\psi)^*(x)&:=\exists u\,\exists v\,
     \bigl(x=c(u,v)\aand\varphi^*(u)\aand\psi^*(v)\bigr), &&\\
  (\varphi\por\psi)^*(x)&:=\exists u\,
     \bigl((x=c(0,u)\aand\varphi^*(u))\aor(x=c(1,u)\aand\psi^*(u))\bigr), &&\\
  (\varphi\pto\psi)^*(x)&:=\forall v\,\bigl(\varphi^*(v)\ato
     \psi^*(U_x(v))\bigr), &&\\
  (\bang A)^*(x)&:=\Prf\bigl(x,\gn{A^*}\bigr). &&
\end{aligned}
\]
\end{definition}

Подчеркнём, что  $A^*$ всегда является замкнутой формулой --- индукцией по построению нетрудно
доказать, что перевод любой задачи --- формула с (не более чем) одной свободной переменной,
а перевод любого высказывания --- замкнутая формула.

\begin{definition}[значение]\label{ru:def:val}
Для высказывания $A$ положим $\val{A}:=\Tru$, если формула $A^*$ истинна в
$\N$, и $\val{A}:=\Fls$ иначе. Для задачи $\varphi$ положим
\[
  \val{\varphi}:=\{n\in\N\ :\ \N\models\varphi^*(n)\}.
\]
\end{definition}

Итак, значение высказывания --- истинностное значение, значение задачи ---
множество натуральных чисел (<<множество решений>>, оно же множество
реализаций). Индекс $*$ в записи $\val{\cdot}$ напоминает, что значение
зависит от интерпретации.

\begin{lemma}[связки ведут себя как обещано]\label{ru:lem:clauses}
Для любой интерпретации $*$:
\begin{enumerate}[leftmargin=2.2em,itemsep=1pt,label=\textup{(\arabic*)}]
\item $\val{\bot}=\Fls$, и связки $\aand,\aor,\ato$ на высказываниях
  действуют классически;
\item $\val{\quest\varphi}=\Tru\iff\val{\varphi}\neq\varnothing$;
\item $\val{\varnothing}=\varnothing$;
\item $\val{\varphi\pand\psi}=\{c(u,v):u\in\val{\varphi},\,v\in\val{\psi}\}$;
\item $\val{\varphi\por\psi}=\{c(0,u):u\in\val{\varphi}\}\cup
      \{c(1,v):v\in\val{\psi}\}$;
\item $\val{\varphi\pto\psi}=\{i:\ \text{для всех}\ u\in\val{\varphi}\
      \text{значение}\ U_i(u)\ \text{определено и лежит в}\ \val{\psi}\}$;
      в частности $\val{\pneg\varphi}=\N$, если $\val{\varphi}=\varnothing$, и
      $\val{\pneg\varphi}=\varnothing$ иначе;
\item $\val{\bang A}=\{n:\ n\ \text{есть номер вывода формулы}\ A^*\
      \text{в}\ \PA\}$.
\end{enumerate}
\end{lemma}

\begin{lemma}[двойное отрицание --- это $\quest$ на стороне задач]\label{ru:lem:dn}
Для всякой задачи $\varphi$ и всякой интерпретации $*$
\[
  \val{\pneg\pneg\varphi}=
  \begin{cases}
    \N, & \text{если } \val{\varphi}\neq\varnothing,\\
    \varnothing, & \text{иначе,}
  \end{cases}
  \qquad\text{то есть}\qquad
  \val{\pneg\pneg\varphi}\neq\varnothing\iff\val{\quest\varphi}=\Tru .
\]
Более того, $\PA$ доказывает $\forall x((\pneg\pneg\varphi)^*(x)\aeq
(\quest\varphi)^*)$, 
и вывод этой эквивалентности строится по
$\gn{*}$.
\end{lemma}

Доказательство --- в приложении.

\remove{
\begin{remark}\label{ru:rem:dn}
Итак, $\pneg\pneg$ --- это тот же оператор <<решение существует>>, только
оставшийся на стороне задач. Отсюда видно, почему $\val{\pneg\pneg\varphi}$
двузначно ($\N$ или $\varnothing$) и почему никакая схема не может извлечь из
$\pneg\pneg\varphi$ решение самой $\varphi$: в $\N$ выбирать нечего.
\end{remark}

\begin{remark}[почему рекурсия должна тащить формулу, а не только значение]
Значение $\val{\bang A}$ зависит не от значения $\val{A}$, а от самой формулы
$A^*$: два высказывания с одним и тем же истинностным значением (скажем, две
разные истинные формулы, из которых одна доказуема, а другая нет) дают разные
$\val{\bang A}$. Поэтому семантика не композициональна на уровне значений, и
определение~\ref{ru:def:tr} обязано определять перевод $*$ параллельно со
значением.
\end{remark}
}

\subsection{Теоретико-множественная семантика: чёрные ящики}\label{ru:sec:valBS}

Наряду с арифметической интерпретацией мы рассматриваем ещё одну,
теоретико-множественную. Вместе с двумя уровнями общезначимости
(\S\ref{ru:sec:validity}) всего получается три семантики: $\sI$ и $\sE$
(арифметические) и $\sB$ (чёрные ящики).

\paragraph{Язык.} Для новой семантики мы вводим в грамматику, задающую
формулы, ещё два правила --- по одному на каждый сорт:
\[
\varphi \leftarrow  \Pi\alpha_i\,\varphi, \qquad
A \leftarrow  \Pi\alpha_i\,A .
\]
Оператор $\Pi\alpha$ (<<равномерно по $\alpha$>>) применяется к
формулам обоих сортов и \emph{связывает}
вхождения переменной $\alpha$ в своей области действия. Кроме того ---
и это существенно --- в семантике $\sB$ связывающей считается и
модальность $\bang$: задача $\bang A$ связывает все задачные
переменные, свободные в $A$. Причина станет ясна из определений:
решения задачи $\bang A$ --- выводы одной формулы, а вывод не может
зависеть от содержимого ящика.

Свободные и связанные вхождения, переименование связанных переменных и
оговорка <<$\theta$ свободна для $\alpha$ в $\varphi$>> понимаются как
обычно для кванторов. Через  $FV(\varphi)$, $FV(A)$ обозначаем множество свободных
задачных переменных формулы (любого сорта).

\begin{definition}[семантика $\sB$]
\label{ru:def:starB}
В семантике $\sB$ \emph{интерпретация} $*$ сопоставляет каждой
переменной сорта высказывание замкнутую формулу языка $\ZF$ (например,
арифметическую); переменные-задачи не интерпретируются --- это
\emph{чёрные ящики}. Интерпретация --- <<прозрачная>> часть семантики
(у неё есть код), ящики --- <<чёрная>> (у ящика кода нет). Перевод
определяется теми же пунктами, что в определении~\ref{ru:def:tr}, с
четырьмя изменениями:
\[
  \alpha^*(x):=(x\in X_\alpha),\qquad
  (\bang A)^*(x):=\Prf_{\ZF}\bigl(x,\gn{\forall\vec X\,A^*}\bigr),
\]
\[
  (\Pi\alpha\,\varphi)^*(x):=\forall X_\alpha{\subseteq}\,\omega\
  \varphi^*(x),\qquad
  (\Pi\alpha\,A)^*:=\forall X_\alpha{\subseteq}\,\omega\ A^* .
\]
Здесь $X_\alpha$ --- множественная переменная языка $\ZF$, своя для
каждой переменной-задачи $\alpha$: ящики представлены в переводе
свободными переменными $X_\alpha$. В пункте для $\bang$ через $\vec X$
обозначены переменные $X_\beta$ для всех задачных переменных $\beta$,
свободных в $A$, то есть связываемых этим знаком $\bang$; гёделев номер
замыкания --- числовой терм, и переменные $X_\beta$ в формуле
$(\bang A)^*$ не встречаются. Замыкание $\forall\vec X\,A^*$ --- это
блок кванторов $\forall X_\beta{\subseteq}\,\omega$ по всем таким
$\beta$ в фиксированном порядке (скажем, по возрастанию индексов
переменных). Иначе говоря, $\forall\vec X\,A^*$ есть перевод
высказывания $\Pi\vec\beta\,A$, где $\vec\beta$ --- свободные задачные
переменные $A$, выписанные в том же порядке, и переводы задач
$\bang A$ и $\bang\Pi\vec\beta\,A$ совпадают. Этим мы воспользуемся при
проверке схемы (A6). \emph{Кодом} интерпретации служит код
одной лишь прозрачной части (определение~\ref{ru:def:code}): у ящика
кода нет.

\emph{Значения.} О значениях формул в семантике $\sB$ мы говорим только
внутри $\ZF$. При данных множествах $X_\alpha\subseteq\omega$
значение задачи $\varphi$ --- это множество
$\val{\varphi}:=\{n\in\omega:\varphi^*(n)\}$, существующее по аксиоме
выделения; слова <<при всяком содержимом ящиков>> означают квантор
всеобщности по $X_\alpha$ внутри $\ZF$. Можно представлять себе, что
рассуждение ведётся в какой-нибудь модели $\ZF$. Но канонической
модели у $\ZF$ нет, поэтому в определение общезначимости
(определение~\ref{ru:def:valid}) значения не входят.
\end{definition}

\begin{remark}[почему $\bang$ связывает]\label{ru:rem:bangbinds}
Индукцией по построению нетрудно доказать, что перевод $A^*$ любого
высказывания $A$ --- формула, свободными переменными которой могут
быть лишь $X_\alpha$ для свободных задачных переменных $\alpha$
высказывания $A$, а перевод $\varphi^*$ любой задачи может иметь сверх
них ещё максимум одну свободную числовую переменную. Решение задачи
$\Pi\alpha\,\varphi$ --- это \emph{одно} решение $\varphi$, годное при
любом содержимом ящика $\alpha$; высказывание $\Pi\alpha\,A$ читается
<<$A$ при любом содержимом ящика $\alpha$>>. 

Пункт для $\bang$ можно было бы записать и короче,
$\Prf_{\ZF}(x,\gn{A^*})$, с открытой формулой $A^*$ под знаком
номера: выводы открытой формулы --- это выводы её замыкания, так что
обе записи задают одно и то же множество решений с точностью до
вычислимого преобразования. Но вторая запись скрывает квантор, и это
опасно. 

Вот пример, иллюстрирующий это явление. Рассмотрим задачу 
$\phi:=\bang\quest\alpha$. Ее интерпретацией является формула
$\phi^*=\Prf(y,\forall X \exists x\in X)$. Формула $\forall X\,\exists x\,(x\in X)$ опровержима в $\ZF$ (ящик
$\varnothing$), и из любого её вывода вместе с фиксированным
опровержением эффективно строится вывод формулы $\arbot $,
причём корректность этого преобразования доказуема в $\ZF$; обратно,
из вывода $\arbot $ строится вывод чего угодно. Таким образом, наша задача
$\phi$ равносильна задаче $\bang\bot$: общезначима формула
$\phi\peq\bang\bot$. 
\remove{
Если бы $\bang$ не считалась связывающей, то задача
$\Pi\alpha(\bang\quest\alpha\pto\varnothing)$ была бы общезначима:
посылка $\bang\quest\alpha$ --- выводы опровержимого замыкания
$\forall X\,\exists y\,(y\in X)$, то есть при $\Con(\ZF)$ пустое
множество. А её экземпляр по закону $(\Pi\mathrm{E})$ с
$\theta:=\pneg\varnothing$ --- задача
$\bang\quest\pneg\varnothing\pto\varnothing$ --- решений не имеет:
формула $\exists x\,(\pneg\varnothing)^*(x)$ доказуема, и выводы у неё
есть. Подстановка внутрь гёделева номера не достаёт, и закон
$(\Pi\mathrm{E})$ ломался бы. Считая $\bang$ связывающей, мы получаем,
что в записи $\Pi\alpha(\bang\quest\alpha\pto\varnothing)$ переменная
$\alpha$ связана знаком $\bang$, квантор $\Pi\alpha$ пуст, и
экземпляров с подстановкой у этой задачи нет.

По той же причине высказывание $\alpha\pto\bang\quest\beta$ не имеет
ничего общего с корефлексией: $\bang\quest\beta$ --- выводы замыкания
$\forall X_\beta\,\exists y\,(y\in X_\beta)$, то есть, по существу,
задача $\bang\bot$, и вся формула равносильна $\pneg\alpha$.}
(Конец замечания.)
\end{remark}

\begin{remark}[об отвергнутой семантике $\sS$]\label{ru:rem:semS}
Наряду с $\sB$ мы рассматривали ещё одну теоретико-множественную
семантику --- \emph{подстановочную} $\sS$ (от слова <<substitute>>),
в двух редакциях. В первой задачные переменные интерпретировались, как
в $\sI$, арифметическими формулами $\alpha^*(x)$, и решатель получал
их код, а квантор $\Pi$ читался подстановочно: <<$\varphi^*$
выполнено при подстановке вместо $\alpha$ любой арифметической
формулы>>. Записать это формулой $\ZF$ можно только через
релятивизацию квантора по ящикам: при вложенных кванторах
подстановочные экземпляры $\varphi^*$ сами содержат предикат
истинности, и по теореме Тарского прямая запись не проходит; поэтому
ящик $X_\alpha$ пробегал арифметически определимые множества, что
выражается формулой $\mathrm{Def}(X)$ языка $\ZF$ (через определимый в
$\ZF$ предикат истинности для арифметики). Ничего принципиально
нового по сравнению с $\sI$ и $\sB$ доказать про эту семантику нам не
удалось, а сочетание кода подставленной формулы со связывающей
модальностью $\bang$ требовало оговорок в схеме $(\Pi\mathrm{E})$ и в
правиле обобщения. Во второй редакции код был убран вовсе: $\sS$ ---
это $\sB$, в которой ящики пробегают только арифметически определимые
множества (кванторы по $X_\alpha$ в пунктах для $\Pi$ и для $\bang$
релятивизованы формулой $\mathrm{Def}$). Такая семантика почти не
отличается от $\sB$: все положительные результаты переносятся
релятивизацией кванторов, а все наши опровержения используют
арифметически определимые ящики. Единственная известная нам точка
возможного расхождения --- предикативность: в $\sS$ значения
$\Pi$-содержащих задач лишь определимы в $\ZF$ и в диапазон квантора
не попадают, поэтому схему $(\Pi\mathrm{E})$ удавалось обосновать
лишь для $\Pi$-свободных $\theta$, а замкнутость относительно
подстановки --- лишь для $\Pi$-свободных подставляемых задач; в $\sB$
квантор не предикативен, и обе оговорки излишни. Платить за
семантику, в доказанном ничем не отличающуюся от $\sB$, оговорками во
всех утверждениях мы не стали и из статьи её исключили. Меняет ли
что-нибудь ограничение ящиков определимыми множествами, остаётся
открытым вопросом (\S\ref{ru:sec:todo}); попутный довод в пользу
определимых ящиков --- переводы высказываний $\quest\Pi\vec\alpha\,\varphi$
становятся арифметическими относительно $0^{(\omega)}$, а не
$\Pi^1_1$-формулами.
(Конец замечания.)
\end{remark}

\begin{remark}\label{ru:rem:clean}
По возможности, мы будем стремиться выполнить 
одно синтаксическое ограничение на формулы для этой семантики: 
\emph{никакая переменная не имеет в формуле одновременно свободных и
связанных вхождений, разные связывающие операторы связывают разные
переменные, а модальность $\bang$ применяется только к высказываниям
без свободных задачных переменных}.
Такие формулы называются \emph{чистыми}. Всякая формула равносильна
чистой: переводы задач $\bang A$ и $\bang\Pi\vec\alpha\,A$ совпадают
(определение~\ref{ru:def:starB}), и достаточно сделать квантор явным, а
затем переименовать связанные переменные.
Мы делаем это, чтобы избежать недоразумений. Например,
в семантике $\sB$
корефлексия $\alpha\pto\bang\quest\alpha$, содержащая как свободное, так и связанное
вхождение переменной $\alpha$, имеет совсем другой смысл, чем в арифметических
семантиках. В семантике $\sB$  корефлексия равносильна чистой формуле
 $\alpha\pto\bang\Pi\beta\,\quest\beta$, в свою очередь равносильной
формуле $\alpha\pto\varnothing$. Последняя совсем не равносильна корефлексии в арифметических семантиках. Смысл же чистых формул примерно одинаков во всех семантиках.
\end{remark}

Лемма~\ref{ru:lem:clauses} переносится на семантику $\sB$ (её утверждения
доказуемы в $\ZF$); дополнительно $\ZF$ доказывает
\[
  \val{\Pi\alpha\,\varphi}
   =\bigcap_{X\subseteq\omega}\val{\varphi}[X_\alpha:=X]
  .
\]

Параметры трёх семантик
собраны в таблице:
\begin{center}\footnotesize
\begin{tabular}{l|l|l|l}
 & $\sI$ & $\sB$ & $\sE$\\
\hline
перем.-задачи & арифм.\ форм. & любые ящики & арифм.\ форм.\\
перем.-выск. & арифм.\ предл. & $\ZF$-предл. & арифм.\ предл.\\
оператор $\Pi$ & нет & $\forall X_\alpha$ & нет\\
$\bang$ связывает перем.-задачи & нет & да & нет\\
внутренняя теория & $\PA$ & $\ZF$ & $\PA$\\
внешняя теория & $\PA$ & $\ZF$ & $\TA$\\
\end{tabular}
\end{center}
Две последние строки вводят два термина, которыми мы будем
пользоваться постоянно. \emph{Внутренняя теория} семантики --- та,
выводы которой служат решениями задач $\bang A$; доказуемость в ней и
выражает модальность $\square=\quest\bang$, а $\Prov$, $\Prf$ и $\Con$ без
уточнения относятся к ней. \emph{Внешняя теория} семантики --- та, в
которой проверяется общезначимость (\S\ref{ru:sec:validity}): формула
считается установленной, если она во внешней теории выводима. Для
$\sE$ внешней теорией служит $\TA$ --- ее аксиомами являются все истинные в
$\N$ арифметические предложения, так что
<<$\TA\vdash\sigma$>> означает просто истинность $\sigma$ в $\N$. Во
внутренних семантиках $\sI$, $\sB$ внутренняя и внешняя теории
совпадают (отсюда и название этих семантик); в $\sE$ внешняя теория
строго сильнее внутренней. Во всех трёх случаях внутренняя теория
содержится во внешней. 

\section{Общезначимость}\label{ru:sec:validity}

\subsection{Определение}

\begin{definition}[код интерпретации]\label{ru:def:code}
Пусть в данной формуле $\varphi$ 
встречаются переменные $A_1,\dots,A_k$ сорта
высказывание и $\alpha_1,\dots,\alpha_m$ сорта задача. \emph{Кодом}
арифметической интерпретации $*$ (относительно формулы $\varphi$) назовём кортеж 
\[
  \gn{*}_\varphi\ :=\ \bigl\langle\, A_1^*,\dots,A_k^*,\
        \alpha_1^*,\dots,\alpha_m^*\,\bigr\rangle .
\]
Это конечный объект: значения $*$ на переменных, не входящих в $\varphi$, в
него не входят. Для семантики $\sB$ вторую часть кортежа
опускаем: задачные переменные там --- ящики, кода у них нет.
\end{definition}

\begin{definition}[общезначимость]\label{ru:def:valid}
Высказывание $A$ \emph{общезначимо} в данной семантике, если при
всякой интерпретации~$*$ формула $A^*$ выводима во внешней теории этой
семантики. 
Задача $\varphi$ \emph{общезначима} в данной семантике, если существует
всюду определённая вычислимая функция $S$ такая, что для всякой
интерпретации $*$ значение $S(\gn{*}_\varphi)$ определено и формула
$\varphi^*\bigl(S(\gn{*})\bigr)$ выводима во внешней теории.\footnote{Точнее, вместо значения функции в формулу нужно
подставлять стандартный терм, задающий это значение в стандартной интерпретации, но
нам не хочется перегружать обозначения.}
В семантике $\sB$ переводы $A^*$ и
$\varphi^*\bigl(S(\gn{*})\bigr)$ могут содержать свободные
множественные переменные $\vec X$ --- ящики свободных задачных
переменных формулы. Выводимостью такой формулы мы называем в $\sB$
выводимость её замыкания $\forall\vec X\,\varphi^*\bigl(S(\gn{*})\bigr)$.
\end{definition}

Общезначимость в семантике $Z$ мы называем также \emph{$Z$-общезначимостью};
общезначимость сразу в нескольких семантиках --- например,
\emph{$\sI,\sB$-общезначимостью}.
Для $\sE$ выводимость в $\TA$ --- это истинность в $\N$, то есть
принадлежность $S(\gn{*})\in\val{\varphi}$; для $\sI$ эта
принадлежность следует из выводимости, поскольку $\PA$ содержится в
$\TA$. Для $\sB$ значения в определение не входят: у $\ZF$
нет канонической модели. Если представлять себе какую-нибудь модель
$M$ теории $\ZF$ и значения в ней (определение~\ref{ru:def:starB}), то
общезначимая задача решается функцией $S$ в $M$ при всякой
интерпретации и всяком содержимом ящиков; но это не часть определения, а его неформальное понимание.

\begin{remark}[метатеоретические допущения]\label{ru:rem:assumptions}
Все доказательства \emph{не}общезначимости в  семантике $\sB$ используют 
предположение о непротиворечивости $\ZF$, что необходимо,
поскольку иначе все формулы общезначимы.

Точнее, чтобы доказать необщезначимость задачи $\phi$ в этих
семантиках, мы должны для любой тотальной вычислимой функции
$S$ предъявить интерпретацию $*$, для которой
формула $\phi^*(S(*))$ не выводима в $\ZF$. 
Вместо этого мы будем доказывать, что отрицание этой формулы
выводимо в $\ZF$. В предположении непротиворечивости
это гарантирует невыводимость самой формулы.
Аналогичным образом, вместо доказательства невыводимости
формулы $A^*$ (для высказывания $A$) мы будем доказывать выводимость
$\aneg A^*$. Точно так же мы будем действовать и для арифметических семантик,
только в этом случае никаких предположений делать не нужно.
В нескольких местах непротиворечивости мало, и мы используем более
сильное допущение --- \emph{$\Sigma_1$-корректность} внешней теории:
всякое доказуемое в ней арифметическое $\Sigma_1$-предложение
истинно; на это всякий раз указывается явно. Других метатеоретических
допущений в статье нет: везде речь идёт либо о непротиворечивости
внешней теории, либо о её $\Sigma_1$-корректности, а для внешних
теорий $\PA$ и $\TA$ оба свойства имеют место.

Оба допущения --- следствие того, что метатеория у нас не
фиксирована и остаётся по возможности слабой. При достаточно сильной
метатеории, например при формализации в системе Lean, они становятся
теоремами: стандартная модель $\N$ и модель $\mathrm{ZFC}$ из
библиотеки Mathlib (кумулятивная иерархия до недостижимого кардинала
$\kappa$) там построены, так что непротиворечивость и
$\Sigma_1$-корректность $\PA$ и $\ZF$ доказуемы. Заодно появляется
аналог семантики $\sE$ для $\ZF$: внешняя теория --- истинность в
этой модели $V_\kappa$, а внутренняя --- по-прежнему $\ZF$.
(Конец замечания.)
\end{remark}

\begin{remark}[другая внутренняя теория]\label{ru:rem:otherth}
В семантике $\sI$ арифметика Пеано играет обе роли --- внутренней и
внешней теории. Вместо неё можно взять любую перечислимую
$\Sigma_1$-корректную теорию $T$, содержащую $\PA$ (или
интерпретирующую арифметику), например $\ZF$: переменные-высказывания
интерпретируются замкнутыми формулами $T$, переменные-задачи ---
формулами $T$ с одной свободной переменной, пробегающей натуральные
числа, решения задачи $\bang A$ --- выводы формулы $A^*$ в $T$, а
общезначимость --- выводимость в $T$. Обозначим эту семантику $\sI_T$,
так что $\sI=\sI_{\PA}$. Всё, что ниже доказано о семантике $\sI$,
переносится на $\sI_T$ вместе с доказательствами: в них используются
лишь $\Sigma_1$-полнота, лемма о неподвижной точке, формализуемость
реализуемости и россеровское предложение --- свойства, общие для всех
таких теорий; для $T=\ZF$ соответствующие выкладки и проведены ниже
при разборе семантики $\sB$. Отметим, что $\sI_{\ZF}$ --- не то же
самое, что $\sB$ без квантора: в $\sB$ модальность $\bang$ связывает
задачные переменные, и, например, высказывание
$\square\quest\alpha\ato\square\bot$ $\sB$-общезначимо
(замечание~\ref{ru:rem:boxqatom}), а в $\sI_{\ZF}$ --- нет (возьмите
$\alpha^*(x):=(x=x)$). Совпадают ли множества общезначимых задач в
$\sI_{\PA}$ и $\sI_{\ZF}$, мы не знаем: переноса доказательств для
этого недостаточно, а полноты на задачах нет. (Конец замечания.)
\end{remark}

\begin{remark}[вывод формулы находится перебором]\label{ru:rem:whycomp}
Заметим, что в определении $\sI,\sB$-общезначимости
выводы формул $A^*$ и $\varphi^*\bigl(S(\gn{*})\bigr)$ можно найти перебором.
Например, для формулы $A^*$: по
коду $\gn{*}$ вычисляется сама формула $A^*$ (перевод $\cdot\mapsto\cdot^*$ из
определения~\ref{ru:def:tr} --- вычислимая операция над кодами),
отношение $\Prf$ разрешимо, а перебор по предположению останавливается при
каждом $*$; значит функция $P\colon\gn{*}\mapsto{}$<<первый найденный вывод
формулы $A^*$>> всюду определена и вычислима. 
\end{remark}


\begin{remark}[три возможных понятия и почему выбрано среднее]\label{ru:rem:three}
В определении общезначимости задачи $\varphi$ можно потребовать
\[
  \text{(а)}\ \ \forall *\ \exists e\ \ e\in\val{\varphi};
  \qquad
  \text{(б)}\ \ \exists\, S\ \forall *\ \ S(\gn{*})\in\val{\varphi};
  \qquad
  \text{(в)}\ \ \exists e\ \forall *\ \ e\in\val{\varphi},
\]
где в~(б) функция $S$ вычислима. Очевидно, (в) влечёт (б), а (б) влечёт~(а).
Крайние понятия не годятся оба.

\emph{При определении \textup{(а)}} 
выбор решения ничем не
стеснён, и фрагмент задач становится классическим: общезначима
$\alpha\por\pneg\alpha$, закон Пирса и вообще все классические тавтологии
(предложение~\ref{ru:prop:noem} ниже показывает, что при~(б) этого
уже нет).

\emph{Понятие \textup{(в)}} --- слишком сильное. Уже задача $\bang(B\ato B)$ не
будет общезначима: её решения при интерпретации $*$ --- это выводы формулы
$B^*\ato B^*$, а одного числа, годного при всех $*$, нет, поскольку у разных
формул разные выводы. Между тем хотелось бы, чтобы все интуиционистски выводимые формулы
были общезначимы. 
При~(б) всё в порядке: вывод формулы
$B^*\ato B^*$ эффективно строится по $B^*$.
(Конец замечания.)
\end{remark}

\remove{
\begin{remark}[семантика интенсиональна]\label{ru:rem:intens}
Функция $S$ получает \emph{синтаксис} интерпретации, а не сами множества:
множество решений --- бесконечный объект, подать его на вход нельзя. Поэтому
на двух кодах логически эквивалентных формул $S$ вправе вести себя
по-разному. Для модальности $\bang$, которая и так смотрит на формулу $A^*$, а
не на её истинностное значение, это естественно.
\end{remark}
}

\begin{remark}[выбор нумерации, спаривания и кодирования не важен]\label{ru:rem:indep}
В определениях~\ref{ru:def:tr}--\ref{ru:def:code} было сделано три произвольных
выбора: гёделева нумерация $U_0,U_1,\dots$, спаривающая функция $c$ и способ
кодирования интерпретации в $\gn{*}$. Ни один из них не влияет на то, какие
формулы общезначимы --- ни в одной из семантик. Обоснование --- в
приложении.
\end{remark}

Мы закончим этот раздел доказательством того, что общезначимость сохраняется
при подстановках формул вместо переменных, а также при применении правила
Modus Ponens (MP).

\begin{lemma}[замкнутость]\label{ru:lem:closure}
В арифметических семантиках множество общезначимых формул 
замкнуто относительно правила \textup{(MP)} и относительно
подстановки: задач --- вместо свободных переменных-задач, высказываний
--- вместо переменных-высказываний. 
\end{lemma}

Для семантики $\sB$ аналог этой леммы верен не для всех подстановок. Вот контрпример:
высказывание $\quest\Pi\alpha(\beta\por\pneg\beta)$
общезначимо в $\sB$: 
формулу $\exists x\,(\beta\por\pneg\beta)^*(x)$ теория $\ZF$ доказывает
разбором случаев: если значение $\beta$ непусто, скажем выполнено $\beta^*(u)$, 
то решением служит
$c(0,u)$, 
а иначе --- $c(1,0)$. Если
подставить $\alpha$ вместо $\beta$
, получится
высказывание $\quest\Pi\alpha(\alpha\por\pneg\alpha)$, которое $\ZF$
\emph{опровергает}: решение вида $c(0,u)$ отвергается ящиком
$\varnothing$, а решение вида $c(1,e)$ --- ящиком $\{0\}$, поскольку
метка $1$ обещает пустоту значения.

Для семантики $\sB$ нам будет достаточно этой леммы для 
так называемых \emph{безопасных подстановок}. По определению, это подстановки
со следующими двумя свойствами:
\begin{enumerate}[leftmargin=2.2em,itemsep=1pt,label=\textup{(\arabic*)}]
\item Свободные задачные переменные подставляемой формулы
  не должны попадать под связывающий квантор $\Pi$. Это явление называется
\emph{захватом переменных},
необходимость запрещения захвата мы только что продемонстрировали контрпримером.
\item Аналогично, свободные задачные переменные подставляемой формулы
  не должны попадать под связывающую модальность $\bang$. Это
  ограничение касается только подстановок высказываний вместо
  переменных-высказываний: при подстановке задачи вместо задачной
  переменной оно выполнено автоматически, поскольку свободные
  вхождения задачной переменной под знаком $\bang$ не лежат. Пример
нарушения: $\bang A\mapsto \bang\quest \alpha$. Необходимость этого ограничения мы не можем продемонстрировать примером:
мы не знаем, необходимо ли оно в самом деле.
\end{enumerate}



\begin{lemma}[замкнутость]\label{ru:lem:closurebs}
В семантике $\sB$ множество общезначимых формул 
замкнуто  относительно правила \textup{(MP)} и относительно
безопасных подстановок.
\end{lemma}

Доказательство лемм --- в приложении.

\subsection{
Схема (A4)}\label{ru:sec:assval}

Единственная схема о высказываниях, которая станет аксиомой, --- это
группа (A4), утверждающая, что модальность $\quest$ --- булев
гомоморфизм; здесь мы проверим её общезначимость. 
\remove{Полное описание
общезначимых $\Pi$-свободных высказываний --- с точностью до
нормальной формы это логика доказуемости $\GL$ Соловея для внутренних
семантик $\sI$, $\sB$ и её внешний вариант $\GLS$ для $\sE$ --- мы
получим в \S\ref{ru:sec:asscompl}, после теоремы о корректности: там оно
достаётся почти даром из того же описания \emph{выводимых}
высказываний. Для арифметических семантик $\sI$ и $\sE$, в языке
которых квантора $\Pi$ нет вовсе, это описание исчерпывает \emph{все}
высказывания; в полном же языке семантики $\sB$ за его пределами
остаются высказывания с квантором $\Pi$, и логика этого
фрагмента --- открытый вопрос (\S\ref{ru:sec:piatoms}, \S\ref{ru:sec:todo}).
}
В этом и двух следующих подразделах все утверждения формулируются и
доказываются сразу для всех трёх семантик.
\remove{
; специфика семантики $\sB$
--- оператор $\Pi$ --- разобрана
отдельно в \S\ref{ru:sec:piatoms}, поэтому здесь рассматриваются
$\Pi$-свободные формулы. Соглашения: внутренняя и внешняя теории каждой
семантики --- как в таблице \S\ref{ru:sec:valBS} (внутренняя: $\PA$ для
$\sI$, $\sE$ и $\ZF$ для $\sB$; внешняя отличается от неё только
в $\sE$, где это $\TA$); $\Con$ обозначает непротиворечивость внутренней
теории. Кодом интерпретации в $\sB$ служит код прозрачной
части (определение~\ref{ru:def:code}); переводы в этой семантике могут
содержать свободные множественные переменные $\vec X$, и выводимость
таких формул понимается как выводимость замыкания
(определение~\ref{ru:def:valid}).
Леммы~\ref{ru:lem:clauses} и~\ref{ru:lem:dn}
используются во всех семантиках: их доказательства схемны и одинаково
проходят в $\PA$ и в $\ZF$, в том числе для переводов со свободными
множественными переменными. Отрицательные утверждения для $\sB$,
как всегда, доказываются в допущениях
замечания~\ref{ru:rem:assumptions}.
}

\begin{lemma}[$\quest$ --- булев гомоморфизм на уровне переводов]\label{ru:lem:qsem}
Для любых задач $\varphi,\psi$ следующие эквивалентности доказуемы в
внутренней теории (а потому и во внешней),
причём выводы вычислимо строятся по коду интерпретации:
\[
\begin{aligned}
  (\quest(\varphi\pto\psi))^* &\aeq\bigl((\quest\varphi)^*\ato(\quest\psi)^*\bigr),&
  (\quest(\varphi\pand\psi))^* &\aeq\bigl((\quest\varphi)^*\aand(\quest\psi)^*\bigr),\\
  (\quest(\varphi\por\psi))^* &\aeq\bigl((\quest\varphi)^*\aor(\quest\psi)^*\bigr),&
  (\quest\varnothing)^* &\aeq\arbot .
\end{aligned}
\]
\end{lemma}

\begin{proof}
Для конъюнкции: формула $(\quest(\varphi\pand\psi))^*$ есть $\exists x\,
\exists u\exists v\,(x=c(u,v)\aand\varphi^*(u)\aand\psi^*(v))$, и внутренняя
теория доказывает её эквивалентность конъюнкции $\exists u\,\varphi^*(u)\aand
\exists v\,\psi^*(v)$: слева направо --- проекции пары, справа налево ---
построение пары. Так же для дизъюнкции, с разбором того, какая из меток
$c(0,u)$, $c(1,v)$ стоит. Для $\varnothing$: $(\quest\varnothing)^*$ есть
$\exists x\,\arbot $.

Для импликации слева направо внутренняя теория доказывает $\exists
x\,(\varphi\pto\psi)^*(x)\ato(\exists y\,\varphi^*(y)\ato\exists z\,\psi^*(z))$
--- это чистая логика. Справа налево нужна работа: надо предъявить индекс.
Внутренняя теория доказывает $(\exists u\,\varphi^*(u)\ato\exists v\,\psi^*(v))\ato\exists
x\,(\varphi\pto\psi)^*(x)$ разбором двух случаев. Если $\aneg\exists
u\,\varphi^*(u)$, то условие в $(\varphi\pto\psi)^*(x)$ пусто и годится любое
$x$. Если же $\exists v\,\psi^*(v)$, то, взяв наименьшее такое $v$, годится
$x:=k(v)$, где $k$ --- всюду определённая вычислимая функция, дающая по $v$ номер
постоянной функции $u\mapsto v$ (теорема о параметризации). Разбор случаев
здесь классический --- это \emph{остаточная классичность} нашей семантики,
к которой мы ещё вернёмся (замечание~\ref{ru:rem:qhom}). Все выкладки
одинаковы в $\PA$ и $\ZF$; в $\sB$ они проводятся поверх свободных
множественных переменных $X_\alpha$, которые ни в разборе случаев, ни в
построении пар не участвуют.
\end{proof}

\subsection{Безмодальные задачи: схемы $\Int$}\label{ru:sec:intval}

Перейдём к задачам и начнём с самого простого их класса --- задач без
модальностей, то есть построенных из переменных и константы $\varnothing$
связками $\pto$, $\pand$, $\por$. Будущие аксиомы (A2) --- схемы
интуиционистского исчисления $\Int$ --- общезначимы.

\begin{proposition}[схемы $\Int$ общезначимы]\label{ru:prop:intval}
Всякая задача, выводимая в интуиционистском исчислении, общезначима во
всех трёх семантиках, причём с постоянной решающей функцией.
\end{proposition}

\begin{proof}
Клиниевские реализации схем аксиом $\Int$ строятся из номеров, доставляемых
теоремой о параметризации, и от интерпретации не зависят; корректность
каждой реализации доказуема во внутренней теории, и вывод строится по коду:
выкладки одни и те же в $\PA$ и $\ZF$, а в $\sB$ они проводятся поверх
свободных множественных переменных. Правило (MP) сохраняет
общезначимость по леммам~\ref{ru:lem:closure} и~\ref{ru:lem:closurebs}.
\end{proof}

Всё остальное, что мы знаем об общезначимых безмодальных задачах ---
в частности, что классических схем среди них нет, а интуиционистскими
они не исчерпываются, --- удобнее доказывать после теоремы
корректности; это сделано в~\S\ref{ru:sec:modalfreeval}.

\subsection{Общезначимые задачи с модальностями}\label{ru:sec:taskval}

О множестве всех общезначимых задач --- с модальностями --- мы знаем гораздо
меньше. Оно (в каждой из трёх семантик) замкнуто относительно (MP) и
подстановки --- в $\sB$ безопасной (леммы~\ref{ru:lem:closure}
и~\ref{ru:lem:closurebs}). Здесь мы проверим общезначимость
тех схем с модальностями, которые станут аксиомами исчислений
\S\ref{ru:sec:calc} для семантик $\sI$ и $\sB$: схемы (A3) и схем (A5) и
(A7) о задачах специального вида. Аксиомы исчисления для $\sE$ ---
рефлексия (A8) и её аналог для задач (A9) --- будут проверены при его
определении. 

\begin{proposition}[схема (A3)]\label{ru:prop:a3}
Для любых высказываний $A$, $B$ задача
$\bang(A\ato B)\pto(\bang A\pto\bang B)$ общезначима во всех трёх
семантиках.
\end{proposition}

\begin{proof}
В $\sI$ и $\sE$: функция, строящая по выводам формул $A^*\ato B^*$ и
$A^*$ во внутренней теории вывод формулы $B^*$ по правилу modus ponens,
вычислима; теорема о параметризации даёт постоянную
решающую функцию $S$, и её корректность доказуема во внешней теории.
В $\sB$ то же, с поправкой на связывание: решения задач
$\bang(A\ato B)$, $\bang A$, $\bang B$ --- выводы замыканий формул
$A^*\ato B^*$, $A^*$, $B^*$ по свободным переменным соответствующих
высказываний. Из вывода
формулы $\forall\vec X\,(A^*\ato B^*)$ и вывода формулы
$\forall\vec X_A\,A^*$ строится вывод формулы $\forall\vec X_B\,B^*$:
переменные, свободные в $A$, но не в $B$, инстанцируем множеством
$\varnothing$, после чего остаются modus ponens и
обобщение. Преобразование вычислимо, решатель постоянен,
его корректность доказуема в $\ZF$.
\end{proof}

\subsubsection{$\Sigma$-задачи}

Выделим класс задач, заданный грамматически; он появится в аксиомах
(A5) и (A7) исчислений (\S\ref{ru:sec:calc}).

\begin{definition}[$\Sigma$-задачи и $\Sigma$-высказывания]\label{ru:def:sigma}
\emph{$\Sigma$-задачами} называем задачи, построенные по грамматике
\[
  \varphi ::= \varnothing \mid \bang A \mid (\varphi\pand\varphi) \mid
             (\varphi\por\varphi) ,
\]
где $A$ --- произвольное высказывание. Иначе говоря, $\Sigma$-задача
собирается из задач вида $\bang A$ и константы $\varnothing$ одними лишь
$\pand$ и $\por$: переменные сорта задача, импликация и, значит, отрицание
сюда не входят. Например, $\bang A$, $\bang A\por\bang B$ и
$(\bang A\pand\bang B)\por\bang C$ --- $\Sigma$-задачи, а $\alpha$, $\pneg\bang A$ и
$\bang A\pto\bang B$ --- нет. В семантике $\sB$ у $\Sigma$-задачи
нет свободных задачных переменных: все её задачные переменные стоят
под $\bang$ и связаны этим знаком.

\emph{$\Sigma$-высказываниями} называем высказывания, построенные по
грамматике
\[
  A ::= \bot \mid \quest\phi \mid (A\aand A) \mid (A\aor A) ,
\]
где $\phi$ пробегает $\Sigma$-задачи. Переменных сорта высказывание,
импликации и отрицания здесь снова нет. В нашем исчислении 
всякое такое
высказывание будет эквивалентно 
комбинации из $\bot$ и формул $\square
B$ при помощи $\aand$ и $\aor$: $\quest\varnothing\aeq\bot$,
$\quest\bang B\aeq\square B$, а $\quest$ проходит сквозь $\pand$ и $\por$.
В частности, $\square B$ --- $\Sigma$-высказывание при любом $B$.
\end{definition}

Проверим две схемы о $\Sigma$-задачах, которые станут аксиомами (A5) и (A7).

\begin{lemma}\label{ru:lem:sigmaenum}
Существует всюду определённая вычислимая функция, которая по $\Sigma$-задаче
$\phi$ и коду $\gn{*}$ выдаёт номер программы, перечисляющей множество
$\val{\phi}$. Более того, внешняя теория доказывает, что эта программа
перечисляет в точности множество, задаваемое формулой $\phi^*$.
\end{lemma}

\begin{proof}
Индукция по построению $\phi$, с опорой на лемму~\ref{ru:lem:clauses}.
Для $\phi=\varnothing$ --- пустой перебор. Для $\phi=\bang A$ множество
$\val{\bang A}$ --- это выводы формулы $A^*$ (в $\sB$ --- её
замыкания), а <<$n$ есть вывод данной формулы>> --- разрешимое отношение,
так что годится перебор всех $n$ с
проверкой. Для $\phi=\psi\pand\chi$ перечисляем пары $c(u,v)$ по
перечислениям $\val{\psi}$ и $\val{\chi}$, для $\phi=\psi\por\chi$ ---
элементы $c(0,u)$ и $c(1,v)$, чередуя два перечисления. В $\sB$
перевод $\phi^*$ свободных множественных переменных не содержит, и рассуждение
то же.
\end{proof}

\begin{proposition}[схемы (A5) и (A7)]\label{ru:prop:a5a7}
Для всякой $\Sigma$-задачи $\phi$ во всех трёх семантиках
общезначимы задачи
\[
  \text{(A5)}\ \ \phi\pto\bang\quest\phi
  \quad\text{(ко-рефлексия)}\qquad\text{и}\qquad
  \text{(A7)}\ \ \pneg\pneg\phi\pto\phi
  \quad\text{(снятие двойного отрицания)}.
\]
\end{proposition}

\begin{proof}
Перевод $\phi^*$ --- $\Sigma_1$-формула: пункт для $\bang A$ даёт
разрешимое отношение $\Prf$, а пункты для $\pand$ и $\por$ навешивают
лишь $\aand$, $\aor$ и $\exists$.

(A5). Решение задачи $(\bang\quest\phi)^*$ --- вывод формулы
$\exists x\,\phi^*(x)$ во внутренней теории. Пусть $e$ --- номер
программы, которая, получив что угодно на вход, перебирает все числа и
выдаёт первое найденное число, являющееся выводом формулы
$\exists x\,\phi^*(x)$ во внутренней теории; положим $S(\gn{*}):=e$.
Если $\val{\phi}=\varnothing$, проверять нечего. Если же
$u\in\val{\phi}$, то $\exists x\,\phi^*(x)$ --- истинное
$\Sigma_1$-предложение, по $\Sigma_1$-полноте внутренней теории оно
выводимо, перебор останавливается, и $U_e(u)$ --- его вывод.
Выводимость во внешней теории: она доказывает
$\phi^*(u)\ato\Prov(\exists x\,\phi^*(x))$ (формализованная
$\Sigma_1$-полнота), а при $\Prov(\exists x\,\phi^*(x))$ перебор
останавливается на первом найденном выводе (формализация
неограниченного поиска, как ниже в (A7)). Итак, внешняя теория
доказывает
$\forall u\,\bigl(\phi^*(u)\ato\Prf(U_e(u),\exists x\,\phi^*(x))\bigr)$.

(A7). По лемме~\ref{ru:lem:sigmaenum} по $\gn{*}$ строится алгоритм,
перечисляющий множество $\val{\phi}$. Пусть $e$ --- номер программы,
которая, получив что угодно на вход, запускает этот перебор и выдаёт
первый найденный элемент; положим $S(\gn{*}):=e$. Если
$\val{\pneg\pneg\phi}=\varnothing$, проверять нечего. Если же
$\val{\pneg\pneg\phi}=\N$, то по лемме~\ref{ru:lem:dn} множество
$\val{\phi}$ непусто, перебор останавливается, и
$U_e(u)\in\val{\phi}$ при всех $u$. Выводимость во внешней теории: по
лемме~\ref{ru:lem:dn} формула $(\pneg\pneg\phi)^*(u)$ доказуемо
эквивалентна $\exists v\,\phi^*(v)$, а по лемме~\ref{ru:lem:sigmaenum}
внешняя теория доказывает, что перечисление отвечает формуле $\phi^*$;
значит она доказывает и то, что при $\exists v\,\phi^*(v)$ перебор
останавливается на первом найденном свидетеле (обычная формализация
ограниченного поиска). Итак, внешняя теория доказывает
$\forall u\,\bigl((\pneg\pneg\phi)^*(u)\ato\phi^*(U_e(u))\bigr)$, и
вывод строится по $\gn{*}$.

Во всех случаях выкладки одинаковы в $\PA$ и $\ZF$; в $\sB$
переводы $\Sigma$-задач свободных множественных переменных не содержат. 
\end{proof}


\subsection{Семантика $\sB$: оператор $\Pi$}
\label{ru:sec:piatoms}

Все факты подразделов \ref{ru:sec:assval}--\ref{ru:sec:taskval} уже доказаны
сразу для всех трёх семантик; новое в этом подразделе --- сам
оператор $\Pi$. Напомним мотивировку из введения. Под знаком $\quest$
конструктивное содержание задачи теряется: высказывание
$\quest(\alpha\por\pneg\alpha)$ общезначимо тривиально ($\PA$ ---
классическая теория и доказывает $\exists
x\,(\alpha\por\pneg\alpha)^*(x)$ разбором случаев при каждой
фиксированной интерпретации), и точно так же тривиально общезначимо
$\quest\rho$ для формулы Роуза (\S\ref{ru:sec:modalfreeval}), хотя сама
задача $\rho$ решается нетривиальным алгоритмом. Квантор $\Pi$ возвращает потерянное: позже
мы докажем (\S\ref{ru:sec:pigain}), что высказывание
$\quest\Pi\alpha(\alpha\por\pneg\alpha)$ опровержимо --- общезначимо
его отрицание, --- а $\quest\Pi\beta\Pi\gamma\rho$ --- содержательное
общезначимое высказывание. Здесь же мы проверим общезначимость законов самого оператора;
они станут аксиомами исчисления для $\sB$.

\begin{proposition}[законы $\Pi$]\label{ru:prop:pilaws}
В $\sB$ общезначимы следующие формулы (в скобках --- реализатор задач;
все реализаторы постоянны).
\begin{enumerate}[leftmargin=1.8em,itemsep=1pt]
\item $(\Pi\mathrm{E})$:
  $\Pi\alpha\,\varphi\pto\varphi[\theta/\alpha]$ и
  $\Pi\alpha\,A\ato A[\theta/\alpha]$, если $\theta$
  свободна для $\alpha$ \ (для задачи реализатор тождественный: перевод
  $\varphi[\theta/\alpha]$ получается из $\varphi^*$ подстановкой
  множества $\{x:\theta^*(x)\}$ вместо $X_\alpha$ --- знаков $\bang$
  подстановка не задевает, поскольку свободные вхождения $\alpha$ под
  ними не лежат, --- а пересечение по всем ящикам содержится в значении
  при этом ящике; для высказывания $\ZF$ доказывает перевод тем же
  инстанцированием квантора $\forall X_\alpha$);
\item $(\Pi\por)$:
  $\Pi\alpha(\varphi\por\psi)\pto\Pi\alpha\,\varphi\por\Pi\alpha\,\psi$
  \ (тождественный; существенно, что метка решения фиксирована: решение
  посылки --- пара из метки и решения $e$, и та же метка с тем же $e$
  годится при всяком ящике);
\item $(\Pi\bang)$: $\bang A\peq\bang\Pi\alpha\,A$ \ (реализатор ---
  преобразование выводов: решения посылки --- выводы замыкания
  $\forall\vec X\,A^*$, решения заключения --- выводы замыкания
  $\forall\vec X'\,\forall X_\alpha\,A^*$, где $\vec X'$ --- ящики
  остальных свободных задачных переменных $A$; эти формулы отличаются
  лишь порядком кванторов (а при $\alpha\notin FV(A)$ --- лишним
  квантором), и если $\alpha$ стоит в $\vec X$ последней, они совпадают;
  перестановка кванторов в выводе --- вычислимое
  преобразование, корректность которого доказуема в $\ZF$).
\end{enumerate}
Все рассуждения проводятся в $\ZF$ поверх свободных множественных
переменных.
\end{proposition}

Кроме этих аксиом исчислению для $\sB$ понадобится одно правило ---
\emph{правило Бернайса}: из $\varphi\pto\psi$ вывести
$\varphi\pto\Pi\alpha\,\psi$ при $\alpha\notin FV(\varphi)$ (и так же
для высказываний). Оно сохраняет общезначимость --- это проверяется
вместе с остальными правилами в теореме~\ref{ru:th:sound}. Остальные
привычные законы квантора --- обобщение, $\Pi\alpha\,\varphi\peq\varphi$
при $\alpha\notin FV(\varphi)$, дистрибутивность относительно
импликации и конъюнкции, обратная к $(\Pi\por)$ импликация --- будут
выведены из $(\Pi\mathrm{E})$ и правила Бернайса
(лемма~\ref{ru:lem:piderived}).

\remove{
Лемма~\ref{ru:lem:qsem} верна и для задач с квантором $\Pi$: её
доказательство не заглядывает внутрь задач $\varphi$ и $\psi$, а
классический разбор случаев проводится в $\ZF$ поверх свободных
множественных переменных. Протолкнуть $\quest$ сквозь сам квантор
$\Pi$ нельзя (замечание после леммы~\ref{ru:lem:semnf}). Исчислениям для
$\sB$ понадобится ещё одно правило --- обобщение по задачной
переменной; его корректность, как и корректность остальных правил, ---
в \S\ref{ru:sec:calc}.
}

\subsection{Сводка: таблица общезначимости}\label{ru:sec:tables}

Таблица~\ref{ru:tab:axioms} сводит общезначимость всех схем, которые в
\S\ref{ru:sec:calc} станут аксиомами исчислений (схемы (A8) и (A9)
исчисления для $\sE$ разобраны в~\S\ref{ru:sec:calcE}). Обозначения: $+$ ---
общезначима; $-$ --- не общезначима (с допущениями
замечания~\ref{ru:rem:assumptions}); $\cdot$ --- не выражается в данной
семантике (нет $\Pi$). 
В последнем
столбце указано, где это доказано; общезначимость схемы (A6) будет
получена в \S\ref{ru:sec:calc} выводом в исчислении (теорема~\ref{ru:th:sound}).

\begin{table}[ht]\centering\small
\caption{Аксиомы исчислений и их общезначимость в трёх семантиках}
\label{ru:tab:axioms}
\medskip
\begin{tabular}{l l c c c l}
 & схема & $\sI$ & $\sB$ & $\sE$ & где доказано\\
\hline
(A1) & классические тавтологии & $+$ & $+$ & $+$ & по определению\\
(A2) & схемы $\Int$ & $+$ & $+$ & $+$ & предл.~\ref{ru:prop:intval}\\
(A3) & $\bang(A\ato B)\pto(\bang A\pto\bang B)$ & $+$ & $+$ & $+$ & предл.~\ref{ru:prop:a3}\\
(A4) & $\quest$ --- булев гомоморфизм & $+$ & $+$ & $+$ & лемма~\ref{ru:lem:qsem}\\
(A5) & $\phi\pto\bang\quest\phi$, $\phi\in\Sigma$ & $+$ & $+$ & $+$ & предл.~\ref{ru:prop:a5a7}\\
(A6) & $\bang(\square A\ato A)\pto\bang A$ & $+$ & $+$ & $+$ & теорема~\ref{ru:th:sound}\\
(A7) & $\pneg\pneg\phi\pto\phi$, $\phi\in\Sigma$
  & $+$ & $+$ & $+$ & предл.~\ref{ru:prop:a5a7}\\
(A8) & рефлексия $\square A\ato A$ & $-$ & $-$ & $+$ & предл.~\ref{ru:prop:hce}(1),(2)\\
(A9) & $\bang\quest\varphi\pto\pneg\pneg\varphi$
  & $-$ & $-$ & $+$ & предл.~\ref{ru:prop:hce}(1),(2)\\
($\Pi$E) & $\Pi\alpha\,\varphi\pto\varphi[\theta/\alpha]$, оба сорта
  & $\cdot$ & $+$ & $\cdot$ & предл.~\ref{ru:prop:pilaws}\\
($\Pi\por$) & $\Pi\alpha(\varphi\por\psi)\pto\Pi\alpha\,\varphi\por\Pi\alpha\,\psi$
  & $\cdot$ & $+$ & $\cdot$ & предл.~\ref{ru:prop:pilaws}\\
($\Pi\bang$) & $\bang A\peq\bang\Pi\alpha\,A$
  & $\cdot$ & $+$ & $\cdot$ & предл.~\ref{ru:prop:pilaws}\\
\hline
\end{tabular}

\end{table}

\remove{
\begin{remark}[что видно из таблицы]\label{ru:rem:read}
Семантики $\sI$, $\sB$ задают \emph{одно и то же} исчисление ---
с двумя оговорками: в $\sI$ $\Pi$-строки отсутствуют в языке, а в
$\sS$ схема ($\Pi$E) ограничена $\Pi$-свободными $\theta$. Семантика
$\sE$ отличается ровно тем, что теряет правило (AtoP) --- оно не
сохраняет $\sE$-общезначимость (теорема~\ref{ru:prop:not}(2)) --- и
приобретает две аксиомы: рефлексию (A8) и рефлексию для задач (A9).
В таблице нет отдельных строк для $\square$-рефлексии
$\square A\ato A$ (это та же (A8), поскольку $\square=\square$),
для формулы Роуза (частный случай (A8)) и для задачи $\pneg\bang\bot$
(частный случай (A9), теорема~\ref{ru:prop:not}(7)).
\end{remark}
}

\section{Исчисления}\label{ru:sec:calc}

Для каждой из трёх
семантик мы определяем исчисление, корректное относительно неё:
$\HCI$, $\HCB$ и $\HCE$. Аксиомы --- в точности схемы из
таблицы~\ref{ru:tab:axioms}, отмеченные плюсом в соответствующем столбце.
Общезначимость всех аксиом исчислений для $\sI$ и $\sB$, кроме (A6), доказана в
\S\ref{ru:sec:validity}, аксиом (A8) и (A9) исчисления для $\sE$ --- в~\S\ref{ru:sec:calcE};
схему (A6) мы не станем проверять напрямую --- она будет выведена из
остальных аксиом и леммы о неподвижной точке (теорема~\ref{ru:th:sound}).

\subsection{Исчисление для семантики $\sI$}\label{ru:sec:calcI}

\begin{definition}[система $\HCI$]\label{ru:def:hck}
Схемы аксиом:
\[
\begin{array}{ll}
  \text{(A1)} & \text{классические тавтологии (сорт высказывание);}\\[2pt]
  \text{(A2)} & \text{схемы аксиом интуиционистской логики (сорт задача);}\\[2pt]
  \text{(A3)} & \bang(A\ato B)\pto(\bang A\pto\bang B);\\[2pt]
  \text{(A4)} & \text{$\quest$ --- булев гомоморфизм:}\\[2pt]
   & \qquad\text{(A4a)}\ \ \quest(\varphi\pto\psi)\aeq(\quest\varphi\ato\quest\psi),\\[2pt]
   & \qquad\text{(A4b)}\ \ \quest(\varphi\pand\psi)\aeq(\quest\varphi\aand\quest\psi),\\[2pt]
   & \qquad\text{(A4c)}\ \ \quest(\varphi\por\psi)\aeq(\quest\varphi\aor\quest\psi),\\[2pt]
   & \qquad\text{(A4d)}\ \ \quest\varnothing\aeq\bot;\\[2pt]
  \text{(A5)} & \phi\pto\bang\quest\phi\quad\text{для $\Sigma$-задач $\phi$ (ко-рефлексия на классе $\Sigma$);}\\[2pt]
  \text{(A6)}& \bang(\square A\ato A)\pto\bang A\ \ \text{(конструктивный вариант теоремы Лёба);}\\[2pt]
  \text{(A7)} & \pneg\pneg\varphi\pto\varphi\quad\text{для $\Sigma$-задач $\varphi$.}
\end{array}
\]
Во всех схемах $A$, $B$ --- произвольные
высказывания, $\varphi$, $\psi$ --- произвольные задачи (а в (A7) ---
произвольные задачи указанного класса); схема обозначает множество всех своих
частных случаев. Все формулы здесь $\Pi$-свободные. Правила вывода:
\[
\begin{array}{lll}
  \text{(MP)} & \dfrac{A,\ A\ato B}{B}\ \ \text{и}\ \
                \dfrac{\varphi,\ \varphi\pto\psi}{\psi}
              & \text{(для обоих сортов);}\\[12pt]
  \text{(AtoP)} & \dfrac{A}{\bang A}
              & \text{из высказывания --- задача;}\\[12pt]
  \text{(PtoA)} & \dfrac{\varphi}{\quest\varphi}
              & \text{из задачи --- высказывание.}
\end{array}
\]
Названия правил читаются как <<assertion to problem>> и <<problem to
assertion>>: (AtoP) по высказыванию $A$ даёт задачу $\bang A$ <<найти
доказательство $A$>>, а (PtoA) по задаче $\varphi$ даёт высказывание
$\quest\varphi$ <<у $\varphi$ есть решение>>.

Частный случай $\sigma:=\bang A$ схемы (A5) --- задача
$\bang A\pto\bang\square A$: по выводу формулы $A$ можно построить
вывод формулы $\square A$.

Через $\HCglm$ обозначаем то же исчисление \emph{без} схем (A6) и (A7). Оно
понадобится потому, что схему Лёба мы не будем проверять напрямую: она будет
выведена в $\HCglm$ из леммы о неподвижной точке. Через $\HCzero$
(\emph{ядро}) обозначаем ещё меньшую систему: схемы (A1), (A2), (A3),
\emph{левая} половина (A4a), (A4b)--(A4d), (A5) и все три правила вывода. Ядро --- это
ровно та часть нашей аксиоматики, которая содержится в логике $\HC$
Мелихова. 
\end{definition}

\subsection{Исчисление для семантики $\sB$}\label{ru:sec:picalc}

\begin{definition}[система $\HCB$]\label{ru:def:hczf}
Формулы --- формулы полного языка с квантором $\Pi$
(\S\ref{ru:sec:valBS}).
Аксиомы системы $\HCB$ --- схемы (A1)--(A7) и
аксиомы оператора $\Pi$:
\[
  (\Pi\mathrm{E})\quad\text{(для обоих сортов)},\qquad
  (\Pi\por),\qquad (\Pi\bang)
\]
(формулировки --- предложение~\ref{ru:prop:pilaws}; в
$(\Pi\mathrm{E})$ подразумевается, что $\theta$ свободна для
$\alpha$). Правила вывода: (MP), (AtoP), (PtoA) и правило Бернайса
для обоих сортов
\[
\begin{array}{lll}
  \text{(Bern)} & \dfrac{\varphi\pto\psi}{\varphi\pto\Pi\alpha\,\psi}
    \ \ \text{и}\ \ \dfrac{A\ato B}{A\ato\Pi\alpha\,B}
  & \text{(для обоих сортов),}
\end{array}
\]
где $\alpha$ не входит свободно в посылку $\varphi$ (соответственно $A$).
\end{definition}

Таким образом, бескванторная часть аксиоматики $\HCB$ совпадает с
$\HCI$; что даёт кванторный слой, обсуждается в
\S\ref{ru:sec:pigain}.

\subsection{Исчисление для семантики $\sE$}\label{ru:sec:calcE}

\begin{definition}[система $\HCE$]\label{ru:def:hcks}
Аксиомы системы $\HCgls$ --- все теоремы $\HCgl$, а также схемы
\[
  \text{(A8)}\ \ \square A\ato A\ \ \text{(рефлексия)},
  \qquad
  \text{(A9)}\ \ \bang\quest\varphi\pto\pneg\pneg\varphi\ \
  \text{(рефлексия для задач)}
\]
($A$ --- произвольное высказывание, $\varphi$ --- произвольная задача).
Правила вывода: (MP) и (PtoA); правила (AtoP) \textbf{нет}. 
\end{definition}

\begin{proposition}[схемы (A8) и (A9) и семантики]\label{ru:prop:hce}
Пусть $A$ --- переменная сорта высказывание, $\alpha$ --- переменная сорта
задача.
\begin{enumerate}[leftmargin=2.2em,itemsep=2pt,label=\textup{(\arabic*)}]
\item Схемы (A8) и (A9) $\sE$-общезначимы.
\item Схемы (A8) и (A9) не $\sI$- и не $\sB$-общезначимы: не общезначимы
  уже их экземпляры $\square\bot\ato\bot$ и
  $\bang\quest\varnothing\pto\pneg\pneg\varnothing$.
\item Задача $\bang(\square\bot\ato\bot)$ не общезначима ни в одной из трёх
  семантик.
\item Ни задача $\bang\quest\alpha\pto\alpha$, ни обратная к ней
  ко-рефлексия $\alpha\pto\bang\quest\alpha$ не общезначимы ни в одной из
  трёх семантик.
\end{enumerate}
\end{proposition}

\begin{proof}
(1) Схема (A8). В $\sE$: если $\val{\quest\bang A}=\Tru$, то $\PA\vdash A^*$, и формула
$A^*$ истинна, поскольку внутренняя теория $\PA$ содержится во внешней
$\TA$.

Схема (A9). В $\sE$ годится $S(\gn{*}):=e$, где $e$ --- номер постоянной функции
$u\mapsto0$; он один и тот же при всех $*$, так что никакого разбора
случаев алгоритм не делает. Проверим, что этот номер подходит. Если
$\val{\bang\quest\varphi}=\varnothing$, то 
$\val{\bang\quest\varphi\pto\pneg\pneg\varphi}=\N$ и годится любой номер.
Если же $\val{\bang\quest\varphi}\neq\varnothing$, то
$\PA\vdash\exists x\,\varphi^*(x)$, и отсюда $\val{\varphi}\neq\varnothing$, и по лемме~\ref{ru:lem:dn}
$\val{\pneg\pneg\varphi}=\N$; значит $U_e$ попадает в
$\val{\pneg\pneg\varphi}$ при любом входе.

(2) Схема (A8). 
Перевод формулы $\square\bot\ato\bot$ рефлексии есть $\Con$, а
непротиворечивость внутренней теории в ней самой (она же внешняя) не
выводима --- вторая теорема Гёделя.

Схема (A9). Во внутренних семантиках не общезначима уже формула 
$\bang\quest\varnothing\pto\pneg\pneg\varnothing$. 
В самом деле, 
задача $\pneg\pneg\varnothing$ доказуемо
эквивалентна $\varnothing$, а $(\quest\varnothing)^*$ ---
формуле $\arbot $ (лемма~\ref{ru:lem:qsem}), так что перевод задачи
$\bang\quest\varnothing\pto\pneg\pneg\varnothing$ доказуемо эквивалентен
переводу задачи $\pneg\bang\bot$, то есть  формуле
$\forall v\,(\Prf(v,\arbot )\ato \arbot )$, доказуемо эквивалентной $\Con$;
а непротиворечивость внутренней теории в ней самой (она же внешняя) не
выводима --- вторая теорема Гёделя.

(3) Допустим, что эта задача общезначима, и обозначим
через $S$ её решатель. 
Применив его к формуле $\arbot $ получим натуральное число $n$,
для которого внешняя теория доказывает
формулу $\Prf(n, \square\bot\ato\bot)$. Это $\Sigma_0$-предложение;
будь оно ложно, его отрицание было бы доказуемо ($\Sigma_1$-полнота)
вопреки непротиворечивости внешней теории. Значит оно истинно: $n$ --- вывод
формулы $\square\bot\ato\bot$ во внутренней теории, то есть внутренняя
теория доказывает собственную непротиворечивость --- вопреки второй
теореме Гёделя.

(4) Для арифметических семантик: 
фиксируем перечислимое неразрешимое множество
$X\subset\N$ и выразим свойство $n\in X$ арифметической
сигма формулой $\exists x A(x,n)$. 
Положим
$\alpha^{*_n}(x):=(x=0\aand\aneg A(x,n))\aor(x=1\aand A(x,n)))$.
Тогда $\PA$ доказывает
$\exists x\,\alpha^{*_n}(x)$ при всяком натуральном $n$, так что
$\val[*_n]{\bang\quest\alpha}\neq\varnothing$, а $\val[*_n]{\alpha}$ есть $\{0\}$ или
$\{1\}$ смотря по тому, принадлежит ли $n$ к $X$; решение выдавало бы
разрешающий алгоритм для $X$.

Эта формула не является
чистой, поэтому в теоретико-множественной семантике доказательство совсем другое.
Равносильная ей чистая формула $\bang\quest\beta\pto\alpha$ равносильна $\bang\bot\pto\alpha$ и не общезначима: её экземпляр при
$\alpha:=\varnothing$ --- задача $\pneg\bang\bot$, перевод которой
доказуемо эквивалентен $\Con$, так что при непротиворечивости $\ZF$
решений у неё нет.

Ко-рефлексия. Для арифметических семантик
допустим, что существует алгоритм $S$, доказуемо доставляющий 
по $*$ решение задачи $(\alpha\pto\bang\quest\alpha)^*$.
Положим $\alpha^*$ равной $\Con$ (без свободных переменных),  то есть, арифметической формуле,
выражающей непротиворечивость внутренней теории. 
Положим $n=S(*)$.
Тогда во внешней теории доказуема формула  
$$
(\alpha\pto\bang\quest\alpha)^*(n)=
\forall x(\Con\ato \Prf(U_n(x),\exists y\,\Con)),
$$
которая доказуемо влечёт 
$$ 
\Con\ato  \exists z \Prf(z,\Con).
$$ 
Между тем,
по второй теореме Гёделя 
во внутренней (а значит и во внешней) теории 
выводимо
 $$ 
\Con\ato \aneg \exists z \Prf(z,\Con).
$$ 
Следовательно, во внешней теории выводима ложная формула
$\aneg\Con$. 

Для семантики $\sB$ сначала заметим, что переменная $\alpha$ в заключении формулы $\alpha\pto\bang\quest\alpha$
связана, поэтому она означает совсем не то, что означает переменная 
$\alpha$ в посылке.
%
Равносильная ей чистая формула $\alpha\pto\bang\quest\beta$ равносильна
формуле $\alpha\pto\bang\bot$ (замечание~\ref{ru:rem:bangbinds}), а та не
общезначима: её экземпляр при $\alpha:=\pneg\varnothing$ равносилен
задаче $\bang\bot$, у которой при непротиворечивости $\ZF$ решений нет.
\end{proof}

Теперь видно, почему в $\HCE$ нет правила (AtoP): высказывание
$\square\bot\ato\bot$ $\sE$-общезначимо (частный случай (A8), пункт~(1)), а
применение к нему правила (AtoP) даёт задачу $\bang(\square\bot\ato\bot)$,
не общезначимую ни в одной семантике (пункт~(3)).

\begin{remark}[зачем (A9) и почему её нельзя усилить]\label{ru:rem:a10}
Схема (A9) --- это в точности рефлексия, перенесённая на сорт задача:
высказывание $\quest\varphi$ и задача $\pneg\pneg\varphi$ имеют, как показывает
лемма~\ref{ru:lem:dn}, один и тот же смысл <<у $\varphi$ есть решение>>, и (A9)
говорит: если доказано, что решение есть, то оно есть. Усилить её до
$\bang\quest\varphi\pto\varphi$ нельзя: уже для переменной эта задача не
общезначима (предложение~\ref{ru:prop:hce}(4)). Двойное отрицание в правой части ---
ровно то ослабление, которое стирает эту информацию: по
лемме~\ref{ru:lem:clauses}(6) множество $\val{\pneg\pneg\varphi}$ есть $\N$ или
$\varnothing$, и выбирать в нём нечего.
\end{remark}

Корректность всех трёх исчислений
--- теорема~\ref{ru:th:sound} (\S\ref{ru:sec:sound}); сначала --- производные
правила и примеры выводов.

Еще один естественный вопрос: почему в (A5) и (A7) взяты именно такие классы
задач. Краткий ответ: мы взяли наибольшие естественные синтаксически определяемые
классы задач, для которых эти схемы общезначимы. Более подробный ответ --- в разделе~\ref{ru:sec:whyclasses} в Приложении.

\subsection{Производные правила вывода}\label{ru:sec:derived}

\begin{proposition}\label{ru:prop:derived}
В исчислениях $\HCglm$, $\HCgl$, $\HCB$ допустимы правила
вывода~(1), (4), (5) и верны лемма о дедукции~(2) и лемма о
переводе выводов~(3). В $\HCE$, где нет правила (AtoP), верны лемма о
дедукции~(2), вторая половина леммы~(3) и правило~(4); правила (1),
(5) и первая половина леммы~(3) опираются на (AtoP).
\begin{enumerate}[leftmargin=2.2em,itemsep=2pt,label=\textup{(\arabic*)}]
\item правило необходимости для $\square$: из $A$ можно вывести $\square A$.
\item Лемма о дедукции. Будем  \emph{выводом из гипотез} $\varphi_1,\dots,\varphi_n$
называть вывод с использованием только правила (MP), в котором в качестве
аксиом используются формулы $\varphi_1,\dots,\varphi_n$ и любые выводимые в
$\HCglm$ формулы. Тогда, если существует вывод задачи $\chi$ из гипотез $\varphi_1,\dots,\varphi_n$, то в $\HCglm$
выводима задача
\[
\varphi_1\pto(\varphi_2\pto(\dots(\varphi_n\pto\chi)\dots)).
\]
Буквы здесь греческие: гипотезы и заключение --- задачи. Аналогичное верно
для высказываний и для $\HCgl$.
\item Лемма о переводе выводов.
Если высказывание $A$ выводимо из гипотез $B_1,\dots,B_n$ в $\HCglm$,
то задача $\bang A$ выводима из гипотез $\bang B_1,\dots,\bang B_n$ в $\HCglm$.
Аналогично, если задача $\chi$ выводима из гипотез
$\varphi_1,\dots,\varphi_n$ в $\HCglm$, то высказывание $\quest\chi$ выводимо из
гипотез $\quest\varphi_1,\dots,\quest\varphi_n$ в $\HCglm$. Аналогичное верно
и для $\HCgl$.
\item правило ($\quest\to$):  
Из $\varphi\pto\psi$ можно вывести 
$\quest\varphi\ato\quest\psi$. 

\item правило ($\bang\to$):  
Из $A\ato B$ можно вывести 
$\bang A\pto\bang B$.

\end{enumerate}
\end{proposition}
\begin{proof}

(1) (AtoP), затем (PtoA).

(2) Годится обычное доказательство леммы о дедукции.

(3) Пусть дан вывод высказывания $A$ из гипотез $B_1,\dots,B_n$ в $\HCglm$.
Заменим в этом выводе каждое высказывание $C$ на задачу $\bang C$.
В результате этой замены все выводимые в $\HCglm$
формулы (принимаемые в качестве аксиом) превратятся по
правилу (AtoP) снова в выводимые формулы. Гипотезы $B_1,\dots,B_n$ превратятся
в данные нам гипотезы $\bang B_1,\dots,\bang B_n$, а последняя формула $A$ ---
в искомую формулу $\bang A$.
Однако полученная последовательность формул перестанет быть выводом
по правилу MP. Но это легко исправить, поскольку если формула $U$ получена из формул $V$ и
$V\ato U$ по правилу MP, то из формул $\bang V$ и
$\bang(V\ato U)$ и аксиомы (A3) можно вывести $\bang U$, применив дважды MP.
Второй перевод обосновывается так же, только вместо (AtoP) применяется
(PtoA), а вместо (A3) --- левая половина (A4a): из $\quest V$ и
$\quest(V\pto U)$ она даёт $\quest U$.

(4) По (PtoA) и левой половине (A4a).

(5) По (AtoP) и (A3).
\end{proof}


Далее мы докажем  лемму о замене на эквивалентную формулу. Она требует аккуратности: наши
формулы двусортны, и заменяемая подформула может стоять под модальностями.
Подчеркнём заранее, что никакие аксиомы <<внутрь модальности>> при этом не
вносятся: замена под знаками $\bang$ и $\quest$ производится подъёмом
\emph{выводимой} эквивалентности сквозь модальность правилами (AtoP) и
(PtoA), а они применимы именно потому, что эквивалентность --- теорема.

\begin{lemma}[о замене на эквивалентную формулу]\label{ru:lem:repl}
Пусть в $\HCglm$ выводима эквивалентность $C\aeq D$ (высказывания) или
$\chi\peq\chi'$ (задачи), и формула $E'$ (любого сорта) получена из формулы
$E$ заменой некоторых вхождений подформулы $C$ на $D$ (соответственно $\chi$
на $\chi'$). Тогда в $\HCglm$ выводима эквивалентность $E\aeq E'$ или
$E\peq E'$, смотря по сорту $E$. То же верно в $\HCgl$, в $\HCgls$ и
--- для формул с квантором $\Pi$ --- в $\HCB$.
\end{lemma}

\begin{proof}
Индукция по построению формулы $E$. Если $E$ --- заменяемое вхождение, то эквивалентность дана по
условию; если же вхождений внутри $E$ нет, то $E=E'$. 

Шаг разбирается по главной связке $E$.
Связки высказываний: если $\vdash E_1\aeq E_1'$ и $\vdash E_2\aeq E_2'$, то
$\vdash(E_1\star E_2)\aeq(E_1'\star E_2')$ для
$\star\in\{\aand,\aor,\ato\}$, поскольку
$(X\aeq X')\ato\bigl((Y\aeq Y')\ato\bigl((X\star Y)\aeq(X'\star
Y')\bigr)\bigr)$ --- классическая тавтология, схема (A1); остаётся (MP).

Связки задач: аналогичные конгруэнц-схемы
$(\varphi\peq\varphi')\pto\bigl((\psi\peq\psi')\pto\bigl((\varphi\star\psi)
\peq(\varphi'\star\psi')\bigr)\bigr)$ для
$\star\in\{\pand,\por,\pto\}$ выводимы в $\Int$, то есть входят в (A2).

Модальность $\bang$: пусть $\vdash A\aeq A'$. Из обеих импликаций по (A1)
получаем $\vdash A\ato A'$ и $\vdash A'\ato A$; правило ($\bang\to$) даёт
$\bang A\pto\bang A'$ и $\bang A'\pto\bang A$. 

Модальность $\quest$: аналогично, только вместо ($\bang\to$)
используем ($\quest\to$)


Квантор $\Pi$ (в $\HCB$): пусть $\vdash\varphi\peq\varphi'$. По
$(\Pi\mathrm{E})$ выводима $\Pi\alpha\,\varphi\pto\varphi$, по (A2) ---
$\varphi\pto\varphi'$, их композиция --- $\Pi\alpha\,\varphi\pto\varphi'$;
переменная $\alpha$ в $\Pi\alpha\,\varphi$ свободно не входит, и правило
Бернайса даёт $\Pi\alpha\,\varphi\pto\Pi\alpha\,\varphi'$; симметрично
обратная импликация. 
Для высказываний
$\Pi\alpha\,A$ то же, с высказывательной $(\Pi\mathrm{E})$, правилом
Бернайса и схемами (A1).
%
\end{proof}

Для исчислений с квантором нужны ещё производные законы самого $\Pi$.

\begin{lemma}[производные законы $\Pi$]\label{ru:lem:piderived}
В $\HCB$ выводимы следующие задачи.
\begin{enumerate}[leftmargin=2.2em,itemsep=1pt,label=\textup{(\arabic*)}]
\item $\Pi\alpha\,\varphi\peq\varphi$, если $\alpha\notin FV(\varphi)$.
\item $(\Pi\pto)$: $\Pi\alpha(\varphi\pto\psi)\peq(\varphi\pto\Pi\alpha\,\psi)$,
  если $\alpha\notin FV(\varphi)$.
\item $\Pi\alpha\Pi\beta\,\varphi\peq\Pi\beta\Pi\alpha\,\varphi$.
\item $(\Pi\pand)$:
  $\Pi\alpha(\varphi\pand\psi)\peq\Pi\alpha\,\varphi\pand\Pi\alpha\,\psi$,
  а также обратная к аксиоме $(\Pi\por)$ импликация
  $\Pi\alpha\,\varphi\por\Pi\alpha\,\psi\pto\Pi\alpha(\varphi\por\psi)$.
\end{enumerate}
Те же утверждения верны и для высказываний:
$\Pi\alpha\,A\aeq A$ при $\alpha\notin FV(A)$;
$\Pi\alpha(A\ato B)\aeq(A\ato\Pi\alpha\,B)$ при $\alpha\notin FV(A)$;
$\Pi\alpha\Pi\beta\,A\aeq\Pi\beta\Pi\alpha\,A$;
$\Pi\alpha(A\aand B)\aeq\Pi\alpha\,A\aand\Pi\alpha\,B$ и
$\Pi\alpha\,A\aor\Pi\alpha\,B\ato\Pi\alpha(A\aor B)$.
\end{lemma}

\begin{proof}
Всюду ниже <<правило Бернайса>> применяется к импликации, в посылку
которой $\alpha$ свободно не входит; это каждый раз видно по виду
посылки.

(1) Импликация $\Pi\alpha\,\varphi\pto\varphi$ --- это $(\Pi\mathrm{E})$
при $\theta:=\alpha$; обратная получается правилом Бернайса из
теоремы $\varphi\pto\varphi$.

(2) Слева направо: из $(\Pi\mathrm{E})$
$\Pi\alpha(\varphi\pto\psi)\pto(\varphi\pto\psi)$ в $\Int$ выводима
$\varphi\pand\Pi\alpha(\varphi\pto\psi)\pto\psi$; правило Бернайса даёт
$\varphi\pand\Pi\alpha(\varphi\pto\psi)\pto\Pi\alpha\,\psi$, и остаётся
снова перейти к импликации в $\Int$. Справа налево: из
$(\Pi\mathrm{E})$ $\Pi\alpha\,\psi\pto\psi$ в $\Int$ выводима
$(\varphi\pto\Pi\alpha\,\psi)\pto(\varphi\pto\psi)$; правило Бернайса
даёт $(\varphi\pto\Pi\alpha\,\psi)\pto\Pi\alpha(\varphi\pto\psi)$.

(3) Из $(\Pi\mathrm{E})$ дважды выводима
$\Pi\alpha\Pi\beta\,\varphi\pto\varphi$; правило Бернайса по $\alpha$,
затем по $\beta$ даёт
$\Pi\alpha\Pi\beta\,\varphi\pto\Pi\beta\Pi\alpha\,\varphi$; обратная
импликация симметрична.

(4) $(\Pi\pand)$ слева направо: из $(\Pi\mathrm{E})$ и проекций
$\Int$ выводимы $\Pi\alpha(\varphi\pand\psi)\pto\varphi$ и
$\Pi\alpha(\varphi\pand\psi)\pto\psi$; правило Бернайса к каждой и
$\pand$-введение. Справа налево: из $(\Pi\mathrm{E})$ дважды
$\Pi\alpha\,\varphi\pand\Pi\alpha\,\psi\pto\varphi\pand\psi$, затем
правило Бернайса. 

Обратная к $(\Pi\por)$ импликация: из
$(\Pi\mathrm{E})$ и $\Int$ выводима $\Pi\alpha\,\varphi\pto\varphi\por\psi$,
правило Бернайса даёт $\Pi\alpha\,\varphi\pto\Pi\alpha(\varphi\por\psi)$;
так же для $\Pi\alpha\,\psi$; остаётся $\por$-исключение в $\Int$.

Для высказываний выводы те же, с высказывательной аксиомой
$(\Pi\mathrm{E})$, правилом Бернайса для высказываний и классическими
тавтологиями вместо теорем $\Int$.
\end{proof}

\remove{
Выводимы также правило обобщения и $(\Pi\mathrm{K})$: из $\varphi$
выводима $\Pi\alpha\,\varphi$, и выводима
$\Pi\alpha(\varphi\pto\psi)\pto(\Pi\alpha\,\varphi\pto\Pi\alpha\,\psi)$
(так же для высказываний). Обобщение: из $\varphi$ в $\Int$ выводима
$(\varnothing\pto\varnothing)\pto\varphi$; правило Бернайса даёт
$(\varnothing\pto\varnothing)\pto\Pi\alpha\,\varphi$, и (MP) с теоремой
$\varnothing\pto\varnothing$ --- $\Pi\alpha\,\varphi$ (для высказываний
вместо $\varnothing\pto\varnothing$ берётся $\aneg\bot$).
$(\Pi\mathrm{K})$: из $(\Pi\mathrm{E})$ дважды в $\Int$ выводима
$\Pi\alpha(\varphi\pto\psi)\pand\Pi\alpha\,\varphi\pto\psi$; правило
Бернайса даёт
$\Pi\alpha(\varphi\pto\psi)\pand\Pi\alpha\,\varphi\pto\Pi\alpha\,\psi$,
и остаётся перейти к импликации.
}

\subsection{Примеры выводов}\label{ru:sec:examples}

\providecommand{\provgl}{\vdash_{\HCgl}}
\providecommand{\provglm}{\vdash_{\HCglm}}
\providecommand{\provgls}{\vdash_{\HCgls}}
\providecommand{\provglb}{\vdash_{\HCB}}

\begin{example}\label{ru:ex:proofs}
Приведем примеры выводов в исчислениях $\HCglm$, $\HCI$, $\HCE$.
Некоторые из них нам понадобятся далее; в каждом пункте указано, в каком
именно исчислении ведётся вывод.
\begin{enumerate}[leftmargin=2.2em,itemsep=2pt,label=\textup{(\arabic*)}]
\item В $\HCglm$ выводима $\square$-версия (A3):
  $\square(A\ato B)\ato(\square A\ato\square B)$. Для вывода применим ($\quest\to$) 
к (A3), затем (A4a) и Лемму о замене на эквивалентную формулу.

\item В $\HCglm$ выводимо $A\ato\square A$ для всякого
$\Sigma$-высказывания $A$ ($\Sigma_1$-полнота на стороне высказываний);
в частности выводима $\square$-версия (A5), схема
$\square A\ato\square\square A$, поскольку $\square A$ ---
$\Sigma$-высказывание. Сначала заметим, что всякое $\Sigma$-высказывание
$A$ выводимо эквивалентно высказыванию $\quest\sigma_A$ для некоторой
$\Sigma$-задачи $\sigma_A$: по (A4d) $\bot\aeq\quest\varnothing$, по (A4b) и
(A4c) $\quest\sigma\aand\quest\tau\aeq\quest(\sigma\pand\tau)$ и
$\quest\sigma\aor\quest\tau\aeq\quest(\sigma\por\tau)$, а задачи
$\varnothing$, $\sigma\pand\tau$, $\sigma\por\tau$ снова лежат в классе
$\Sigma$; остаётся индукция по построению $A$ с леммой~\ref{ru:lem:repl}.
Теперь аксиома (A5) даёт $\sigma_A\pto\bang\quest\sigma_A$, правило
$(\quest\to)$ --- $\quest\sigma_A\ato\square\quest\sigma_A$, и
лемма~\ref{ru:lem:repl} заменяет $\quest\sigma_A$ на $A$ в обеих частях.


\item В $\HCgl$ выводимо
\[   
 \text{(D4)}\  \square (\square A\ato A)\ato\square A\ \text{(теорема Лёба)}
  \]
Применим ($\quest\to$) к (A6).

\item  В $\HCE$ выводимы формулы $\pneg\bang\bot$ и $\aneg\square\bot$.
Сначала применим правило $(\bang\to)$ к формуле $\bot\ato\quest\varnothing$, классической тавтологии (A1). 
Мы получим формулу $\bang\bot\pto\bang\quest\varnothing$.

 Далее, схема (A9) при
$\varphi=\varnothing$ даёт $\bang\quest\varnothing\pto\pneg\pneg\varnothing$, а
$\pneg\pneg\varnothing\pto\varnothing$ --- интуиционистская теорема
. Итого $\bang\bot\pto\varnothing$, то
есть $\pneg\bang\bot$. Теперь применим к  этой формуле ($\quest\to$), 
получим $\square\bot\pto\quest\varnothing$. Наконец, воспользовавшись (A4d),
получим $\square\bot\pto\bot$.
(Заметим в скобках, что последняя формула является частным случаем аксиомы A8,
так что ее можно было вывести короче.)

\item Этот пример на выводимость в $\HCzero$ нам понадобится в доказательстве общезначимости аксиомы (A6).
 Пусть $P$ и $S$ --- произвольные высказывания. Тогда  в $\HCglm$ выводима задача
\begin{equation}\tag{L}\label{ru:eq:lobimpl}
  \bang\bigl(P\aeq(\square P\ato S)\bigr)\ \pto\
  \bigl(\bang(\square S\ato S)\pto\bang S\bigr).
\end{equation}
Сначала мы воспользуемся леммой о переводе выводов, а затем леммой о дедукции. Для этого 
примем в качестве гипотез высказывания
$$
\square S\ato S, \qquad
\square (P\ato (\square P\ato S)), \qquad (\square P\ato S)\ato P
$$
и выведем формулы $\square P\ato S$ и $P$. 

\begin{enumerate}[leftmargin=2.6em,itemsep=5pt,label=\textup{(\arabic*)}]

\item $\square P\ato\square(\square P\ato S)$\\
  \emph{эта формула выводится по MP из второй гипотезы и первого примера выводимых формул
}
\item $\square(\square P\ato S)\ato(\square\square P\ato\square S)$\\
  \emph{это формула из нашего первого примера выводимых формул при $A:=\square P$, $B:=S$.}
\item $\square P\ato\square\square P$\\
  \emph{это формула из второго примера выводимых формул.}
\item $\square P\ato\square S$\\
  \emph{Выводится из  (1), (2), (3) в классическом ИВ.}
\item $\square P\ato S$\\
  \emph{Вывод из (4) и первой гипотезы в классическом ИВ.}
\item $P$\\
\emph{из третьей гипотезы и  (5) по MP.}

\item  Заметим, что весь вывод (1)--(6) является выводом из трёх принятых гипотез.
Поэтому можно применить лемму о переводе выводов. По этой лемме
имеется вывод формул $\bang(\square P\ato S)$ и $\bang P$ из гипотез
$$
\bang(\square S\ato S),
\quad \bang\square (P\ato(\square P\ato S)),
\quad \bang((\square P\ato S)\ato P).
$$
Докажем, что этот вывод можно продолжить до вывода задачи $\bang S$.
\item  Из тех же гипотез можно вывести ещё и формулу $\bang\square P$
(приписав к выводу аксиому (A5) при $\sigma:=\bang P$ и применив MP).
\item Приписав к выводу аксиому (A3)
и применив MP, можно далее вывести задачу $\bang\square P\pto\bang S$. 
\item  Наконец,  можно вывести формулу $\bang S$, ещё раз
применив MP.

Итак, из этих трёх гипотез выводится формула $\bang S$.

\item Теперь покажем, что вторую и третью гипотезу можно заменить на одну гипотезу: 
$$
F:=\bang\bigl(P\aeq(\square P\ato S)\bigr),
$$ 
являющуюся первой посылкой исходной формулы.
Во-первых, из $F$ выводится третья гипотеза: высказывание
$(\square P\ato S)\ato P$ выводится из $P\aeq(\square P\ato S)$ в
классическом исчислении, а лемма о переводе выводов переносит этот вывод
под знак~$\bang$.
Аналогичным образом из $F$ выводится формула
$\bang(P\ato(\square P\ato S))$. Она отличается от второй гипотезы только
наличием бокса после восклицательного знака, а аксиома (A5) при
$\sigma:=\bang(P\ato(\square P\ato S))$ как раз позволяет по MP вывести из неё вторую
гипотезу.

Итак, из гипотез $F$ и $\bang(\square S\ato S)$
можно вывести формулу $\bang S$, причём одними лишь применениями MP. Осталось
применить лемму о дедукции и получить формулу~\eqref{ru:eq:lobimpl}.

Отметим, что вывод этого пункта был формализован и проверен системой
Metamath, причём в полностью развёрнутом виде: лемма о дедукции и лемма о
переводе выводов в формализации не используются, а раскрыты
в явный вывод из аксиом (около семисот элементарных шагов). Проверены и
выводы пунктов (1)--(4).
\end{enumerate} 
\end{enumerate}
\end{example}

\begin{remark}[о классе в (A5)]
Аксиома (A5) записана сразу для всего класса $\Sigma$, хотя её
частного случая $\sigma:=\bang A$ было бы достаточно. В самом деле,
индукцией по построению $\Sigma$-задачи $\sigma$ выводится
$\sigma\pto\bang\quest\sigma$: для $\varnothing$ --- по $\Int$, для
$\bang A$ --- частный случай; если же $\sigma=\sigma_1\pand\sigma_2$, то
из $\sigma_i\pto\bang\quest\sigma_i$ получаем
$\sigma\pto\bang\quest\sigma_1\pand\bang\quest\sigma_2$, далее
$\sigma\pto\bang(\quest\sigma_1\aand\quest\sigma_2)$ (из тавтологии
$A\ato(B\ato A\aand B)$ правилом (AtoP) и схемой (A3)), а
$\quest\sigma_1\aand\quest\sigma_2\aeq\quest\sigma$ по (A4b), и
лемма~\ref{ru:lem:repl} заканчивает; для $\por$ так же, с правилом
$(\bang\to)$ из тавтологий $A\ato A\aor B$, $B\ato A\aor B$ и с (A4c).
Это то же наблюдение, что в пункте~(2) примера~\ref{ru:ex:proofs}: схемы
(A4) склеивают $\quest$-образы $\Sigma$-задач в $\quest$-образ одной
$\Sigma$-задачи. Записать (A5) сразу для класса $\Sigma$
естественнее: семантически это эффективная $\Sigma_1$-полнота
\emph{со свидетелем на входе} (предложение~\ref{ru:prop:a5a7}), и
ограничение на $\bang A$ в ней случайно.

Аксиома (A5) и пункт~(2) примера~\ref{ru:ex:proofs} --- одно и то же
утверждение, увиденное с двух сторон: эффективная $\Sigma_1$-полнота в
$\bang$-форме порождает и ко-рефлексию на $\Sigma$-задачах, и
$\Sigma_1$-полноту на $\Sigma$-высказываниях, а $\Sigma$-высказывания
--- это с точностью до выводимой эквивалентности высказывания
$\quest\sigma$ для $\Sigma$-задач $\sigma$. Границы у обоих классов
поэтому одинаковые. Переменную в $\Sigma$-высказывания не пустить: при
$A^*:=\Con$ схема $A\ato\square A$ даёт $\Con\ato\Prov(\Con)$, что ложно по
второй теореме Гёделя. Импликацию тоже: $\square\bot\ato\bot$ есть $\Con$, и
получается то же самое.
(Конец замечания.)
\end{remark}

\remove{
\begin{remark}[почему в $\HCgls$ нет (AtoP)]\label{ru:rem:lobfork}
Рефлексия и необходимость несовместимы. Из (A8) применением (AtoP) и
(PtoA) получается формула $\square(\square A\ato A)$, схема Лёба (D4) даёт $\square A$, а рефлексия --- формулу $A$. При этом
$A$ здесь может быть любым высказыванием. Именно
поэтому в $\HCE$ остаются только (MP) и (PtoA); непротиворечивость
$\HCE$ следует из теоремы~\ref{ru:th:sound}.
\end{remark}

\begin{remark}[о $\square$-фрагментах]\label{ru:rem:names}
В $\HCI$ модальность $\square$ удовлетворяет схеме K, схеме
$\square A\ato\square\square A$ (пример~\ref{ru:ex:proofs}(1),(2)) и
правилу необходимости (предложение~\ref{ru:prop:derived}(1)), то есть
аксиоматике $\GL$ без схемы Лёба, --- это $\mathrm{K4}$; добавление
(A6) даёт $\GL$. Система $\HCE$ добавляет рефлексию при единственных
правилах (MP) и (PtoA) --- так же, как $\GLS$ получается из $\GL$.
$\square$-фрагменты и совпадают с этими логиками в точности ---
теорема~\ref{ru:th:assval} ниже.
\end{remark}
}

\subsection{Корректность исчислений}\label{ru:sec:sound}

\begin{theorem}[корректность]\label{ru:th:sound}
Для каждой из трёх семантик $Z\in\{\sI,\sB,\sE\}$ всякая теорема
исчисления $\HC_Z$ $Z$-общезначима.
\end{theorem}

\begin{proof}
Индукция по выводу: нужно проверить, что каждая аксиома общезначима в
соответствующей семантике и что каждое правило сохраняет
общезначимость.

(1) \emph{Аксиомы, кроме (A6).} Общезначимость во всех трёх
семантиках схем (A1) видна прямо: перевод --- подстановочный экземпляр
классической тавтологии, и его вывод вычислимо строится по коду
интерпретации. Общезначимость схем (A2) --- это
предложение~\ref{ru:prop:intval}, схемы (A3) ---
предложение~\ref{ru:prop:a3}, группы (A4) --- лемма~\ref{ru:lem:qsem}
(переводы доказуемы во внутренней теории, выводы строятся по коду),
схем (A5) и (A7) --- предложение~\ref{ru:prop:a5a7}. $\sE$-общезначимость схем (A8) и (A9) ---
предложение~\ref{ru:prop:hce}(1). Аксиомы $\Pi$ в $\sB$ --- предложение~\ref{ru:prop:pilaws}. Всё это и сведено в
таблице~\ref{ru:tab:axioms}.

(2) \emph{Правила.} Всюду ниже, говоря о выводимости во внешней
теории, мы имеем в виду в $\sB$ выводимость замыкания из
определения~\ref{ru:def:valid}.

(MP) сохраняет общезначимость по леммам~\ref{ru:lem:closure}
и~\ref{ru:lem:closurebs}.

(PtoA): пусть $S$ решает задачу $\varphi$, то есть внешняя теория
доказывает $\varphi^*(S(\gn{*}))$ при всякой интерпретации $*$
(её замыкание $\overline{\varphi^*}(S(\gn{*}))$ в интерпретациях $\sB$). 
Припишем к её выводу  одну аксиому $\varphi^*(S(\gn{*}))\ato \exists x\,\varphi^*(x)$
исчисления предикатов и формулу $\exists x\,\varphi^*(x)$ --- получится вывод
формулы $(\quest\varphi)^*$.  В $\sB$ 
мы выведем формулу $\exists x\,\overline{\varphi^*}(x)$, из которой выводится 
искомая формула $\overline{\exists x\,\varphi^*(x)}$ стандартными средствами исчисления
предикатов. Итак, высказывание $\quest\varphi$
общезначимо.

(AtoP) в $\sI$, $\sB$: пусть высказывание $A$ общезначимо, то
есть внешняя теория, совпадающая здесь с внутренней, доказывает $A^*$
(в $\sB$ --- замыкание $\overline{A^*}$). Это ровно та формула,
выводы которой суть решения задачи $\bang A$ (определения~\ref{ru:def:tr},
\ref{ru:def:starB}). По коду интерпретации 
некоторый вывод $P(*)$ этой формулы находится перебором (замечание~\ref{ru:rem:whycomp}). 
Утверждение $\Prf(P(\gn{*}),A^*)$  (или $\Prf(P(\gn{*}),\overline{A^*})$) --- истинное
$\Sigma_0$-предложение, и внешняя теория его доказывает. Значит, $P$
решает задачу $\bang A$.

Правило Бернайса в $\sB$: пусть $S$ решает $\varphi\pto\psi$, где
$\alpha\notin FV(\varphi)$, то есть $\ZF$ доказывает замыкание
$\forall\vec X\,\forall X_\alpha\,\forall u\,\bigl(\varphi^*(u)\ato
\psi^*(U_{S(\gn{*})}(u))\bigr)$ по всем ящикам, в том числе по
$X_\alpha$. Перевод $(\varphi\pto\Pi\alpha\,\psi)^*(S(\gn{*}))$ --- это
$\forall u\,\bigl(\varphi^*(u)\ato\forall X_\alpha\,
\psi^*(U_{S(\gn{*})}(u))\bigr)$; переменная $X_\alpha$ в $\varphi^*$ не
встречается, и замыкание этой формулы получается из предыдущего чистой
логикой --- переносом квантора $\forall X_\alpha$ через посылку. Тот же
$S$ решает $\varphi\pto\Pi\alpha\,\psi$. Для высказываний так же: из
$\ZF\vdash\forall\vec X\,\forall X_\alpha\,(A^*\ato B^*)$ при
$X_\alpha$, не входящей в $A^*$, следует
$\ZF\vdash\forall\vec X\,(A^*\ato\forall X_\alpha\,B^*)$.

(3) \emph{Схема (A6).} 
Сначала заметим, что в семантике $\sB$ достаточно разобрать
высказывания $A$ без свободных задачных переменных: экземпляр схемы
для произвольного $A$ \emph{выводится} в $\HCB$ из экземпляра для
высказывания $\bar A:=\Pi\vec\alpha\,A$, где $\vec\alpha$ --- свободные
задачные переменные $A$, с помощью аксиом и правил, уже проверенных в
пунктах~(1) и~(2); поэтому из общезначимости экземпляра для $\bar A$
следует общезначимость экземпляра для $A$. 
Для этого достаточно доказать, что задачи
$$
\bang(\square A\ato A)\pto\bang A,  \qquad
\bang(\square \bar A\ato \bar A)\pto\bang \bar A
$$
выводимо эквивалентны в $\HCB$ (без (A6)).
По лемме~\ref{ru:lem:repl} достаточно установить эквивалентность их посылок
и заключений.

По аксиоме $(\Pi\bang)$
выводима эквивалентность заключений
$\bang A\peq\bang\bar A$, а из неё по правилу $(\quest\to)$ ---
$\square A\aeq\square\bar A$.
Для посылок та же аксиома $(\Pi\bang)$ даёт
$$
\bang(\square A\ato A)\peq\bang\Pi\vec\alpha(\square A\ato A)
$$
(свободные задачные переменные высказывания $\square A\ato A$ --- те же
$\vec\alpha$, поскольку у $\square A=\quest\bang A$ их нет).
По лемме~\ref{ru:lem:piderived}(2) для высказываний
($\vec\alpha$ в $\square A$ не входят) выводима эквивалентность
$$
\Pi\vec\alpha(\square A\ato A)\aeq(\square A\ato\Pi\vec\alpha\,A),
$$
а по лемме~\ref{ru:lem:repl} в её правой части $\square A$ можно заменить
на $\square\bar A$; получаем
$\Pi\vec\alpha(\square A\ato A)\aeq(\square\bar A\ato\bar A)$.
Наконец, ещё раз применим лемму~\ref{ru:lem:repl}, теперь под знаком
$\bang$, и получим $\bang\Pi\vec\alpha(\square A\ato A)\peq
\bang(\square\bar A\ato\bar A)$, а с ней и требуемую эквивалентность
посылок $\bang(\square A\ato A)\peq\bang(\square\bar A\ato\bar A)$.

Теперь, когда общезначимость остальных аксиом и
корректность правил установлены, схему Лёба можно не проверять
напрямую, а \emph{вывести} из них и леммы о неподвижной точке для
внутренней теории. Нам нужна лемма о неподвижной точке в эффективной
форме:
\begin{quote}
по замкнутой формуле $\sigma$ (арифметической --- для $\PA$; языка
$\ZF$ --- для $\ZF$) можно вычислить замкнутую формулу $\pi$ и вывод во
внутренней теории эквивалентности $\pi\aeq(\Prov(\pi)\ato\sigma)$.
\end{quote}

Вывод
разбивается на две части. Сначала --- вывод в $\HCglm$ формулы
\begin{equation}\tag{L$'$}
  \bang\bigl(P_0\aeq(\square P_0\ato A)\bigr)\ \pto\
  \bigl(\bang(\square A\ato A)\pto\bang A\bigr),
\end{equation}
где $P_0$ --- любая переменная сорта высказывание, не входящая в $A$:
это пример~\ref{ru:ex:proofs}(5) при $P:=P_0$, $S:=A$. Поскольку
(L$'$) выведена без (A6), по уже доказанному она общезначима --- в
$\sI$, а также в $\sB$, где $\HCglm$ содержится в $\HCB$. Пусть $A$ --- высказывание (в $\sB$ --- по доказанному, без
свободных задачных переменных), для которого нужно доказать
общезначимость задачи $\bang(\square A\ato A)\pto\bang A$. Нам нужно
построить алгоритм, который по интерпретации $*$ находит такое $x$,
для которого внешняя теория доказывает
$\bigl(\bang(\square A\ato A)\pto\bang A\bigr)^*(x)$. При этом мы можем
считать, что $*$ присваивает значения только переменным, входящим в
нашу формулу, в частности не присваивает значения переменной $P_0$.
Поскольку формула (L$'$) общезначима, у нас имеется алгоритм, который
доставляет нужное нам $x$, если ему дать на вход любые $y,*'$, для
которых во внешней теории выводимо
$\bigl(\bang(P_0\aeq(\square P_0\ato A))\bigr)^{*'}(y)$ и $*'$
продолжает $*$ на переменную $P_0$. Такие $y,*'$ и доставляются леммой
о неподвижной точке: положим $*'(P_0)=\pi$, для которого во
внутренней теории доказуема эквивалентность
$\pi\aeq(\Prov(\pi)\ato A^*)$; $y$ --- её вывод, найденный перебором
всех доказательств (формула $\Prf(y,\dots)$ --- истинное
$\Sigma_0$-предложение, и внешняя теория его доказывает). Именно здесь
нужна замкнутость формулы $A^*$: в $\sB$ решение задачи
$\bang(P_0\aeq(\square P_0\ato A))$ --- вывод замыкания эквивалентности
по свободным задачным переменным $A$, а лемма о неподвижной точке даёт
эквивалентность лишь для замкнутой $A^*$; сама же формула (L$'$)
общезначима при любом $A$.

(4) \emph{Семантика $\sE$.} Аксиомы $\HCE$ --- теоремы $\HCI$ и схемы
(A8), (A9). Всякая теорема $\HCI$ $\sI$-общезначима по доказанному, а
$\sI$-общезначимость влечёт $\sE$-общезначимость, поскольку внутренняя
теория $\PA$ содержится во внешней $\TA$. Схемы (A8) и (A9)
$\sE$-общезначимы по предложению~\ref{ru:prop:hce}(1), а правила (MP)
и (PtoA) сохраняют $\sE$-общезначимость по пункту~(2).
\end{proof}

\subsection{Что даёт квантор $\Pi$}\label{ru:sec:pigain}

Напомним мотивировку квантора (\S\ref{ru:sec:piatoms}): под знаком
$\quest$ конструктивное содержание задачи теряется, и высказывания
$\quest(\alpha\por\pneg\alpha)$, $\quest\rho$ общезначимы тривиально.
Теперь, когда есть исчисления и теорема о корректности, обещанное
можно доказать.

\begin{theorem}[три задачи с квантором $\Pi$, равносильные $\varnothing$]\label{ru:prop:piatoms}
В $\sB$ общезначимы эквивалентности
\[
  \Pi\alpha\,\alpha\peq\varnothing,\qquad
  \Pi\alpha(\alpha\por\pneg\alpha)\peq\varnothing,\qquad
  \Pi\alpha(\pneg\pneg\alpha\pto\alpha)\peq\varnothing
\]
и отрицания высказываний
\[
  \quest\Pi\alpha\,\alpha,\qquad
  \quest\Pi\alpha(\alpha\por\pneg\alpha),\qquad
  \quest\Pi\alpha(\pneg\pneg\alpha\pto\alpha).
\]
Сами эти три высказывания, следовательно, не общезначимы (в
предположении непротиворечивости $\ZF$).
\end{theorem}

\begin{proof}Выведем все шесть формул  в исчислении $\HCB$.
Импликации справа налево для первых трех формул
являются аксиомами $\Int$. Теперь выведем импликации в обратную сторону.

%
(1) 
Формула  $\Pi\alpha\,\alpha\pto\varnothing$ является аксиомой
$(\Pi\mathrm{E})$ при $\theta:=\varnothing$.

(2) 
Действительно, $(\Pi\mathrm{E})$ при
$\theta:=\varnothing$ даёт $\Pi\alpha\,\alpha\pto\varnothing$.
Аналогично, $(\Pi\mathrm{E})$ при
$\theta:=\pneg\varnothing$ даёт $\Pi\alpha\,\pneg\alpha\pto\pneg\pneg\varnothing$,
а $\pneg\pneg\varnothing\pto\varnothing$ выводима в $\Int$.
Поэтому выводимо $\Pi\alpha\,\pneg\alpha\pto\varnothing$.
Вместе с
$(\Pi\por)$ это даёт $\Pi\alpha(\alpha\por\pneg\alpha)\pto\varnothing$.

(3) 
Действительно,
для $H:=\Pi\alpha(\pneg\pneg\alpha\pto\alpha)$: $(\Pi\mathrm{E})$ при
$\theta:=\alpha\por\pneg\alpha$ 
даёт
$$
H\pto  (\pneg\pneg  (\alpha\por\pneg\alpha) \pto (\alpha\por\pneg\alpha)).
$$
Теорема $\Int$
$\pneg\pneg(\alpha\por\pneg\alpha)$ даёт $H\pto(\alpha\por\pneg\alpha)$. Далее
правило Бернайса ($\alpha$ в $H$ свободно не входит) даёт
$H\pto\Pi\alpha(\alpha\por\pneg\alpha)$,
а отсюда $H\pto\varnothing$ по предыдущему пункту. 

Для доказательства второго утверждения (о выводимости отрицаний) 
 мы к каждой из выведенных импликаций применяем ($\quest\to$) и затем (A4d).
 %
%
%
\end{proof}

Итак, $\quest\Pi\alpha(\alpha\por\pneg\alpha)$ опровержимо, тогда как
без квантора $\quest(\alpha\por\pneg\alpha)$ --- теорема (по (A4) оно
равносильно тавтологии $\quest\alpha\aor\aneg\quest\alpha$)
: сквозь квантор классический разбор случаев не проходит, и теперь и на уровне
высказываний становится видна разница между законом исключённого третьего
и его конструктивной версией.

\begin{example}\label{ru:ex:pider}
Пусть $\rho$ --- формула Роуза, то есть задача \eqref{ru:rose-formula} для
$\varphi=\pneg\beta\por\pneg\gamma$, где $\beta,\gamma$ --- задачные
переменные (\S\ref{ru:sec:modalfreeval} ниже). В $\sB$ общезначимо высказывание
$\quest\Pi\beta\Pi\gamma\rho$: решатель формулы Роуза постоянен и на
содержимое ящиков не смотрит (предложение~\ref{ru:cor:rose} ниже), так что
$\ZF$ доказывает $\forall X_\beta\forall X_\gamma\,\rho^*(r)$ для его
кода $r$, то есть $(\Pi\beta\Pi\gamma\rho)^*(r)$, а с ним и
$(\quest\Pi\beta\Pi\gamma\rho)^*$.
Без квантора высказывание $\quest\rho$ по (A4) равносильно тавтологии и
ничего не говорит; с квантором оно утверждает, что формула Роуза
решается одним алгоритмом при любом содержимом ящиков. Это первый
случай, когда открытие на задачной стороне порождает новое
общезначимое высказывание. Выводимо ли оно в $\HCB$, мы не знаем:
сама $\rho$ в $\HCI$ не выводима (\S\ref{ru:sec:indep}), а
консервативность $\HCB$ над $\HCI$ на чистых задачах открыта
(замечание~\ref{ru:rem:qfree}).
\end{example}

Второй вопрос --- насколько новы задачи и высказывания с квантором.
Во всех нетривиальных случаях, которые нам удалось разобрать, квантор
устраняется: задача с $\Pi$ оказывается равносильной задаче без
$\Pi$. Соберём эти случаи вместе.

\begin{proposition}[устранение квантора]\label{ru:prop:pienc}
Пусть $\varphi$ --- задача, в которую переменная $\beta$ свободно не
входит, а $\alpha,\beta,\gamma$ --- переменные. Следующие
эквивалентности $\sB$-общезначимы:
\begin{enumerate}[leftmargin=2.2em,itemsep=2pt,label=\textup{(\arabic*)}]
\item $\Pi\beta\,(\varphi\pto\beta)\ \peq\ \pneg\varphi$;
\item $\Pi\beta\,\bigl((\varphi\pto\beta)\pto\beta\bigr)\ \peq\ \varphi$;
\item $\Pi\beta\,\bigl(((\alpha\pto\beta)\pto\alpha)\pto\alpha\bigr)
  \ \peq\ (\pneg\pneg\alpha\pto\alpha)$;
\item $\Pi\alpha\,\bigl(((\alpha\pto\beta)\pto\alpha)\pto\alpha\bigr)
  \ \peq\ \beta$;
\item $\Pi\beta\,\bigl(((\alpha\pto\beta)\pto\beta)\pto(\alpha\por\beta)\bigr)
  \ \peq\ \alpha\por\pneg\alpha$;
\item $\Pi\alpha\,\bigl(((\alpha\pto\beta)\pto\beta)\pto(\alpha\por\beta)\bigr)
  \ \peq\ \beta$;
\item $\Pi\alpha\,\bigl((\alpha\pto\beta\por\gamma)\pto
  ((\alpha\pto\beta)\por(\alpha\pto\gamma))\bigr)
  \ \peq\ (\gamma\pto\beta)\por(\beta\pto\gamma)$;
\item $\Pi\alpha\,\bigl((\beta\pand\gamma\pto\alpha)\pto
  ((\beta\pto\alpha)\por(\gamma\pto\alpha))\bigr)
  \ \peq\ (\beta\pto\gamma)\por(\gamma\pto\beta)$.
\end{enumerate}
Все они, кроме пунктов~(4) и~(6), выводимы в $\HCB$; в пунктах~(4)
и~(6) выводима импликация справа налево. Вместе с каждой
эквивалентностью задач $\Pi\vec\alpha\,\chi\peq\psi$ выводима
(соответственно общезначима) и эквивалентность
высказываний $\quest\Pi\vec\alpha\,\chi\aeq\quest\psi$.
\end{proposition}

\begin{proof}
Последнее утверждение: из задачи $\Pi\vec\alpha\,\chi\peq\psi$ правило
$(\quest\to)$ даёт обе импликации между $\quest\Pi\vec\alpha\,\chi$ и
$\quest\psi$; семантический вариант --- лемма~\ref{ru:lem:qsem} и
теорема~\ref{ru:th:sound}.

Во всех задачах импликация справа налево доказывается одинаково. 
Пусть $\Pi\delta\,\lambda$ обозначает левую часть, а $\mu$ --- правую часть формулы.
Тогда импликация $\mu\pto\lambda$ выводима в Int, что проверяется непосредственно
в каждом из пунктов. А затем мы применяем правило Бернайса.

Импликация слева направо во всех пунктах тоже доказывается одинаково, за исключением 
пунктов (4) и (6). А именно, мы применяем аксиому $(\Pi\mathrm{E})$, подставляя в формулу $\lambda$
некоторую формулу 
$\theta$ вместо связанной переменной так, что для полученной формулы $\lambda(\theta)$ формула 
$\lambda(\theta)\pto\mu$ выводима в Int. Укажем, как выбирать формулу $\theta$ в каждом пункте:
(1) $\theta:=\varnothing$, 
(2) $\theta:=\varphi$, (3) $\theta:=\varnothing$,
(5) $\theta:=\pneg\alpha$, (7)  $\theta:=\beta\por\gamma$;
(8) $\theta:=\beta\pand\gamma$.

Пункты (4) и (6) разберём отдельно. 
\remove{
Направление справа налево --- правило Бернайса из
$\Int$-теорем $\beta\pto(((\alpha\pto\beta)\pto\alpha)\pto\alpha)$ и
$\beta\pto(((\alpha\pto\beta)\pto\beta)\pto(\alpha\por\beta))$.
Направление слева направо докажем семантически,} 
Начнём с (4). Пишем
$Y$ вместо $X_\beta$; рассуждаем в $\ZF$ при свободной $Y$. Для
множеств $X,Y$ положим $X{\to}Y:=\{q:\forall x\in X\ U_q(x){\downarrow}\in Y\}$
и 
$$
\Phi(X):
=\{p:\forall q\in X{\to}Y\ U_p(q){\downarrow}\in X\}.
$$
 Тогда
$e\in\val{\Pi\alpha(((\alpha\pto\beta)\pto\alpha)\pto\alpha)}$ означает:
для всех $X$ и всех $p\in\Phi(X)$ выполнено $U_e(p){\downarrow}\in X$.
Зафиксируем код $I$ тождественной программы ($U_I(q)=q$) и код $a$
программы, тождественно равной нулю. Решатель такой:
\[
  f(e):=U_n(a),\qquad\text{где } n:=U_e(I).
\]
Покажем, что для всякого $Y$ и всякого такого $e$ значение $f(e)$
определено и лежит в $Y$.

\emph{Шаг 1: $n$ определено.} При $X:=\omega$ имеем $I\in\Phi(\omega)$,
поскольку $U_I(q)=q\in\omega$ при любом $q$; значит $n=U_e(I)$
определено.

\emph{Шаг 2: $U_n$ отображает $\omega\setminus\{n\}$ в $Y$.} Допустим
противное: есть $x\neq n$, для которого $U_n(x)$ не определено или не
лежит в $Y$. Положим $X:=\omega\setminus\{n\}$. Тогда $I\in\Phi(X)$:
если $q\in X{\to}Y$, то $q\neq n$ (из-за точки $x\in X$ число $n$ не
принадлежит $X{\to}Y$), то есть $U_I(q)=q\in X$. По определению $e$
отсюда $n=U_e(I)\in X=\omega\setminus\{n\}$ --- противоречие.

\emph{Шаг 3: $f(e)\in Y$.} Если $a\neq n$, то по шагу~2 значение
$U_n(a)$ определено и лежит в $Y$. Если $a=n$, то $U_n$ тождественно
равна нулю, по шагу~2 число $0=U_n(x)$ при любом $x\neq n$ лежит в $Y$,
и $f(e)=0\in Y$.

Всё это --- рассуждение в $\ZF$ при свободной $Y$. 

Для (6) рассуждение похожее. Будем использовать обозначения из рассуждения
для случая (4). Положим $k:=\lambda q.\,U_q(q)$, то есть, $U_k(q)=U_q(q)$. Число $k$ принадлежит 
$\val{(\omega\pto Y)\pto Y}$, 
так что $U_e(k)$ определено и лежит в
$\omega\sqcup Y$. Если это $\langle1,y\rangle$, то $y\in Y$. Если это
$\langle0,n\rangle$, то, как и
в шаге~2, $U_n$
отображает $\omega\setminus\{n\}$ в $Y$. В самом деле, иначе
$k$ является решением задачи $(\omega\setminus\{n\} \pto Y)\pto Y$
 (если $q$ --- решение задачи $\omega\setminus\{n\}\pto Y$, то $q\ne n$ по предположению, следовательно, 
 $U_q(q)\in Y$, поскольку $U_q(x)\in Y$ для всех $x\ne n$). Следовательно, $U_n(a)\in Y$, что доказывается так же, как и в шаге~3 выше.
 
 В первом случае решатель выдает то $y$, для которого $U_e(k)=c(1,y)$, а во втором --- $U_n(a)$.
\end{proof}

Эквивалентность~(2) --- известное непредикативное представление
задачи квантором второго порядка: <<решить $\varphi$>> равносильно
<<уметь решить всякую задачу $\beta$, сводящуюся к $\varphi$>>.

\begin{remark}[о пунктах (4) и (6)]\label{ru:rem:peirce}
Выводимы ли в $\HCB$ импликации слева направо в пунктах~(4) и~(6)
предложения~\ref{ru:prop:pienc}, мы не знаем. Выводимо их ослабление с
$\pneg\pneg\beta$ в заключении: при гипотезе $\pneg\beta$ импликация
$\alpha\pto\beta$ равносильна $\pneg\alpha$, закон Пирса превращается в
$\pneg\pneg\alpha\pto\alpha$, правило Бернайса даёт
$\Pi\alpha(\pneg\pneg\alpha\pto\alpha)$, а эта задача равносильна
$\varnothing$ (теорема~\ref{ru:prop:piatoms}). 

Невыводимость же импликаций справа налево, как есть (с формулой $\beta$), нам удалось доказать только 
в исчислении HCI минус $(\Pi\por)$.  
В самом деле возьмём двухточечную шкалу Крипке
$w_0<w_1$, в которой $\beta$ истинна только в $w_1$, задачи --- конусы,
$\Pi$ --- пересечение по всем конусам, $\bang A$ --- верх, а
$\quest\varphi$ истинно, когда $\varphi$ истинна в $w_1$. Все аксиомы и
правила $\HCB$, кроме $(\Pi\por)$, в ней верны (ср.\
теорему~\ref{ru:prop:modalfree}), и в $w_0$ верна задача
$\Pi\alpha(((\alpha\pto\beta)\pto\alpha)\pto\alpha)$: при любом конусе
$X$ посылка $X\pto\beta$ истинна в $w_0$ ровно когда $X$ не содержит
$w_0$, и тогда $(X\pto\beta)\pto X$ в $w_0$ ложна. А $\beta$ в $w_0$
ложна. Аксиому $(\Pi\por)$ эта шкала не выполняет: в $w_1$ посылка
$\Pi\alpha(\alpha\por\pneg\alpha)$ истинна, а заключение
$\Pi\alpha\,\alpha\por\Pi\alpha\,\pneg\alpha$ ложно (оба дизъюнкта пусты:
среди конусов есть пустой и полный), так что этот экземпляр
$(\Pi\por)$ ложен и в $w_1$, и в $w_0$. Значит, возможный вывод обязан
использовать $(\Pi\por)$ существенно. В семантическом
доказательстве решающим было инстанцирование квантора ящиком
$\omega\setminus\{n\}$, зависящим от ответа самого решателя; значением
формулы такое множество не бывает. Пункты~(4) и~(6) --- кандидаты на
примеры неполноты $\HCB$ (следствие~\ref{ru:cor:incomplete}).
(Конец замечания.)
\end{remark}

Примера \emph{неустранимости} квантора --- задачи с $\Pi$, не
равносильной никакой задаче без $\Pi$, --- нам найти не удалось. Это,
конечно, не означает, что квантор не нужен: устранение в
предложении~\ref{ru:prop:pienc} само требует квантора, а в пунктах~(4)
и~(6) ещё и выходит за пределы исчисления.

\begin{remark}[бескванторная часть]\label{ru:rem:qfree}
Аксиоматика $\HCB$ устроена как $\HCI$ плюс кванторный слой:
выбрасывая аксиомы $\Pi$ и правило Бернайса, получаем в точности $\HCI$.
Цель
расширения --- \emph{увеличить запас доказуемых формул}. В расширенном
языке рост заведомый: выводимы новые содержательные высказывания вроде
$\quest\Pi\beta\Pi\gamma\rho$ (пример~\ref{ru:ex:pider}), невыразимые без
квантора. Рост есть и на \emph{бескванторных} формулах:
$\square\quest\alpha\ato\square\bot$ выводимо в $\HCB$ и не выводимо в
$\HCI$ (замечание~\ref{ru:rem:boxqatom}). Но этот рост --- не новое
знание о задачах, а следствие того, что в $\sB$ знак $\bang$ связывает
задачные переменные: в $\HCB$ формула означает
$\square\Pi\alpha\,\quest\alpha\ato\square\bot$ и не является чистой. Поэтому
вопрос о консервативности правильно ставить для \emph{чистых}
бескванторных формул (замечание~\ref{ru:rem:clean}): совпадает ли
множество чистых бескванторных теорем $\HCB$ с множеством теорем
$\HCI$? Для высказываний ответ известен: в нормальной форме чистого
$\Pi$-свободного высказывания атомы $\quest\alpha$ стоят вне знаков
$\square$, и при $\Sigma_1$-корректности $\ZF$ его выводимость в $\HCB$
равносильна выводимости в $\HCI$ (предложение~\ref{ru:prop:nf},
теорема~\ref{ru:th:assval}(1),(3)). Для чистых задач вопрос открыт
(\S\ref{ru:sec:todo}).
\end{remark}

\remove{
\begin{remark}[лёбовская механика не задета]\label{ru:rem:lobpi}
Вывод задачи (L), пример~\ref{ru:ex:proofs} о $\square$-версиях, вывод
схемы (A6) из эффективной леммы о неподвижной точке и
теорема~\ref{ru:th:lob} переносятся в $\Pi$-расширение дословно: в них не
участвуют ни ящики, ни $\Pi$, а условия выводимости и диагональная
лемма верны для обеих внутренних теорий, $\PA$ и $\ZF$.
\end{remark}
}

\subsection{Нормальная форма и полнота на стороне высказываний}
\label{ru:sec:asscompl}

Здесь мы опишем два семейства $\Pi$-свободных высказываний:
\emph{выводимые} в каждом из исчислений и \emph{общезначимые} в
каждой из семантик. Оба описания получаются через одну и ту же
\emph{нормальную форму}, в которой модальность $\quest$ устранена до
сочетаний $\quest\alpha$ и $\square B$, и оба отсылают к известным
объектам: с точностью до нормальной формы это логика доказуемости
$\GL$ Соловея (для $\HCI$ и для внутренних семантик $\sI$, $\sB$) и её
внешний вариант $\GLS$ (для $\HCE$ и для $\sE$). Доказательство идёт по
кругу: из выводимости в нашем исчислении --- общезначимость (теорема о
корректности), из общезначимости --- выводимость в $\GL$ или $\GLS$
(теорема Соловея), а из выводимости в $\GL$ или $\GLS$ --- выводимость
в нашем исчислении (прямая проверка). В итоге на стороне высказываний
исчисления $\HCI$ и $\HCE$ полны.

Нормальная форма строится проталкиванием $\quest$ внутрь по схемам (A4):
к переменным и к задачам вида $\bang B$.

\begin{definition}[перевод в нормальную форму]\label{ru:def:nf}
Каждой $\Pi$-свободной формуле $\chi$ сопоставим её нормальную форму
$\nf{\chi}$. Для высказываний: переменные и константа $\bot$ переходят
в себя, классические связки сохраняются
($\nf{(A\aand B)}=\nf{A}\aand\nf{B}$ и так далее для $\aor$, $\ato$),
а $\nf{(\quest\varphi)}:=\varphi^\quest$, где вспомогательный перевод
$\varphi^\quest$ задач в высказывания задаётся правилами
\[
\begin{aligned}
  \alpha^\quest&:=\quest\alpha, &
  \varnothing^\quest&:=\bot, &
  (\bang B)^\quest&:=\square\,\nf{B},\\
  (\varphi\pto\psi)^\quest&:=\varphi^\quest\ato\psi^\quest, &
  (\varphi\pand\psi)^\quest&:=\varphi^\quest\aand\psi^\quest, &
  (\varphi\por\psi)^\quest&:=\varphi^\quest\aor\psi^\quest.
\end{aligned}
\]
В дальнейшем мы будем назвать $\varphi^\quest$ \emph{булевой свёрткой} задачи $\varphi$.
Для задач: $\nf{\alpha}=\alpha$, $\nf{\varnothing}=\varnothing$,
задачные связки сохраняются, а $\nf{(\bang B)}:=\bang\,\nf{B}$ (при этом
высказывание $B$ нормализуется по правилам выше).
\end{definition}

Нетрудно видеть, что в нормальной форме всякого $\Pi$-свободного
высказывания $A$ модальность встречается только в сочетаниях
$\quest\alpha$ ($\alpha$ --- переменная сорта задача) и
$\quest\bang B=\square B$.

\begin{proposition}[нормальная форма]\label{ru:prop:nf}
Для всякой $\Pi$-свободной формулы в $\HCglm$ (а значит, и в $\HCgl$,
$\HCB$, $\HCgls$) выводима эквивалентность
\[
  A\aeq\nf{A}\quad\text{для высказываний},\qquad
  \varphi\peq\nf{\varphi}\quad\text{для задач}.
\]
Следовательно, при всякой интерпретации внешняя теория каждой из трёх
семантик доказывает $A^*\aeq\nf{A}^*$ (в $\sB$ --- замыкание этой
эквивалентности), и высказывание $A$ общезначимо в данной семантике
тогда и только тогда, когда в ней общезначимо $\nf{A}$.
\end{proposition}

\begin{proof}
Совместная индукция по построению формулы. Для переменной и $\bot$
очевидно, связки проходят по лемме~\ref{ru:lem:repl}. Для $A=\quest\varphi$
надо проверить $\provglm\quest\varphi\aeq\varphi^\quest$ индукцией по
$\varphi$: при $\varphi=\alpha$ это тождество
($\alpha^\quest=\quest\alpha$); при $\varphi=\varnothing$ --- схема
(A4d); при $\varphi=\bang B$ требуется
$\provglm\square B\aeq\square\,\nf{B}$, что получается из внешнего
предположения индукции $\provglm B\aeq\nf{B}$ леммой~\ref{ru:lem:repl};
случаи $\pto$, $\pand$, $\por$ --- это в точности три эквивалентности
(A4a), (A4b), (A4c) вместе с предположением индукции и
леммой~\ref{ru:lem:repl}. Здесь и работает булевость $\quest$. Для задач:
переменная и $\varnothing$ очевидны, связки --- по
лемме~\ref{ru:lem:repl}, а $\nf{(\bang B)}=\bang\,\nf{B}$ получается из
$\provglm B\aeq\nf{B}$ той же леммой. Использованы только (A1), (A2),
(A3), группа (A4) и три правила.

Второе утверждение: по теореме~\ref{ru:th:sound} выведенная
эквивалентность $\sI$-общезначима, а потому и $\sE$-общезначима
($\PA\subseteq\TA$), и $\sB$-общезначима, поскольку
$\HCglm\subseteq\HCB$. Множество теорем внешней теории замкнуто по
(MP), поэтому $A^*$ выводима во внешней теории тогда и только тогда,
когда выводима $\nf{A}^*$.
\end{proof}

Оговорка о $\Pi$-свободности существенна: схемы (A4) проталкивают
$\quest$ сквозь связки $\pto,\pand,\por$ и константу $\varnothing$, но не
сквозь квантор $\Pi$, и в полном языке семантики $\sB$ нормализация
обрывается на высказываниях вида $\quest\Pi\alpha\,\varphi$. Аналога
схем (A4) для $\Pi$ нет и быть не может: высказывание
$\quest\Pi\alpha(\alpha\por\pneg\alpha)$ \emph{опровержимо} ---
общезначимо его отрицание (теорема~\ref{ru:prop:piatoms}), тогда как без
квантора $\quest(\alpha\por\pneg\alpha)\aeq\quest\alpha\aor\aneg\quest\alpha$
--- общезначимая классическая тавтология: сквозь $\Pi$ классический
разбор случаев не проходит. Законы квантора $\Pi$ позволяют
упростить и такие высказывания; что при этом получается, сказано в
\S\ref{ru:sec:pinf}.

\begin{corollary}[фрагмент высказываний]\label{ru:cor:assfrag}
В нормальной форме высказывания строятся классическими связками и
модальностью $\square$ из переменных $A_i$ и из формул $\quest\alpha_j$. При этом
формулы $\quest\alpha_j$ ведут себя как \emph{ещё одни независимые переменные
сорта высказывание}: взяв $\alpha^*:=\sigma$ для произвольного замкнутого
$\sigma$, получаем $(\quest\alpha)^*=\exists x\,\sigma$, что логически эквивалентно
$\sigma$. Поэтому вопрос о фрагменте высказываний целиком сводится к вопросу о
$\square$-фрагменте, которым мы дальше и займёмся.
\end{corollary}

\emph{$\square$-формулой} называется высказывание, построенное
классическими связками $\bot,\aor,\aand,\ato$ из \emph{атомов} двух
видов --- переменных сорта высказывание и высказываний вида
$\quest\beta$, где $\beta$ --- переменная сорта задача, --- и модальности
$\square=\quest\bang$, применяемой к $\square$-формулам; иных вхождений
$\bang$ и $\quest$ нет. При выводимости в $\GL$ и $\GLS$ атомы обоих видов
играют роль пропозициональных переменных модального языка.
В нормальной форме $\Pi$-свободное высказывание построено классическими связками из
переменных $A_i$, высказываний $\quest\alpha_j$ и высказываний $\square B$, где
$B$ --- снова высказывания в нормальной форме. Иными словами,
\emph{нормальная форма $\Pi$-свободного высказывания --- $\square$-формула}, и отдельный перевод в модальный язык не нужен. $\square$-формулой является и $\quest$-свёртка $\varphi^\quest$ всякой $\Pi$-свободной
задачи $\varphi$ (определение~\ref{ru:def:nf}).

Напомним определения логик $\GL$ и $\GLS$
(см.~\cite[гл.~1]{ru:Boolos93}); обе формулируются в языке $\square$-формул. 
\begin{itemize}[leftmargin=1.6em,itemsep=1pt]
\item $\GL$: аксиомы --- все классические тавтологии (в
  модальном языке), схема K $\square(M\ato N)\ato(\square M\ato\square N)$,
  схема $\square M\ato\square\square M$ (по традиции называемая
  схемой~4) и схема Лёба $\square(\square M\ato M)\ato\square M$;
  правила --- (MP) и \emph{правило необходимости} $M/\square M$.
  (Схема~4 выводима из остальных аксиом, но мы оставляем её в списке.)
\item $\GLS$: множество формул, выводимых по правилу (MP) из теорем
  $\GL$ и всех экземпляров схемы рефлексии $\square M\ato M$.
  Правила необходимости в $\GLS$ нет, и множество её теорем
  относительно него не замкнуто: $\GLS\vdash\square\bot\ato\bot$, но
  $\GLS\nvdash\square(\square\bot\ato\bot)$, поскольку эта формула
  $\GL$-эквивалентна $\square\bot$.
\end{itemize}
Теоремы Соловея говорят, что $\GL$ --- это в точности $\square$-формулы, все арифметические реализации которых доказуемы в $\PA$, а
$\GLS$ --- в точности формулы, все реализации которых истинны; в наших
терминах это общезначимость в $\sI$ и в $\sE$.

Для семантик $\sI$ и $\sE$ мы воспользуемся теоремами
Соловея~\cite{ru:Solovay76} для $\PA$; для $\sB$ нужен их
вариант для $\ZF$, и здесь требуется второе допущение
замечания~\ref{ru:rem:assumptions} --- $\Sigma_1$-корректность $\ZF$.
Обойтись одной
непротиворечивостью нельзя уже на замкнутых формулах:
необщезначимость формулы $\square \bot$ означает, что
$\ZF$ не доказывает $\Prov_{\ZF}(\bot)$.
Непротиворечивость не гарантирует, а $\Sigma_1$-корректность гарантирует это.
Для упрощения формулировки воспользуемся тем, что
(1) язык $\GL$ является частью нашего языка и (2) для
семантик $\sI$ и $\sE$ наше определение общезначимости 
$\square$-формул есть определение общезначимости по Соловею: интерпретациями
атомов обоих видов служат произвольные замкнутые арифметические
формулы. В семантике $\sB$ так читаются только атомы первого
вида (им отвечают произвольные замкнутые формулы языка $\ZF$); атом
$\quest\beta$ там --- предложение $\exists y\,(y\in X_\beta)$ о ящике, а
под знаком $\square$ он к тому же читается замыканием
(замечание~\ref{ru:rem:bangbinds}). Поэтому для $\sB$ лемма
формулируется для $\square$-формул без атомов второго вида.

\begin{lemma}[теорема  Соловея]\label{ru:lem:solzf}
\begin{enumerate}[leftmargin=2.2em,itemsep=1pt,label=\textup{(\arabic*)}]
\item Для всякой $\square$-формулы $M$: $\GL\vdash M$ тогда и только
  тогда, когда $M$ общезначима в $\sI$; $\GLS\vdash M$ тогда и только
  тогда, когда $M$ общезначима в $\sE$.
\item Пусть $\ZF$ $\Sigma_1$-корректна. Для всякой $\square$-формулы
  $M$, не содержащей атомов вида $\quest\beta$, следующие утверждения
  равносильны: $\GL\vdash M$ и $M$ общезначима в $\sB$.
\end{enumerate}
\end{lemma}

\begin{remark}
Пункт (1) --- это первая и вторая теоремы Соловея для $\PA$, пункт
(2) --- теорема Соловея для $\ZF$. Метатеоретическое допущение нужно
для направления <<общезначима $\Rightarrow$ выводима>>: в конструкции
Соловея по невыводимой в $\GL$ формуле строится реализация, невыводимость
которой в $\ZF$ устанавливается в предположении
$1$-непротиворечивости, то есть $\Sigma_1$-корректности. В пункте~(2)
интерпретации переменных-высказываний --- произвольные замкнутые
формулы языка $\ZF$, а теорема Соловея доказана для арифметических
реализаций; это не мешает: направление <<$\GL\vdash M\Rightarrow
\ZF\vdash M^*$>> верно при любых подставленных предложениях (условия
выводимости и формализованная теорема Лёба от них не зависят), а
обратное достаточно иметь для арифметических реализаций, которые суть
частный случай. В $\sB$ перевод формулы без атомов
$\quest\beta$ не содержит множественных переменных, и оба замыкания в
определении~\ref{ru:def:valid} пусты. Подробности о конструкции Соловея
--- в~\cite{ru:Solovay76} и~\cite[гл.~9]{ru:Boolos93}.
\end{remark}

\begin{theorem}[выводимые и общезначимые $\square$-формулы --- это $\GL$ и $\GLS$]
\label{ru:th:assval}
Пусть $M$ --- $\square$-формула.
\begin{enumerate}[leftmargin=2.2em,itemsep=1pt,label=\textup{(\arabic*)}]
\item Следующие три утверждения равносильны: $\provgl M$;\
  $M$ общезначима в $\sI$;\ $\GL\vdash M$.
\item Следующие три утверждения равносильны: $\provgls M$;\
  $M$ общезначима в $\sE$;\ $\GLS\vdash M$.
\item Пусть атомы вида $\quest\alpha$ входят в $M$ только вне знаков
  $\square$ (под знаком $\square$ в семантике $\sB$ задачная переменная
  была бы связана и читалась бы замыканием ---
  замечание~\ref{ru:rem:bangbinds}). Тогда из $\provglb M$ следует
  общезначимость $M$ в $\sB$, из $\GL\vdash M$ следует $\provglb M$, а
  в предположении $\Sigma_1$-корректности $\ZF$ из общезначимости в
  $\sB$ следует $\GL\vdash M$, так что все три утверждения равносильны.
\end{enumerate}
\end{theorem}

\begin{proof}
Каждый пункт доказывается по кругу: из выводимости в нашем исчислении
--- общезначимость (теорема о корректности), из общезначимости ---
выводимость в $\GL$ или $\GLS$ (теорема Соловея), из выводимости в
$\GL$ или $\GLS$ --- выводимость в нашем исчислении (прямая проверка).
Последний шаг общий для всех пунктов, с него и начнём.

\emph{Из $\GL\vdash M$ следует $\provgl M$.} Всякий $\GL$-вывод,
прочитанный на наших атомах (переменных-высказываниях и формулах
$\quest\alpha$), превращается в вывод в $\HCgl$: схема K --- это
пункт~(1) примера~\ref{ru:ex:proofs}, схема~4 --- пункт~(2), схема Лёба
--- пункт~(3); классические тавтологии дают схемы (A1); правило
необходимости допустимо, поскольку из $\provgl A$ по (AtoP) следует
$\provgl\bang A$, а по (PtoA) --- $\provgl\quest\bang A$, то есть
$\provgl\square A$. \emph{Из $\GLS\vdash M$ следует $\provgls M$:}
$\GLS$ --- это теоремы $\GL$ вместе со всеми частными случаями
рефлексии $\square N\ato N$, замкнутые относительно (MP); теоремы $\GL$
выводимы в $\HCgl\subseteq\HCgls$, рефлексия --- схема (A8), а правило
необходимости здесь не используется. \emph{Из $\GL\vdash M$ следует
$\provglb M$}, поскольку $\HCgl\subseteq\HCB$.

\emph{Пункты (1) и (2).} Из $\provgl M$ следует $\sI$-общезначимость
$M$, а из $\provgls M$ --- $\sE$-общезначимость (теорема~\ref{ru:th:sound}).
Из $\sI$-общезначимости следует $\GL\vdash M$, а из
$\sE$-общезначимости --- $\GLS\vdash M$ (лемма~\ref{ru:lem:solzf}(1)).
Круг замкнут.

\emph{Пункт (3).} Из $\provglb M$ следует $\sB$-общезначимость по
теореме~\ref{ru:th:sound}. Остаётся вывести $\GL\vdash M$ из
$\sB$-общезначимости $M$ в предположении $\Sigma_1$-корректности.
Формула $M$ --- булева комбинация переменных $A_i$, атомов
$\quest\alpha_1,\dots,\quest\alpha_m$ и высказываний $\square B$, причём
по условию атомы $\quest\alpha_j$ стоят \emph{вне} знаков $\square$, а
аргументы $B$ задачных переменных не содержат. Зафиксируем
интерпретацию переменных-высказываний. Перевод $M^*$ --- булева
комбинация предложений $A_i^*$, атомов $\exists y\,(y\in X_j)$ и
предложений $\Prov(B^*)$ (формулы $B^*$ замкнуты). Общезначимость $M$
означает $\ZF\vdash\forall\vec X\,M^*$.

Покажем, что это условие равносильно следующему: 
\[
  \text{($\dagger$)}\qquad \ZF\vdash M^*[\vec\sigma]\ \ \text{для каждого
  набора}\ \vec\sigma\in\{0{=}0,\ \arbot \}^m,
\]
где $M^*[\vec\sigma]$ получается из $M^*$ заменой каждого атома
$\exists y\,(y\in X_j)$ на предложение $\sigma_j$. Из общезначимости
следует ($\dagger$): подставим вместо $X_j$ множество $\varnothing$
при $\sigma_j=(\arbot )$ и $\omega$ при $\sigma_j=(0{=}0)$; после
подстановки атом
$\exists y\,(y\in X_j)$ $\ZF$-эквивалентен $\sigma_j$, и замена
эквивалентных в булевом остове даёт $M^*[\vec\sigma]$. Обратно, пусть
выполнено ($\dagger$). Рассуждаем в $\ZF$: по классической логике
$\bigwedge_j\bigl(\exists y\,(y\in X_j)\aor\aneg\exists y\,(y\in X_j)\bigr)$,
что даёт разбор $2^m$ случаев; в каждом случае каждый атом
$\ZF$-эквивалентен своему $\sigma_j$, замена эквивалентных в булевом
остове даёт $M^*[\vec\sigma]$, а она доказуема по ($\dagger$). Итак,
$\ZF\vdash M^*$, а с ним и замыкание.

Осталось связать ($\dagger$) с $\GL$. Пусть $M[\vec\sigma]$ ---
результат замены в $M$ каждого атома $\quest\alpha_j$ на $\aneg\bot$ при
$\sigma_j=(0{=}0)$ и на $\bot$ при $\sigma_j=(\arbot )$. Это 
$\square$-формула без атомов $\quest\beta$, и $M^*[\vec\sigma]$ $\ZF$-эквивалентна
её переводу. По лемме~\ref{ru:lem:solzf}(2) условие ($\dagger$)
равносильно тому, что $\GL\vdash M[\vec\sigma]$ при всех $\vec\sigma$.
А это равносильно $\GL\vdash M$: атомы $\quest\alpha_j$ стоят вне
$\square$, поэтому формула
$M\aeq\bigvee_{\vec\sigma}\bigl(\bigwedge_j(\quest\alpha_j\aeq\sigma_j)
\aand M[\vec\sigma]\bigr)$ --- классическая тавтология, 
и из всех $\GL\vdash M[\vec\sigma]$ следует $\GL\vdash M$; а
обратно $M[\vec\sigma]$ --- подстановочные экземпляры $M$.
\end{proof}

\begin{remark}[атомы $\quest\alpha$ под знаком $\square$]\label{ru:rem:boxqatom}
Оговорка пункта~(3) существенна: если атомы $\quest\alpha$ стоят под
знаком $\square$, общезначимость в $\sB$ выводимости в $\GL$ не
влечёт. Пример: $\square\quest\alpha\ato\square\bot$
общезначимо в $\sB$, поскольку $\ZF$ доказывает импликацию
$(\square\quest\alpha)^*\ato(\square\bot)^*$: аргумент внутреннего
$\square$ читается замыканием $\forall X\,\exists y\,(y\in X)$,
которое опровержимо, и из его
вывода вместе с фиксированным опровержением получается вывод $\bot$. В
$\GL$ же формула $\square p\ato\square\bot$ не выводима, и в $\sI$, $\sE$
высказывание $\square\quest\alpha\ato\square\bot$ не общезначимо: при
$\alpha^*:=(0{=}0)$ оно превращается в формулу, выводимо эквивалентную
невыводимой формуле $\square\bot$. В $\HCB$ же это высказывание
выводимо, и вывод показывает, как исчисление учитывает скрытый
квантор: высказывательная схема $(\Pi\mathrm{E})$ при
$\theta:=\varnothing$ даёт $\Pi\alpha\,\quest\alpha\ato\quest\varnothing$, то
есть $\Pi\alpha\,\quest\alpha\ato\bot$ по (A4d); правило $(\bang\to)$
--- задачу $\bang\Pi\alpha\,\quest\alpha\pto\bang\bot$; аксиома
$(\Pi\bang)$ --- $\bang\quest\alpha\pto\bang\Pi\alpha\,\quest\alpha$; итого
выводима задача $\bang\quest\alpha\pto\bang\bot$, а правило $(\quest\to)$
даёт $\square\quest\alpha\ato\square\bot$. В $\HCI$ эта $\Pi$-свободная
формула не выводима: по пункту~(1) теоремы она была бы теоремой
$\GL$. Формула, конечно, не чистая (замечание~\ref{ru:rem:clean}):
переменная $\alpha$ стоит под знаком $\bang$; к вопросу о
консервативности $\Pi$-расширения мы вернёмся в
замечании~\ref{ru:rem:qfree}.
\end{remark}

\begin{corollary}[полнота на стороне высказываний]\label{ru:cor:solovay}
Для всякого $\Pi$-свободного высказывания $A$  следующие утверждения равносильны:
\begin{itemize}
\item $\provgl A$,
\item $A$ $\sI$-общезначимо, 
\item $\GL\vdash\nf{A}$.
\end{itemize}
Равносильны также следующие утверждения:
\begin{itemize}
\item 
$\provgls A$,
\item  $A$
$\sE$-общезначимо, 
\item $\GLS\vdash\nf{A}$.
\end{itemize}
В частности, $\square$-фрагмент $\HCgl$ есть в точности $\GL$, а
$\square$-фрагмент $\HCgls$ --- в точности $\GLS$.
\end{corollary}

\begin{proof}
По предложению~\ref{ru:prop:nf} выводимость и общезначимость $A$
равносильны выводимости и общезначимости $\square$-формулы $\nf{A}$;
остаётся теорема~\ref{ru:th:assval}.
\end{proof}

Итак, на стороне $\Pi$-свободных высказываний наша семантика не даёт
ничего нового по сравнению с логикой доказуемости: двусортность и
модальности не добавляют к $\GL$ и $\GLS$ ни одной новой схемы. Вся
новизна сосредоточена на стороне задач, а для $\sB$ --- ещё и в
высказываниях с квантором $\Pi$ (\S\ref{ru:sec:pigain}, \S\ref{ru:sec:pinf}).

\subsection{Нормальная форма в $\Pi$-расширении}\label{ru:sec:pinf}

В $\Pi$-свободном языке нормальная форма
(предложение~\ref{ru:prop:nf}) строится из переменных $A_i$ и
высказываний $\quest\alpha$ и $\square B$. В языке с квантором $\Pi$
приведение по (A4) обрывается на высказываниях
$\quest\Pi\vec\alpha\,\chi$ (замечание после предложения~\ref{ru:prop:nf}), но
законы $\Pi$ позволяют упростить и их. \emph{Атомом} мы называем
формулу, на которой такие упрощения заканчиваются. Это понятие
синтаксическое: оно определяется выписанным ниже набором
эквивалентностей, и атом может оказаться равносилен более простой
формуле. В этом подразделе рассматриваются формулы, в которых квантор
$\Pi$ применяется только к задачам; высказывания вида $\Pi\alpha\,A$ мы
здесь не упрощаем.

\begin{definition}[$\Pi$-атом]\label{ru:def:piatom}
\emph{$\Pi$-атомом} называем высказывание
$\quest\Pi\alpha_1\ldots\Pi\alpha_k\,\chi$, $k\geqslant1$, в котором каждая
переменная $\alpha_i$ входит в задачу $\chi$ свободно, а задача $\chi$ --- либо
переменная, либо задача вида $\bang A$, либо импликация
$\varphi\pto\psi$, в посылку $\varphi$ которой свободно входит каждая
$\alpha_i$.
\end{definition}

\begin{proposition}[нормальная форма в $\Pi$-расширении]\label{ru:prop:pinf}
Всякая формула, в которой квантор $\Pi$ применяется только к
задачам, эквивалентна в $\HCB$ формуле, в
которой $\quest$ встречается только в сочетаниях $\quest\alpha$ и
$\square B$ и в $\Pi$-атомах.
\end{proposition}

\begin{proof}
Как в предложении~\ref{ru:prop:nf}, но к эквивалентностям (A4a)--(A4d)
добавляются эквивалентности $(\Pi\pand)$, $(\Pi\por)$ (аксиома и
лемма~\ref{ru:lem:piderived}(4)) и эквивалентности
леммы~\ref{ru:lem:piderived}; все замены законны по лемме~\ref{ru:lem:repl}. Пусть
$\quest\Pi\vec\alpha\,\chi$ --- вхождение, не являющееся $\Pi$-атомом.
Квантор по переменной, не входящей в $\chi$ свободно, снимаем
(лемма~\ref{ru:lem:piderived}(1)). Если $\chi$ --- конъюнкция или дизъюнкция,
вносим кванторы в члены по $(\Pi\pand)$, $(\Pi\por)$, и $\quest$
проходит связку по (A4b) или (A4c). Если $\chi=\varphi\pto\psi$ и
некоторая $\alpha_i$ не входит в посылку, переставляем её квантор
внутрь блока (лемма~\ref{ru:lem:piderived}(3)) и вносим его в заключение
(лемма~\ref{ru:lem:piderived}(2)); если кванторов над импликацией не осталось,
$\quest$ проходит её по (A4a). Процесс обрывается: каждый шаг, кроме
перестановки кванторов, уменьшает в лексикографическом порядке тройку
из числа задачных связок, стоящих под знаками $\quest$, суммарного
числа задачных связок в областях действия кванторов $\Pi$, стоящих под
знаками $\quest$, и числа таких кванторов; а перестановка делается лишь
непосредственно перед шагом, который тройку уменьшает.
\end{proof}

Атомы $\Pi$-свободной нормальной формы, переменные $A_i$ и
высказывания $\quest\alpha$, неупрощаемы и в семантическом смысле: им
можно независимо приписать любые истинностные значения
(следствие~\ref{ru:cor:assfrag}). Для $\Pi$-атомов это не так.

\begin{definition}[устранимый $\Pi$-атом]\label{ru:def:elim}
$\Pi$-атом $\theta$ \emph{устраним}, если для некоторого
$\Pi$-свободного высказывания $A$ эквивалентность $\theta\aeq A$
$\sB$-общезначима.
\end{definition}

Отметим сразу, что устранение квантора в задачах
(предложение~\ref{ru:prop:pienc}) даёт и устранимость атомов: если
задача $\Pi\vec\alpha\,\chi$ $\sB$-равносильна $\Pi$-свободной задаче
$\psi$, то высказывание $\quest\Pi\vec\alpha\,\chi$ равносильно
$\quest\psi$, а то приводится по (A4) к $\Pi$-свободной нормальной
форме; так устранимы многие атомы, на первый взгляд неприводимые,
например $\quest\Pi\alpha((\alpha\pto\beta\por\gamma)\pto
((\alpha\pto\beta)\por(\alpha\pto\gamma)))$, равносильное
$(\quest\gamma\ato\quest\beta)\aor(\quest\beta\ato\quest\gamma)$, то есть
классической тавтологии.

Все $\Pi$-атомы, встречавшиеся до сих пор, устранимы:
$\quest\Pi\alpha\,\alpha$ и $\quest\Pi\alpha(\pneg\pneg\alpha\pto\alpha)$
равносильны $\bot$ (теорема~\ref{ru:prop:piatoms}),
$\quest\Pi\beta\bigl((\bang A\pto\beta)\pto\beta\bigr)$ равносильно
$\square A$ (предложение~\ref{ru:prop:pienc}(2) при $\varphi:=\bang A$), а $\quest\Pi\beta\Pi\gamma\rho$
общезначимо (пример~\ref{ru:ex:pider}) и потому равносильно $\aneg\bot$.
Устраним и атом $\quest\Pi\alpha(((\alpha\pto\beta)\pto\alpha)\pto\alpha)$
--- он $\sB$-равносилен $\quest\beta$, --- но эта эквивалентность доказана
лишь семантически, и выводима ли она в $\HCB$, неизвестно
(предложение~\ref{ru:prop:pienc}(4), замечание~\ref{ru:rem:peirce}). Существуют ли
неустранимые $\Pi$-атомы, мы не знаем (\S\ref{ru:sec:todo}).

\subsection{Безмодальные $\Pi$-свободные общезначимые задачи}\label{ru:sec:modalfreeval}

Вернёмся к самому простому классу задач --- задачам без модальностей,
то есть построенным из переменных и константы $\varnothing$
связками $\pto$, $\pand$, $\por$. Для каждой семантики
$Z\in\{\sI,\sB,\sE\}$ обозначим через $\Lp_Z$ множество безмодальных
задач, общезначимых в $Z$ (отметим, что по своей грамматике
безмодальные задачи не содержат и квантора $\Pi$).
По
леммам~\ref{ru:lem:closure} и~\ref{ru:lem:closurebs} все три множества
замкнуты относительно (MP) и подстановки безмодальных задач (такие
подстановки не затрагивают ни $\bang$, ни $\Pi$ и потому безопасны), то
есть являются
\emph{суперинтуиционистскими логиками} в языке
$\pto,\pand,\por,\varnothing$ (константа $\varnothing$ играет роль лжи):
$\Int$ они содержат по предложению~\ref{ru:prop:intval}. Теперь, когда доказана
теорема корректности, всякий факт о необщезначимости даёт и факт о
невыводимости; на это мы будем ссылаться в~\S\ref{ru:sec:indep}.

Классических схем ни в одной из трёх логик при этом нет.

\begin{proposition}[задачи остались интуиционистскими]\label{ru:prop:noem}
Пусть $\alpha$, $\beta$ --- переменные сорта задача. Задачи
$\alpha\por\pneg\alpha$, $\pneg\alpha\por\pneg\pneg\alpha$,
$\pneg\pneg\alpha\pto\alpha$ и $((\alpha\pto\beta)\pto\alpha)\pto\alpha$ не
общезначимы ни в одной из трёх семантик.
\end{proposition}

\begin{proof}
Достаточно доказать необщезначимость второй задачи --- слабого закона
исключённого третьего $\pneg\alpha\por\pneg\pneg\alpha$. В самом деле,
первая задача влечёт вторую в $\Int$, а третья и четвёртая --- схемы, каждая
из которых вместе с $\Int$ даёт всю классическую логику; поскольку
общезначимость сохраняется подстановкой задач вместо задачных переменных и
правилом (MP) (леммы~\ref{ru:lem:closure} и~\ref{ru:lem:closurebs}), из
общезначимости любой из них следовала бы общезначимость всех классических
тавтологий, в том числе второй задачи.

В семантиках $\sI,\sE$: фиксируем перечислимое неразрешимое множество
$X\subset\N$ и выразим свойство $n\in X$ арифметической
$\Sigma$-формулой $\exists x A(x,n)$. 
Рассуждая от противного,  
допустим существует решение
$S$ задачи $\pneg\alpha\por\pneg\pneg\alpha$.
Для произвольного $n\in\N$  рассмотрим интерпретацию $*_n$, интерпретирующую переменную $\alpha$ формулой $A(n,x)$.
Тогда решениями задачи $(\pneg\alpha\por\pneg\pneg\alpha)^{*_n}$ являются все числа вида 
$c(0,m)$ (и только они),
если $n\notin X$ и все числа $c(1,m)$ (и только они),
если $n\in X$. И этот факт доказуем во внешней теории.
Таким образом, внешняя теория доказывает
$n\in X\aeq \exists m ( S(*_n)=c(1,m))$, а значит это утверждение истинно для всех $n$.
Поэтому с помощью алгоритма $S$, применённого к коду интерпретации
$*_n$, можно выяснить принадлежность $n$ к $X$. Противоречие.

\remove{
--- её решение. По
лемме~\ref{ru:lem:clauses}(5),(6) элемент множества
$\val{\alpha\por\pneg\alpha}$ имеет вид $c(0,u)$ с $u\in\val{\alpha}$, если
$\val{\alpha}\neq\varnothing$, и вид $c(1,v)$, если $\val{\alpha}=\varnothing$
(при непустом $\val{\alpha}$ множество $\val{\pneg\alpha}$ пусто). В
семантиках $\sI$, $\sE$ по
метке числа $S(\gn{\alpha^*})$ вычисляется истинностное значение предложения
$\exists x\,\alpha^*(x)$. Но уже для $\Sigma_1$-формул $\alpha^*$ это
невозможно: получился бы разрешающий алгоритм для проблемы остановки.}
В
семантике $\sB$ ещё проще: кода нет, решение --- одно фиксированное
число $c(i,u)$, и при $i=0$ его отвергает ящик $\varnothing$, а при
$i=1$ --- ящик $\{0\}$, и $\ZF$ доказывает отрицание
перевода при таком ящике.
\end{proof}

Итак, все три логики $\Lp_Z$ лежат между $\Int$ и
классической логикой, не совпадая с последней. С $\Int$ они тоже не
совпадают, и разделяющий пример был найден Роузом~\cite{ru:Rose53} --- как
пропозициональная формула, реализуемая при всякой подстановке, но не
выводимая интуиционистски; ею было впервые показано, что пропозициональная
логика реализуемости строго шире $\Int$. В нашей семантике она получается как
частный случай общего утверждения об \emph{обозримых} задачах (теорема~\ref{ru:prop:obz}
ниже).

\subsubsection{Обозримые задачи}
Зафиксируем интерпретацию $*$ (в семантике $\sB$ рассуждаем
внутри $\ZF$ при данном содержимом ящиков). \emph{Программами выбора}
для задачи $\varphi$ при данной интерпретации будем называть элементы множества
$\val{\pneg\pneg\varphi\pto\varphi}$, то есть решения задачи
$\pneg\pneg\varphi\pto\varphi$ при этой интерпретации. По
лемме~\ref{ru:lem:clauses}(6) они устроены так:
если $\val{\varphi}\neq\varnothing$, то $\val{\pneg\pneg\varphi}=\N$, вход
ничего не сообщает, и программа выбора --- это всюду определённая программа со
значениями в $\val{\varphi}$; если же $\val{\varphi}=\varnothing$, то
$\val{\pneg\pneg\varphi}=\varnothing$ и программой выбора является
\emph{любая} программа, в том числе нигде не определённая.

\begin{definition}[обозримая задача]\label{ru:def:obz}
Задача $\varphi$ называется \emph{обозримой}, если по коду $\gn{*}$ вычислим
конечный список программ $q_1,\dots,q_k$ (и их число, и они сами могут
зависеть от кода) со свойством: если $\val{\varphi}\neq\varnothing$, то хотя
бы одна из $q_1,\dots,q_k$ всюду определена и все её значения лежат в
$\val{\varphi}$. (Напомним, что в $\sB$ кодом служит код прозрачной
части, а свойство списка требуется при любом содержимом ящиков.)
\end{definition}

Иными словами, обозримость означает, что кандидатов в программы выбора
конечное число и они выписываются заранее по коду интерпретации. Например,
задача $\pneg\beta\por\pneg\gamma$ 
обозрима с
постоянным списком из двух программ --- <<выдать $c(0,0)$>> и <<выдать
$c(1,0)$>>.

\begin{theorem}[схема Роуза для обозримых задач]\label{ru:prop:obz}
Пусть задача $\varphi$ обозрима и $\psi:=\pneg\pneg\varphi\por\pneg\varphi$.
Тогда задача
\begin{align}\tag{R}\label{ru:rose-formula}
  \bigl((\pneg\pneg\varphi\pto\varphi)\pto \psi\bigr)\pto \psi
\end{align}
общезначима в данной семантике, как только определяющее свойство
списка из определения~\ref{ru:def:obz} при каждом $*$ доказуемо в её
внешней теории. В $\sE$ это условие выполнено всегда (внешняя теория
--- $\TA$, а свойство списка истинно по определению), так что там
задача~\eqref{ru:rose-formula} общезначима безусловно.
\end{theorem}

Перед доказательством напомним, что 
формула 
Роуза $\rho$ --- это задача \eqref{ru:rose-formula} для (обозримой) формулы $\phi=\pneg\beta\por\pneg\gamma$.

\begin{proof}
Разберёмся сперва, что такое решение задачи $\psi$. По леммам~\ref{ru:lem:clauses}(5)
и~\ref{ru:lem:dn} множество $\val{\pneg\pneg\varphi}$ равно $\N$ при
$\val{\varphi}\neq\varnothing$ и пусто иначе, а $\val{\pneg\varphi}$ ---
наоборот. Поэтому решение задачи $\psi$ --- это \emph{одна лишь метка}, и притом
верная: при $\val{\varphi}\neq\varnothing$ решения --- это в точности числа с
меткой~$0$, а при $\val{\varphi}=\varnothing$ --- в точности числа с
меткой~$1$. В частности, годится $c(0,0)$ в первом случае и $c(1,0)$ во втором.

Решение посылки --- это программа $p$, которая по всякой программе выбора для
$\varphi$ останавливается и выдаёт решение задачи~$\psi$. Построим по ней решение
заключения.

Вычислим по коду $\gn{*}$ список $q_1,\dots,q_k$ из
определения~\ref{ru:def:obz} и запустим $p$ параллельно на всех $q_i$. Далее:
\begin{itemize}[leftmargin=1.6em,itemsep=1pt]
\item как только какой-нибудь запуск остановился с ответом, имеющим метку~$0$,
  останавливаемся и выдаём $c(0,0)$;
\item если все $k$ запусков остановились с ответами, имеющими метку~$1$,
  выдаём $c(1,0)$.
\end{itemize}

\emph{Проверка.} Пусть $\val{\varphi}=\varnothing$. Тогда всякая программа
является программой выбора, в частности все $q_i$; значит $p$ на каждой из них
останавливается, и все ответы имеют метку~$1$ --- других решений задачи $\psi$
при пустом $\val{\varphi}$ нет. Поэтому первый пункт не сработает, все $k$
запусков кончатся, и мы выдадим $c(1,0)$ --- верное решение.

Пусть теперь $\val{\varphi}\neq\varnothing$. По определению обозримости
какая-то $q_i$ является программой выбора, значит $p$ на ней останавливается и
выдаёт решение задачи~$\psi$; при непустом $\val{\varphi}$ это число с меткой~$0$.
Стало быть первый пункт рано или поздно сработает, и мы выдадим $c(0,0)$ ---
снова верное решение.

Список $q_1,\dots,q_k$ вычисляется по коду, а всё остальное построение от кода
не зависит, так что решающая функция вычислима. Выводимость формулы
$\varphi^*(S(\gn{*}))$ во внешней теории
получается тем же рассуждением, проведённым внутри неё: разбор
случая по
$\exists x\,\varphi^*(x)$ классический, а всё остальное ---
$\Sigma_1$-рассуждение о работе машины; в $\sB$ рассуждение проводится
поверх свободных множественных переменных, а список вычисляется по коду
прозрачной части. Единственное место, где требуется
дополнительное условие, --- переход <<если $\exists x\,\varphi^*(x)$, то
какая-то $q_i$ --- программа выбора>>: внутри внешней теории им можно
пользоваться, только если оно в ней доказано (для $\TA$ это и есть
истинность, для $\PA$ и $\ZF$ истинности недостаточно).
\end{proof}

Класс обозримых задач достаточно широк; чтобы описать его явно, выделим
ещё один класс, заданный грамматически.

\begin{definition}[$\mathcal O$-задачи]\label{ru:def:obzgram}
\emph{$\mathcal O$-задачами} называем задачи, построенные по грамматике
\[
  \varphi ::= \sigma \mid \pneg\psi \mid (\varphi\pand\varphi) \mid
              (\varphi\por\varphi) ,
\]
где $\sigma$ пробегает $\Sigma$-задачи, а $\psi$ --- произвольные задачи.
Иначе говоря, к $\Sigma$-задачам добавлены всевозможные отрицания, и всё это
замкнуто относительно $\pand$ и $\por$. Всякая $\Sigma$-задача является
$\mathcal O$-задачей; обратное неверно уже для $\pneg\bang A$.
\end{definition}

\begin{lemma}[запас обозримых задач]\label{ru:lem:obz}
\begin{enumerate}[leftmargin=2.2em,itemsep=1pt,label=\textup{(\arabic*)}]
\item Всякая задача $\phi$, для которой общезначима задача $\pneg\pneg\varphi\pto\varphi$,
  обозрима. В частности обозримы все $\Sigma$-задачи и все задачи вида
  $\pneg\psi$.
\item Если $\varphi$ и $\psi$ обозримы, то обозримы $\varphi\pand\psi$ и
  $\varphi\por\psi$.
\end{enumerate}
Всё сказанное относится к каждой из трёх семантик, и во всех случаях
определяющее свойство списка доказуемо во внешней теории.
\end{lemma}

\begin{proof}
(1) Пусть $S$ --- вычислимая функция, решающая
$\pneg\pneg\varphi\pto\varphi$. Возьмём список из одной программы:
$q:=$ <<на любом входе выдать $U_{S(\gn{*})}(0)$>>. Если
$\val{\varphi}\neq\varnothing$, то $U_{S(\gn{*})}$ всюду определена со
значениями в $\val{\varphi}$, значит $q$ всюду определена и её значения лежат
в $\val{\varphi}$. Для $\Sigma$-задач схема
общезначима во всех трёх семантиках --- это схема (A7),
предложение~\ref{ru:prop:a5a7};
для $\pneg\psi$ она выводима уже в $\Int$ и общезначима по
предложению~\ref{ru:prop:intval}.

Доказуемость во внешней теории: для $\Sigma$-задач она входит в
общезначимость (A7) по определению~\ref{ru:def:valid} --- внешняя теория
доказывает
$\forall u\,((\pneg\pneg\varphi)^*(u)\ato\varphi^*(U_{S(\gn{*})}(u)))$, а
посылка эквивалентна $\exists x\,\varphi^*(x)$ по лемме~\ref{ru:lem:dn};
для $\pneg\psi$ --- в общезначимость схемы $\Int$.

(2) Пусть $q_1,\dots,q_k$ --- список для $\varphi$, а $q'_1,\dots,q'_m$ ---
для $\psi$. Для $\varphi\por\psi$ возьмём список из $k+m$ программ: <<запустить
$q_i$ и выдать $c(0,\text{результат})$>> и <<запустить $q'_j$ и выдать
$c(1,\text{результат})$>>. Если $\val{\varphi\por\psi}\neq\varnothing$, то
непусто $\val{\varphi}$ или $\val{\psi}$, и соответствующая обёрнутая
программа всюду определена со значениями в $\val{\varphi\por\psi}$. Для
$\varphi\pand\psi$ возьмём $km$ программ <<запустить $q_i$ и $q'_j$ и выдать
$c(\text{результаты})$>>; если $\val{\varphi\pand\psi}\neq\varnothing$, то
непусты оба сомножителя, и годная пара найдётся. Оба рассуждения дословно
формализуются во внешней теории, так что доказуемость свойства списка
сохраняется.
\end{proof}

\begin{proposition}[схема Роуза общезначима на классе $\mathcal O$]\label{ru:cor:rose}
Для всякой $\mathcal O$-задачи $\varphi$ задача~\eqref{ru:rose-formula}
общезначима  во всех трёх семантиках. 
\end{proposition}

\begin{proof}
Индукция по грамматике из определения~\ref{ru:def:obzgram} показывает, что
всякая $\mathcal O$-задача обозрима, причём определяющее свойство списка
доказуемо во внешней теории: для $\Sigma$-задач и задач вида $\pneg\psi$ это
пункт~(1) леммы~\ref{ru:lem:obz}, а связки $\pand$ и $\por$ сохраняют
обозримость по её пункту~(2). Остаётся применить
теорему~\ref{ru:prop:obz}.
\end{proof}

Схему Роуза мы \emph{не} включили в число аксиом исчислений, 
хотя она
общезначима во всех трёх семантиках. Исчисления составляются так,
чтобы облегчать доказательства общезначимости, а для этого аксиомы
должны быть естественными принципами с простой формулировкой --- таковы аксиомы
(A1)--(A7).
Схема Роуза --- не
принцип, а готовый результат: она длинна, класс $\mathcal O$-задач
определяется непросто, полноты на задачах она, по всей видимости, не
даёт, и трудно представить, чтобы кто-нибудь воспользовался ею для
доказательства общезначимости чего-то другого. Поэтому мы 
записываем её в примеры \emph{неполноты}: формула Роуза общезначима,
но не выводима в $\Int$, а потому и в $\HCI$  (\S\ref{ru:sec:indep}).

Обозримость в теореме~\ref{ru:prop:obz} существенна: если $\varphi$ ---
переменная, схема~\eqref{ru:rose-formula} не общезначима. В $\Int$ она
равносильна схеме
\[
  \bigl((\pneg\pneg\varphi\pto\varphi)\pto(\varphi\por\pneg\varphi)\bigr)\pto
  (\pneg\pneg\varphi\por\pneg\varphi),
\]
отличающейся тем, что в заключении посылки стоит
$\varphi\por\pneg\varphi$, а не $\pneg\pneg\varphi\por\pneg\varphi$:
посылки равносильны. В самом деле,
$(\varphi\por\pneg\varphi)\pto(\pneg\pneg\varphi\por\pneg\varphi)$, так как
$\varphi\pto\pneg\pneg\varphi$; обратно, в предположении
$\pneg\pneg\varphi\pto\varphi$ из $\pneg\pneg\varphi\por\pneg\varphi$
следует $\varphi\por\pneg\varphi$. Эту схему называют \emph{аксиомой
Скотта}: Крайзел и Патнэм~\cite{ru:KP57} упоминают, что Д.~Скотт доказал
дизъюнктивное свойство логики, получающейся добавлением её к $\Int$.
Формула Роуза --- её частный случай при $\varphi:=\pneg\beta\por\pneg\gamma$.
Мы будем называть эту задачу \emph{принципом Скотта}. Янков~\cite{ru:Yankov63}
доказал, что принцип Скотта с переменной не реализуем по Клини;
следующая теорема --- аналог этого результата для наших семантик (в
доказательстве, как и у Янкова, используется диагонализация по числу
шагов вычисления).

\begin{theorem}[принцип Скотта не общезначим]\label{ru:prop:scott}
Задача
\[\phi:=
  \bigl((\pneg\pneg\alpha\pto\alpha)\pto(\pneg\pneg\alpha\por\pneg\alpha)\bigr)\pto
  (\pneg\pneg\alpha\por\pneg\alpha),
\]
где $\alpha$ --- переменная сорта задача (то есть схема~\eqref{ru:rose-formula}
при $\varphi:=\alpha$), не общезначима ни в одной из трёх семантик.
\end{theorem}


\begin{proof}
Пусть $S$ --- произвольный алгоритм. Нам нужно определить 
интерпретацию $*$, для которой невыводима или ложна формула $\phi^*(S(*))$.
Как уже говорилось в Замечании~\ref{ru:rem:assumptions},
мы докажем выводимость отрицания этой формулы во внешней теории.
Буквой $X$ будем обозначать значение $\val{\alpha}$: в семантике $\sB$ это
содержимое ящика, в остальных --- множество, определяемое формулой
$\alpha^*$. 

\emph{Семантика $\sB$.} 
В ней доказательство проще: переменных-высказываний в задаче $\phi$
нет, код интерпретации пуст, и решатель --- одно число. Поэтому с неё
и начинаем.

Опишем, что такое решения входящих в задачу подзадач.

По лемме~\ref{ru:lem:clauses}(6) решение задачи $\pneg\pneg\alpha\pto\alpha$ ---
это, при $X\neq\varnothing$, номер всюду определённой функции со значениями
в~$X$ (множество $\val{\pneg\pneg\alpha}$ равно тогда $\N$, так что вход
ничего не сообщает), а при $X=\varnothing$ --- любое число. Назовём такой
номер $i$ \emph{$X$-корректным}.

Решение задачи $\pneg\pneg\alpha\por\pneg\alpha$ --- это число $c(0,v)$
(если $X\neq\varnothing$) или число $c(1,v)$ (если $X=\varnothing$), где $v$
произвольно: значение каждого из дизъюнктов, когда оно непусто, есть
$\N$. Значит, решение --- это по существу одна лишь метка, то есть ответ
на вопрос, пусто ли $X$; будем говорить, что программа отвечает <<не
пусто>> или <<пусто>>.

Итак, решение посылки --- это программа $p$, которая по всякому
$X$-корректному номеру $i$ останавливается и правильно отвечает, пусто ли
$X$.

Код 
интерпретации $*$ пуст,  значит
$S(\gn{*})$ --- некоторая  
фиксированная программа, не зависящая
от  ящика. Положим $r:=S(\gn{*})$ 

\medskip
\emph{Диагональная конструкция.} 
По
теореме Клини о неподвижной точке существует программа $p$, которая на входе
$i$ делает следующее:
\begin{enumerate}[leftmargin=2.2em,itemsep=1pt,label=\textup{(\arabic*)}]
\item моделирует работу программы $r$ на входе $p$, то есть на
  собственном номере, ровно $i$ шагов;
\item если за это время $r$ остановилась с ответом <<пусто>>, то отвечает
  <<не пусто>>;
\item иначе отвечает <<пусто>>.
\end{enumerate}
Ограничение числом шагов существенно дважды. Во-первых, моделирование конечно,
и никакого зацикливания не возникает, хотя $r$ по дороге может моделировать
саму~$p$. Во-вторых --- и это главное --- если $r$ останавливается с
ответом <<пусто>>, то входы, на которых $p$ отвечает <<пусто>>, образуют
\emph{конечный} начальный отрезок. (Сам вход $i$ программа $p$ использует
только как число шагов.)

Свяжем с $r$ множество $X_r$: если $r$ на входе $p$ останавливается на шаге
$s$ с ответом <<пусто>>, положим $X_r:=\{\nu(s)\}$, где $\nu(s)$ ---
наименьшее число, не встречающееся среди $U_i(0)$ при $i<s$ (таких чисел
конечное число, так что $\nu(s)$ определено; вычислимость $\nu$ не
понадобится, достаточно арифметической определимости); во всех остальных
случаях положим $X_r:=\varnothing$.

\medskip
\emph{Разбор случаев.} Проверим: если $X=X_r$, то $p$ --- решение посылки,
а $r$ на входе $p$ расходится или выдаёт неверную метку.

\emph{Случай 1: $r$ на входе $p$ расходится.} Тогда при всяком $i$ шаг~(2) не
срабатывает, и $p(i)$ отвечает <<пусто>>. Значит $p$ всюду определена и всегда
отвечает <<пусто>>. При этом $X=\varnothing$, всякий номер $X$-корректен, и
ответ <<пусто>> --- правильное решение задачи $\pneg\pneg\alpha\por\pneg\alpha$. Стало
быть $p$ --- решение посылки, а $r$ на нём расходится.

\emph{Случай 2: $r$ на входе $p$ останавливается с ответом <<не пусто>>.}
Шаг~(2) снова не срабатывает ни при каком $i$, так что $p$ по-прежнему всюду
отвечает <<пусто>> и является решением посылки при $X=\varnothing$. Но тогда
правильный ответ --- <<пусто>>, а $r$ ответила иначе.

\emph{Случай 3: $r$ на входе $p$ останавливается на шаге $s$ с ответом
<<пусто>>.} Тогда $X=\{\nu(s)\}$ --- непустое множество. Проверим, что $p$ ---
решение посылки. Пусть $i$ --- $X$-корректный номер; тогда $U_i$ всюду
определена и $U_i(0)=\nu(s)$. Если бы было $i<s$, то $\nu(s)$ встречалось бы
среди $U_i(0)$ при $i<s$, чего быть не может по выбору $\nu(s)$. Значит
$i\geq s$, при моделировании $i\geq s$ шагов программа $p$ видит остановку $r$
с ответом <<пусто>>, срабатывает шаг~(2), и $p$ отвечает <<не пусто>>
--- правильно, поскольку $X$ непусто. Итак,
$p$ --- решение посылки, множество $X$ непусто, а $r$ ответила <<пусто>>.

Во всех трёх случаях $p$ --- решение посылки при $X=X_r$, а $r$ на нём
ошибается или расходится.


Возьмём в качестве содержимого ящика множество $X_r$: по
доказанному $p$ решает посылку, а $r$ на ней ошибается или расходится, то
есть $r$ не принадлежит значению задачи при этом содержимом ящика.
Это заканчивает доказательство отрицания формулы
$\phi^*(S(*))$ в теории множеств.

\emph{Семантики $\sE$ и $\sI$.} 
На этот раз нам нужно начать с построения формулы $\alpha^*(x)$,
которую мы определим с помощью теоремы о неподвижной точке для PA.
Опять
положим $r:=S(\gn{\alpha^*})$. Теперь программа $r$ строится алгоритмом $S$ по 
формуле $\alpha^*$, которую мы должны предъявить, а предъявить её нужно так,
чтобы значение $\val{\alpha}$ совпало с $X_r$. Получается порочный круг,
и существование $\alpha^*$ обеспечивается леммой о неподвижной точке для
$\PA$: существует формула $\alpha^*(x)$, для которой $\PA$ доказывает
\[
  \alpha^*(x)\ \aeq\ \exists s\,\bigl(\text{<<$
  S(\gn{\alpha^*})
  $ на входе
  $p$ останавливается на шаге $s$ с ответом <<пусто>>>>}\ \aand\ x=\nu(s)\bigr),
\]
где $p$ строится по $r=S(\gn{\alpha^*})$ вышеприведённой Диагональной конструкцией. Лемма применима: условие в
скобках арифметично, а $\nu$ арифметически определима. 
По построению формула $\alpha^*(x)$  равносильна в PA
формуле $x\in X_r$, так что разбор
случаев опровергает формулу  $\phi^*(S(*))$ во внешней теории. 
\end{proof}

Отметим попутно, что класс $\mathcal O$ в схеме Роуза ни расширить до
всех задач, ни существенно сузить нельзя: с переменной вместо
$\varphi$ схема влечёт принцип Скотта, который не общезначим
(теорема~\ref{ru:prop:scott}); на $\Sigma$-задачах она выводима из (A7)
одним (MP), так что всё её содержание держится на отрицаниях.



\subsubsection{Позитивный безмодальный фрагмент}\label{ru:sec:posfrag}

\emph{Позитивной} назовём безмодальную задачу, не содержащую
$\varnothing$ (а значит, и $\pneg$). \emph{Позитивный фрагмент} логики --- множество её позитивных
формул; логика \emph{имеет интуиционистский позитивный фрагмент}, если её
позитивный фрагмент совпадает с позитивным фрагментом $\Int$. Мы докажем,
что позитивные фрагменты всех трёх логик $\Lp_Z$
интуиционистские: вся разница этих логик с $\Int$ живёт на формулах
с~$\varnothing$.

\begin{definition}[Медведев]\label{ru:def:crit}
\emph{Критическая импликация} --- позитивная формула вида
$\bigwedge_i((P_i\pto Q_i)\pto Q_i)\pto R$, где $P_i$ --- конъюнкции
переменных, $Q_i$ и $R$ --- дизъюнкции переменных, причём при каждом $i$
формулы $P_i$ и $Q_i$ не имеют общих переменных, а каждая из $P_i$, $Q_i$, $R$
содержит хотя бы одну переменную. Для каждого $n>0$ положим
\[
  J_n:=\bigwedge_{\varnothing\neq E\subsetneq\{1,\dots,n\}}
    \Bigl(\bigl(\bigwedge_{j\notin E}\alpha_j\pto\bigvee_{i\in E}\alpha_i\bigr)
      \pto\bigvee_{i\in E}\alpha_i\Bigr)\pto\bigvee_{i=1}^n \alpha_i
\]
--- самую слабую критическую импликацию в переменных-задачах
$\alpha_1,\dots,\alpha_n$; в частности $J_1=\alpha_1$.
\end{definition}

\begin{theorem}[Медведев, 1962;~\cite{ru:Medvedev62}]\label{ru:th:medv}
Суперинтуиционистская логика имеет интуиционистский позитивный фрагмент в том
и только том случае, когда в ней не выводима ни одна критическая импликация,
и это равносильно тому, что в ней не выводима ни одна из формул $J_n$.
\end{theorem}

\begin{theorem}[Янков, 1968;~\cite{ru:Yankov68}]\label{ru:th:jan}
Суперинтуиционистская логика имеет интуиционистский позитивный фрагмент в том
и только том случае, когда она содержится в $\J:=\Int+(\pneg\alpha\por
\pneg\pneg\alpha)$ --- логике слабого закона исключённого третьего.
\end{theorem}

Критерий Медведева удобен для логик, заданных семантически (а наша задана
именно так), критерий Янкова --- для заданных аксиоматически. Мы пользуемся
ими в связке: первый даёт интуиционистский позитивный фрагмент, второй
превращает это во включение $\Lp_Z\subseteq\J$.

Логика реализуемых пропозициональных формул устроена так: переменные
заменяются арифметическими формулами, связки --- клиниевскими операциями, и
формула считается принадлежащей логике, если при всякой подстановке у неё есть
реализация. Различают три варианта~\cite{ru:Plisko09}: $\R$ --- реализация просто
существует; $\Reff$ --- реализацию можно \emph{эффективно найти} по
подставляемым формулам; $\Rcnst$ --- одна и та же реализация годится при всех
подстановках. Очевидно $\Rcnst\subseteq\Reff\subseteq\R$.

Три уровня из замечания~\ref{ru:rem:three} --- это в точности они. Наше
понятие общезначимости --- средний: решение выдаётся вычислимой функцией по
коду интерпретации. Разница с $\Reff$ только в том, что у нас значением
переменной служит произвольное арифметически определимое множество, а там ---
множество реализаций арифметической формулы.

Для логик реализуемости $\R$, $\Reff$, $\Rcnst$ интуиционистичность
позитивного фрагмента известна (Медведев, Плиско;~\cite{ru:Plisko09}), но
сослаться на это здесь нельзя. Казалось бы, достаточно
включения $\Lp_{\sE}\subseteq\R$ --- условие $\exists S\,\forall{*}$ сильнее, чем
$\forall{*}\,\exists$, --- и тогда позитивный фрагмент $\Lp_{\sE}$ оказался бы зажат
между совпадающими позитивными фрагментами $\Int$ и $\R$. Но семантики не
совпадают: у нас значением переменной служит произвольное арифметически
определимое множество, а в определении $\R$ --- множество реализаций
арифметической формулы, и \emph{не всякое арифметическое множество является
множеством реализаций}. Поэтому мы пойдём прямым путём: по
теореме~\ref{ru:th:medv} достаточно проверить, что ни одна формула $J_n$ не
общезначима. Нужная конструкция есть в~\cite{ru:VSSCh}; операции над множествами
там ровно те же, что в нашем определении~\ref{ru:def:val} (метки $0$ и~$1$ у
дизъюнкции, пары у конъюнкции, $e(x)\in Y$ у импликации), так что
перенос сводится к тому, чтобы замкнуть круговую зависимость решателя от кода
интерпретации.

Скажем, что программа \emph{коперечисляет} множество, если она перечисляет его
дополнение.

\begin{lemma}[\cite{ru:VSSCh}]\label{ru:lem:vssch}
Пусть даны натуральное $m\ge 1$, программа $q$, перечисляющая множество $M$ из
не более чем $2^m$ чисел, и критическая импликация $J_n$ от $n$ переменных.
Тогда можно эффективно (то есть некоторой всюду определённой вычислимой
функцией от $m$, $n$ и $q$) найти программы $a_1,\dots,a_n$, коперечисляющие
непустые множества $A_1,\dots,A_n$, для которых
\[
  M\cap J_n(A_1,\dots,A_n)=\varnothing .
\]
При этом все элементы множеств $A_j$ не превосходят $2^{Cm}$, где константа $C$
зависит только от $n$.
\end{lemma}

Доказательство этой леммы не выходит за рамки PA.
Конструкция в двух словах: множества $A_j$ берутся вначале
конечными отрезками подходящих размеров, а затем из них по шагам удаляются
элементы --- всякий раз, когда очередная программа $e\in M$ обнаруживает себя
на подходящих аргументах, её ответ удаляется из соответствующего $A_j$.
Одновременно по теореме о неподвижной точке строятся вспомогательные
программы, свидетельствующие о принадлежности посылок
$((P_i\pto Q_i)\pto Q_i)$, и главная работа --- оценка числа удалений, которая
и показывает, что процесс останавливается, а множества не пустеют. Именно
здесь используется условие из определения критической импликации, что $P_i$ и
$Q_i$ не имеют общих переменных.

\begin{theorem}[позитивный фрагмент интуиционистский]\label{ru:th:posfrag}
Позитивная задача общезначима (в любой из трёх семантик) тогда и
только тогда, когда она выводима в~$\Int$.
\end{theorem}

\begin{proof}
Одна сторона --- предложение~\ref{ru:prop:intval}. Обратно, по
теореме~\ref{ru:th:medv}
достаточно проверить, что ни одна формула $J_n$ не общезначима --- ни в
одной из трёх семантик.

\emph{Семантика $\sB$ --- самый простой случай.} Высказывательных
переменных в $J_n$ нет, поэтому код прозрачной части пуст и решатель
сводится к одному числу $e_0$. Применим лемму~\ref{ru:lem:vssch} при
$m:=1$ к программе <<напечатать $e_0$ и остановиться>>: получим непустые
коперечислимые множества $A_1,\dots,A_n$ с
$e_0\notin J_n(A_1,\dots,A_n)$. Положим содержимое ящиков
$X_{\alpha_j}:=A_j$; по лемме~\ref{ru:lem:clauses} значения связок --- это
в точности операции из леммы~\ref{ru:lem:vssch}, так что
$\val{J_n}=J_n(A_1,\dots,A_n)$ и утверждение $e_0\notin\val{J_n}$ доказуемо в PA, а тем более в ZF.

\emph{Семантики $\sE$ и $\sI$.} Здесь решатель получает код, и
нужна самоссылка. 
Пусть $S$ --- всюду определённая вычислимая функция, предположительно решающая $J_n$.
Это 
означает:  формула $J_n^*(S(\gn{*}))$ доказуема или истинна  при всякой интерпретации~$*$. 

\emph{Построение по коду.} Возьмём $m:=1$. Для произвольного кортежа $c$  из $n$ формул
обозначим через $q_c$ программу <<напечатать $S(c)$ и остановиться>>; её номер
вычислим по $c$. Пусть $a_1(c),\dots,a_n(c)$ --- программы, которые выдаёт по
$q_c$ алгоритм из леммы~\ref{ru:lem:vssch}. Этот алгоритм есть всюду определённая
вычислимая функция от номера входной программы, поэтому $a_j(c)$ определены
при \emph{всяком} $c$; положим
\[
  A_j(c):=\{x:\ \text{программа } a_j(c)\ \text{не печатает } x\}.
\]
Условие справа --- это $\forall i\ \aneg\bigl(U(a_j(c),i)\!\downarrow=x\bigr)$,
то есть $\Pi_1$-условие, причём равномерно по~$c$. Значит есть арифметические
формулы $\theta_j(c,x)$, задающие $A_j(c)$ при каждом~$c$.

Утверждение же леммы гарантирует, что для всех $c$ 
множества $A_j(c)$ непусты и что
\[
  S(c)\notin J_n\bigl(A_1(c),\dots,A_n(c)\bigr) \tag{$\ast$},
\]
причем это утверждение доказуемо в PA. Заметим, что 
это утверждение есть отрицание формулы $(J_n)^{*}(S(c))$ для интерпретации 
$*$ (зависящей от $c$), определенной равенствами $\alpha_j^*:=\theta_j(c,x)$. 
Если бы нам удалось выбрать $c$ так, чтобы $c$ совпадало с $\langle \theta_1(c,x),\dots,\theta_n(c,x)\rangle$,
мы достигли бы цели: формула $\aneg(J_n)^{*}(S(c))$ была бы
доказуема в PA для  соответствующей интерпретации  $*$. Но ясно, что это безнадёжно. 


\emph{
Порочный круг разрывается}
совместной леммой о неподвижной точке (см., например,
\cite[гл.~1]{ru:Boolos93}): для любых формул $\theta_1(c,x),\dots,\theta_n(c,x)$
существуют формулы $\chi_1(x),\dots,\chi_n(x)$, для которых
\[
  \PA\vdash \chi_j(x)\aeq\theta_j\bigl(\langle
   \gn{\chi_1},\dots,\gn{\chi_n}\rangle,\,x\bigr),\qquad j=1,\dots,n .
\]
Теперь мы можем определить интерпретацию $*$, положив $\alpha_j^*:=\chi_j$ и $c:=\langle
   \gn{\chi_1},\dots,\gn{\chi_n}\rangle$. 
По построению мы можем доказать в PA равносильность $\theta_j(c,x)$ и $\chi_j(x)$,
а значит и равносильность $J_n\bigl(A_1(c),\dots,A_n(c)\bigr)$
и $J_n\bigl(\chi_1,\dots,\chi_n\bigr)$. Значит формула $S(c)\notin J_n\bigl(\chi_1,\dots,\chi_n\bigr)$
так же доказуема в PA. А эта формула и есть $\aneg(J_n)^{*}(S(c))$.

Как обычно, для семантики $\sE$ мы просто можем сказать, что из доказуемости следует истинность,
поэтому $J_n$ не общезначима. А для семантики $\sI$, как и выше для
$\sB$, воспользоваться непротиворечивостью внешней теории.
\end{proof}

Сведём вместе обе границы. 

\begin{corollary}\label{ru:cor:sandwich}
Для каждой из трёх семантик $Z$
\[
  \Int\ \subsetneq\ \Lp_Z\ \subsetneq\ \J .
\]
Кроме того, $\Lp_{\sI}\subseteq\Lp_{\sE}$.
\end{corollary}

\begin{proof}
Первое включение и его строгость --- предложения~\ref{ru:prop:intval}
и~\ref{ru:cor:rose} (формула Роуза общезначима во всех трёх семантиках,
причём с постоянным решателем). Далее: по теореме~\ref{ru:th:posfrag}
каждая логика $\Lp_Z$ имеет интуиционистский позитивный фрагмент, а
тогда по критерию Янкова (теорема~\ref{ru:th:jan}) $\Lp_Z\subseteq\J$.
Строгость: слабый закон исключённого третьего принадлежит $\J$ по
определению, но не общезначим ни в одной семантике
(предложение~\ref{ru:prop:noem}). Включение
$\Lp_{\sI}\subseteq\Lp_{\sE}$: внутренняя теория $\PA$ содержится во
внешней $\TA$.
\end{proof}

\subsection{Независимость аксиом}\label{ru:sec:indep}

Нетривиальные факты о выводимости --- это прежде всего факты о
\emph{не}выводимости. Соберём их здесь.

\paragraph{Неполнота на задачах.} Семейства \emph{общезначимых} безмодальных задач --- логики $\Lp_Z$ ---
описаны в~\S\ref{ru:sec:modalfreeval} с точностью до сэндвича. Здесь речь пойдёт о
\emph{выводимых} безмодальных задачах, причём о выводимости в $\HCI$.
Про $\HCE$ скажем сразу то
немногое, что знаем: его безмодальный фрагмент содержит фрагмент $\HCI$
(всякая теорема $\HCI$ --- аксиома $\HCE$) и содержится в $\Lp_{\sE}$ (корректность,
теорема~\ref{ru:th:sound}), но совпадает ли он с фрагментом
$\HCI$ --- неизвестно: из (A8) и (A9) по (MP) в принципе могли бы
вывестись новые безмодальные задачи, и модель
теоремы~\ref{ru:prop:modalfree} ниже этого не исключает --- схемы (A8) и
(A9) в ней ложны (см.\ также \S\ref{ru:sec:todo}).

Про $\HCI$ известно следующее.

\begin{theorem}[безмодальный фрагмент --- это в точности $\Int$]
\label{ru:prop:modalfree}
Пусть $\varphi$ --- задача, не содержащая ни $\bang$, ни $\quest$. Тогда
$\varphi$ выводима в $\HCI$ в том и только том случае, когда она
выводима в $\Int$.
\end{theorem}

Доказательство --- в приложении.

\begin{corollary}[неполнота на задачах]\label{ru:cor:incomplete}
Исчисление $\HCI$ неполно относительно каждой из семантик $\sI$, $\sE$:
формула Роуза общезначима в них (предложение~\ref{ru:cor:rose}), но
безмодальна и в $\Int$ не выводима по теореме Роуза~\cite{ru:Rose53}, а
потому не выводима и в $\HCI$. 
\end{corollary}
Выводима ли формула Роуза в $\HCB$, мы не знаем
(замечание~\ref{ru:rem:qfree}); ещё два кандидата на пример неполноты
$\HCB$ --- задачи
\[
  \Pi\alpha\bigl(((\alpha\pto\beta)\pto\alpha)\pto\alpha\bigr)\pto\beta
  \qquad\text{и}\qquad
  \Pi\alpha\bigl(((\alpha\pto\beta)\pto\beta)\pto(\alpha\por\beta)\bigr)\pto\beta,
\]
$\sB$-общезначимые по предложению~\ref{ru:prop:pienc}(4),(6); их
выводимость в $\HCB$ неизвестна (замечание~\ref{ru:rem:peirce}).

\paragraph{Схема (A7).}

\begin{theorem}\label{ru:prop:a7main}
Схема (A7) не выводима в $\HCI$ из остальных аксиом.
\end{theorem}

Доказательство --- двухточечная модель Крипке, в которой верны все схемы
(A1)--(A6) и все три правила вывода, а задача
$\pneg\pneg\bang\bot\pto\bang\bot$ --- нет; она построена в
приложении~\ref{ru:sec:a7indep} (предложение~\ref{ru:prop:a7indep}). Для $\HCE$
вопрос открыт, и моделями его так не решить (замечание ниже).

\paragraph{Схемы (A8) и (A9).} Они не просто не выводимы в $\HCI$, а
несовместимы с ним.

\begin{remark}
\label{ru:rem:a10fork}
Ни (A8), ни (A9) нельзя добавить к $\HCgl$: \emph{оба исчисления
$\HCgl+$(A8) и $\HCgl+$(A9) противоречивы}. Достаточно одного экземпляра
(A8) --- высказывания $\square\bot\ato\bot$: правило (AtoP) даёт
$\bang(\square\bot\ato\bot)$, схема Лёба (A6) при $A:=\bot$ --- задачу
$\bang\bot$, правило (PtoA) --- высказывание $\quest\bang\bot$, то есть
$\square\bot$, и $\square\bot\ato\bot$ даёт $\bot$.

Схема (A9) этот экземпляр доставляет: в примере~\ref{ru:ex:proofs}(4) из
(A9) выведена задача $\pneg\bang\bot$, а из неё правилом (PtoA) и
схемами (A4a), (A4d) --- высказывание $\quest\bang\bot\ato\bot$, то есть
$\square\bot\ato\bot$.

Иначе говоря, (A9) относится к сорту задача ровно так же, как рефлексия (A8)
--- к сорту высказывание: и та, и другая
несовместимы с парой <<необходимость $+$ Лёб>>. Поэтому обе стоят в
$\HCE$, где правила (AtoP) нет, и поэтому же $\HCE$, вопреки виду,
\emph{не} замкнуто относительно (AtoP): его непротиворечивость даёт только
теорема~\ref{ru:th:sound}.

Практическое следствие: никакая одна алгебраическая модель не может
удовлетворять всем аксиомам и правилам $\HCgl$ и вдобавок всем частным
случаям (A9). Независимость чего-либо от $\HCgls$ моделями так не
доказывается --- нужна техника вроде той, какой пользуются для $\GLS$.
(Конец замечания.)
\end{remark}

Обратный вопрос --- выводится ли (A8) из (A9) --- решается
отрицательно, хотя частные случаи (A8) из (A9) получаются: так в
примере~\ref{ru:ex:proofs}(4) выведено $\square\bot\ato\bot$.

\begin{proposition}[(A8) не выводима из (A9)]\label{ru:prop:a8indep}
Пусть $A$ --- переменная сорта высказывание. Высказывание
$\square A\ato A$ не выводимо из теорем $\HCI$ и схемы (A9) правилами
(MP) и (PtoA). Выводимы же из них все экземпляры рефлексии
$\square\quest\varphi\ato\quest\varphi$.
\end{proposition}

\begin{proof}
Выводимость: применим к экземпляру (A9) правило $(\quest\to)$ и схемы
(A4a), (A4d): $\quest\bang\quest\varphi\ato\quest\pneg\pneg\varphi$, а
$\quest\pneg\pneg\varphi\aeq\aneg\aneg\quest\varphi\aeq\quest\varphi$.

Невыводимость. Пусть $C_1,\dots,C_m$ --- вывод высказывания $C_m$ из
теорем $\HCI$ и экземпляров (A9) правилами (MP) и (PtoA). Сопоставим
каждой формуле вывода высказывание $C_i'$: если $C_i$ --- высказывание,
то $C_i':=C_i$, а если задача, то $C_i':=\quest C_i$ 
(внутрь формул мы не заглядываем). 

Покажем, что
$\nf{C_m}$ выводится из теорем $\GL$ и конечного числа высказываний
вида $\square N\ato N$ одним правилом (MP), причём $N$ здесь ---
$\square$-формулы, в которых все переменные-высказывания стоят под
$\square$ (свободно стоят только атомы $\quest\alpha$); в частности,
$N=A$ среди них нет. Если $C_i$ --- теорема
$\HCI$, то и $C_i'$ --- теорема $\HCI$ (при необходимости применяем
(PtoA)), и $\GL\vdash\nf{C_i'}$ по следствию~\ref{ru:cor:solovay}. Если
$C_i$ --- экземпляр (A9) $\bang\quest\varphi\pto\pneg\pneg\varphi$, то
$\nf{C_i'}=\square N\ato\aneg\aneg N$, где $N=\varphi^\quest$; она
получается по (MP) из $\square N\ato N$ и теоремы классической логики.
Формула $N$ --- указанного вида, поскольку $\varphi$ --- задача:
высказывание входит в задачу только в виде $\bang B$, а
$(\bang B)^\quest=\square\,\nf{B}$. Если $C_i$ получена по (MP) из задач $C_j$ и
$C_l=C_j\pto C_i$, то $\nf{C_i'}$ получается по (MP) из $\nf{C_j'}$ и
$\nf{C_l'}=\nf{C_j'}\ato\nf{C_i'}$; правило (MP) на высказываниях
переходит в себя. Наконец, если $C_i=\quest C_j$ получена по (PtoA),
то $C_i'=C_j'$. 

Итак, $\nf{C_m}$ получается из теорем $\GL$ и конечного
числа высказываний $\square N_i\ato N_i$ указанного вида
одними применениями (MP), и потому
\[
  \GL\vdash\textstyle\bigwedge_i(\square N_i\ato N_i)\ato\nf{C_m}.
\]
Для $C_m=\square A\ato A$ (здесь $\nf{C_m}=C_m$) это невозможно. Пусть $k$
больше модальной глубины всех $N_i$, и рассмотрим $\GL$-модель на цепи
$w_0Rw_1R\dots Rw_k$ (транзитивной и иррефлексивной), в которой $A$
ложна в $w_0$ и истинна в остальных точках, а все прочие атомы
постоянны. Тогда $w_0\Vdash\square A$ и $w_0\nVdash A$. С другой
стороны, в точках $w_1,\dots,w_k$ все атомы постоянны, и индукцией по
модальной глубине $d$ произвольной $\square$-формулы $X$ (в ней $A$
может стоять и вне модальностей) проверяется, что во всех точках
$w_1,\dots,w_{k-d}$ истинность $X$ одна и та же. База ($d=0$): $X$
--- булева комбинация атомов, постоянных в этих точках. Шаг индукции: пусть
$X=\square Y$ имеет глубину $d+1$ и $1\le i\le k-d-1$; преемники точки
$w_i$ --- это $w_{i+1},\dots,w_{k-d}$, в которых $Y$ по предположению
индукции имеет одно и то же значение, и $w_{k-d+1},\dots,w_k$, общие
для всех таких $w_i$. Аналогично, для формулы $N_i$, в которой $A$
встречается только под $\square$, истинность $N_i$ в $w_0$ совпадает с
истинностью в $w_1$. Действительно, $N_i$ --- булева комбинация атомов,
отличных от $A$, и формул вида $\square X$; атомы, отличные от $A$,
в $w_0$ и $w_1$ имеют одинаковые значения, а $\square X$ в $w_0$
отличается от $\square X$ в $w_1$ лишь значением $X$ в $w_1$,
которое совпадает со значением $X$ в $w_2$ по доказанному (глубина
$X$ не превосходит $k-2$). Значит, из $w_0\Vdash\square N_i$ следует $w_1\Vdash N_i$, а с
ним и $w_0\Vdash N_i$; все посылки $\square N_i\ato N_i$ в $w_0$
верны, а $\square A\ato A$ --- нет.
\end{proof}

Таким образом, схемы (A8) и (A9) в $\HCE$ независимы хотя бы в одну
сторону; обратная независимость --- (A9) от $\HCI+$(A8) --- остаётся
открытой (\S\ref{ru:sec:todo}).

\paragraph{Схемы, не ставшие аксиомами.} Аксиомы исчислений отбирались по
общезначимости, и некоторые естественные кандидаты на этом пути были
отвергнуты. Два из них стоит назвать --- как следы поиска аксиом:
распределение $\bang$ по дизъюнкции (и его $\square$-версия) и попытка
обратить схему (A3), получив из преобразования выводов вывод импликации
между доказуемостями.

\begin{theorem}\label{ru:prop:notaxioms}
Пусть $A$, $B$ --- переменные сорта высказывание. Следующие формулы не
общезначимы ни в одной из трёх семантик, а потому (теорема~\ref{ru:th:sound})
не выводимы ни в одном из трёх исчислений:
\begin{enumerate}[leftmargin=2.2em,itemsep=2pt,label=\textup{(\arabic*)}]
\item 
$\square(A\aor B)\ato(\square A\aor  \square B)$;
\item $(\bang A\pto\bang B)\pto\bang(\square A\ato\square B)$.
\end{enumerate}
\end{theorem}

\begin{proof}
(1) 
Пусть $\sigma$ --- россеровское предложение для внутренней теории:
из непротиворечивости теории следует, что ни $\sigma$, ни $\aneg\sigma$ в
ней не доказуемы, причём это рассуждение формализуется в самой теории
(см., напр.,~\cite[гл.~2]{ru:Boolos93}). Возьмём $A^*:=\sigma$,
$B^*:=\aneg\sigma$. Дизъюнкция $A^*\aor B^*$ --- тавтология и потому
доказуема, а ни $A^*$, ни $B^*$ не доказуемы.
В семантике $\sE$ уже всё ясно:
посылка истинна, а заключение ложно. Во
внутренних семантиках надо опровергнуть выводимость во внешней теории
$T$ (она же внутренняя): посылка --- истинное $\Sigma_1$-предложение и
потому доказуема, так что из доказуемости высказывания следовало бы
$T\vdash\Prov(\sigma)\aor\Prov(\aneg\sigma)$, а вместе с
формализованной теоремой Россера
\[
  T\vdash\Con\ato\aneg\Prov(\sigma)\aand\aneg\Prov(\aneg\sigma)
\]
это дало бы $T\vdash\aneg\Con$. Последнее невозможно при
$\Sigma_1$-корректности $T$: $\aneg\Con$ --- ложное $\Sigma_1$-предложение.

(2) Возьмём $A^*:=\aneg\Con$ (то есть $\square\bot$) и $B^*:=\arbot $;
обозначим через $L$ перевод высказывания $\square A\ato\square B$, то
есть формулу $\square\square\bot\ato\square\bot$, а через $T$ --- внешнюю
теорию. Здесь понадобится $\Sigma_1$-корректность $T$
(замечание~\ref{ru:rem:assumptions});
поскольку внутренняя теория содержится во внешней, она тогда тоже
$\Sigma_1$-корректна. Надо показать, что $T\nvdash\varphi^*(n)$ при
любом $n$, где 
$$
\varphi=(\bang\square\bot\pto\bang\bot)\ato \bang(\square\square\bot\pto\square\bot)
$$
 (для $\sE$ это означает,
что значение задачи пусто). Формула $\varphi^*(n)$ имеет
вид $$\forall v\,\bigl(P(v)\ato\Prf(U_n(v),L)\bigr),$$ где $P(v)$ ---
перевод посылки $\bang A\pto\bang B$, то есть
$$\forall w\,\bigl(\Prf(w,\aneg\Con)\ato\Prf(U_v(w),\arbot )\bigr).
$$ 
Допустим
$T\vdash\varphi^*(n)$. Рассуждаем в $T$: если $\aneg\square\square\bot$,
то есть выводов формулы $\aneg\Con$ нет, то $P(v)$ выполнено при всяком
$v$, так как ее посылка ложна, поэтому $\Prf(U_n(0),L)$ и тем более $\square L$. Итак,
$T\vdash\aneg\square\square\bot\ato\square L$. С другой стороны,
формализованная теорема Лёба $\square(\square C\ato C)\ato\square C$ ---
теорема $\PA$, а потому и $T$ --- при
$C:=\square\bot$ даёт $T\vdash\square L\ato\square\square\bot$. Вместе
получаем $T\vdash\square\square\bot$. Это $\Sigma_1$-предложение, и по
$\Sigma_1$-корректности $T$ оно истинно: внутренняя теория доказывает
$\square\bot$, то есть $\aneg\Con$, --- вопреки её
$\Sigma_1$-корректности. 
\end{proof}


\subsection{Выводимые задачи с модальностями: что известно}\label{ru:sec:taskder}

Полного описания семейства выводимых задач у нас нет --- как нет и описания
семейства общезначимых (\S\ref{ru:sec:taskval}); неизвестно даже, полна ли
система $\HCgl$ на задачах (\S\ref{ru:sec:todo}). Здесь собрано то, что известно о
задачах специального вида --- о $\triangledown$-фрагменте.

Для $\triangledown=\bang\quest$ выводимы схемы K, 4 и Лёба. Первая схема
$\triangledown(\varphi\pto\psi)\pto(\triangledown\varphi\pto\triangledown\psi)$
получается из (A3) при $A:=\quest\varphi$, $B:=\quest\psi$, поскольку
$\bang\quest(\varphi\pto\psi)\peq\bang(\quest\varphi\ato\quest\psi)$ по (A4a) и
лемме~\ref{ru:lem:repl}. Вторая схема $\triangledown\varphi\pto\triangledown\triangledown\varphi$
--- это (A5) при $\sigma:=\bang\quest\varphi$. Наконец,
$\triangledown(\triangledown\varphi\pto\varphi)\pto\triangledown\varphi$
получается из (A6) при $A:=\quest\varphi$, поскольку
$\bang\quest(\triangledown\varphi\pto\varphi)\peq\bang(\square\quest\varphi\ato\quest\varphi)$
по (A4a) и лемме~\ref{ru:lem:repl}. Правило необходимости для
$\triangledown$ --- это (PtoA), затем (AtoP). Не выводима ко-рефлексия
$\varphi\pto\triangledown\varphi$: она не общезначима
(предложение~\ref{ru:prop:hce}(4)). Обратная схема
$\triangledown\triangledown\varphi\pto\triangledown\varphi$
$\sE$-общезначима, но не $\sI$-общезначима и, стало быть, в $\HCgl$ не
выводима: она потребовала бы $\PA\vdash\square\square A\ato\square A$.

Отдельно отметим уже доказанное: схема $\pneg\pneg\triangledown\varphi\pto\triangledown
\varphi$ выводима как частный случай (A7), а ко-рефлексия
$\varphi\pto\triangledown\varphi$ выводима на всём классе $\Sigma$-задач и
только на нём ($\triangledown\varphi=\bang\quest\varphi$ ---
$\Sigma$-задача по определению~\ref{ru:def:sigma}; аксиома (A5) и
предложение~\ref{ru:prop:hce}(4)). Главное же --- следующее.

\begin{theorem}\label{ru:prop:nabla}
Если булева свёртка задачи $\chi$ по схемам (A4) --- классическая
тавтология, то задача $\triangledown\chi$ выводима в $\HCgl$. В частности,
выводимы
\[
  \triangledown(\varphi\por\pneg\varphi),\qquad
  \triangledown(\pneg\pneg\varphi\pto\varphi),\qquad
  \triangledown\bigl(((\varphi\pto\psi)\pto\varphi)\pto\varphi\bigr).
\]
\end{theorem}

\begin{proof}
По схемам (A4) высказывание $\quest\chi$ выводимо эквивалентно
своей булевой свёртке (предложение~\ref{ru:prop:nf}), а та ---
классическая тавтология, схема (A1); значит $\quest\chi$ --- теорема
$\HCgl$, и, применяя (AtoP), получаем задачу
$\bang\quest\chi=\triangledown\chi$.
\end{proof}

\remove{
\begin{remark}\label{ru:rem:nabla}
Причина расхождения структурная, а не случайная. Модальность
$\triangledown=\bang\quest$ пропускает свой аргумент через $\quest$, а $\quest$ по (A4)
--- булев гомоморфизм; поэтому под знаком $\triangledown$ интуиционистская
природа задач стирается, и туда попадает вся классическая логика.
\end{remark}
}

\section{Сравнение с близкими логиками}\label{ru:sec:hc}

Из систем, перечисленных во введении, только $\HC$ имеет тот же синтаксис, и
только с ней сравнение требует отдельной работы; ей посвящена большая часть
раздела. Остальные сравнения опираются на уже доказанное и здесь лишь
названы --- в конце раздела.

\subsection{Логика $\HC$ Мелихова}

Напомним аксиоматику $\HC$: классические схемы (высказывания),
интуиционистские схемы (задачи), схемы
\[
  \bang(A\ato B)\pto(\bang A\pto\bang B),\quad
  \quest(\varphi\pto\psi)\ato(\quest\varphi\ato\quest\psi),\quad
  \quest\bang A\ato A,\quad
  \varphi\pto\bang\quest\varphi,\quad
  \pneg\bang\bot
\]
и правила (MP), (AtoP), (PtoA).

Сравнение двусторонне: какие теоремы $\HC$ общезначимы в нашей семантике и
какие из наших аксиом выводимы в $\HC$. Первую половину дают
предложение~\ref{ru:prop:hce} и следующая теорема.

\begin{proposition}\label{ru:prop:hcval}
Задача $\pneg\bang\bot$ $\sE$-общезначима, но не $\sI$- и не
$\sB$-общезначима.
\end{proposition}

\begin{proof}
В $\sE$: эта задача выводима в $\HCE$ (Пример 1(4)), следовательно, общезначима.
Во внутренних
семантиках: формула $(\pneg\bang\bot)^*(x)$ доказуемо эквивалентна $\Con$
--- снова вторая теорема Гёделя.
\end{proof}

Первые две выписанные схемы --- это наши (A3) и
левая половина (A4a); они $\sI$-общезначимы (теорема~\ref{ru:th:sound}).
Рефлексия $\quest\bang A\ato A$ и схема $\pneg\bang\bot$ только $\sE$-общезначимы
(предложение~\ref{ru:prop:hce}(1),(2) и теорема~\ref{ru:prop:hcval}), а ко-рефлексия
$\varphi\pto\bang\quest\varphi$ не общезначима вовсе
(предложение~\ref{ru:prop:hce}(4)). Правила: (MP) и (PtoA) сохраняют и $\sI$-, и
$\sE$-общезначимость (лемма~\ref{ru:lem:closure} и теорема~\ref{ru:th:sound}), а (AtoP)
--- только $\sI$-общезначимость: по предложению~\ref{ru:prop:hce}(3) применение (AtoP) к
$\sE$-общезначимому высказыванию $\square\bot\ato\bot$ даёт необщезначимую
задачу. Не допустимо в нашей семантике и правило $\quest\varphi/\varphi$,
обратное к (PtoA): высказывание $\quest(\alpha\por\pneg\alpha)$ общезначимо,
поскольку $\val{\alpha\por\pneg\alpha}\neq\varnothing$ при всякой
интерпретации, а задача $\alpha\por\pneg\alpha$ не общезначима
(предложение~\ref{ru:prop:noem}).

Итог удобно назвать заранее: общая часть двух исчислений --- это \emph{ядро}
$\HCzero$, аксиомы (A1), (A2), (A3), левая половина (A4a), (A4b)--(A4d) и (A5) при всех
трёх правилах; всё, что сверх ядра, расходится, причём в обе стороны.

Осталась вторая половина --- что из нашей аксиоматики выводимо в $\HC$.
Начнём с группы (A4): она разбирается поштучно.

\begin{remark}[группа (A4) в $\HC$; остаточная классичность]\label{ru:rem:qhom}
Схема (A4d) выводима в $\HC$ целиком. Направление $\bot\ato\quest\varnothing$
--- классическая логика; для обратного берём интуиционистскую теорему
$\varnothing\pto\bang\bot$, применяем (PtoA) и левую половину (A4a), что даёт
$\quest\varnothing\ato\square\bot$, и остаётся рефлексия $\square\bot\ato\bot$.

Схема (A4c) тоже выводима целиком. Лёгкое направление
$(\quest\varphi\aor\quest\psi)\ato\quest(\varphi\por\psi)$ получается из вложений
$\varphi\pto\varphi\por\psi$ тем же приёмом: (PtoA), затем левая половина
(A4a). Содержательное направление проходит через ко-рефлексию: из
$\varphi\pto\bang\quest\varphi$ и $\bang$-версии тавтологии
$\quest\varphi\ato(\quest\varphi\aor\quest\psi)$ (правило (AtoP)) по (A3) получаем
$\varphi\pto\bang(\quest\varphi\aor\quest\psi)$, так же для $\psi$, а значит и
$\varphi\por\psi\pto\bang(\quest\varphi\aor\quest\psi)$; применяя (PtoA) и левую
половину (A4a), получаем
$\quest(\varphi\por\psi)\ato\square(\quest\varphi\aor\quest\psi)$, и остаётся
рефлексия. Содержательно это ожидаемо: в $\HC$ оператор $\quest$ --- левый
сопряжённый к $\bang$, а левые сопряжённые сохраняют дизъюнкции.

У схемы (A4b) в $\HC$ выводима половина
$\quest(\varphi\pand\psi)\ato(\quest\varphi\aand\quest\psi)$ --- из проекций, снова
(PtoA) плюс левая половина (A4a). В обратную сторону: к интуиционистской аксиоме
$\phi\ato(\psi\ato(\phi\aand\psi))$ применяем левую половину (A4a), получая
$\quest\phi\ato\quest(\psi\pto(\phi\pand\psi))$. С помощью аксиомы 
$\quest(\psi\ato(\phi\aand\psi))\ato(\quest\psi\ato\quest(\phi\aand\psi))$
и MP выводим $\quest\phi\ato(\quest\psi\ato\quest(\phi\aand\psi))$, что классически то же самое,
что искомая формула.


Существенно новыми, таким образом, остаются правая половина (A4a),
\[
  (\quest\varphi\ato\quest\psi)\ato\quest(\varphi\pto\psi),
\]
и производная от неё $\aneg\quest\varphi\ato\quest\pneg\varphi$ (получается из
(A4a) и (A4d)). В $\HC$ их нет: в шкалах Крипке в посылке возникает квантор
по будущим мирам. [TODO Добавить ссылку на работу Мелихова, в которой это доказано.]
(Конец замечания.)
\end{remark}

\begin{proposition}\label{ru:prop:incl}\phantom{l}
\begin{enumerate}[leftmargin=2.2em,itemsep=2pt,label=\textup{(\arabic*)}]
\item $\HCzero\subseteq\HC$.
\item $\HC\nsubseteq\HCgls$; тем более $\HC\nsubseteq\HCgl$.
\item $\HCglm\nsubseteq\HC$; тем более $\HCgl\nsubseteq\HC$ и
  $\HCgls\nsubseteq\HC$.
\item Исчисление $\HC+$(A6) противоречиво.
\end{enumerate}
\end{proposition}

\begin{proof}
(1) Схемы (A1), (A2), (A3), левая половина (A4a) и правила --- это в точности
часть аксиоматики $\HC$, схемы (A4b)--(A4d) выводимы в $\HC$
(замечание~\ref{ru:rem:qhom}), а (A5) --- частный случай ко-рефлексии
$\varphi\pto\bang\quest\varphi$ (для $\Sigma$-задач $\varphi$).

(2) В $\HC$ выводима ко-рефлексия, которая по предложению~\ref{ru:prop:hce}(4)
не $\sE$-общезначима, а $\HCE$ корректна относительно $\sE$ (теорема~\ref{ru:th:sound});
теоремы $\HCgl$ входят в $\HCgls$.

(3) Уже правая половина (A4a) в $\HC$ не выводима (замечание~\ref{ru:rem:qhom});
не выводима и (A7): её частный случай $\pneg\pneg\bang A\pto\bang A$ --- марковский
принцип для $\bang$, а задачи в $\HC$ интуиционистские, без марковских добавок.

(4) В $\HC$ есть рефлексия $\square A\ato A$ и правило необходимости для
$\square$ (применяем (AtoP), затем (PtoA)); из (A6) мы получили бы
$\square(\square A\ato A)$, затем $\square A$ при любом $A$, а при $A=\bot$
--- противоречие с рефлексией. Между тем $\HC$ непротиворечива: она корректна
относительно шкал Крипке с отмеченными мирами~\cite{ru:Onoprienko20}.
\end{proof}

Итак, внутри $\HC$ из нашего исчисления живёт ровно ядро, а расхождение идёт
по нескольким независимым линиям: классичность модальности $\quest$ (правая
половина (A4a) с производной $\aneg\quest\varphi\ato\quest\pneg\varphi$),
марковский принцип (A7) и лёбовская линия --- схема (A6), которая не просто
отсутствует в $\HC$, а несовместима с ней. 
В обратную сторону линия одна ---
ко-рефлексия.
Существенно, что вывод формулы (L) проходит уже в ядре
(теорема~\ref{ru:th:lob}): поэтому следствие~\ref{ru:cor:hcfp} о
несовместимости $\HC$ с неподвижными точками опирается только на общую часть
двух исчислений.

Несовместимость $\HC$ с чтением $\bang$ как формальной доказуемости
видна и без неподвижных точек. В $\HC$ выводима задача
$\bang(\square\bot\ato\bot)$: применяем правило
(AtoP) к частному случаю рефлексии $\square\bot\ato\bot$.
Если решения задачи $\bang A$ --- выводы $A$ в непротиворечивой теории
$T$, доказывающей условия выводимости Гильберта--Бернайса, 
то решение этой задачи --- вывод
$\mathrm{Con}(T)$ в $T$, которого нет по второй теореме Гёделя. Это
двусортный вариант наблюдения Гёделя~\cite{ru:Goedel33}: $\square$-фрагмент
$\HC$ есть $\mathrm{S4}$~\cite{ru:Melikhov13}, а в $\mathrm{S4}$ выводимо
$\square\aneg\square\bot$, поэтому $\square$ в $\mathrm{S4}$ нельзя
читать как формальную доказуемость. Исключается именно доказуемость в
одной фиксированной формальной теории; семантика с явными
доказательствами в духе логики доказательств
Артёмова~\cite{ru:Artemov01} под этот довод не подпадает: там
$\aneg(t{:}\bot)$ выводимо для каждого терма $t$, но утверждения
<<никакого доказательства $\bot$ нет>> в языке нет.
\remove{
Наконец, всё это сравнение держится на выбранном понятии общезначимости. При
слабом варианте~(а) замечания~\ref{ru:rem:three} фрагмент задач классичен, схему
(A2) пришлось бы усиливать до классической, и связь с интуиционистской частью
$\HC$ пропала бы вовсе; содержательным сопоставление двух исчислений делает
именно требование предъявлять решение вычислимо.
}

\subsubsection*{Неподвижные точки несовместимы с $\HC$}

Из вывода формулы (L) в примере~\ref{ru:ex:proofs}(5) можно извлечь больше, чем
корректность (A6). Внимательный осмотр показывает, что вывод проходит уже в
общей части двух исчислений --- в ядре $\HCzero$.

\begin{theorem}\label{ru:th:lob}
Уже в ядре $\HCzero$ (а значит, и в $\HCglm$, и в $\HC$) выводима задача
\[
  \bang\bigl(P\aeq(\square P\ato S)\bigr)\ \pto\
  \bigl(\bang(\square S\ato S)\pto\bang S\bigr)
\]
при любых высказываниях $P$ и $S$.
\end{theorem}

\begin{proof}
Это пункт~(5) примера~\ref{ru:ex:proofs}. Достаточно заметить, что в нём
используются только (A1), (A2), (A3), (A5), правила (MP), (AtoP), (PtoA),
лемма о дедукции и лемма о переводе выводов
(предложение~\ref{ru:prop:derived}(2),(3)) и левая половина (A4a) --- последняя
лишь в пунктах~(1) и~(2) того же примера, при выводе $\square$-версий (A3)
и (A5). Ни правая половина (A4a), ни (A6), ни (A7) не нужны.
\end{proof}
\remove{
\begin{remark}[как устроен вывод и почему посылку удаётся снять]
Посылкой служит не сама эквивалентность $P\aeq(\square P\ato S)$, а
задача $\bang(\dots)$ <<найти её вывод>>: (A3) умеет проносить $\bang$ через
импликацию, и потому работать можно только с выводами.

Сам вывод устроен в два этажа. Внизу, на сорте высказывание, всё делается
классически: шаги~(1)--(6) пункта~(5) примера~\ref{ru:ex:proofs} --- это обычный
модальный вывод, использующий лишь $\square$-версии (A3) и (A5). Наверх, под
знак~$\bang$, он переносится целиком леммой о переводе выводов
(предложение~\ref{ru:prop:derived}(3)); правило (AtoP) при этом применяется
только к уже доказанным формулам, а к гипотезам --- никогда, и потому лемму о
дедукции в конце применить можно.

Схема (A5) работает дважды, и оба раза для того, чтобы приписать бокс под
знаком~$\bang$: она даёт $\bang\square P$ из $\bang P$ и вторую гипотезу из
$\bang(P\ato(\square P\ato S))$. Она занимает то самое место, которое в
стандартном доказательстве теоремы Лёба занимает схема
$\square A\ato\square\square A$.
(Конец замечания.)
\end{remark}
}
\begin{corollary}[неподвижные точки несовместимы с $\HC$]\label{ru:cor:hcfp}
Пусть к $\HC$ добавлена неподвижная точка для лжи: то есть, аксиома 
$\bang\bigl(P\aeq(\square P\ato\bot)\bigr)$,
где $P$ --- переменная сорта высказывание. 
Тогда полученная
система противоречива --- в ней выводима задача $\varnothing$, а значит и
любая задача.
\end{corollary}

\begin{proof}
По предложению~\ref{ru:prop:incl}(1) ядро $\HCzero$ содержится в
$\HC$, поэтому теорема~\ref{ru:th:lob} даёт в $\HC$ импликацию
\[
  \bang\bigl(P\aeq(\square P\ato\bot)\bigr)\pto
  \bigl(\bang(\square\bot\ato\bot)\pto\bang\bot\bigr) .
\]
Посылка выводима по условию. Далее, $\square\bot\ato\bot$ --- частный случай
рефлексии, то есть теорема $\HC$; по (AtoP) выводимо
$\bang(\square\bot\ato\bot)$, и двумя применениями (MP) получаем $\bang\bot$.
Вместе с аксиомой $\pneg\bang\bot$ это даёт $\varnothing$.
\end{proof}
\remove{
Это синтаксическая форма развилки Лёба: $\HC$ несовместима не с какой-то одной
арифметической конструкцией, а с самим наличием неподвижных точек --- а в
арифметической интерпретации они есть всегда. Система $\HCgl$, наоборот,
неподвижные точки спокойно выдерживает и даже извлекает из них схему Лёба.
}

\subsection{Остальные системы}

Здесь только перечень, с указанием мест, где соответствующие утверждения
доказаны.

\emph{Логика доказуемости.} $\square$-фрагмент $\HCgl$ есть в точности $\GL$,
$\square$-фрагмент $\HCgls$ --- в точности $\GLS$; вместе с теоремами
Соловея~\cite{ru:Solovay76} отсюда следует полнота на стороне высказываний:
следствие~\ref{ru:cor:solovay}.

\emph{Интуиционистские логики доказуемости.} Напомним
определения~\cite{ru:VZ19,ru:DvdGM24}: $\mathrm{iK}$ --- интуиционистские
тавтологии плюс схема $\square(A\ato B)\ato(\square A\ato\square B)$ при
правилах (MP) и необходимости; $\mathrm{iGL}=\mathrm{iK}+{}$схема Лёба
$\square(\square A\ato A)\ato\square A$; наконец
$\mathrm{iGLC}=\mathrm{iGL}+{}$\emph{принцип полноты} $A\ato\square A$.
Семантика Крипке для $\mathrm{iGL}$ --- бирелационные шкалы
$\langle W,\preccurlyeq,\sqsubset\rangle$, где $\preccurlyeq$ ---
частичный порядок (интуиционистское отношение), а $\sqsubset$ ---
транзитивное обратно фундированное модальное отношение, согласованное
с порядком в обе стороны: из $x\preccurlyeq y\sqsubset z$ следует
$x\sqsubset z$, и из $x\sqsubset y\preccurlyeq z$ следует $x\sqsubset z$
(такие шкалы Иемхоф называет \emph{brilliant}). Оценки переменных
монотонны по $\preccurlyeq$, а $\square A$ вынуждено в точке $x$, если
$A$ вынуждено во всех её $\sqsubset$-преемниках. Логика $\mathrm{iGL}$
полна относительно таких шкал~\cite[предложение~9]{ru:Iemhoff99}; см.\
также~\cite{ru:DvdGM24}. Сравнивать с ними будем
$\triangledown$-фрагмент $\HCgl$ (\S\ref{ru:sec:taskder}), читая $\square$ как
$\triangledown$.

\begin{proposition}\label{ru:prop:igl}
Логика $\mathrm{iGL}$ содержится в $\triangledown$-фрагменте $\HCgl$,
причём строго: всякая теорема $\mathrm{iGL}$, прочитанная с
$\square:=\triangledown$, выводима в $\HCgl$, а задача
$\triangledown(\varphi\por\pneg\varphi)$ выводима в $\HCgl$, тогда как
формула $\square(p\aor\aneg p)$ в $\mathrm{iGL}$ не выводима.
\end{proposition}

\begin{proof}
Включение: схема K, схема Лёба и правило необходимости для
$\triangledown$ выведены в \S\ref{ru:sec:taskder}. Строгость даёт задача
$\triangledown(\varphi\por\pneg\varphi)$ (теорема~\ref{ru:prop:nabla}): покажем,
что $\square(p\aor\aneg p)$ не выводима в $\mathrm{iGL}$.

Возьмём шкалу $\langle W,\preccurlyeq,\sqsubset\rangle$ с $W=\{w,v,u\}$, где
$\preccurlyeq$ --- рефлексивное отношение вместе с $v\preccurlyeq u$, а
$\sqsubset$ состоит из пар $(w,v)$ и $(w,u)$. Это шкала для
$\mathrm{iGL}$: отношение $\sqsubset$ транзитивно и обратно фундировано
(цепочек длины два нет), а оба условия согласования проверяются
перебором --- единственный нетривиальный случай $w\sqsubset v\preccurlyeq u$
даёт $w\sqsubset u$. Пусть переменная $p$ истинна только в $u$;
оценка монотонна, поскольку $u$ максимальна. Тогда $v\nVdash p$ и
$v\nVdash\aneg p$ (так как $v\preccurlyeq u$ и $u\Vdash p$), значит
$v\nVdash p\aor\aneg p$, а из $w\sqsubset v$ следует
$w\nVdash\square(p\aor\aneg p)$. По полноте $\mathrm{iGL}$ относительно
таких шкал формула $\square(p\aor\aneg p)$ в $\mathrm{iGL}$ не выводима.
\end{proof}

\begin{remark}\label{ru:rem:igl}
Причина расхождения --- 
под знаком
$\triangledown$ оказывается вся классическая логика, чего никакая
интуиционистская модальная логика дать не может, и $\mathrm{iGL}$ --- не
исключение. Это же показывает, где искать точную характеризацию:
$\triangledown$-фрагмент --- это, по-видимому, не какая-нибудь
интуиционистская модальная логика, а $\mathrm{iGL}$ вместе со всеми
формулами $\triangledown\chi$, у которых булева свёртка $\chi$ ---
тавтология. Совпадает ли фрагмент в точности с этим пополнением ---
открытый вопрос (\S\ref{ru:sec:todo}).

Принцип полноты логики $\mathrm{iGLC}$ --- это ко-рефлексия,
опровергнутая предложением~\ref{ru:prop:hce}(4); но ограничение полноты на
$\Sigma_1$, принятое в~\cite{ru:VZ19}, отвечает нашей ко-рефлексии на
$\Sigma$-задачах (аксиома (A5)). Наконец, в логике $\HC$ Мелихова для
$\triangledown$ верны схема K, схема~4 и ко-рефлексия, но не схема Лёба
--- это интуиционистская модальная логика $\mathrm{H4}$, она же
$\mathrm{IEL}^+$ Артёмова и Протопопеску~\cite{ru:Protopopescu16}; у нас
ко-рефлексию вытесняет схема Лёба, и $\triangledown$-фрагмент $\HCgl$ ---
<<лёбовский>> аналог $\mathrm{H4}$ над интуиционистской логикой задач.
(Конец замечания.)
\end{remark}

\emph{Логики реализуемости.} Наше понятие общезначимости --- это $\Reff$ из
классификации $\Rcnst\subseteq\Reff\subseteq\R$~\cite{ru:VSSCh}; три уровня из
замечания~\ref{ru:rem:three} суть в точности эти три логики. Отсюда результаты
\S\ref{ru:sec:modalfreeval} о фрагменте задач: позитивная часть интуиционистская
(теорема~\ref{ru:th:posfrag}), а сам фрагмент лежит строго между $\Int$ и логикой
слабого закона исключённого третьего (следствие~\ref{ru:cor:sandwich}).

\section{Открытые вопросы}\label{ru:sec:todo}

\begin{enumerate}[leftmargin=1.8em,itemsep=3pt]
\item \textbf{Что за логика получилась на задачах?} Известно, что
  $\Int\subsetneq\Lp_Z\subsetneq\J$ для всех трёх семантик
  (следствие~\ref{ru:cor:sandwich}) и что
  позитивный фрагмент интуиционистский (теорема~\ref{ru:th:posfrag}), так что вся
  разница с $\Int$ живёт на формулах с $\varnothing$. Слева нас отделяет
  формула Роуза (предложение~\ref{ru:cor:rose}), справа --- принцип Скотта
  (теорема~\ref{ru:prop:scott}). Значит
  схемы (A2) на стороне задач заведомо не хватает, и $\HCI$ на задачах
  неполна (следствие~\ref{ru:cor:incomplete}); закрывать брешь готовой
  схемой Роуза мы не стали (\S\ref{ru:sec:modalfreeval}), а как выглядит
  полное исчисление --- неизвестно. Постановка
  близка к пропозициональной логике реализуемости Роуза и
  Плиско~\cite{ru:Rose53,ru:Plisko09} --- переменные там тоже заменяются
  арифметическими формулами, --- и оттуда же стоит брать следующие пробные
  схемы.

  Стоит отдавать себе отчёт в том, насколько разные это две стороны. На
  стороне высказываний мы имеем дело с известным объектом: фрагменты суть в
  точности $\GL$ и $\GLS$, полнота получается
  из теорем Соловея (следствие~\ref{ru:cor:solovay}), и эта сторона (в
  $\Pi$-свободном языке) уже не изменится, какие бы общезначимые
  схемы о задачах мы ни добавляли: новые аксиомы дают лишь
  общезначимые теоремы, а все общезначимые высказывания уже выводимы. На стороне задач известного объекта нет: логики
  $\Lp_Z$ лежат строго между $\Int$ и $\J$ и ни с одной из них не
  совпадают. Безмодальный фрагмент исчисления $\HCI$ есть в точности
  $\Int$, и само исчисление читается как склейка двух готовых систем ---
  $\Int$ на задачах и $\GL$ на высказываниях, соединённых схемами
  (A3)--(A5) и (A7). Семантика же склейкой не является: с одной стороны
  стоит классическая логика доказуемости, с другой --- промежуточная
  логика, которую ещё предстоит опознать.

  Не исключено, что опознать её в виде исчисления не удастся вовсе. Ближайшие
  родственники логик $\Lp_Z$ --- логики реализуемости, а там дело обстоит так: для
  логики реализуемых пропозициональных формул неизвестна ни конечная, ни даже
  рекурсивная аксиоматизируемость, а логика финитных задач Медведева,
  как недавно показано, $\Pi^0_1$-полна и, следовательно, не перечислима:
  никакого перечислимого полного исчисления для неё
  нет~\cite{ru:Pawlowski26,ru:AlmeidaKnudstorp26}. Наше понятие
  общезначимости --- это $\Reff$, то есть ровно из этого семейства. Поэтому и
  то, что общезначимая схема Роуза (предложение~\ref{ru:cor:rose}) есть серия
  по числу дизъюнктов, скорее свойство предмета, чем случайность. И
  поэтому же отбор аксиом по общезначимости не
  может служить признаком полноты: новые схемы могут перестать находиться не
  оттого, что список полон.

\item \textbf{Насколько сильны $\bang$-формы.} Из (A5) и (A6) их
  $\square$-версии выводятся (пример~\ref{ru:ex:proofs}(2),(3)); верно ли
  обратное? $\square$-версии говорят лишь о существовании выводов, $\bang$-версии
  --- о вычислимом преобразовании выводов, так что обратный переход
  неочевиден. От ответа зависит, насколько содержательно
  теорема~\ref{ru:th:lob} по сравнению с его модальным аналогом
  (приложение~\ref{ru:sec:novelty}).

\item \textbf{Нужна ли (A5) для вывода теоремы Лёба}, то есть выводима ли
  (A6) из $\HCglm$ без (A5) (точнее, без её частного случая
  $\bang A\pto\bang\square A$) вместе с посылкой
  $\bang(P\aeq(\square P\ato S))$? В доказательстве
  теоремы~\ref{ru:th:sound} схема (A5) занимает место схемы
  $\square A\ato\square\square A$, без которой рассуждение с неподвижной точкой
  не проходит и в модальном случае; но доказательства невыводимости нет.
  (Что (A5) не выводима из остальных аксиом $\HCgls$ вместе с правилом
  <<из $\bang A\pto\bang B$ вывести $\bang(\square A\ato\square B)$>> и
  всеми задачами вида $\bang T$ для классических тавтологий $T$,
  видно из рефлексивной нетранзитивной трёхмировой шкалы Крипке, в которой
  $\bang$ читается как $\square$, а $\quest$ --- тождественно.)

\item \textbf{Полнота на задачах.} На стороне высказываний вопрос закрыт
  (следствие~\ref{ru:cor:solovay}), а предложение~\ref{ru:prop:nf} показывает, что
  схемами про $\quest$ система (A4) исчерпана: $\quest$ устраняется до $\quest\alpha$ и
  $\square B$. Значит всё нерешённое собралось на стороне задач. Известные
  бреши --- формула Роуза (следствие~\ref{ru:cor:incomplete}) и, возможно,
  закон Пирса с квантором (предложение~\ref{ru:prop:pienc}(4),
  замечание~\ref{ru:rem:peirce}). Ближайшая проверка:
  не осталось ли общезначимых схем
  того же рода, где решателю нужно угадать метку дизъюнкции, а список
  кандидатов не конечен; см. раздел~\ref{ru:sec:whyclasses}.

\item \textbf{Допустимо ли в $\HCB$ правило <<из $\bang A$ вывести
  $A$>>?} Это правило сохраняет $\sB$-общезначимость без всяких
  допущений: если $\ZF\vdash\mathrm{Prf}(n,\gn{\forall\vec X A^*})$, то
  при непротиворечивости $\ZF$ эта $\Delta_0$-формула истинна и
  $\ZF\vdash\forall\vec X A^*$, а при противоречивости общезначимо всё.
  В $\HCE$ правило производное ((PtoA), затем (A8)), в $\HCI$ допустимо:
  из $\provgl\bang A$ по теореме~\ref{ru:th:sound} следует
  $\sI$-общезначимость $A$, а по следствию~\ref{ru:cor:solovay} --- его
  выводимость. Верно ли, что из $\HCB\vdash\bang A$ следует
  $\HCB\vdash A$? Новых задач добавление правила к $\HCB$ не даёт
  (индукцией по выводу: для всякого высказывания $A$, выводимого с этим
  правилом, задача $\bang A$ выводима в $\HCB$), так что вопрос касается
  лишь высказываний с квантором $\Pi$. Отрицательный ответ означал бы,
  что в $\HCB$ не хватает какого-то закона: под знаком $\bang$
  исчисление выводит больше, чем снаружи.

\remove{
\item \textbf{Внутренний уровень:} совпадают ли внутренняя и внешняя
  общезначимость на каких-нибудь естественных фрагментах (для задач вида
  $\bang A$ совпадают, замечание~\ref{ru:rem:whycomp}).
}

\item \textbf{Безмодальный фрагмент $\HCgls$.} Совпадает ли множество
  безмодальных теорем-задач $\HCE$ с таковым для $\HCI$? Известно лишь
  $\HCI\subseteq\HCE$ и что оба фрагмента лежат в $\Lp_{\sE}$
  (\S\ref{ru:sec:indep}); шкала с отмеченной максимальной точкой вопроса
  не решает --- схемы (A8) и (A9) в ней ложны.

\item \textbf{Независимость (A7) в $\HCE$.} Для $\HCI$
  вопрос решён: предложение~\ref{ru:prop:a7indep} (приложение~\ref{ru:sec:a7indep})
  даёт двухточечную модель Крипке, в
  которой верны (A1)--(A6) и все правила, а
  $\pneg\pneg\bang\bot\pto\bang\bot$ --- нет. Для $\HCgls$ моделями так действовать нельзя
  (замечание~\ref{ru:rem:a10fork}). Заметим, что синтаксически задача сводится к
  следующему: теоремы-задачи $\HCgls$ --- это замыкание по (MP) множества
  теорем-задач $\HCgl$ вместе со всеми частными случаями (A9), поскольку из
  высказывания в задачу без (AtoP) хода нет. Значит нужна структура, в
  которой все теоремы-задачи $\HCgl$ получают $1$ (но правило (AtoP) в целом
  не сохраняется) и все частные случаи (A9) --- тоже.

\item \textbf{Независимость (A9).} В обратную сторону вопрос решён:
  (A8) из (A9) не выводится (предложение~\ref{ru:prop:a8indep}), хотя
  (A9) даёт рефлексию для всех $\square$-формул с помодализованными
  переменными-высказываниями. Трудность с самой (A9) в том, что её
  $\quest$-образ есть частный случай (A8) (предложение~\ref{ru:prop:nf}):
  отличить (A9) от (A8) может только конструкция, различающая задачу и
  её $\quest$-образ.

\item \textbf{Как далеко тянется (A7).} Она принята для $\Sigma$-задач;
  приложение~\ref{ru:sec:whyclasses} показывает, что ни переменные, ни
  дизъюнкции произвольных задач добавить нельзя. Точный класс задач $\varphi$,
  для которых $\pneg\pneg\varphi\pto\varphi$ общезначима, не описан. Заметим, что
  этот класс зависит от синтаксиса: по $\gn{*}$ надо уметь построить
  перечисление $\val{\varphi}$, а для импликации оно неоткуда взяться.

\item \textbf{Совпадает ли $\triangledown$-фрагмент с $\mathrm{iGL}$?}
  (Определения логик $\mathrm{iK}$, $\mathrm{iGL}$, $\mathrm{iGLC}$ --- в
  \S\ref{ru:sec:hc}.) Вопрос о совпадении решён отрицательно: по
  предложению~\ref{ru:prop:igl} $\mathrm{iGL}$ содержится в
  $\triangledown$-фрагменте $\HCgl$ \emph{строго} --- уже
  $\triangledown(\varphi\por\pneg\varphi)$ выводима у нас и не выводима в
  $\mathrm{iGL}$. Осталось найти точную характеризацию; правдоподобная догадка
  выписана в замечании~\ref{ru:rem:igl}.

  А вот $\mathrm{iGLC}$ нашим фрагментом \emph{не} является: принцип полноты
  $\varphi\pto\triangledown\varphi$ --- это ко-рефлексия, и она опровергнута
  (предложение~\ref{ru:prop:hce}(4)). Тем не менее сходство с арифметической
  частью~\cite{ru:VZ19} есть, и содержательное: там полнота требуется не от всех
  формул, а от $\Sigma_1$ (условие $\square\triangle S\ato\square S$ при
  $S\in\Sigma_1$), --- у нас ровно так же ко-рефлексия верна в точности на
  $\Sigma$-задачах (аксиома (A5)). Стоит понять,
  случайное ли это совпадение.

  \emph{Оговорка.} Речь идёт о \emph{выводимых} формулах. Множество
  \emph{общезначимых} $\triangledown$-формул строго шире $\mathrm{iGL}$ уже
  потому, что общезначима формула Роуза (предложение~\ref{ru:cor:rose}), а она
  не выводима даже в $\Int$.

\remove{
\item \textbf{Не следует ли (A6) из теоремы реализации.} Схема (A6) говорит о
  вычислимом преобразовании выводов; в логике доказательств такому отвечает
  явный аналог $\GL$ с теоремой реализации~\cite{ru:Artemov01,ru:Sidon97,ru:Nogina14}.
  Стоит проверить, нельзя ли получить (A6) прямо оттуда, а не выводом из
  неподвижной точки, как при выводе (A6) в теореме~\ref{ru:th:sound}.
}

\item \textbf{Интуиционистский вариант семантики.} Раз расхождение с
  $\mathrm{iGL}$ создаёт булевость $\quest$,  
  напрашивается заменить $\PA$ на $\mathrm{HA}$. Полезно заранее видеть, что
  этого мало: у булевости два независимых источника. \emph{Внутренний} ---
  классичность $\PA$: правая половина (A4a) выведена разбором случаев
  <<либо $\aneg\exists u\,\varphi^*(u)$, либо $\exists v\,\psi^*(v)$>>, и
  $\mathrm{HA}$ такого разбора не даёт, так что здесь замена сработает.
  \emph{Внешний} --- двузначность самой стандартной модели: внешняя
  общезначимость определена через истинность в~$\N$, а
  $\val{\quest\varphi}=\Tru$ равносильно $\val{\varphi}\neq\varnothing$ ---
  двузначный факт о~$\N$ при любой теории под знаком~$\bang$. Чтобы снять и
  его, придётся считать значения не в стандартной модели, а в шкалах Крипке
  или в гейтинговозначных моделях; но тогда сорт высказываний перестанет быть
  классическим, то есть речь пойдёт уже не о $\HC$, а о её интуиционистской
  родне. Стоит иметь в виду и цену: при $\mathrm{HA}$ под знаком~$\bang$
  модальность $\square$ перестаёт отвечать $\GL$ и становится логикой
  доказуемости $\mathrm{HA}$~\cite{ru:AM19}, а из дизъюнктивного свойства
  $\mathrm{HA}$ добавляется схема $\square(A\aor B)\ato\square A\aor\square B$
  для замкнутых $A$, $B$.

\item \textbf{Определимые ящики.} Меняется ли что-нибудь, если ящики
  пробегают только арифметически определимые множества
  (замечание~\ref{ru:rem:semS})? На всех разобранных фактах --- нет;
  единственный кандидат на расхождение --- полная схема
  $(\Pi\mathrm{E})$, которую при определимых ящиках мы умеем
  обосновать лишь для $\Pi$-свободных $\theta$, и родственный ей
  вопрос о подстановке $\Pi$-содержащих задач
  (лемма~\ref{ru:lem:closurebs}). Сюда же ---
  консервативность $\Pi$-расширения над $\HCI$ на чистых
  бескванторных задачах: для чистых высказываний она доказана при
  $\Sigma_1$-корректности $\ZF$, а на нечистых формулах нарушается
  (замечание~\ref{ru:rem:qfree}). Существуют ли
  \emph{неустранимые} $\Pi$-атомы (определение~\ref{ru:def:elim})? Все
  встретившиеся нам приводятся к $\Pi$-свободным высказываниям
  (\S\ref{ru:sec:pinf}), но для одного из них --- закона Пирса с
  квантором по $\alpha$ --- приведение доказано только семантически, и
  выводимость задачи
  $\Pi\alpha(((\alpha\pto\beta)\pto\alpha)\pto\alpha)\pto\beta$ в
  $\HCB$ --- отдельный открытый вопрос (предложение~\ref{ru:prop:pienc}(4),
  замечание~\ref{ru:rem:peirce}). Для существования неустранимых атомов
  достаточно, чтобы безмодальная
  чёрноящичная логика не была рекурсивно перечислимой; её сложность
  --- отдельный вопрос. Наконец, \emph{возвращение к $\PA$}: для каких
  $\varphi$ перевод $(\quest\Pi\vec\alpha\,\varphi)^*$ эквивалентен
  арифметической формуле? Это доказано для класса задач, в которых все
  вхождения переменных лежат под $\pneg$ (зависимость от ящиков
  факторизуется через конечный вектор битов пустоты); общий случай
  открыт; при определимых ящиках перевод арифметичен относительно
  $0^{(\omega)}$.

\item \textbf{Двухтеорийный уровень.} Провести полный аудит пары
  $(\PA,\ZF)$ (модальность --- доказуемость в $\PA$, общезначимость ---
  выводимость в $\ZF$) --- она возвращает
  исчислению $\HCE$ модельно-независимую синтаксическую семантику, ---
  сопоставить её с теоретико-множественной интерпретацией Соловея
  (логика $SOL$,~\cite[гл.~13]{ru:Boolos93},~\cite{ru:EG08}) и подобрать
  верификационную теорию для $\sB$ (например,
  $\ZF$ плюс недостижимый кардинал). Сюда же примыкает вопрос о
  двойственном операторе $\Sigma\alpha$ (значение
  $\bigcup_{A}$): нужен ли он для пренексной нормализации
  $\Pi$-атомов и какова его логика.

\item \textbf{Логика $\Pi$-задач без модальностей.} Самое
  неожиданное из найденного происходит на стороне задач без всякого
  участия высказываний: фрагмент языка $\sB$ из задач, построенных из
  переменных и $\varnothing$ связками $\pto,\pand,\por$ и квантором
  $\Pi$, и его общезначимость в $\sB$ образуют логику, которую мы не
  опознали. Синтаксически это второпорядковая интуиционистская
  пропозициональная логика с одним квантором всеобщности,
  семантически --- реализуемость по Клини с равномерной квантификацией
  по множествам, та же, какой интерпретируют второпорядковую арифметику
  $\mathrm{HAS}$; аксиома $(\Pi\por)$ --- пропозициональная форма
  принципа равномерности Трулстры. На $\Pi$-свободных задачах это
  $\Lp_{\sB}$. Квантор во всех разобранных нами случаях (все задачи
  глубины не более трёх от двух переменных) устраним: например,
  $\Pi\alpha\bigl((\alpha\pto\beta\por\gamma)\pto((\alpha\pto\beta)\por
  (\alpha\pto\gamma))\bigr)$ равносильна аксиоме линейности
  $(\gamma\pto\beta)\por(\beta\pto\gamma)$, а
  $\Pi\beta\bigl(((\alpha\pto\beta)\pto\beta)\pto(\alpha\por\beta)\bigr)$
  --- закону $\alpha\por\pneg\alpha$; это родственно равномерной
  интерполяции Питтса для $\Int$, но устранение иногда опирается на
  ящики вроде $\omega\setminus\{n\}$, которые значением формулы не
  бывают (предложение~\ref{ru:prop:pienc}(4), замечание~\ref{ru:rem:peirce}). Вопросы: устраним ли квантор
  всегда; рекурсивно ли перечислима эта логика (для предикатной
  логики реализуемости Плиско доказал неарифметичность); какие
  принципы сверх $(\Pi\mathrm{E})$, правила Бернайса и $(\Pi\por)$ в
  ней верны --- первый кандидат
  $\Pi\alpha(((\alpha\pto\beta)\pto\alpha)\pto\alpha)\peq\beta$.
\end{enumerate}

\appendix

\section{Отложенные доказательства}\label{ru:sec:app}

\begin{proof}[Доказательство леммы~\ref{ru:lem:dn}]
По лемме~\ref{ru:lem:clauses}(6) множество $\val{\pneg\varphi}$ есть $\N$ при
$\val{\varphi}=\varnothing$ и $\varnothing$ иначе; применяя то же ещё раз,
получаем требуемое. Для второго утверждения раскроем определения:
$(\pneg\varphi)^*(u)$ есть $\forall v\,(\varphi^*(v)\ato \arbot )$ ---
формула, не зависящая от $u$ и эквивалентная $\aneg\exists v\,\varphi^*(v)$; поэтому
$(\pneg\pneg\varphi)^*(x)$ эквивалентна $\aneg\aneg\exists v\,\varphi^*(v)$, а
$\PA$ классична.
\end{proof}

\begin{proof}[Обоснование замечания~\ref{ru:rem:indep}]
Для $\gn{*}$ это очевидно: любые два эффективных кодирования одного и того же
конечного набора формул переводятся друг в друга вычислимо, и переводчик можно
приписать к $S$ спереди.

Для $U$ и $c$ рассуждение чуть длиннее. Пусть $\val{\cdot}$ и
$\val[*]{\cdot}'$ --- значения, посчитанные по двум наборам выборов, а
$\varphi^*$ и $\varphi^{*\prime}$ --- соответствующие переводы. Одновременной
индукцией по построению формулы строится \emph{пара} вычислимых
(равномерно по $\varphi$ и по $\gn{*}$) функций $h_\varphi$ и $g_\varphi$,
переводящих решения в обе стороны:
\[
  u\in\val{\varphi}\ \Longrightarrow\ h_\varphi(u)\in\val[*]{\varphi}' ,
  \qquad
  v\in\val[*]{\varphi}'\ \Longrightarrow\ g_\varphi(v)\in\val{\varphi} ,
\]
причём само это утверждение доказуемо в $\PA$, и его вывод строится по
$\gn{*}$. Для конъюнкции и дизъюнкции обе функции переставляют метки, заменяя
$c$ на~$c'$ и обратно; для импликации $h_{\varphi\pto\psi}$ по индексу $i$
выдаёт индекс программы $v\mapsto h_\psi\bigl(U_{i}(g_\varphi(v))\bigr)$ в
новой нумерации (теорема о параметризации), и симметрично для
$g_{\varphi\pto\psi}$. Здесь и видно, зачем нужны обе функции сразу: при
импликации направление перевода меняется на посылке. Случай $\bang A$ разбирается так: по
предположению индукции $\PA$ доказывает
$\exists x\,\varphi^*(x)\aeq\exists x\,\varphi^{*\prime}(x)$ для
всех задач $\varphi$, входящих в $A$ под знаком~$\quest$, а потому и
$A^*\aeq A^{*\prime}$; вывод этой эквивалентности вычисляется по
$\gn{*}$, и вывод формулы $A^*$ перерабатывается в вывод формулы
$A^{*\prime}$.

Значит если $S$ --- вычислимое сечение при первом наборе выборов, то
$\gn{*}\mapsto h_\varphi(S(\gn{*}))$ --- вычислимое сечение при втором. То же
рассуждение, проведённое внутри $\PA$, даёт независимость и для
$\sI$-общезначимости.

Тот же довод проходит и при замене формулы $\Prf$ на любую другую стандартную
формулу доказуемости для $\PA$: выводы одной системы перерабатываются в выводы
другой вычислимо, и корректность перевода доказуема в $\PA$. Существенно
только то, что $\Prf$ задаёт \emph{разрешимое} отношение и что $\PA$ доказывает
формализованные условия выводимости.
\end{proof}

\begin{proof}[Доказательство лемм~\ref{ru:lem:closure} и~\ref{ru:lem:closurebs}]
\emph{(MP)}. Для высказываний правило очевидно: выводимость в $\PA$ или $\ZF$ и
истинность в $\N$ замкнуты относительно modus ponens; в $\sB$ --- вместе с
замыканием (из $\ZF\vdash\forall\vec X\,A^*$ и
$\ZF\vdash\forall\vec X\,(A^*\ato B^*)$ следует
$\ZF\vdash\forall\vec X\,B^*$). Для задач: если $S_1$ решает $\varphi$, а
$S_2$ решает $\varphi\pto\psi$, то $S(c):=U_{S_2(c)}(S_1(c))$ вычислима и
по лемме~\ref{ru:lem:clauses}(6) даёт элемент $\val{\psi}$. В $\sB$ решения
собираются так же (от ящиков они не зависят), а вывод замыкания для
$\psi$ получается из выводов замыканий для посылок по той же схеме.

\emph{Подстановка.} Семантики $\sI$ и $\sE$:
 пусть $\varphi(\alpha)$ общезначима с функцией $S$, и мы
подставляем вместо $\alpha$ формулу $\chi$. Перевод $(\varphi(\chi))^*$
совпадает с переводом $\varphi(\alpha)^{*'}$, где интерпретация $*'$ совпадает с $*$ на
прочих переменных, а переменной $\alpha$ сопоставляет формулу $\chi^*$. Код
$\gn{*'}$ вычисляется по $\gn{*}$, поскольку перевод $\chi\mapsto\chi^*$ ---
вычислимая операция над кодами (определение~\ref{ru:def:tr}). Значит
$S'(\gn{*}):=S(\gn{*'})$ --- искомая функция.

Подстановка высказывания $C$ вместо переменной-высказывания $A$ в
семантиках $\sI$ и $\sE$ разбирается той же композицией интерпретаций
($A^{*'}:=C^*$).

В семантике $\sB$ подстановка задачи $\chi$ вместо свободной
переменной $\alpha$ --- это инстанцирование квантора. Свободные
вхождения $\alpha$ под знаком $\bang$ не лежат (они были бы связаны),
поэтому перевод обращается к $\alpha$ только атомами $x\in X_\alpha$.
Из $\ZF\vdash\forall\vec X\,\forall X_\alpha\,
\varphi^*\bigl(\overline{S(\gn{*})}\bigr)$ подстановкой множества
$\{x:\chi^*(x)\}$ вместо $X_\alpha$ получается замыкание перевода
$\varphi[\chi/\alpha]$: свободные множественные переменные формулы
$\chi^*$ уходят под общее замыкание $\forall\vec X$, а захвата нет по
свойству (1) безопасной подстановки. Решатель менять не нужно, поскольку ящики в код не
входят. В $\sB$ множество $\{x:\chi^*(x)\}$ --- законное содержимое
ящика при любой $\chi$, в том числе с квантором $\Pi$: ящики
пробегают \emph{все} множества (непредикативность).

Подстановка высказывания $C$ вместо переменной-высказывания $A$
в формулу
разбирается двумя способами. Если $C$ не имеет свободных
задачных переменных, то $C^*$ --- замкнутая формула языка $\ZF$, то
есть допустимая интерпретация, и работает композиция интерпретаций
$A^{*\prime}:=C^*$, в том числе для вхождений $A$ под знаком $\bang$.

Если же у $C$ есть свободные задачные переменные, то по свойству (2)
все вхождения $A$ лежат вне $\bang$. Такое возможно, только если формула
--- высказывание $B(A)$ (в задачи высказывания входят лишь через
$\bang$), и по той же причине вхождения $A$ не лежат и под знаком
$\quest$. Поэтому $(B(C))^*$ получается из $B^*$ заменой вхождений $A^*$
на $C^*$, и ни одно из них не стоит внутри гёделева номера. Так как
$B(A)$ общезначимо, а $(0=0)$ и $\arbot $ --- допустимые значения
переменной $A$, $\ZF$ доказывает замыкания $B^*(0=0)$ и $B^*\arbot $.
Рассуждая в $\ZF$ при произвольном содержимом ящиков, разберём случаи:
при  истинном $C^*$ формула $C^*$ эквивалентна $(0=0)$, при ложном $C^*$ --- формуле
$\arbot $. По замене эквивалентных $B^*(C^*)$ эквивалентна $B^*(0=0)$ или
$B^*\arbot $. Замена законна и под кванторами $\forall X_\beta$, отвечающими
$\Pi\beta$: по свойству (1) $X_\beta$ в $C^*$ не входит. Следовательно,
$\ZF$ доказывает замыкание $B^*(C^*)=(B(C))^*$.
\end{proof}

\begin{proof}[Доказательство теоремы~\ref{ru:prop:modalfree}]
Одна сторона --- это схема (A2). Обратная получается моделью.

Возьмём конечную шкалу Крипке, то есть конечное частично упорядоченное
множество $\langle P,\le\rangle$, и зафиксируем в ней \emph{максимальную}
точку~$p$. Построим модель $K$. Переменные-задачи оцениваются
произвольными открытыми вверх подмножествами $P$ (как в моделях Крипке
для $\Int$), переменные-высказывания --- произвольными истинностными
значениями. Вынуждение $x\Vdash\varphi$ задач в точках $x\in P$
определяется обычными пунктами Крипке для $\varnothing$, $\pand$, $\por$,
$\pto$; высказывания получают истинностные значения, их связки
понимаются классически. Модальности:
\[
  x\Vdash\bang A\ \ \text{при всех $x$ и всех $A$},
  \qquad
  \quest\varphi\ \text{истинно}\iff p\Vdash\varphi .
\]
Определение корректно: $\bang A$ вынуждена всюду независимо от $A$,
поэтому никакой циклической зависимости между сортами нет; вынуждение
монотонно по $\le$ (индукция по построению, для $\bang A$ тривиально).
Задача \emph{верна в $K$}, если вынуждена во всех точках при всякой
оценке переменных; высказывание --- если истинно при всякой оценке.

Проверим аксиомы. (A1) --- высказывания двузначны и связки классические;
(A2) --- схемы $\Int$ верны во всякой модели Крипке. Схемы (A3), (A5) и (A6)
тривиальны: заключение каждой из них --- задача вида $\bang X$ (у (A3)
--- после снятия обеих посылок), вынужденная всюду, а импликация с
вынужденным всюду заключением сама вынуждена всюду.

Схема (A4) --- единственное место, где работает выбор точки. Так как $p$
максимальна, единственная точка над нею --- она сама, поэтому вынуждение
в~$p$ классично:
\[
  p\Vdash\varphi\pto\psi\iff
  \bigl(p\Vdash\varphi\Rightarrow p\Vdash\psi\bigr),
\]
и точно так же $p\Vdash\varphi\pand\psi$ равносильно конъюнкции, а
$p\Vdash\varphi\por\psi$ --- дизъюнкции; наконец $p\nVdash\varnothing$.
Это и есть (A4a)--(A4d): отображение $\varphi\mapsto\quest\varphi$ --- булев
гомоморфизм.

Схема (A7): $\Sigma$-задача собрана из $\varnothing$ и задач вида
$\bang X$ связками $\pand$ и $\por$, поэтому (индукция по построению)
она либо вынуждена всюду, либо не вынуждена нигде. Для таких задач
$\pneg\pneg\sigma\pto\sigma$ вынуждена всюду: если $\sigma$ вынуждена
всюду, это ясно; если нигде, то $\pneg\sigma$ вынуждена всюду,
$\pneg\pneg\sigma$ --- нигде, и импликация вынуждена вхолостую.

Правила. (MP): если в точке вынуждены $\varphi$ и $\varphi\pto\psi$, то
вынуждена и $\psi$; для высказываний --- классика. (AtoP): $\bang A$
вынуждена всюду при любом $A$, тем более при выводимом. (PtoA): если
$\varphi$ вынуждена всюду, то и в $p$, значит $\quest\varphi$ истинно.

Итак, все аксиомы и правила $\HCI$ в этой модели верны. (Формулу Роуза
модель, разумеется, не выполняет --- она не выводима в $\Int$; именно
этим модель и полезна.) Пусть теперь задача
$\varphi$ не содержит модальностей. Тогда её вынуждение --- это обычное
вынуждение в шкале Крипке~$P$.
Если $\varphi$ выводима в $\HCI$, то она верна в каждой такой модели, то есть
в каждой конечной шкале Крипке, а значит выводима в $\Int$ по теореме о
полноте интуиционистской логики относительно конечных шкал.
\end{proof}
\section
{Почему в (A5) и (A7) взяты именно такие классы задач}
\label{ru:sec:whyclasses}

\label{ru:rem:whysigma}
Для схемы (A5) ответ короткий: шире класса $\Sigma$ ко-рефлексия
$\varphi\pto\bang\quest\varphi$ неверна уже для переменной
(предложение~\ref{ru:prop:hce}(4), где $\alpha^*:=\Con$), а семантически
(A5) --- это эффективная $\Sigma_1$-полнота со свидетелем на входе
(предложение~\ref{ru:prop:a5a7}), которая держится на том, что переводы
$\Sigma$-задач суть $\Sigma_1$-формулы. Со схемой (A7) сложнее.
Класс $\Sigma$-задач нельзя естественным образом расширить. Например, его
нельзя заменить на класс $\mathcal O$ из схемы Роуза: там схема (A7) неверна
уже на $\bang A\por\pneg\bang A$, как мы увидим ниже в этом разделе.

Прежде всего, схема (A7) перестаёт быть общезначимой на переменных: задача
$\pneg\pneg\alpha\pto\alpha$ с переменной $\alpha$ не общезначима
(предложение~\ref{ru:prop:noem}).

Для дизъюнкций не общезначима даже схема
$\pneg\pneg(\varphi\por\psi)\pto(\pneg\pneg\varphi\por\pneg\pneg\psi)$.
Слово <<даже>> оправдано вот чем. Прямым расширением (A7) на дизъюнкции
произвольных задач была бы схема $\pneg\pneg(\varphi\por\psi)\pto
(\varphi\por\psi)$, а у выписанной заключение слабее: в $\Int$ верно
$\varphi\pto\pneg\pneg\varphi$, значит $(\varphi\por\psi)\pto
(\pneg\pneg\varphi\por\pneg\pneg\psi)$. От решения выписанной схемы требуется
не решение одного из дизъюнктов, а лишь \emph{метка} --- указание, у какого
из них значение непусто (лемма~\ref{ru:lem:dn}). Тем не менее не общезначима и
она, а тогда по леммам~\ref{ru:lem:closure} и~\ref{ru:lem:closurebs} не
общезначимо и прямое расширение.
Возьмём переменные $\alpha$, $\beta$ и положим
$\alpha^*(x):=\Prf(x,\sigma)$, $\beta^*:=\aneg\Prov(\sigma)$. Тогда
$\val{\alpha\por\beta}\neq\varnothing$ при всяком $\sigma$, значит
$\val{\pneg\pneg(\alpha\por\beta)}=\N$ по лемме~\ref{ru:lem:dn}; а по
лемме~\ref{ru:lem:clauses}(5),(6) элемент множества
$\val{\pneg\pneg\alpha\por\pneg\pneg\beta}$ несёт метку, различающую случаи
$\val{\alpha}\neq\varnothing$ и $\val{\beta}\neq\varnothing$, то есть
разрешающую доказуемость~$\sigma$.

Видно и то, где именно ломается доказательство: перечисление, на котором
держится лемма~\ref{ru:lem:sigmaenum}, для переменной и для импликации взять
неоткуда. Для переменной
$\alpha$ множество $\val{\alpha}=\{n:\N\models\alpha^*(n)\}$ произвольно:
формула $\alpha^*$ может быть какой угодно, и перечислять тут нечего. Для
импликации условие <<$x\in\val{\psi\pto\chi}$>> по
лемме~\ref{ru:lem:clauses}(6) есть <<$U_x$ определена всюду на $\val{\psi}$ и
попадает в $\val{\chi}$>> --- условие с квантором всеобщности, тоже не
перечислимое.

Теперь про обещанный класс $\mathcal O$ из схемы Роуза
(предложение~\ref{ru:cor:rose}): на нём схема (A7)
тоже перестаёт быть общезначимой. Достаточно взять
$\varphi:=\bang A\por\pneg\bang A$ --- первый дизъюнкт есть $\Sigma$-задача, второй
имеет вид $\pneg\psi$, так что $\varphi$ есть $\mathcal O$-задача. При всякой
интерпретации $\val{\varphi}\neq\varnothing$: если $A^*$ доказуема, непусто
$\val{\bang A}$, а если нет, то $\val{\pneg\bang A}=\N$. Значит
$\val{\pneg\pneg\varphi}=\N$ всегда, и схема (A7) для этой $\varphi$
требовала бы вычислять по коду решение $\varphi$ --- то есть метку, отвечающую
на вопрос, доказуема ли $A^*$. Это в точности разобранный выше контрпример с
дизъюнкцией, метка которой разрешала бы доказуемость, --- записанный внутри
класса $\mathcal O$. Схема (A7) и схема Роуза поэтому не случайно заданы на разных
классах: схема Роуза требует лишь \emph{конечного списка кандидатов} --- их
выписывают заранее, а какой из них годен, узнают у собеседника, решающего
посылку, --- тогда как (A7) требует \emph{поиска}, находящего свидетеля, то
есть перечислимости значения. Класс $\mathcal O$ даёт первое и не даёт
второго: как только в дизъюнкцию входит отрицание, перечислимость пропадает.
Внутри же $\Sigma$ дизъюнкция перечислимость сохраняет --- потому там и
работает (A7), а схема Роуза, наоборот, вырождается: на $\Sigma$-задачах
она выводима из (A7) одним (MP).

Спросим теперь наоборот: нельзя ли класс сузить, ничего не потеряв? Иначе
говоря, какие из двух допущенных связок действительно усиливают схему (A7)?
Оказывается, только дизъюнкция.

Задачу $\varphi$ называют \emph{стабильной}, если для неё верна схема снятия
двойного отрицания $\pneg\pneg\varphi\pto\varphi$; в этих словах схема (A7)
говорит, что все $\Sigma$-задачи стабильны. Про связку говорят, что она
\emph{наследует} стабильность, если из стабильности аргументов следует
стабильность результата; в $\Int$ стабильность наследуют $\pand$ и $\pto$,
но не $\por$.

Конъюнкция поэтому схему (A7) не усиливает: случай
$\varphi=\psi\pand\chi$ выводится из случаев $\psi$ и
$\chi$, поскольку уже в $\Int$ из $\pneg\pneg(\psi\pand\chi)$ следует
$\pneg\pneg\psi\pand\pneg\pneg\chi$. Дизъюнкция усиливает: схема
$\pneg\pneg(\bang A\por\bang B)\pto(\bang A\por\bang B)$ из $\pneg\pneg\bang A\pto\bang A$ и
$\pneg\pneg\bang B\pto\bang B$ в $\Int$ уже не выводится: даже
$\pneg\pneg\alpha\pto\alpha$ и $\pneg\pneg\beta\pto\beta$ не дают
$\pneg\pneg(\alpha\por\beta)\pto(\alpha\por\beta)$ (шкала Крипке с корнем и
двумя несравнимыми потомками, в которых истинны соответственно $\alpha$
и~$\beta$). Между тем
лемма~\ref{ru:lem:sigmaenum} эту схему даёт: параллельный перебор выводов обеих
формул. Значит всё содержание (A7), сверх марковского принципа
$\pneg\pneg\bang A\pto\bang A$, состоит в том, что \emph{внутри класса $\Sigma$
стабильность объявляется наследуемой дизъюнкцией}.

Конъюнкцию, однако, из грамматики убрать всё равно нельзя: без неё не
построить дизъюнкций вроде $(\bang A\pand\bang B)\por\bang C$, для которых схема (A7)
по только что сказанному уже не выводится.

%

\section{Независимость схемы (A7)}\label{ru:sec:a7indep}

Схема (A7) не выводится из остальных аксиом $\HCgl$. Доказательство ---
модель Крипке из двух точек.

\begin{definition}[двухточечная модель $M$]\label{ru:def:model3}
Шкала --- двухэлементная цепь $a<b$. Переменные-задачи оцениваются
произвольными открытыми вверх множествами ($\varnothing$, $\{b\}$,
$\{a,b\}$), переменные-высказывания --- произвольными истинностными
значениями; связки задач понимаются по пунктам Крипке, связки
высказываний --- классически, $\varnothing$ не вынуждена нигде, $\bot$
ложно. Модальности:
\[
  x\Vdash\bang A\iff\bigl(x=b\ \text{ или }\ A\ \text{истинно}\bigr),
  \qquad
  \quest\varphi\ \text{истинно}\iff b\Vdash\varphi .
\]
Иначе говоря, $\bang A$ вынуждена в верхней точке всегда, а в нижней ---
тогда и только тогда, когда $A$ истинно (вынуждение монотонно); $\quest\varphi$
означает, что $\varphi$ вынуждена хотя бы в одной точке, поскольку всякое
непустое открытое вверх множество содержит $b$. Формула \emph{верна в
$M$}, если при всякой оценке переменных она вынуждена в обеих точках (для
задач) или истинна (для высказываний).
\end{definition}

\begin{proposition}\label{ru:prop:a7indep}
Все схемы (A1)--(A6) верны в $M$, и все три правила (MP), (AtoP), (PtoA)
сохраняют верность в $M$. Между тем задача $\pneg\pneg\bang\bot\pto\bang\bot$ ---
частный случай (A7) --- в $M$ не верна: она не вынуждена в точке $a$.
Следовательно, схема \textup{(A7)} не выводима в $\HCgl$ из остальных аксиом.
\end{proposition}

\begin{proof}
Заметим сразу: $\bang A$ вынуждена в $b$ при любом $A$, поэтому
\[
  \square A=\quest\bang A\quad\text{истинно при всех }A .
\]
В этой модели <<доказуемо всё>>; именно поэтому схема Лёба в ней верна даром.

Два общих наблюдения. Во-первых, точка $b$ максимальна, поэтому
вынуждение в ней классично: $b\Vdash\varphi\pto\psi$ равносильно
$(b\Vdash\varphi\Rightarrow b\Vdash\psi)$, и так же для $\pand$, $\por$.
Во-вторых, для задач вида $\bang X$ вынуждение в $a$ --- это истинность
$X$, а в $b$ оно автоматическое; поэтому для задач, составленных из
задач $\bang X$ импликациями, вынуждение в $a$ сводится к импликации
истинностей: $a\Vdash\bang X\pto\bang Y$ равносильно
$(X\ \text{истинно}\Rightarrow Y\ \text{истинно})$, и так же для
вложенных импликаций.

(A1) Высказывания двузначны, связки классические.

(A2) Схемы $\Int$ верны во всякой модели Крипке.

(A3) В $b$ всё вынуждено. В $a$ по второму наблюдению нужно: если
$A\ato B$ и $A$ истинны, то истинно $B$.

(A4) По первому наблюдению $\quest\varphi$ истинно тогда и только тогда,
когда $b\Vdash\varphi$, а вынуждение в $b$ классично; отсюда (A4a)--(A4c),
а (A4d) --- это $b\nVdash\varnothing$.

(A5) Если $\quest\sigma$ истинно, то $\bang\quest\sigma$ вынуждена всюду
и импликация $\sigma\pto\bang\quest\sigma$ вынуждена всюду. Если
$\quest\sigma$ ложно, то $b\nVdash\sigma$, а тогда по монотонности
$\sigma$ не вынуждена нигде, и импликация вынуждена вхолостую.

(A6) $\square A$ истинно, поэтому $\square A\ato A$ истинно тогда и
только тогда, когда истинно $A$; по второму наблюдению
$a\Vdash\bang(\square A\ato A)\pto\bang A$, а в $b$ всё вынуждено.

Правила. (MP) --- обычное свойство вынуждения и классической импликации.
(AtoP): если $A$ истинно, то $\bang A$ вынуждена в обеих точках. (PtoA):
если $\varphi$ вынуждена в обеих точках, то и в $b$, и $\quest\varphi$
истинно.

Наконец, $\bot$ ложно, поэтому $\bang\bot$ вынуждена только в $b$. Тогда
$\pneg\bang\bot$ не вынуждена нигде: в $b$ --- потому что $b\Vdash\bang\bot$,
в $a$ --- потому что над $a$ лежит $b$. Значит, $\pneg\pneg\bang\bot$
вынуждена всюду, и в точке $a$ посылка задачи
$\pneg\pneg\bang\bot\pto\bang\bot$ вынуждена, а заключение --- нет.
\end{proof}

\section{Что ново в выводе схемы Лёба}\label{ru:sec:novelty}

Здесь мы разбираем, есть ли что-либо новое 
в нашем доказательстве схемы Лёба (пункт~(3) доказательства
теоремы~\ref{ru:th:sound}).
В модальном языке такой вывод известен и стандартен: в
$\GL$ без схемы Лёба неподвижная точка $\psi\aeq(\square\psi\ato
\sigma)$ даёт формулу Лёба. Действительно, из $\psi\ato(\square\psi\ato\sigma)$
правилом необходимости и схемой K получаем
$\square\psi\ato\square(\square\psi\ato\sigma)$, далее по K и схеме
$\square\psi\ato\square\square\psi$ --- импликацию $\square\psi\ato\square\sigma$;
отсюда $(\square\sigma\ato\sigma)\ato(\square\psi\ato\sigma)$, а значит и
$\square(\square\sigma\ato\sigma)\ato\square(\square\psi\ato\sigma)$, что вместе с
$\square(\square\psi\ato\sigma)\ato\square\psi$ (вторая половина неподвижной
точки) даёт $\square(\square\sigma\ato\sigma)\ato\square\sigma$.

Поэтому в теореме~\ref{ru:th:lob} ново не то, что схема Лёба выводится из
неподвижной точки, а то, \emph{в какой форме} она выводится. Схема (A6)
записана через $\bang$: она утверждает существование \emph{вычислимого
преобразования выводов}, а не существование выводов, и в модальном языке
невыразима. Из (A6) $\square$-версия получается применением (PtoA) и (A4a);
обратный переход неочевиден и остаётся открытым (см.~\S\ref{ru:sec:todo}).

Содержательно получается вот что: \emph{конструктивная теорема Лёба выводится
из эффективной $\Sigma_1$-полноты}. Обычно её доказывают, разбирая по шагам
доказательство теоремы Лёба и следя за тем, что все преобразования выводов
вычислимы; здесь же достаточно эффективной $\Sigma_1$-полноты --- схемы (A5) --- и
эффективности леммы о неподвижной точке, а остальное делает формальный вывод.
Функцию, перерабатывающую вывод формулы $\Prov(\sigma)\ato\sigma$ в вывод
формулы $\sigma$, строить руками не нужно: она читается из этого вывода.


\begin{thebibliography}{9}

\bibitem{Kolm32} A.~Kolmogoroff. Zur Deutung der intuitionistischen Logik.
  \emph{Mathematische Zeitschrift}, Bd.~35, H.~1, 1932, S.~58--65.

\bibitem{Kolm85} A.\,N.~Kolmogorov. On the papers on intuitionistic logic.
  In: \emph{Selected Works of A.\,N.~Kolmogorov}, vol.~1, Mathematics and
  its Applications (Soviet Series), vol.~25, Kluwer, Dordrecht, 1991,
  pp.~451--452. (Russian original: Избранные труды. Математика и
  механика. М.: Наука, 1985, с.~393.)

\bibitem{Goedel33} K.~G\"odel. Eine Interpretation des intuitionistischen
  Aussagenkalk\"uls. \emph{Ergebnisse eines mathematischen Kolloquiums}
  4 (1933), 39--40.

\bibitem{Melikhov13} S.\,A.~Melikhov. A Galois connection between classical
  and intuitionistic logics. I: Syntax. arXiv:1312.2575, 2013--2017.
  \emph{(The logic $\HC$ is defined here.)}

\bibitem{Melikhov15} S.\,A.~Melikhov. A Galois connection between classical
  and intuitionistic logics. II: Semantics. arXiv:1504.03379, 2015--2018.

\bibitem{Melikhov24} S.\,A.~Melikhov. A joint logic of problems and
  propositions. \emph{Doklady Mathematics} 109(2):130--139, 2024.

\bibitem{Onoprienko20} A.\,A.~Onoprienko. Kripke semantics for the logic of
  problems and propositions. \emph{Sbornik: Mathematics} 211(5):709--732, 2020.

\bibitem{Kleene45} S.\,C.~Kleene. On the interpretation of intuitionistic
  number theory. \emph{Journal of Symbolic Logic} 10(4):109--124, 1945.

\bibitem{Rose53} G.\,F.~Rose. Propositional calculus and realizability.
  \emph{Transactions of the American Mathematical Society} 75(1):1--19, 1953.

\bibitem{Solovay76} R.\,M.~Solovay. Provability interpretations of modal
  logic. \emph{Israel Journal of Mathematics} 25(3--4):287--304, 1976.

\bibitem{Medvedev62} Yu.\,T.~Medvedev. Finite problems. \emph{Soviet
  Mathematics Doklady} 3 (1962), 227--230. (Russian original: Финитные
  задачи. \emph{Доклады АН СССР}, т.~142, вып.~5, 1962, с.~1015--1018.)

\bibitem{Yankov63} V.\,A.~Yankov. On realizable formulas of propositional
  logic. \emph{Soviet Mathematics Doklady} 4 (1963), 1146--1148. (Russian
  original: О реализуемых формулах логики высказываний. \emph{Доклады АН
  СССР}, т.~151, вып.~5, 1963, с.~1035--1037.)

\bibitem{KP57} G.~Kreisel, H.~Putnam. Eine Unableitbarkeitsbeweismethode
  f\"ur den intuitionistischen Aussagenkalk\"ul. \emph{Archiv f\"ur
  mathematische Logik und Grundlagenforschung} 3(3--4):74--78, 1957.

\bibitem{Yankov68} V.\,A.~Yankov. The calculus of the weak ``law of excluded
  middle''. \emph{Mathematics of the USSR---Izvestiya} 2(5):997--1004, 1968.
  (Russian original: Об исчислении слабого закона исключённого третьего.
  \emph{Известия АН СССР}, сер.~матем., т.~32, вып.~5, 1968, с.~1044--1051.)

\bibitem{VSSCh} N.\,K.~Vereshchagin, D.\,P.~Skvortsov, E.\,Z.~Skvortsova,
  A.\,V.~Chernov. Variants of realizability for propositional formulas and
  the logic of weak excluded middle. \emph{Proc. Steklov Inst. Math.}
  242 (2003), 67--85. (Russian original: \emph{Труды МИАН}, т.~242, 2003,
  с.~77--97.)

\bibitem{Pawlowski26} P.~Pawlowski. Medvedev logic is not decidable. It is
  $\Pi^0_1$-complete. Who would have guessed? arXiv:2609.11576, 2026.

\bibitem{AlmeidaKnudstorp26} R.\,N.~Almeida, S.\,B.~Knudstorp. Medvedev
  logic is undecidable. arXiv:2609.13359, 2026.

\bibitem{Plisko09} V.\,E.~Plisko. A survey of propositional realizability
  logic. \emph{Bulletin of Symbolic Logic} 15(1):1--42, 2009.

\bibitem{Artemov01} S.\,N.~Artemov. Explicit provability and constructive
  semantics. \emph{Bulletin of Symbolic Logic} 7(1):1--36, 2001.

\bibitem{Sidon97} T.~Sidon. Provability logic with operations on proofs.
  \emph{Proc. LFCS'97}, Lecture Notes in Computer Science 1234, 342--353,
  Springer, 1997.

\bibitem{Nogina14} E.~Nogina. On logic of formal provability and explicit
  proofs. arXiv:1405.2559, 2014.

\bibitem{Protopopescu16} T.~Protopopescu. An arithmetical interpretation of
  verification and intuitionistic knowledge. \emph{Proc. LFCS 2016}, Lecture
  Notes in Computer Science 9537, 317--330, Springer, 2016; arXiv:1601.03059.

\bibitem{AM19} M.~Ardeshir, M.~Mojtahedi. The $\Sigma_1$-provability logic of
  $\mathrm{HA}^*$. \emph{Journal of Symbolic Logic} 84(3):1118--1135, 2019.

\bibitem{VZ19} A.~Visser, J.~Zoethout. Provability logic and the completeness
  principle. \emph{Annals of Pure and Applied Logic} 170(6), 2019;
  arXiv:1804.09451.

\bibitem{Iemhoff99} R.~Iemhoff. A modal analysis of some principles of
  the provability logic of Heyting Arithmetic. Logic Group Preprint
  Series, Utrecht University, 1999;
  \texttt{https://dspace.library.uu.nl/handle/1874/27770}.

\bibitem{DvdGM24} A.~Das, I.~van der Giessen, S.~Marin. Intuitionistic
  G\"odel--L\"ob logic, \`a la Simpson: labelled systems and birelational
  semantics. \emph{CSL 2024}, LIPIcs 288, article~22; arXiv:2309.00532.

\bibitem{Boolos93} G.~Boolos. \emph{The Logic of Provability.} Cambridge
  University Press, 1993.

\bibitem{EG08} L.~Esakia, R.~Grigolia. Formulas of one propositional
  variable in intuitionistic logic with the Solovay modality.
  \emph{Logic and Logical Philosophy}, 17 (2008), 111--127.

\end{thebibliography}

\begin{thebibliography}{9}

\bibitem{ru:Kolm32} A.~Kolmogoroff. Zur Deutung der intuitionistischen Logik.
  \emph{Mathematische Zeitschrift}, Bd.~35, H.~1, 1932, S.~58--65.

\bibitem{ru:Kolm85}
Колмогоров A.H. К работам no интуиционистской логике. Избранные труды // Математи-
ка и механика. М.: Наука, 1985. С. 393; English
transl., On the papers on intuitionistic logic. Selected
Works of A.N. Kolmogorov // Mathematics and
its Applications. 1991. V. 1. P. 451–452. Kluwer,
Dordrecht: Soviet Series. V. 25.


\bibitem{ru:Goedel33} K. G\"odel. Eine Interpretation des intuitionistischen
  Aussagenkalk\"uls. \emph{Ergebnisse eines mathematischen Kolloquiums}
  4 (1933), 39--40.

\bibitem{ru:Melikhov13} S.\,A. Melikhov. A Galois connection between classical
  and intuitionistic logics. I: Syntax. arXiv:1312.2575, 2013--2017.
  \emph{(Здесь определена логика $\HC$.)}

\bibitem{ru:Melikhov15} S.\,A. Melikhov. A Galois connection between classical
  and intuitionistic logics. II: Semantics. arXiv:1504.03379, 2015--2018.

\bibitem{ru:Melikhov24} S.\,A. Melikhov. A joint logic of problems and
  propositions. \emph{Doklady Mathematics} 109(2):130--139, 2024.

\bibitem{ru:Onoprienko20} A.\,A. Onoprienko. Kripke semantics for the logic of
  problems and propositions. \emph{Sbornik: Mathematics} 211(5):709--732, 2020.

\bibitem{ru:Kleene45} S.\,C. Kleene. On the interpretation of intuitionistic
  number theory. \emph{Journal of Symbolic Logic} 10(4):109--124, 1945.

\bibitem{ru:Rose53} G.\,F. Rose. Propositional calculus and realizability.
  \emph{Transactions of the American Mathematical Society} 75(1):1--19, 1953.

\bibitem{ru:Solovay76} R.\,M. Solovay. Provability interpretations of modal
  logic. \emph{Israel Journal of Mathematics} 25(3--4):287--304, 1976.

\bibitem{ru:Medvedev62} Ю.\,Т.~Медведев. Финитные задачи.
\emph{Доклады АН СССР}, т.~142, вып.~5, 1962, с.~1015--1018.

\bibitem{ru:Yankov63} В.\,А.~Янков. О реализуемых формулах логики
высказываний. \emph{Доклады АН СССР}, т.~151, вып.~5, 1963, с.~1035--1037.

\bibitem{ru:KP57} G.~Kreisel, H.~Putnam. Eine Unableitbarkeitsbeweismethode
für den intuitionistischen Aussagenkalkül. \emph{Archiv für mathematische
Logik und Grundlagenforschung} 3(3--4):74--78, 1957.

\bibitem{ru:Yankov68} В.\,А.~Янков. Об исчислении слабого закона исключённого
третьего. \emph{Известия АН СССР}, сер.~матем., т.~32, вып.~5, 1968,
с.~1044--1051.

\bibitem{ru:VSSCh} Н.\,К.~Верещагин, Д.\,П.~Скворцов, Е.\,З.~Скворцова,
А.\,В.~Чернов. Варианты понятия реализуемости для пропозициональных формул,
приводящие к логике слабого закона исключенного третьего.
\emph{Труды МИАН}, т.~242, 2003, с.~77--97. Английский перевод:
N.\,K.~Vereshchagin, D.\,P.~Skvortsov, E.\,Z.~Skvortsova, A.\,V.~Chernov.
Variants of realizability for propositional formulas and the logic of weak
excluded middle. \emph{Proc. Steklov Inst. Math.}, 242 (2003), 67--85.

\bibitem{ru:Pawlowski26} P.~Pawlowski. Medvedev logic is not decidable. It is
  $\Pi^0_1$-complete. Who would have guessed? arXiv:2609.11576, 2026.

\bibitem{ru:AlmeidaKnudstorp26} R.\,N. Almeida, S.\,B. Knudstorp. Medvedev
  logic is undecidable. arXiv:2609.13359, 2026.

\bibitem{ru:Plisko09} V.\,E. Plisko. A survey of propositional realizability
  logic. \emph{Bulletin of Symbolic Logic} 15(1):1--42, 2009.

\bibitem{ru:Artemov01} S.\,N. Artemov. Explicit provability and constructive
  semantics. \emph{Bulletin of Symbolic Logic} 7(1):1--36, 2001.

\bibitem{ru:Sidon97} T.~Sidon. Provability logic with operations on proofs.
  \emph{Proc. LFCS'97}, Lecture Notes in Computer Science 1234, 342--353,
  Springer, 1997.

\bibitem{ru:Nogina14} E.~Nogina. On logic of formal provability and explicit
  proofs. arXiv:1405.2559, 2014.

\bibitem{ru:Protopopescu16} T.~Protopopescu. An arithmetical interpretation of
  verification and intuitionistic knowledge. \emph{Proc. LFCS 2016}, Lecture
  Notes in Computer Science 9537, 317--330, Springer, 2016; arXiv:1601.03059.

\bibitem{ru:AM19} M.~Ardeshir, M.~Mojtahedi. The $\Sigma_1$-provability logic of
  $\mathrm{HA}^*$. \emph{Journal of Symbolic Logic} 84(3):1118--1135, 2019.

\bibitem{ru:VZ19} A.~Visser, J.~Zoethout. Provability logic and the completeness
  principle. \emph{Annals of Pure and Applied Logic} 170(6), 2019;
  arXiv:1804.09451.

\bibitem{ru:Iemhoff99} R.~Iemhoff. A modal analysis of some principles of
  the provability logic of Heyting Arithmetic. Logic Group Preprint
  Series, Utrecht University, 1999;
  \texttt{https://dspace.library.uu.nl/handle/1874/27770}.

\bibitem{ru:DvdGM24} A.~Das, I.~van der Giessen, S.~Marin. Intuitionistic
  G\"odel--L\"ob logic, \`a la Simpson: labelled systems and birelational
  semantics. \emph{CSL 2024}, LIPIcs 288, article~22; arXiv:2309.00532.

\bibitem{ru:Boolos93} G. Boolos. \emph{The Logic of Provability.} Cambridge
  University Press, 1993.

\bibitem{ru:EG08} L.~Esakia, R.~Grigolia. Formulas of one propositional
  variable in intuitionistic logic with the Solovay modality.
  \emph{Logic and Logical Philosophy}, 17 (2008), 111--127.

\end{thebibliography}
\end{document}